\documentclass[12pt,oneside,english,reqno,a4paper]{article}
\usepackage[utf8]{inputenc}
\usepackage{amsmath,amsthm}
\usepackage{lmodern}
\usepackage[T1]{fontenc}
\usepackage{fixltx2e,babel}
\usepackage[numbers]{natbib}
\usepackage{graphicx}
\usepackage{url}
\usepackage{ifpdf}
\ifpdf
  \pdfoutput=1
    \usepackage[bitstream-charter]{mathdesign}
    \usepackage[pdfusetitle,colorlinks=true,linktoc=all]{hyperref}
\else
    \usepackage{times}
    \usepackage[dvips,ps2pdf,linktoc=all]{hyperref}
\fi
\ifdefined\hypersetup
  \hypersetup{
    pdfkeywords={}, linkcolor=black, citecolor=black, filecolor=black, urlcolor=black,
  }
\fi
\usepackage{fancyhdr}
\usepackage{blindtext}
\usepackage{fancyhdr}
\usepackage{algpseudocode}
\usepackage{tikz}
\usetikzlibrary{arrows.meta}

\newtheorem{theorem}{Theorem}[section]
\newtheorem{lemma}[theorem]{Lemma}
\newtheorem{definition}[theorem]{Definition}
\newtheorem{remark}[theorem]{Remark}
\newtheorem{example}[theorem]{Example}
\newtheorem{corollary}[theorem]{Corollary}

\begin{document}
\pagenumbering{Alph}
\renewcommand{\thepage}{C-\Roman{page}}
\begin{titlepage}
      \textsf{
    \begin{center}
      \vspace*{1cm}
    \huge \textbf{$p$-Laplacian Operators on Hypergraphs} \\
    \vspace{2cm}
    \LARGE\textbf{Master's Thesis}\\[5mm]
    \normalsize
    by \textbf{Ariane Fazeny}
    \vfill
    Advisors: MSc. Piero Deidda and Dr. Daniel Tenbrinck\\[5mm]
    Examiner: Prof. Dr. Martin Burger
    \\\vspace{2cm}
    This work is licensed under a Creative Commons Attribution-NonCommercial-ShareAlike 4.0 International License.\\[3em]
    \end{center}
    \begin{minipage}{0.36\textwidth}
    \vspace*{10mm}
    \includegraphics[width=\textwidth]{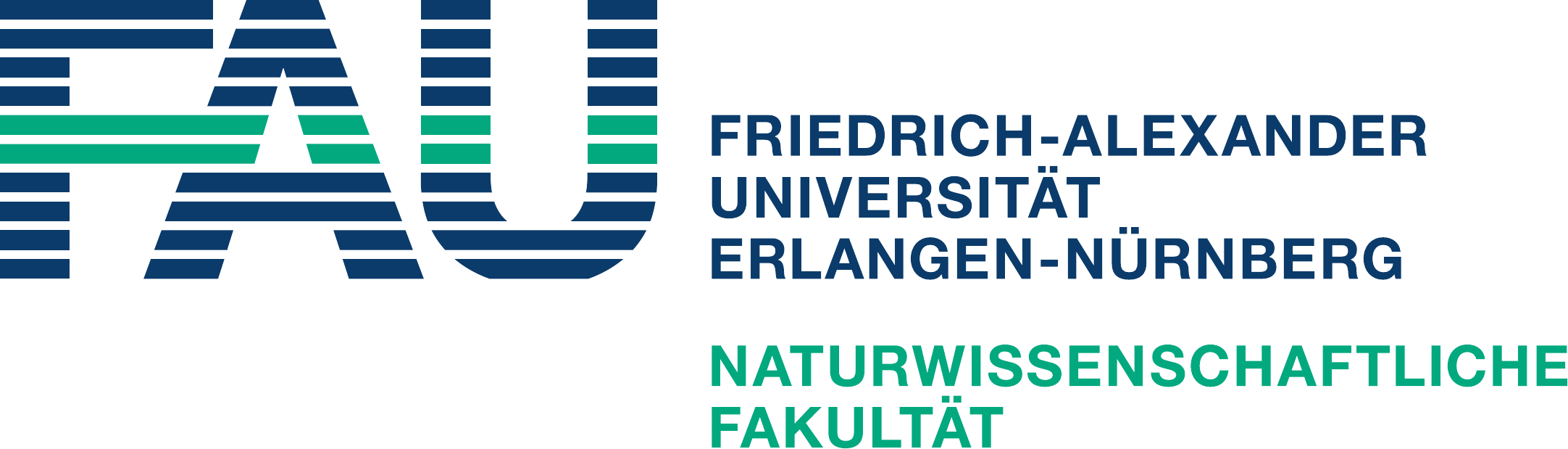}
    \end{minipage}
    \hfill
    \begin{minipage}{0.3\textwidth}
    \includegraphics[width=\textwidth]{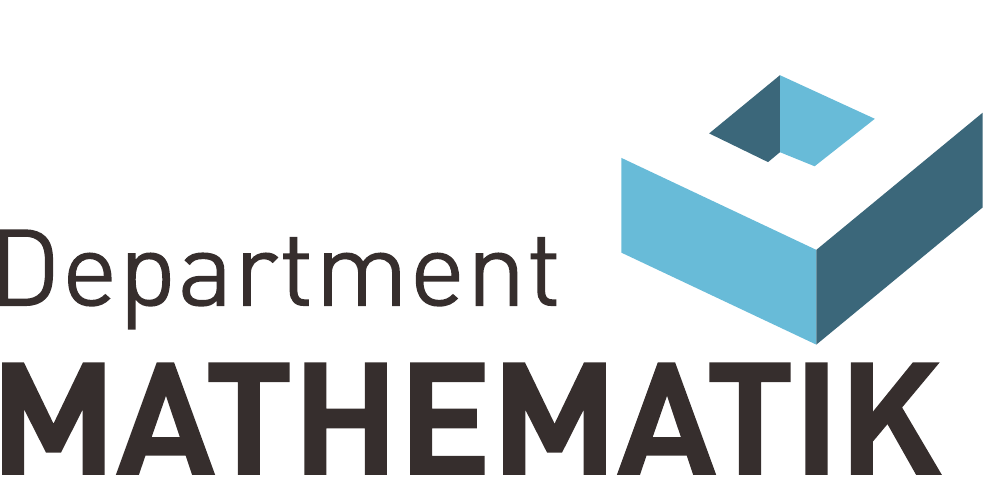} \\[3mm]
    \end{minipage}
}
\end{titlepage}

\pagenumbering{Roman}
\clearpage
\addcontentsline{toc}{section}{Abstract}
\section*{Abstract}

Literature about $p$-Laplacian operators on hypergraphs already exists, however the introduced concepts are mostly tailored to applications such as signal processing or chemical reaction networks. Therefore, this thesis introduces generalized definitions for the gradient, the adjoint and the $p$-Laplacian operators on hypergraphs and focuses on the theoretical background and properties of the described concepts. Firstly this thesis introduces normal graphs and hypergraphs (not oriented and oriented) and shows, by critically analyzing two normal graph representations of hypergraphs, why the extension to hypergraphs is justified and even necessary in order to describe higher-level connections between multiple objects. Secondly, this thesis generalizes the already existing definitions for gradients, adjoints and $p$-Laplacians (for vertices and arcs) in the normal graph setting. Lastly, this thesis extends the concepts of the normal graph case to the hypergraph case and thereby introduces generally applicable definitions for gradients, adjoints and $p$-Laplacians for both vertices and hyperarcs.
\clearpage
\addcontentsline{toc}{section}{Acknowledgements}
\section*{Acknowledgments}

I would like to thank Prof. Dr. Martin Burger for welcoming me at his chair with such open arms and for sharing his expertise about this topic with me. I also want to thank him for taking the time to discuss the general direction of this thesis and for giving valuable input about how to further develop my ideas.\\

Moreover, I want to say thank you to my advisors MSc. Piero Deidda and Dr. Daniel Tenbrinck for supporting me through all the steps of this thesis, for providing so much useful feedback and for encouraging me during the more difficult times. Daniel, thank you very much for even meeting up in the evening to work on the most tricky parts of my thesis and Piero, thank you very much for joining us online even though you are currently busy enough with your PhD thesis. I really appreciate your help and all the time you dedicated to this thesis.\\

I am really looking forward to continuing to work with you!
\clearpage
\addcontentsline{toc}{section}{Introduction}
\section*{Introduction}

The aim of this thesis is to generalize already existing definitions of gradients, adjoints, divergences, Laplacians and $p$-Laplacians both on oriented normal graphs and oriented hypergraphs. Currently, most publications only introduce application-specific concepts tailored to the corresponding field such as chemistry, biology, sociology or electrical engineering. Instead of introducing another application-focused version of calculus on oriented normal graphs and hypergraphs, this thesis provides very versatile concepts, which match the specific definitions in related publications.\\

Since hypergraphs are a rather new topic, the thesis starts with a detailed introduction to hypergraphs, especially highlighting how hypergraphs can be represented algebraically by different adjacency and incidence matrices and adjacency tensors. Furthermore, the thesis also includes an extensive analysis on why hypergraphs are a natural extension of the typical normal graph and it includes proof that every normal graph is a special case of a hypergraph.\\

Moreover, this thesis examines two different normal graph representations of hypergraphs, which could undermine the necessity of introducing hypergraphs in the first place. The two normal graph representations, one using bipartite normal graphs and one using normal graphs with complete subgraphs, however show major drawbacks compared to hypergraphs. The first disadvantage is the generally non-unique mapping back to the underlying hypergraph and the second drawback is a possibly strongly increasing number of vertices and $\backslash$ or edges and arcs in the normal graph representation. This implies that using the more complex concept of hypergraphs, instead of larger normal graphs, is justified.\\

The second part of the thesis then focuses on introducing gradient, adjoint, divergence, Laplacian and $p$-Laplacian operators on oriented normal graphs and hypergraphs. Compared to the already existing definitions in other publications, these operators are more general and can be individually adapted to different use cases by choosing different parameters and weight functions.\\

Section (\ref{1}) gives all fundamental definitions for normal graphs together with describing their algebraic representation and section (\ref{2}) then introduces hypergraphs as well as their matrix and tensor representation. Section (\ref{3}) lists some further definitions on normal graphs and hypergraphs, which are needed in later parts of this thesis. Furthermore, section (\ref{4}) focuses on the connection between normal graphs and hypergraphs by proving that every normal graph is a hypergraph and by critically analyzing two normal graph representations of hypergraphs.\\

\clearpage
Section (\ref{5}) then introduces the concept of real functions being defined on the domain of vertices, arcs or edges of normal graphs, for which in section (\ref{6}) vertex and arc gradient operators as well as vertex and arc adjoint operators are defined. Section (\ref{7}) uses these operators to calculate the corresponding divergence, Laplacian and $p$-Laplacian operators on normal graphs. Section (\ref{8}) generalizes real vertex, arc and edge functions on normal graphs to real vertex, hyperarc and hyperedge functions on hypergraphs and section (\ref{9}) uses these generalized functions on hypergraphs to give definitions for the vertex gradient, the hyperarc gradient, the vertex adjoint and the hyperarc adjoint on oriented hypergraphs. Lastly, section (\ref{10}) presents the divergence, Laplacian and $p$-Laplacian operators both for vertices and hyperarcs on hypergraphs as a extension from the normal graph case.\\

\clearpage
\addcontentsline{toc}{section}{Related work}
\section*{Related work}

The following publications discuss the topic of gradients, adjoints, Laplacians or $p$-Laplacians on graphs or hypergraphs. Furthermore, most of the sources below include useful applications of the newly introduced concepts.

\begin{itemize}
	\item Elmoataz, Toutain, Tenbrinck (2015): On the $p$-Laplacian and $\infty$-Laplacian on Graphs with Applications in Image and Data Processing \cite{elmoataz2015p}\\
	
	This paper introduces $p$-Laplacian and $\infty$-Laplacian operators on normal graphs as well as the necessary gradient, adjoint and divergence definitions, which match the definitions of the continuum setting. The described applications include image restoration, image simplification, interpolation, active contours, semisupervised segmentation, classification and nonlocal image inpainting.
	
	The gradient, adjoint, divergence, Laplacian and $p$-Laplacian operators for oriented normal graphs serve as a starting point for this thesis to generalize the corresponding operators in the normal graph case and to extend the generalized definitions in a second step to the hypergraph case.\\
	
	\item Jost, Mulas, Zhang (2021): $p$-Laplace Operators for Oriented Hypergraphs \cite{jost2021plaplaceoperators}\\
	
	This publication includes both, a vertex and a hyperedge $p$-Laplacian, for a general oriented hypergraph (Note: hyperedges of \cite{jost2021plaplaceoperators} are called hyperarcs in this thesis). Moreover, the applications discussed in this paper include vertex and hyperedge partition problems as well as smallest, largest and smallest nonzero eigenvalue problems of the $1$-Laplacians.
	
	The definitions for the vertex and hyperedge $p$-Laplacians in \cite{jost2021plaplaceoperators} are a special case of the more general definitions for the vertex and hyperarc $p$-Laplacian introduced in this thesis.\\
	
	\item Jost, Mulas (2019): Hypergraph Laplace Operators for Chemical Reaction Networks \cite{jost2019hypergraph}\\
	
	This paper introduces a special case of an oriented hypergraph called chemical hypergraph, which is specifically tailored for the application of hypergraphs to chemical reaction networks. Hence, each hyperedge (Note: hyperedges of \cite{jost2019hypergraph} are called hyperarcs in this thesis) of the hypergraph encodes a chemical reaction of the elements, which are encoded as vertices. Moreover, this paper includes a definition of the Laplacian operator and analyses the corresponding eigenvalues with a focus on the application to chemical reactions.\\
	
	\item Mulas, Kuehn, Böhle, Jost (2021): Random walks and Laplacians on hypergraphs, When do they match? \cite{mulas2022random}\\
	
	This publication introduces the concept of random walks on not oriented hypergraphs. The corresponding random walk Laplacians for hypergraphs are compared to the normalized hypergraph Laplacians from \cite{jost2021plaplaceoperators}, which are based on applications in biology or chemistry instead of random walks. Furthermore, this paper also includes a spectral analysis of the two different hypergraph Laplacians.\\	
	
	\item Ding, Cui (2020): Introducing Hypergraph Signal Processing, Theoretical Foundation and Practical Applications \cite{zhang2019introducing}\\
	
	This paper introduces a new framework for hypergraph signal processing and new methods for data analysis applications, which rely on the tensor representation of not oriented hypergraphs in order to capture higher-order interactions. The presented results show a significant improvement in accuracy when using the newly introduced hypergraph signal processing framework compared to traditional signal processing solutions. Their applications include data compression, spectral clustering, classification and denoising, and Internet of Things, social networks and natural language processing are listed as further possible applications for the hypergraph signal processing framework.
	
	Especially their definition of the tensor representation of a not oriented hypergraph is relevant for this thesis. The publication also includes a definition of the Laplacian on hypergraphs as $L = D - A$ based on the definition in the normal graph case, but with tensors instead of matrices.\\
\end{itemize}

\clearpage
\tableofcontents
\clearpage

\pagenumbering{arabic}

\section{Normal graphs}\label{1}

Given a data set consisting of $N \in \mathbb{N} \backslash \left\{0, 1\right\}$ data points $v_1, v_2, \dots v_N$, the aim is to capture and describe complex relationships between these data points using first normal graphs and later on hypergraphs. In both cases, the given data points $v_1, v_2, \dots v_N$ are represented by the vertices of graphs and the relationships within the data points are described by edges, arcs, hyperedges or hyperarcs \big(depending on the type of graph\big).\\

Subsection (\ref{1.1}) introduces the definitions of the not oriented and the oriented normal graph. These definitions are commonly known as graph and directed graph in standard graph theory, where vertices are linked through pairwise connections called edges or arcs. Furthermore, it discusses the difference between the orientation and the direction of graphs and gives a definition of the oriented normal graph with switched orientation. Finally, subsection (\ref{1.2}) gives definitions for adjacency matrices and incidence matrices, which can be used to describe normal graphs algebraically.\\

\subsection{General definitions for normal graphs}\label{1.1}

In order to highlight the enhanced possibilities of hypergraphs in comparison to the standard normal graphs, the definitions of the not oriented and oriented normal graphs are given below. They are based on the definition of simple graphs and simple digraphs in \cite{williamson2010lists}.\\

\begin{definition}[\textbf{Not oriented normal graph $NG$}]\label{NG} (Page 148 in \cite{williamson2010lists}: Simple graph) \ \\
	A not oriented normal graph $NG = \left(\mathcal{V}, \mathcal{E}_G\right)$ is a tuple consisting of a finite set of vertices $\mathcal{V} = \left\{v_1, v_2, \dots v_N\right\}$ and a set of edges $\mathcal{E}_G \subseteq \left\{\left\{v_i, v_j\right\} ~ \middle| ~ v_i, v_j \in \mathcal{V}, v_i \neq v_j\right\}$.\\
\end{definition}

\begin{example}[\textbf{Not oriented normal graph}]\ \\
	Given a set of vertices
	\begin{equation*}
		\mathcal{V} = \left\{v_1, v_2, v_3, v_4, v_5, v_6\right\}
	\end{equation*} together with a set of edges
	\begin{equation*}
		\mathcal{E}_G = \left\{\left\{v_2, v_4\right\}, \left\{v_2, v_5\right\}, \left\{v_3, v_6\right\}, \left\{v_4, v_5\right\}\right\},
	\end{equation*}
	then the not oriented normal graph $NG = \left(\mathcal{V}, \mathcal{E}_G\right)$ can be visualized in the following way:
	
	\begin{center}\begin{tikzpicture}
		\begin{scope}[every node/.style={circle,thick,draw}]
			\node (v1) at (0,0) {$v_1$};
			\node (v2) at (3,0) {$v_2$};
			\node (v3) at (6,0) {$v_3$};
			\node (v4) at (0,-3) {$v_4$};
			\node (v5) at (3,-3) {$v_5$};
			\node (v6) at (6,-3) {$v_6$};
		\end{scope}
		
		\begin{scope}[>={Stealth[black]},
			every node/.style={fill=white,circle},
			every edge/.style={draw=gray,very thick}]
			\path [-] (v2) edge node[text = gray, above left] {$e_1$} (v4);
			\path [-] (v2) edge node[text = gray, right] {$e_2$} (v5);
			\path [-] (v3) edge node[text = gray, left] {$e_3$} (v6);
			\path [-] (v4) edge node[text = gray, below] {$e_4$} (v5);
		\end{scope}
	\end{tikzpicture}\end{center}\vspace{0em}
\end{example}

In the case of a not oriented normal graph, the edge $\left\{v_i, v_j\right\}$ is equal to the edge $\left\{v_j, v_i\right\}$ for all vertices $v_i, v_j \in \mathcal{V}$ and hence the graph is called ``not oriented'' or ``undirected''. In contrast to this, there also exist normal graphs for which the orientation \big(the order of the vertices in an edge\big) matters.\\

\begin{definition}[\textbf{Oriented normal graph $OG$}]\label{OG} (Page 161 in \cite{williamson2010lists}: Simple digraph) \ \\
	An oriented normal graph $OG = \left(\mathcal{V}, \mathcal{A}_G\right)$ is defined by the same finite set of vertices $\mathcal{V} = \left\{v_1, v_2, \dots v_N\right\}$ as the not oriented normal graph, however instead of a set of edges $\mathcal{E}_G$ it is composed of a set of arcs $\mathcal{A}_G \subseteq \left\{\left(v_i, v_j\right) ~ \middle| ~ v_i, v_j \in \mathcal{V}, v_i \neq v_j\right\} \subset \mathcal{V} \times \mathcal{V}$.
	
	For an arc $a_q = \left(v_i, v_j\right) \in \mathcal{A}_G$, the vertex $v_i$ is called output vertex of arc $a_q$ and respectively vertex $v_j$ is called input vertex of arc $a_q$.\\
\end{definition}

\begin{example}[\textbf{Oriented normal graph}] \ \\
	Given a set of vertices
	\begin{equation*}
		\mathcal{V} = \left\{v_1, v_2, v_3, v_4, v_5, v_6\right\}
	\end{equation*} together with a set of arcs
	\begin{equation*}
		\mathcal{A}_G = \left\{\left(v_2, v_5\right), \left(v_3, v_6\right), \left(v_4, v_2\right), \left(v_4, v_5\right), \left(v_6, v_3\right)\right\},
	\end{equation*}
	then the oriented normal graph $OG = \left(\mathcal{V}, \mathcal{A}_G\right)$ can be visualized in the following way:
	
	\begin{center}\begin{tikzpicture}
		\begin{scope}[every node/.style={circle,thick,draw}]
			\node (v1) at (0,0) {$v_1$};
			\node (v2) at (3,0) {$v_2$};
			\node (v3) at (6,0) {$v_3$};
			\node (v4) at (0,-3) {$v_4$};
			\node (v5) at (3,-3) {$v_5$};
			\node (v6) at (6,-3) {$v_6$};
		\end{scope}
		
		\begin{scope}[>={Stealth[gray]},
			every node/.style={circle},
			every edge/.style={draw=gray,very thick}]
			\path [->] (v2) edge node[text = gray, right] {$a_1$} (v5);
			\path [->] (v3) edge[bend left=20] node[text = gray, right] {$a_2$} (v6);
			\path [->] (v4) edge node[text = gray, above left] {$a_3$} (v2);
			\path [->] (v4) edge node[text = gray, below] {$a_4$} (v5);
			\path [->] (v6) edge[bend left=20] node[text = gray, left] {$a_5$} (v3);
		\end{scope}
	\end{tikzpicture}\end{center}\vspace{0em}
\end{example}

\begin{remark}[\textbf{Finiteness of vertex, edge and arc sets}]\label{Gfinite} \ \\
	For any given not oriented normal graph $NG = \left(\mathcal{V}, \mathcal{E}_G\right)$ or oriented normal graph $OG = \left(\mathcal{V}, \mathcal{A}_G\right)$, the vertex, the edge and the arc sets are assumed to be finite:
	\begin{equation}
		\left\lvert\mathcal{V}\right\rvert, \left\lvert\mathcal{E}_G\right\rvert, \left\lvert\mathcal{A}_G\right\rvert < \infty.
	\end{equation}
	Furthermore, in order to ensure clear identification of edges and arcs, it is assumed that each edge in the set of edges $\mathcal{E}_G$, and each arc in the set of arcs $\mathcal{A}_G$ is unique and thus occurs only once.\\
\end{remark}

In an oriented normal graph, changing the order of the vertices in an arc $\left(v_i, v_j\right) \in \mathcal{A}_G$ describes a different arc, which means $\left(v_i, v_j\right) \neq \left(v_j, v_i\right)$ for vertices $v_i, v_j \in \mathcal{V}$ with $v_i \neq v_j$.\\

\begin{remark}[\textbf{Orientation versus direction}]\label{Gorient} \ \\
	As proposed in Remark 3.1 in \cite{mulas2022random}, this thesis differentiates between oriented and not oriented graphs instead of directed and undirected graphs, because for every oriented arc there is only one orientation but two possible directions if interpreted as edge. Changing the orientation of a graph would hence mean exchanging the output and input vertices of each arc. The same reasoning also holds true for the later introduced hypergraphs.\\
\end{remark}

Based on this remark, it is possible to define an oriented normal graph with switched orientation for any given oriented normal graph.\\

\begin{definition}[\textbf{Oriented normal graph with switched orientation $\widetilde{OG}$}]\label{switchedOG} \ \\
	In order to switch the orientation of an oriented normal graph $OG = \left(\mathcal{V}, \mathcal{A}_G\right)$, the vertex set $\mathcal{V}$ remains the same, however for every arc $a_q \in \mathcal{A}_G$ the output and input vertex are exchanged:
	\begin{equation}
		\widetilde{OG} = \left(\mathcal{V}, \widetilde{\mathcal{A}_G}\right)
	\end{equation}
	with
	\begin{equation}
		\widetilde{\mathcal{A}_G} = \left\{\left(v_j, v_i\right) ~ \middle| ~ \left(v_i, v_j\right) \in \mathcal{A}_G\right\}.
	\end{equation}\vspace{0em}
\end{definition}

\begin{example}[\textbf{Oriented normal graph with switched orientation}]\ \\
	Given the previous oriented normal graph $OG = \left(\mathcal{V}, \mathcal{A}_G\right)$ with $\mathcal{V} = \left\{v_1, v_2, v_3, v_4, v_5, v_6\right\}$ and $\mathcal{A}_G = \left\{\left(v_2, v_5\right), \left(v_3, v_6\right), \left(v_4, v_2\right), \left(v_4, v_5\right), \left(v_6, v_3\right)\right\}$, then the corresponding oriented normal graph with switched orientation $\widetilde{OG} = \left(\mathcal{V}, \widetilde{\mathcal{A}_G}\right)$ with
	\begin{equation*}
		\widetilde{\mathcal{A}_G} = \left\{\left(v_5, v_2\right), \left(v_6, v_3\right), \left(v_2, v_4\right), \left(v_5, v_4\right), \left(v_3, v_6\right)\right\}
	\end{equation*}
	can be visualized in the following way:\\
	
	\begin{minipage}{.5\textwidth}$OG = \left(\mathcal{V}, \mathcal{A}_G\right)$\\
			
		\begin{tikzpicture}
			\begin{scope}[every node/.style={circle,thick,draw}]
				\node (v1) at (0,0) {$v_1$};
				\node (v2) at (2.5,0) {$v_2$};
				\node (v3) at (5,0) {$v_3$};
				\node (v4) at (0,-2.5) {$v_4$};
				\node (v5) at (2.5,-2.5) {$v_5$};
				\node (v6) at (5,-2.5) {$v_6$};
			\end{scope}
		
			\begin{scope}[>={Stealth[gray]},
				every node/.style={circle},
				every edge/.style={draw=gray,very thick}]
				\path [->] (v2) edge node[text = gray, right] {$a_1$} (v5);
				\path [->] (v3) edge[bend left=20] node[text = gray, right] {$a_2$} (v6);
				\path [->] (v4) edge node[text = gray, above left] {$a_3$} (v2);
				\path [->] (v4) edge node[text = gray, below] {$a_4$} (v5);
				\path [->] (v6) edge[bend left=20] node[text = gray, left] {$a_5$} (v3);
			\end{scope}
	\end{tikzpicture}\end{minipage}
	\begin{minipage}{.5\textwidth} $\widetilde{OG} = \left(\mathcal{V}, \widetilde{\mathcal{A}_G}\right)$\\
				
		\begin{tikzpicture}
			\begin{scope}[every node/.style={circle,thick,draw}]
				\node (v1) at (0,0) {$v_1$};
				\node (v2) at (2.5,0) {$v_2$};
				\node (v3) at (5,0) {$v_3$};
				\node (v4) at (0,-2.5) {$v_4$};
				\node (v5) at (2.5,-2.5) {$v_5$};
				\node (v6) at (5,-2.5) {$v_6$};
			\end{scope}
			
			\begin{scope}[>={Stealth[gray]},
				every node/.style={circle},
				every edge/.style={draw=gray,very thick}]
				\path [->] (v5) edge node[text = gray, right] {$\widetilde{a_1}$} (v2);
				\path [->] (v6) edge[bend right=20] node[text = gray, right] {$\widetilde{a_2}$} (v3);
				\path [->] (v2) edge node[text = gray, above left] {$\widetilde{a_3}$} (v4);
				\path [->] (v5) edge node[text = gray, below] {$\widetilde{a_4}$} (v4);
				\path [->] (v3) edge[bend right=20] node[text = gray, left] {$\widetilde{a_5}$} (v6);
			\end{scope}
	\end{tikzpicture}\end{minipage}\vspace{0em}
\end{example}
\subsection{Representation of normal graphs}\label{1.2}

In order to represent normal graphs not only visually, but also algebraically, matrices are typically used, for example in \cite{ijiri1965generalized}, \cite{jeurissen1981incidence} and \cite{singh2012role}.\\

\begin{definition}[\textbf{Adjacency matrix $A_{NG}$}]\label{ANG} (Page 179 in \cite{singh2012role}: Adjacency matrix) \ \\
	The adjacency matrix $A_{NG} \in \left\{0, 1\right\}^{\lvert\mathcal{V}\rvert \times \lvert\mathcal{V}\rvert} = \left\{0, 1\right\}^{N \times N}$ of a not oriented normal graph $NG = \left(\mathcal{V}, \mathcal{E}_G\right)$ is defined as:
	\begin{equation}
		\left(A_{NG}\right)_{i j} = \left\{\begin{array}{ll}
			1 & \quad \left\{v_i, v_j\right\} \in \mathcal{E}_G\\
			0 & \quad \text{otherwise}\end{array}\right. .
	\end{equation}
	
	Thus, each entry $\left(A_{NG}\right)_{i j}$ of the adjacency matrix  indicates whether two vertices $v_i, v_j \in \mathcal{V}$ are connected by an edge $\left\{v_i, v_j\right\} \in \mathcal{E}_G$ \big(or equivalently $\left\{v_j, v_i\right\} \in \mathcal{E}_G$\big). Due to this property, the adjacency matrix $A_{NG}$ is symmetric.\\
\end{definition}

The adjacency matrix can also be used to represent oriented normal graphs. In this case, it is not only necessary to capture whether two vertices are connected, but also the orientation of the connecting arc.\\

\begin{definition}[\textbf{Adjacency matrix $A_{OG}$}]\label{AOG} (Page 179 in \cite{singh2012role}: Adjacency matrix) \ \\
	The adjacency matrix $A_{OG} \in \left\{0, 1\right\}^{\lvert\mathcal{V}\rvert \times \lvert\mathcal{V}\rvert} = \left\{0, 1\right\}^{N \times N}$ of an oriented normal graph $OG = \left(\mathcal{V}, \mathcal{A}_G\right)$ is defined as:
	\begin{equation}
		\left(A_{OG}\right)_{i j} = \left\{\begin{array}{ll}
			1 & \quad \left(v_i, v_j\right) \in \mathcal{A}_G\\
			0 & \quad \text{otherwise}\end{array}\right. .
	\end{equation}
	
	In contrast to the previously introduced adjacency matrix $A_{NG}$ for not oriented normal graphs, this adjacency matrix is not necessarily symmetric, since $\left(v_i, v_j\right) \in \mathcal{A}_G$ does not imply $\left(v_j, v_i\right) \in \mathcal{A}_G$, hence $\left(A_{NG}\right)_{i j} = 1$ but $\left(A_{NG}\right)_{j i} = 0$ is possible.\\
\end{definition}

\begin{example}[\textbf{Adjacency matrices}]\ \\
	Given the previous not oriented normal graph $NG = \left(\mathcal{V}, \mathcal{E}_G\right)$ with $\mathcal{V} = \left\{v_1, v_2, v_3, v_4, v_5, v_6\right\}$ and $\mathcal{E}_G = \left\{\left\{v_2, v_4\right\}, \left\{v_2, v_5\right\}, \left\{v_3, v_6\right\}, \left\{v_4, v_5\right\}\right\}$, then the corresponding adjacency matrix is a binary $6 \times 6$ matrix:
	
	\begin{minipage}{.5\textwidth}\begin{tikzpicture}
		\begin{scope}[every node/.style={circle,thick,draw}]
			\node (v1) at (0,0) {$v_1$};
			\node (v2) at (2.5,0) {$v_2$};
			\node (v3) at (5,0) {$v_3$};
			\node (v4) at (0,-2.5) {$v_4$};
			\node (v5) at (2.5,-2.5) {$v_5$};
			\node (v6) at (5,-2.5) {$v_6$} ;
		\end{scope}
		
		\begin{scope}[>={Stealth[black]},
			every node/.style={fill=white,circle},
			every edge/.style={draw=gray,very thick}]
			\path [-] (v2) edge node[text = gray, above left] {$e_1$} (v4);
			\path [-] (v2) edge node[text = gray, right] {$e_2$} (v5);
			\path [-] (v3) edge node[text = gray, left] {$e_3$} (v6);
			\path [-] (v4) edge node[text = gray, below] {$e_4$} (v5);
		\end{scope}
	\end{tikzpicture}\end{minipage}
	\begin{minipage}{.5\textwidth} 
		$A_{NG} = \begin{pmatrix}
			0 & 0 & 0 & 0 & 0 & 0\\
			0 & 0 & 0 & 1 & 1 & 0\\
			0 & 0 & 0 & 0 & 0 & 1\\
			0 & 1 & 0 & 0 & 1 & 0\\
			0 & 1 & 0 & 1 & 0 & 0\\
			0 & 0 & 1 & 0 & 0 & 0
		\end{pmatrix}$
	\end{minipage}
	
	Similarly, for the previous oriented normal graph $OG = \left(\mathcal{V}, \mathcal{A}_G\right)$ with the same vertex set $\mathcal{V} = \left\{v_1, v_2, v_3, v_4, v_5, v_6\right\}$ and arc set $\mathcal{A}_G = \left\{\left(v_2, v_5\right), \left(v_3, v_6\right), \left(v_4, v_2\right), \left(v_4, v_5\right), \left(v_6, v_3\right)\right\}$, the adjacency matrix is given by:
	
	\begin{minipage}{.5\textwidth}\begin{tikzpicture}
		\begin{scope}[every node/.style={circle,thick,draw}]
			\node (v1) at (0,0) {$v_1$};
			\node (v2) at (2.5,0) {$v_2$};
			\node (v3) at (5,0) {$v_3$};
			\node (v4) at (0,-2.5) {$v_4$};
			\node (v5) at (2.5,-2.5) {$v_5$};
			\node (v6) at (5,-2.5) {$v_6$} ;
		\end{scope}
		
		\begin{scope}[>={Stealth[gray]},
			every node/.style={circle},
			every edge/.style={draw=gray,very thick}]
			\path [->] (v2) edge node[text = gray, right] {$a_1$} (v5);
			\path [->] (v3) edge[bend left=20] node[text = gray, right] {$a_2$} (v6);
			\path [->] (v4) edge node[text = gray, above left] {$a_3$} (v2);
			\path [->] (v4) edge node[text = gray, below] {$a_4$} (v5);
			\path [->] (v6) edge[bend left=20] node[text = gray, left] {$a_5$} (v3);
		\end{scope}
	\end{tikzpicture}\end{minipage}
	\begin{minipage}{.5\textwidth} 
		$A_{OG} = \begin{pmatrix}
			0 & 0 & 0 & 0 & 0 & 0\\
			0 & 0 & 0 & 0 & 1 & 0\\
			0 & 0 & 0 & 0 & 0 & 1\\
			0 & 1 & 0 & 0 & 1 & 0\\
			0 & 0 & 0 & 0 & 0 & 0\\
			0 & 0 & 1 & 0 & 0 & 0
		\end{pmatrix}$
	\end{minipage}\vspace{0em}
\end{example}

In addition to the adjacency matrices, which express whether two vertices are connected or not, incidence matrices are introduced. Incidence matrices encode the information whether a vertex $v_i \in \mathcal{V}$ is part of an edge $e_q \in \mathcal{E}_G$ or of an arc $a_q \in \mathcal{A}_G$ in the $i$-th row and $q$-th column.\\

\begin{definition}[\textbf{Incidence matrix $I_{NG}$}]\label{ING} (Page 291 in \cite{jeurissen1981incidence}: Incidence matrix) \ \\
	The incidence matrix $I_{NG} \in \left\{0, 1\right\}^{\lvert\mathcal{V}\rvert \times \lvert\mathcal{E}_G\rvert} = \left\{0, 1\right\}^{N \times M}$ of a not oriented normal graph $NG = \left(\mathcal{V}, \mathcal{E}_G\right)$ is defined as:
	\begin{equation}
		\left(I_{NG}\right)_{i q} = \left\{\begin{array}{ll}
			1 & \quad \left\{v_i, v_j\right\} = e_q ~ \text{for some} ~ v_j \in \mathcal{V} \backslash \left\{v_i\right\}\\
			0 & \quad \text{otherwise}\end{array}\right. .
	\end{equation}
	
	Therefore, a positive entry $\left(I_{NG}\right)_{i q}$ in the incidence matrix indicates that vertex $v_i$ is part of edge $e_q \in \mathcal{E}_G$.\\
\end{definition}

For oriented normal graphs, it is also necessary to differentiate between the output and the input vertex of each arc $a_q \in \mathcal{A}_G$.\\

\begin{definition}[\textbf{Incidence matrix $I_{OG}$}]\label{IOG} (Page 827 in \cite{ijiri1965generalized}: Incidence matrix) \ \\
	The incidence matrix $I_{OG} \in \left\{-1, 0, 1\right\}^{\lvert\mathcal{V}\rvert \times \lvert\mathcal{A}_G\rvert} = \left\{-1, 0, 1\right\}^{N \times M}$ of an oriented normal graph $OG = \left(\mathcal{V}, \mathcal{A}_G\right)$ is defined as:
	\begin{equation}
		\left(I_{OG}\right)_{i q} = \left\{\begin{array}{ll}
			1 & \quad \left(v_i, v_j\right) = a_q ~ \text{for some} ~ v_j \in \mathcal{V} \backslash \left\{v_i\right\}\\
			-1 & \quad \left(v_j, v_i\right) = a_q ~ \text{for some} ~ v_j \in \mathcal{V} \backslash \left\{v_i\right\}\\
			0 & \quad \text{otherwise}\end{array}\right. .
	\end{equation}
	
	Hence, a positive entry $\left(I_{OG}\right)_{i q}$ in the incidence matrix indicates that vertex $v_i$ is the output vertex of arc $a_q \in \mathcal{A}_G$, and analogously an entry $\left(I_{OG}\right)_{i q}$ has a negative sign if $v_i$ is the input vertex of arc $a_q$.\\
\end{definition}

\begin{example}[\textbf{Incidence matrices}]\ \\
	Given the not oriented normal graph $NG = \left(\mathcal{V}, \mathcal{E}_G\right)$ with $\mathcal{V} = \left\{v_1, v_2, v_3, v_4, v_5, v_6\right\}$ and $\mathcal{E}_G = \left\{\left\{v_2, v_4\right\}, \left\{v_2, v_5\right\}, \left\{v_3, v_6\right\}, \left\{v_4, v_5\right\}\right\}$, then the corresponding incidence matrix is a binary $6 \times 4$ matrix:
	
	\begin{minipage}{.5\textwidth}\begin{tikzpicture}
		\begin{scope}[every node/.style={circle,thick,draw}]
			\node (v1) at (0,0) {$v_1$};
			\node (v2) at (2.5,0) {$v_2$};
			\node (v3) at (5,0) {$v_3$};
			\node (v4) at (0,-2.5) {$v_4$};
			\node (v5) at (2.5,-2.5) {$v_5$};
			\node (v6) at (5,-2.5) {$v_6$} ;
		\end{scope}
		
		\begin{scope}[>={Stealth[black]},
			every node/.style={fill=white,circle},
			every edge/.style={draw=gray,very thick}]
			\path [-] (v2) edge node[text = gray, above left] {$e_1$} (v4);
			\path [-] (v2) edge node[text = gray, right] {$e_2$} (v5);
			\path [-] (v3) edge node[text = gray, left] {$e_3$} (v6);
			\path [-] (v4) edge node[text = gray, below] {$e_4$} (v5);
		\end{scope}
	\end{tikzpicture}\end{minipage}
	\begin{minipage}{.5\textwidth} 
		$I_{NG} = \begin{pmatrix}
			0 & 0 & 0 & 0\\
			1 & 1 & 0 & 0\\
			0 & 0 & 1 & 0\\
			1 & 0 & 0 & 1\\
			0 & 1 & 0 & 1\\
			0 & 0 & 1 & 0
		\end{pmatrix}$
	\end{minipage}
	
	For the oriented normal graph $OG = \left(\mathcal{V}, \mathcal{A}_G\right)$ with $\mathcal{V} = \left\{v_1, v_2, v_3, v_4, v_5, v_6\right\}$ and $\mathcal{A}_G = \left\{\left(v_2, v_5\right), \left(v_3, v_6\right), \left(v_4, v_2\right), \left(v_4, v_5\right), \left(v_6, v_3\right)\right\}$, the corresponding incidence matrix is a $6 \times 5$ matrix:
	
	\begin{minipage}{.5\textwidth}\begin{tikzpicture}
				\begin{scope}[every node/.style={circle,thick,draw}]
					\node (v1) at (0,0) {$v_1$};
					\node (v2) at (2.5,0) {$v_2$};
					\node (v3) at (5,0) {$v_3$};
					\node (v4) at (0,-2.5) {$v_4$};
					\node (v5) at (2.5,-2.5) {$v_5$};
					\node (v6) at (5,-2.5) {$v_6$} ;
				\end{scope}
				
				\begin{scope}[>={Stealth[gray]},
					every node/.style={circle},
					every edge/.style={draw=gray,very thick}]
					\path [->] (v2) edge node[text = gray, right] {$a_1$} (v5);
					\path [->] (v3) edge[bend left=20] node[text = gray, right] {$a_2$} (v6);
					\path [->] (v4) edge node[text = gray, above left] {$a_3$} (v2);
					\path [->] (v4) edge node[text = gray, below] {$a_4$} (v5);
					\path [->] (v6) edge[bend left=20] node[text = gray, left] {$a_5$} (v3);
				\end{scope}
	\end{tikzpicture}\end{minipage}
	\begin{minipage}{.5\textwidth} 
		$I_{OG} = \begin{pmatrix}
			0 & 0 & 0 & 0 & 0 \\
			1 & 0 & -1 & 0 & 0 \\
			0 & 1 & 0 & 0 & -1 \\
			0 & 0 & 1 & 1 & 0 \\
			-1 & 0 & 0 & -1 & 0 \\
			0 & -1 & 0 & 0 & 1 
		\end{pmatrix}$
	\end{minipage}\vspace{0em}
\end{example}
\clearpage
\section{Hypergraphs}\label{2}

In the previously introduced normal graphs, edges or arcs only model a pairwise interaction or connection between the data points $v_1, v_2, \dots v_N$, which makes the representation of higher-order relationships ineffective or even impossible. Such higher-order connections can for example occur in social networks, since not only two people can be friends, but there can be friend groups with more than two people as well.\\

Thus, this section fills this gap by introducing the concept of not oriented and oriented hypergraphs in subsection (\ref{2.1}) as a possible extension of normal graphs. By using hyperedges and hyperarcs, a link between multiple vertices within a graph can be described mathematically. Furthermore, based on the definitions in the normal graph case, the difference between the orientation and the direction of hypergraphs is discussed and a definition for the oriented hypergraph with switched orientation is given.\\

Subsection (\ref{2.2}) then generalizes the concept of adjacency and incidence matrices such that they are also able to describe hypergraphs. For the generalized adjacency matrix of oriented hypergraphs, two different versions are introduced: a non-unique version with a dimension matching the normal graph case and the not oriented hypergraph case and a unique version with bigger dimension. Since the generalized adjacency matrices for oriented hypergraphs are not straight forward to interpret, a detailed explanation on how to retrieve the underlying oriented hypergraph from its adjacency tensor is given as well.\\

\subsection{General definitions for hypergraphs}\label{2.1}

Hypergraphs are introduced in order to capture more complex relationships between the data points $v_1, v_2, \dots v_N$ elegantly, while simultaneously being a natural extension of normal graphs. The definitions of the not oriented and oriented hypergraphs are based on \cite{jost2021plaplaceoperators} and \cite{zhang2019introducing}.\\

\begin{definition}[\textbf{Not oriented hypergraph $NH$}]\label{NH} (Page 641 in \cite{zhang2019introducing}: Hypergraph) \ \\
	Definition (\ref{NG}) of the not oriented normal graph $NG = \left(\mathcal{V}, \mathcal{E}_G\right)$ can be extended to a not oriented hypergraph. A not oriented hypergraph is defined as a pair $NH = \left(\mathcal{V}, \mathcal{E}_H\right)$ with the same set of vertices $\mathcal{V} = \left\{v_1, v_2, \dots v_N\right\}$ as the normal graph, but a more general set of so-called hyperedges $\mathcal{E}_H$. The set of hyperedges $\mathcal{E}_H \subseteq \left\{e_q ~ \middle| ~ e_q \subseteq \mathcal{V} ~ \text{with} ~ \left\lvert e_q\right\rvert \geq 2\right\} \subset 2^\mathcal{V}$ is a subset of the power set of the vertices $2^{\mathcal{V}}$, where each hyperedge can contain a different number of vertices \big(at least two\big).\\
\end{definition}

\clearpage
\begin{example}[\textbf{Not oriented hypergraph}]\ \\
	Given a set of vertices
	\begin{equation*}
		\mathcal{V} = \left\{v_1, v_2, v_3, v_4, v_5, v_6, v_7, v_8\right\}
	\end{equation*} and a set of hyperedges
	\begin{equation*}
		\mathcal{E}_H = \left\{\left\{v_1, v_2, v_5\right\}, \left\{v_2, v_3, v_7, v_8\right\}, \left\{v_6, v_7\right\}\right\},
	\end{equation*}
	then the not oriented hypergraph $NH = \left(\mathcal{V}, \mathcal{E}_H\right)$ can be visualized in the following way:
	
	\begin{center}\begin{tikzpicture}
			\tikzstyle{vertex} = [fill,shape=circle,node distance=80pt]
			\tikzstyle{edge} = [fill,opacity=.5,fill opacity=.5,line cap=round, line join=round, line width=40pt]
			\tikzstyle{elabel} =  [fill,shape=circle,node distance=40pt]
			
			\pgfdeclarelayer{background}
			\pgfsetlayers{background,main}
			
			\begin{scope}[every node/.style={circle,thick,draw}]
				\node (v1) at (0,0) {$v_1$};
				\node (v2) at (3,0) {$v_2$};
				\node (v3) at (6,0) {$v_3$};
				\node (v4) at (9, 0) {$v_4$};
				\node (v5) at (0,-3) {$v_5$};
				\node (v6) at (3,-3) {$v_6$};
				\node (v7) at (6,-3) {$v_7$};
				\node (v8) at (9,-3) {$v_8$};
			\end{scope}
			
			\begin{pgfonlayer}{background}
				\begin{scope}[transparency group,opacity=.5]
					\draw[edge,opacity=1,color=blue] (v1.center) -- (v2.center) -- (v5.center) -- (v1.center);
					\fill[edge,opacity=1,color=blue] (v1.center) -- (v2.center) -- (v5.center) -- (v1.center);
				\end{scope}
				\begin{scope}[transparency group,opacity=.5]
					\draw[edge,opacity=1,color=green] (v2.center) -- (v3.center) -- (v8.center) -- (v7.center) -- (v2.center);
					\fill[edge,opacity=1,color=green] (v2.center) -- (v3.center) -- (v8.center) -- (v7.center) -- (v2.center);
				\end{scope}
				\begin{scope}[transparency group,opacity=.5]
					\draw[edge,opacity=1,color=lightgray] (v6.center) -- (v7.center) -- (v6.center);
					\fill[edge,opacity=1,color=lightgray] (v6.center) -- (v7.center) -- (v6.center);
				\end{scope}
			\end{pgfonlayer}
		
			\node[text=blue] (e1) at (-1,-1.5) {$e_1$};
			\node[text=green] (e2) at (8.21,-0.79) {$e_2$};
			\node[text=lightgray] (e3) at (4.5,-4) {$e_3$};
	\end{tikzpicture}\end{center}\vspace{0em}
\end{example}

As in the case of normal graphs, there also exists an oriented version of hypergraphs.\\

\begin{definition}[\textbf{Oriented hypergraph $OH$}]\label{OH} (Page 31 in \cite{jost2021plaplaceoperators}: Chemical hypergraph) \ \\
	An oriented hypergraph $OH = \left(\mathcal{V}, \mathcal{A}_H\right)$ is defined as the generalization of an oriented normal graph $OG = \left(\mathcal{V}, \mathcal{A}_G\right)$ from definition (\ref{OG}), consisting of the same set of vertices $\mathcal{V}$ and a different set of so-called hyperarcs $\mathcal{A}_H$. In contrast to the hyperedges in a not oriented hypergraph, the hyperarcs of an oriented hypergraph are oriented, which means that each of them is made up of two disjoint sets of vertices
	\begin{equation}
		a_q = \left(a_q^{out}, a_q^{in}\right)
	\end{equation}
	with $\emptyset \subset a_q^{out}, a_q^{in} \subset \mathcal{V}$, $a_q^{out} \cap a_q^{in} = \emptyset$, $a_q^{out}$ being the set of all output vertices and $a_q^{in}$ being the set of all input vertices of hyperarc $a_q$.\\
\end{definition}

\begin{example}[\textbf{Oriented hypergraph}]\ \\
	Given a set of vertices
	\begin{equation*}
		\mathcal{V} = \left\{v_1, v_2, v_3, v_4, v_5, v_6, v_7, v_8\right\}
	\end{equation*} and a set of hyperarcs
	\begin{equation*}
		\mathcal{A}_H = \left\{\left(\left\{v_1, v_2\right\}, \left\{v_5\right\}\right), \left(\left\{v_3\right\}, \left\{v_2, v_7, v_8\right\}\right), \left(\left\{v_6\right\}, \left\{v_7\right\}\right)\right\},
	\end{equation*}
	then the oriented hypergraph $OH = \left(\mathcal{V}, \mathcal{A}_H\right)$ can be visualized in the following way:
	
	\begin{center}\begin{tikzpicture}
			\tikzstyle{vertex} = [fill,shape=circle,node distance=80pt]
			\tikzstyle{edge} = [fill,opacity=.5,fill opacity=.5,line cap=round, line join=round, line width=40pt]
			\tikzstyle{elabel} =  [fill,shape=circle,node distance=40pt]
			
			\pgfdeclarelayer{background}
			\pgfsetlayers{background,main}
			
			\begin{scope}[every node/.style={circle,thick,draw}]
				\node (v1) at (0,0) {$v_1$};
				\node (v2) at (3,0) {$v_2$};
				\node (v3) at (6,0) {$v_3$};
				\node (v4) at (9, 0) {$v_4$};
				\node (v5) at (0,-3) {$v_5$};
				\node (v6) at (3,-3) {$v_6$};
				\node (v7) at (6,-3) {$v_7$};
				\node (v8) at (9,-3) {$v_8$};
			\end{scope}
			
			\begin{pgfonlayer}{background}
				\begin{scope}[transparency group,opacity=.2]
					\draw[edge,opacity=1,color=blue] (v1.center) -- (v2.center) -- (v5.center) -- (v1.center);
					\fill[edge,opacity=1,color=blue] (v1.center) -- (v2.center) -- (v5.center) -- (v1.center);
				\end{scope}
				\begin{scope}[transparency group,opacity=.2]
					\draw[edge,opacity=1,color=green] (v2.center) -- (v3.center) -- (v8.center) -- (v7.center) -- (v2.center);
					\fill[edge,opacity=1,color=green] (v2.center) -- (v3.center) -- (v8.center) -- (v7.center) -- (v2.center);
				\end{scope}
				\begin{scope}[transparency group,opacity=.2]
					\draw[edge,opacity=1,color=lightgray] (v6.center) -- (v7.center) -- (v6.center);
					\fill[edge,opacity=1,color=lightgray] (v6.center) -- (v7.center) -- (v6.center);
				\end{scope}
			\end{pgfonlayer}
			
			\node[text=blue] (v1a1) at (0.5, -0.5) {$out$};
			\node[text=blue] (v2a1) at (2.5, -0.5) {$out$};
			\node[text=blue] (v5a1) at (0.5, -2.5) {$in$};
			\node[text=green] (v3a2) at (6, -0.75) {$out$};
			\node[text=green] (v2a2) at (3.5, -0.5) {$in$};
			\node[text=green] (v7a2) at (6, -2.25) {$in$};
			\node[text=green] (v8a2) at (8.5, -2.5) {$in$};
			\node[text=lightgray] (v6a3) at (3.75, -3) {$out$};
			\node[text=lightgray] (v7a3) at (5.25, -3) {$in$};
			\node[text=blue] (a1) at (-1,-1.5) {$a_1$};
			\node[text=green] (a2) at (8.21,-0.79) {$a_2$};
			\node[text=lightgray] (a3) at (4.5,-4) {$a_3$};
	\end{tikzpicture}\end{center}
	
	Alternatively, hyperarcs can also be visualized similarly to arcs in normal graphs:
	
	\begin{center}\begin{tikzpicture}
			\tikzstyle{vertex} = [fill,shape=circle,node distance=80pt]
			\tikzstyle{edge} = [fill,opacity=.5,fill opacity=.5,line cap=round, line join=round, line width=40pt]
			\tikzstyle{elabel} =  [fill,shape=circle,node distance=40pt]
			
			\pgfdeclarelayer{background}
			\pgfsetlayers{background,main}
			
			\begin{scope}[every node/.style={circle,thick,draw}]
				\node (v1) at (0,0) {$v_1$};
				\node (v2) at (3,0) {$v_2$};
				\node (v3) at (6,0) {$v_3$};
				\node (v4) at (9, 0) {$v_4$};
				\node (v5) at (0,-3) {$v_5$};
				\node (v6) at (3,-3) {$v_6$};
				\node (v7) at (6,-3) {$v_7$};
				\node (v8) at (9,-3) {$v_8$};
			\end{scope}
			
			\begin{scope}[>={Stealth[blue]},
				every node/.style={circle},
				every edge/.style={draw=blue,very thick}]
				\path [->] (0.75, -1.5) edge node[] {} (v5);
				\path [-] (v2) edge node[] {} (0.75, -1.5);
				\path [-] (v1) edge node[] {} (0.75, -1.5);
			\end{scope}
			\begin{scope}[>={Stealth[green]},
				every node/.style={circle},
				every edge/.style={draw=green,very thick}]
				\path [->] (v3) edge node[] {} (v7);
				\path [->] (6, -1.5) edge node[] {} (v2);
				\path [->] (6, -1.5) edge node[] {} (v8);
			\end{scope}
			\begin{scope}[>={Stealth[lightgray]},
				every node/.style={circle},
				every edge/.style={draw=lightgray,very thick}]
				\path [->] (4.5,-3) edge node[] {} (v7);
				\path [-] (v6) edge node[] {} (4.5,-3);
			\end{scope}
			\node[text=blue, left] (a1) at (0.75, -1.5) {$a_1$};
			\node[text=green, above right] (a2) at (6, -1.5) {$a_2$};
			\node[text=lightgray, below] (a3) at (4.5,-3) {$a_3$};
	\end{tikzpicture}\end{center}\vspace{0em}
\end{example}

\begin{remark}[\textbf{Finiteness of vertex, hyperedge and hyperarc sets}]\label{Hfinite} \ \\
	For any given not oriented hypergraph $NH = \left(\mathcal{V}, \mathcal{E}_H\right)$ or oriented hypergraph $OH = \left(\mathcal{V}, \mathcal{A}_H\right)$, the vertex, the hyperedge and the hyperarc sets are assumed to be finite:
	\begin{equation}
		\left\lvert\mathcal{V}\right\rvert, \left\lvert\mathcal{E}_H\right\rvert, \left\lvert\mathcal{A}_H\right\rvert < \infty.
	\end{equation}
	Moreover, in order to ensure clear identification of hyperedges and hyperarcs, it is assumed, that each hyperedge and each hyperarc occurs in the set of hyperedges $\mathcal{E}_H$ or in the set of hyperarcs $\mathcal{A}_H$ at most once and can thus be clearly identified.\\
\end{remark}

The orientation of oriented hypergraphs can also be switched just like in the case of oriented normal graphs.\\

\begin{definition}[\textbf{Oriented hypergraph with switched orientation $\widetilde{OH}$}]\label{switchedOH} \ \\
	In order to switch the orientation of an oriented hypergraph $OH = \left(\mathcal{V}, \mathcal{A}_H\right)$, the vertex set $\mathcal{V}$ remains the same, however for every hyperarc $a_q \in \mathcal{A}_G$ the output and input vertex sets $a_q^{out}$ and $a_q^{in}$ are exchanged, which results in
	\begin{equation}
		\widetilde{OH} = \left(\mathcal{V}, \widetilde{\mathcal{A}_H}\right)
	\end{equation}
	with
	\begin{equation}
		\widetilde{\mathcal{A}_H} = \left\{\widetilde{a_q} = \left(\widetilde{a_q^{out}}, \widetilde{a_q^{in}}\right) ~ \middle| ~ a_q \in \mathcal{A}_H ~ \text{and} ~ \widetilde{a_q^{out}} = a_q^{in}, \widetilde{a_q^{in}} = a_q^{out}\right\}.
	\end{equation}\vspace{0em}
\end{definition}

\begin{example}[\textbf{Oriented hypergraph with switched orientation}]\ \\
	Given the oriented hypergraph $OH = \left(\mathcal{V}, \mathcal{A}_H\right)$ with $\mathcal{V} = \left\{v_1, v_2, v_3, v_4, v_5, v_6, v_7, v_8\right\}$ and $\mathcal{A}_H = \left\{\left(\left\{v_1, v_2\right\}, \left\{v_5\right\}\right), \left(\left\{v_3\right\}, \left\{v_2, v_7, v_8\right\}\right), \left(\left\{v_6\right\}, \left\{v_7\right\}\right)\right\}$, then the corresponding oriented hypergraph with switched orientation $\widetilde{OH} = \left(\mathcal{V}, \widetilde{\mathcal{A}_H}\right)$ with
	\begin{equation*}
		\widetilde{\mathcal{A}_H} = \left\{\left(\left\{v_5\right\}, \left\{v_1, v_2\right\}\right), \left(\left\{v_2, v_7, v_8\right\}, \left\{v_3\right\}\right), \left(\left\{v_7\right\}, \left\{v_6\right\}\right)\right\}
	\end{equation*}
	can be visualized in the following two ways:\\
	
	\begin{minipage}{.5\textwidth}$OH = \left(\mathcal{V}, \mathcal{A}_H\right)$
		
		\begin{flushleft}\begin{tikzpicture}
				\tikzstyle{vertex} = [fill,shape=circle,node distance=80pt]
				\tikzstyle{edge} = [fill,opacity=.5,fill opacity=.5,line cap=round, line join=round, line width=30pt]
				\tikzstyle{elabel} =  [fill,shape=circle,node distance=40pt]
				
				\pgfdeclarelayer{background}
				\pgfsetlayers{background,main}
				
				\begin{scope}[every node/.style={circle,thick,draw}]
					\node (v1) at (0,0) {$v_1$};
					\node (v2) at (2,0) {$v_2$};
					\node (v3) at (4,0) {$v_3$};
					\node (v4) at (6, 0) {$v_4$};
					\node (v5) at (0,-2) {$v_5$};
					\node (v6) at (2,-2) {$v_6$};
					\node (v7) at (4,-2) {$v_7$};
					\node (v8) at (6,-2) {$v_8$};
				\end{scope}
				
				\begin{pgfonlayer}{background}
					\begin{scope}[transparency group,opacity=.2]
						\draw[edge,opacity=1,color=blue] (v1.center) -- (v2.center) -- (v5.center) -- (v1.center);
						\fill[edge,opacity=1,color=blue] (v1.center) -- (v2.center) -- (v5.center) -- (v1.center);
					\end{scope}
					\begin{scope}[transparency group,opacity=.2]
						\draw[edge,opacity=1,color=green] (v2.center) -- (v3.center) -- (v8.center) -- (v7.center) -- (v2.center);
						\fill[edge,opacity=1,color=green] (v2.center) -- (v3.center) -- (v8.center) -- (v7.center) -- (v2.center);
					\end{scope}
					\begin{scope}[transparency group,opacity=.2]
						\draw[edge,opacity=1,color=lightgray] (v6.center) -- (v7.center) -- (v6.center);
						\fill[edge,opacity=1,color=lightgray] (v6.center) -- (v7.center) -- (v6.center);
					\end{scope}
				\end{pgfonlayer}
				
				\node[text=blue] (v1a1) at (0.5, -0.5) {\footnotesize{$out$}};
				\node[text=blue] (v2a1) at (1.5, -0.5) {\footnotesize{$out$}};
				\node[text=blue] (v5a1) at (0.5, -1.5) {\footnotesize{$in$}};
				\node[text=green] (v3a2) at (4, -0.7) {\footnotesize{$out$}};
				\node[text=green] (v2a2) at (2.5, -0.5) {\footnotesize{$in$}};
				\node[text=green] (v7a2) at (4, -1.3) {\footnotesize{$in$}};
				\node[text=green] (v8a2) at (5.5, -1.5) {\footnotesize{$in$}};
				\node[text=lightgray] (v6a3) at (2.75, -2) {\footnotesize{$out$}};
				\node[text=lightgray] (v7a3) at (3.35, -2) {\footnotesize{$in$}};
				\node[text=blue] (a1) at (1,0.75) {$a_1$};
				\node[text=green] (a2) at (6,-1) {$a_2$};
				\node[text=lightgray] (a3) at (3,-2.75) {$a_3$};
	\end{tikzpicture}\end{flushleft}\end{minipage}
	\begin{minipage}{.5\textwidth} $~ ~ ~ ~ \widetilde{OH} = \left(\mathcal{V}, \widetilde{\mathcal{A}_H}\right)$
		
		\begin{flushright}\begin{tikzpicture}
				\tikzstyle{vertex} = [fill,shape=circle,node distance=80pt]
				\tikzstyle{edge} = [fill,opacity=.5,fill opacity=.5,line cap=round, line join=round, line width=30pt]
				\tikzstyle{elabel} =  [fill,shape=circle,node distance=40pt]
				
				\pgfdeclarelayer{background}
				\pgfsetlayers{background,main}
				
				\begin{scope}[every node/.style={circle,thick,draw}]
					\node (v1) at (0,0) {$v_1$};
					\node (v2) at (2,0) {$v_2$};
					\node (v3) at (4,0) {$v_3$};
					\node (v4) at (6, 0) {$v_4$};
					\node (v5) at (0,-2) {$v_5$};
					\node (v6) at (2,-2) {$v_6$};
					\node (v7) at (4,-2) {$v_7$};
					\node (v8) at (6,-2) {$v_8$};
				\end{scope}
				
				\begin{pgfonlayer}{background}
					\begin{scope}[transparency group,opacity=.2]
						\draw[edge,opacity=1,color=blue] (v1.center) -- (v2.center) -- (v5.center) -- (v1.center);
						\fill[edge,opacity=1,color=blue] (v1.center) -- (v2.center) -- (v5.center) -- (v1.center);
					\end{scope}
					\begin{scope}[transparency group,opacity=.2]
						\draw[edge,opacity=1,color=green] (v2.center) -- (v3.center) -- (v8.center) -- (v7.center) -- (v2.center);
						\fill[edge,opacity=1,color=green] (v2.center) -- (v3.center) -- (v8.center) -- (v7.center) -- (v2.center);
					\end{scope}
					\begin{scope}[transparency group,opacity=.2]
						\draw[edge,opacity=1,color=lightgray] (v6.center) -- (v7.center) -- (v6.center);
						\fill[edge,opacity=1,color=lightgray] (v6.center) -- (v7.center) -- (v6.center);
					\end{scope}
				\end{pgfonlayer}
				
				\node[text=blue] (v1a1) at (0.5, -0.5) {\footnotesize{$in$}};
				\node[text=blue] (v2a1) at (1.5, -0.5) {\footnotesize{$in$}};
				\node[text=blue] (v5a1) at (0.5, -1.5) {\footnotesize{$out$}};
				\node[text=green] (v3a2) at (4, -0.7) {\footnotesize{$in$}};
				\node[text=green] (v2a2) at (2.5, -0.5) {\footnotesize{$out$}};
				\node[text=green] (v7a2) at (4, -1.3) {\footnotesize{$out$}};
				\node[text=green] (v8a2) at (5.5, -1.5) {\footnotesize{$out$}};
				\node[text=lightgray] (v6a3) at (2.65, -2) {\footnotesize{$in$}};
				\node[text=lightgray] (v7a3) at (3.25, -2) {\footnotesize{$out$}};
				\node[text=blue] (a1) at (1,0.75) {$\widetilde{a_1}$};
				\node[text=green] (a2) at (6,-1) {$\widetilde{a_2}$};
				\node[text=lightgray] (a3) at (3,-2.75) {$\widetilde{a_3}$};
	\end{tikzpicture}\end{flushright}\end{minipage}

	\begin{minipage}{.5\textwidth}$OH = \left(\mathcal{V}, \mathcal{A}_H\right)$
		
		\begin{flushleft}\begin{tikzpicture}
				\tikzstyle{vertex} = [fill,shape=circle,node distance=80pt]
				\tikzstyle{edge} = [fill,opacity=.5,fill opacity=.5,line cap=round, line join=round, line width=40pt]
				\tikzstyle{elabel} =  [fill,shape=circle,node distance=40pt]
				
				\pgfdeclarelayer{background}
				\pgfsetlayers{background,main}
				
				\begin{scope}[every node/.style={circle,thick,draw}]
					\node (v1) at (0,0) {$v_1$};
					\node (v2) at (2,0) {$v_2$};
					\node (v3) at (4,0) {$v_3$};
					\node (v4) at (6, 0) {$v_4$};
					\node (v5) at (0,-2) {$v_5$};
					\node (v6) at (2,-2) {$v_6$};
					\node (v7) at (4,-2) {$v_7$};
					\node (v8) at (6,-2) {$v_8$};
				\end{scope}
				
				\begin{scope}[>={Stealth[blue]},
					every node/.style={circle},
					every edge/.style={draw=blue,very thick}]
					\path [->] (0.5, -1) edge node[] {} (v5);
					\path [-] (v2) edge node[] {} (0.5, -1);
					\path [-] (v1) edge node[] {} (0.5, -1);
				\end{scope}
				\begin{scope}[>={Stealth[green]},
					every node/.style={circle},
					every edge/.style={draw=green,very thick}]
					\path [->] (v3) edge node[] {} (v7);
					\path [->] (4, -1) edge node[] {} (v2);
					\path [->] (4, -1) edge node[] {} (v8);
				\end{scope}
				\begin{scope}[>={Stealth[lightgray]},
					every node/.style={circle},
					every edge/.style={draw=lightgray,very thick}]
					\path [->] (3,-2) edge node[] {} (v7);
					\path [-] (v6) edge node[] {} (3,-2);
				\end{scope}
				\node[text=blue, left] (a1) at (0.5, -1) {$a_1$};
				\node[text=green, above right] (a2) at (4, -1) {$a_2$};
				\node[text=lightgray, below] (a3) at (3,-2) {$a_3$};
	\end{tikzpicture}\end{flushleft}\end{minipage}
	\begin{minipage}{.5\textwidth} $~ ~ ~ ~ \widetilde{OH} = \left(\mathcal{V}, \widetilde{\mathcal{A}_H}\right)$
		
		\begin{flushright}\begin{tikzpicture}
								\tikzstyle{vertex} = [fill,shape=circle,node distance=80pt]
				\tikzstyle{edge} = [fill,opacity=.5,fill opacity=.5,line cap=round, line join=round, line width=40pt]
				\tikzstyle{elabel} =  [fill,shape=circle,node distance=40pt]
				
				\pgfdeclarelayer{background}
				\pgfsetlayers{background,main}
				
				\begin{scope}[every node/.style={circle,thick,draw}]
					\node (v1) at (0,0) {$v_1$};
					\node (v2) at (2,0) {$v_2$};
					\node (v3) at (4,0) {$v_3$};
					\node (v4) at (6, 0) {$v_4$};
					\node (v5) at (0,-2) {$v_5$};
					\node (v6) at (2,-2) {$v_6$};
					\node (v7) at (4,-2) {$v_7$};
					\node (v8) at (6,-2) {$v_8$};
				\end{scope}
				
				\begin{scope}[>={Stealth[blue]},
					every node/.style={circle},
					every edge/.style={draw=blue,very thick}]
					\path [-] (0.5, -1) edge node[] {} (v5);
					\path [<-] (v2) edge node[] {} (0.5, -1);
					\path [<-] (v1) edge node[] {} (0.5, -1);
				\end{scope}
				\begin{scope}[>={Stealth[green]},
					every node/.style={circle},
					every edge/.style={draw=green,very thick}]
					\path [<-] (v3) edge node[] {} (v7);
					\path [-] (4, -1) edge node[] {} (v2);
					\path [-] (4, -1) edge node[] {} (v8);
				\end{scope}
				\begin{scope}[>={Stealth[lightgray]},
					every node/.style={circle},
					every edge/.style={draw=lightgray,very thick}]
					\path [-] (3,-2) edge node[] {} (v7);
					\path [<-] (v6) edge node[] {} (3,-2);
				\end{scope}
				\node[text=blue, left] (a1) at (0.5, -1) {$\widetilde{a_1}$};
				\node[text=green, above right] (a2) at (4, -1) {$\widetilde{a_2}$};
				\node[text=lightgray, below] (a3) at (3,-2) {$\widetilde{a_3}$};
	\end{tikzpicture}\end{flushright}\end{minipage}\vspace{0em}
\end{example}
\subsection{Representation of hypergraphs}\label{2.2}

This subsection generalizes the representation of not oriented hypergraphs, as introduced in \cite{zhang2019introducing}, to also represent oriented hypergraphs.\\

The following definitions already use the maximum hyperedge cardinality $max_e$ of a not oriented hypergraph and the maximum hyperarc cardinality $max_a$ of an oriented hypergraph, which are formally introduced in definition (\ref{maxea}). The maximum hyperedge cardinality $max_e$ is the maximum number of vertices contained in one hyperedge of a not oriented hypergraph. Analogously, the maximum hyperarc cardinality $max_a$ is the maximum number of vertices contained in one hyperarc \big(in the union of the output and input vertices of one hyperarc\big) of an oriented hypergraph.\\

For describing not oriented hypergraphs algebraically, a generalized version of the incidence matrix $I_{NG}$ from definition (\ref{ING}) can be used.\\

\begin{definition}[\textbf{Incidence matrix $I_{NH}$}]\label{INH} (Page 642 in \cite{zhang2019introducing})\ \\
	The incidence matrix $I_{NH} \in \left\{0, 1\right\}^{\lvert\mathcal{V}\rvert \times \lvert\mathcal{E}_H\rvert} = \left\{0, 1\right\}^{N \times M}$ of a not oriented hypergraph $NH = \left(\mathcal{V}, \mathcal{E}_H\right)$ is defined as:
	\begin{equation}
		\left(I_{NH}\right)_{i q} = \left\{\begin{array}{ll}
			1 & \quad v_i \in e_q\\
			0 & \quad \text{otherwise}\end{array}\right. .
	\end{equation}

	In contrast to the incidence matrix $I_{NG}$ of a not oriented normal graph, the number of nonzero entries in each column of the incidence matrix $I_{NH}$ \big(describing a hyperedge $e_q \in \mathcal{E}_H$\big) is not limited by two, but by $\lvert\mathcal{V}\rvert$, or more specifically the maximum hyperedge cardinality $max_e$ of the hypergraph.\\
\end{definition}

An oriented hypergraph $OH = \left(\mathcal{V}, \mathcal{A}_H\right)$ can also be described by an incidence matrix $I_{OH}$, similar to the incidence matrix $I_{OG}$ of a not oriented normal graph from definition (\ref{IOG}).\\

\begin{definition}[\textbf{Incidence matrix $I_{OH}$}] \ \\
	The incidence matrix $I_{OH} \in \left\{-1, 0, 1\right\}^{\lvert\mathcal{V}\rvert \times \lvert\mathcal{A}_H\rvert} = \left\{-1, 0, 1\right\}^{N \times M}$ of an oriented hypergraph $OH = \left(\mathcal{V}, \mathcal{A}_H\right)$ is defined as:
	\begin{equation}
		\left(I_{OH}\right)_{i q} = \left\{\begin{array}{ll}
			1 & \quad v_i \in a_q^{out}\\
			-1 & \quad v_i \in a_q^{in}\\
			0 & \quad \text{otherwise}\end{array}\right. .
	\end{equation}
	Due to set of all output vertices and the set of all input vertices of hyperarc $a_q \in \mathcal{A}_H$ being disjoint \big($a_q^{out} \cap a_q^{in} = \emptyset$\big) by definition, the incidence matrix $I_{OH}$ is well-defined.\\
\end{definition}

\begin{example}[\textbf{Incidence matrices}]\ \\
	Given the previous not oriented hypergraph $NH = \left(\mathcal{V}, \mathcal{E}_H\right)$ with the vertex set $\mathcal{V} = \left\{v_1, v_2, v_3, v_4, v_5, v_6, v_7, v_8\right\}$ and hyperedges $\mathcal{E}_H = \left\{\left\{v_1, v_2, v_5\right\}, \left\{v_2, v_3, v_7, v_8\right\}, \left\{v_6, v_7\right\}\right\}$, then the corresponding incidence matrix is a binary $8 \times 3$ matrix:
	
	\begin{minipage}{.5\textwidth}\begin{tikzpicture}
				\tikzstyle{vertex} = [fill,shape=circle,node distance=80pt]
				\tikzstyle{edge} = [fill,opacity=.5,fill opacity=.5,line cap=round, line join=round, line width=30pt]
				\tikzstyle{elabel} =  [fill,shape=circle,node distance=40pt]
				
				\pgfdeclarelayer{background}
				\pgfsetlayers{background,main}
				
				\begin{scope}[every node/.style={circle,thick,draw}]
					\node (v1) at (0,0) {$v_1$};
					\node (v2) at (2,0) {$v_2$};
					\node (v3) at (4,0) {$v_3$};
					\node (v4) at (6, 0) {$v_4$};
					\node (v5) at (0,-2) {$v_5$};
					\node (v6) at (2,-2) {$v_6$};
					\node (v7) at (4,-2) {$v_7$};
					\node (v8) at (6,-2) {$v_8$};
				\end{scope}
				
				\begin{pgfonlayer}{background}
					\begin{scope}[transparency group,opacity=.5]
						\draw[edge,opacity=1,color=blue] (v1.center) -- (v2.center) -- (v5.center) -- (v1.center);
						\fill[edge,opacity=1,color=blue] (v1.center) -- (v2.center) -- (v5.center) -- (v1.center);
					\end{scope}
					\begin{scope}[transparency group,opacity=.5]
						\draw[edge,opacity=1,color=green] (v2.center) -- (v3.center) -- (v8.center) -- (v7.center) -- (v2.center);
						\fill[edge,opacity=1,color=green] (v2.center) -- (v3.center) -- (v8.center) -- (v7.center) -- (v2.center);
					\end{scope}
					\begin{scope}[transparency group,opacity=.5]
						\draw[edge,opacity=1,color=lightgray] (v6.center) -- (v7.center) -- (v6.center);
						\fill[edge,opacity=1,color=lightgray] (v6.center) -- (v7.center) -- (v6.center);
					\end{scope}
				\end{pgfonlayer}
				
				\node[text=blue] (a1) at (1,0.75) {$e_1$};
				\node[text=green] (a2) at (6,-1) {$e_2$};
				\node[text=lightgray] (a3) at (3,-2.75) {$e_3$};
	\end{tikzpicture}\end{minipage}
	\begin{minipage}{.5\textwidth} 
		$I_{NH} = \begin{pmatrix}
			1 & 0 & 0\\
			1 & 1 & 0\\
			0 & 1 & 0\\
			0 & 0 & 0\\
			1 & 0 & 0\\
			0 & 0 & 1\\
			0 & 1 & 1\\
			0 & 1 & 0\\
		\end{pmatrix}$
	\end{minipage}
	
	\clearpage
	Similarly, for the previous oriented hypergraph $OH = \left(\mathcal{V}, \mathcal{A}_H\right)$ with vertex set $\mathcal{V} = \left\{v_1, v_2, v_3, v_4, v_5, v_6, v_7, v_8\right\}$ and the set of hyperarcs $\mathcal{A}_H =$\\ $\left\{\left(\left\{v_1, v_2\right\}, \left\{v_5\right\}\right), \left(\left\{v_3\right\}, \left\{v_2, v_7, v_8\right\}\right), \left(\left\{v_6\right\}, \left\{v_7\right\}\right)\right\}$, the incidence matrix is given by:
	
	\begin{minipage}{.5\textwidth}\begin{tikzpicture}
				\tikzstyle{vertex} = [fill,shape=circle,node distance=80pt]
				\tikzstyle{edge} = [fill,opacity=.5,fill opacity=.5,line cap=round, line join=round, line width=30pt]
				\tikzstyle{elabel} =  [fill,shape=circle,node distance=40pt]
				
				\pgfdeclarelayer{background}
				\pgfsetlayers{background,main}
				
				\begin{scope}[every node/.style={circle,thick,draw}]
					\node (v1) at (0,0) {$v_1$};
					\node (v2) at (2,0) {$v_2$};
					\node (v3) at (4,0) {$v_3$};
					\node (v4) at (6, 0) {$v_4$};
					\node (v5) at (0,-2) {$v_5$};
					\node (v6) at (2,-2) {$v_6$};
					\node (v7) at (4,-2) {$v_7$};
					\node (v8) at (6,-2) {$v_8$};
				\end{scope}
				
				\begin{pgfonlayer}{background}
					\begin{scope}[transparency group,opacity=.2]
						\draw[edge,opacity=1,color=blue] (v1.center) -- (v2.center) -- (v5.center) -- (v1.center);
						\fill[edge,opacity=1,color=blue] (v1.center) -- (v2.center) -- (v5.center) -- (v1.center);
					\end{scope}
					\begin{scope}[transparency group,opacity=.2]
						\draw[edge,opacity=1,color=green] (v2.center) -- (v3.center) -- (v8.center) -- (v7.center) -- (v2.center);
						\fill[edge,opacity=1,color=green] (v2.center) -- (v3.center) -- (v8.center) -- (v7.center) -- (v2.center);
					\end{scope}
					\begin{scope}[transparency group,opacity=.2]
						\draw[edge,opacity=1,color=lightgray] (v6.center) -- (v7.center) -- (v6.center);
						\fill[edge,opacity=1,color=lightgray] (v6.center) -- (v7.center) -- (v6.center);
					\end{scope}
				\end{pgfonlayer}
				
				\node[text=blue] (v1a1) at (0.5, -0.5) {\footnotesize{$out$}};
				\node[text=blue] (v2a1) at (1.5, -0.5) {\footnotesize{$out$}};
				\node[text=blue] (v5a1) at (0.5, -1.5) {\footnotesize{$in$}};
				\node[text=green] (v3a2) at (4, -0.7) {\footnotesize{$out$}};
				\node[text=green] (v2a2) at (2.5, -0.5) {\footnotesize{$in$}};
				\node[text=green] (v7a2) at (4, -1.3) {\footnotesize{$in$}};
				\node[text=green] (v8a2) at (5.5, -1.5) {\footnotesize{$in$}};
				\node[text=lightgray] (v6a3) at (2.75, -2) {\footnotesize{$out$}};
				\node[text=lightgray] (v7a3) at (3.35, -2) {\footnotesize{$in$}};
				\node[text=blue] (a1) at (1,0.75) {$a_1$};
				\node[text=green] (a2) at (6,-1) {$a_2$};
				\node[text=lightgray] (a3) at (3,-2.75) {$a_3$};
	\end{tikzpicture}\end{minipage}
	\begin{minipage}{.5\textwidth} 
		$I_{OH} = \begin{pmatrix}
			1 & 0 & 0\\
			1 & -1 & 0\\
			0 & 1 & 0\\
			0 & 0 & 0\\
			-1 & 0 & 0\\
			0 & 0 & 1\\
			0 & -1 & -1\\
			0 & -1 & 0\\
		\end{pmatrix}$
	\end{minipage}\vspace{0em}
\end{example}

The description of hypergraphs by the incidence matrices $I_{NH}$ or $I_{OH}$ is easily understandable, however at the same time it limits the possible mathematical operations on the hypergraph. Thus, in addition to the matrix representation, the adjacency tensor representation for not oriented hypergraphs was introduced in \cite{zhang2019introducing}. Similar to the adjacency matrix $A_{NG}$ for normal graphs from definition (\ref{ANG}), this representation indicates whether a subset of vertices is connected by a hyperedge.\\

\begin{definition}[\textbf{Adjacency tensor $T_{NH}$}]\label{TNH} (Page 643 in \cite{zhang2019introducing}: Adjacency tensor) \ \\
	A not oriented hypergraph $NH = \left(\mathcal{V}, \mathcal{E}_H\right)$ can be described by a ${max_e}^{\text{th}}$-order $N$-dimensional tensor:
	\begin{equation*}
		T_{NH} \in \mathbb{R}^{\tiny{\underbrace{N \times N \times \dots \times N}_{max_e-\text{times}}}} .
	\end{equation*}
	Every hyperedge $e_q = \left\{v_{i_1}, v_{i_2}, \dots v_{i_n}\right\} \in \mathcal{E}_H$ consisting of $2 \leq n \leq max_e$ vertices is represented by all entries $\left(T_{NH}\right)_{m_1, m_2, \dots m_{max_e}}$ for which it holds true that
	\begin{equation*}
		\left\{i_1, i_2, \dots i_n\right\} = \left\{m_1, m_2, \dots m_{max_e}\right\},
	\end{equation*}
	which means that $n$ indices of $\left(T_{NH}\right)_{m_1, m_2, \dots m_{max_e}}$ are exactly the same as $\left\{i_1, i_2, \dots i_n\right\}$ and the remaining $max_e - n$ indices are arbitrarily picked from $\left\{i_1, i_2, \dots i_n\right\}$. All such entries $\left(T_{NH}\right)_{m_1, m_2, \dots m_{max_e}}$ fulfilling these two conditions, and therefore representing a hyperedge $e_q \in \mathcal{E}_H$, are defined as
	\begin{equation} \label{TNHdenom}
		\left(T_{NH}\right)_{m_1, m_2, \dots m_{max_e}} = \frac{n}{\sum_{\left\{l_1, l_2, \dots l_n \in \mathbb{N} \backslash \left\{0\right\}: ~ \sum_{k = 1}^{n} l_k = max_e\right\}} \frac{max_e!}{l_1! ~ l_2! ~ \cdots ~ l_n!}} > 0,
	\end{equation}
	where the sum in the denominator iterates over all possible combinations of $n$ positive integers summing up to $max_e$. All entries of $T_{NH}$, which do not correspond to a hyperedge $e_q \in \mathcal{E}_H$, are set to $0$.\\
\end{definition}

The adjacency tensor $T_{NH}$ for not oriented hypergraphs can be adapted in order to also describe oriented hypergraphs and to constitute a feasible generalization of the adjacency matrix $A_{OG}$ for oriented normal graphs from definition (\ref{AOG}).\\

\begin{definition}[\textbf{Adjacency tensor $T_{OH}$}]\label{TOH} \ \\
	An oriented hypergraph $OH = \left(\mathcal{V}, \mathcal{A}_H\right)$ can be described by a ${max_a}^{\text{th}}$-order $N$-dimensional tensor:
	\begin{equation*}
		T_{OH} \in \mathbb{R}^{\tiny{\underbrace{N \times N \times \dots \times N}_{max_a-\text{times}}}} .
	\end{equation*}
	Every hyperarc $a_q = \left(a_q^{out}, a_q^{in}\right) = \left(\left\{v_{i_1}, v_{i_2}, \dots v_{i_r}\right\}, \left\{v_{i_{r + 1}}, v_{i_{r + 2}}, \dots v_{i_n}\right\}\right) \in \mathcal{A}_H$ consisting of $2 \leq n \leq max_a$ vertices with $1 \leq r \leq n - 1$ is represented by all entries of the tensor $\left(T_{OH}\right)_{m_1, m_2, \dots m_s, m_{s + 1}, \dots m_{max_a}}$  with $r \leq s \leq r + max_a - n$, where
	\begin{equation*}
		\left\{i_1, i_2, \dots i_r\right\} = \left\{m_1, m_2, \dots m_s\right\},
	\end{equation*}
	\begin{equation*}
		\left\{i_{r + 1}, i_{r + 2}, \dots i_n\right\} = \left\{m_{s + 1}, m_{s + 2}, \dots m_{max_a}\right\},
	\end{equation*}
	which means that the first $s$ indices of $\left(T_{OH}\right)_{m_1, m_2, \dots m_{max_a}}$ only consist of $\left\{i_1, i_2, \dots i_r\right\}$, with each of the vertex indices $\left\{i_1, i_2, \dots i_r\right\}$ being represented in $m_1, m_2, \dots m_s$ at least once. The remaining $max_a - s$ indices $m_{s + 1}, m_{s + 2}, \dots m_{max_a}$ similarly only consist of the vertex indices $\left\{i_{r + 1}, i_{r + 2}, \dots i_n\right\}$. For all such indices $m_1, m_2, \dots m_{max_a}$ representing a hyperarc $a_q \in \mathcal{A}_H$, define:
	\begin{equation} \label{TOHdenom}
		\left(T_{OH}\right)_{m_1, m_2, \dots m_{max_a}}^{a_q} = \frac{n}{\sum_{\left\{l_1, l_2, \dots l_n \in \mathbb{N} \backslash \left\{0\right\}: ~ \sum_{k = 1}^{n} l_k = max_a\right\}} \frac{max_a!}{l_1! ~ l_2! ~ \cdots ~ l_n!}} > 0.
	\end{equation}
	Since indices $m_1, m_2, \dots m_{max_a}$ can represent more than one hyperarc at the same time, the entry of the adjacency tensor $T_{OH}$ is set to:
	\begin{equation}
		\left(T_{OH}\right)_{m_1, m_2, \dots m_{max_a}} = \sum_{a_q \in \mathcal{A}_H ~ \text{represented by} ~ m_1, m_2, \dots m_{max_a}} \left(T_{OH}\right)_{m_1, m_2, \dots m_{max_a}}^{a_q}.
	\end{equation}
	All entries of $T_{OH}$, which do not correspond to any hyperarc $a_q \in \mathcal{A}_H$, are set to $0$.\\
\end{definition}

\begin{remark}[\textbf{Simplification of the adjacency tensors}]\label{Tsimple} \ \\
	Since the denominators in the tensor entries  $\left(T_{NH}\right)_{m_1, m_2, \dots m_{max_e}}$ (\ref{TNHdenom}) and $\left(T_{OH}\right)_{m_1, m_2, \dots m_{max_a}}^{a_q}$ (\ref{TOHdenom}) do not play a role in identifying the hyperedges or hyperarcs encoded by the corresponding indices and are costly to compute, they can be simplified to
	\begin{equation}
		\left(T_{NH}\right)_{m_1, m_2, \dots m_{max_e}} = \frac{n}{2} > 0
	\end{equation}
	for any hyperedge $e_q = \left\{v_{i_1}, v_{i_2}, \dots v_{i_n}\right\} \in \mathcal{E}_H$, and
	\begin{equation}
		\left(T_{OH}\right)_{m_1, m_2, \dots m_{max_a}}^{a_q} = \frac{n}{2} > 0
	\end{equation}
	for any hyperarc $a_q = \left(\left\{v_{i_1}, v_{i_2}, \dots v_{i_r}\right\}, \left\{v_{i_{r + 1}}, v_{i_{r + 2}}, \dots v_{i_n}\right\}\right) \in \mathcal{A}_H$ without any loss of information.\\
	
	The more complicated definition was introduced in \cite{zhang2019introducing} with the aim of applying not oriented hyperarcs and their adjacency tensor to signal processing, however those specific properties are not necessary for the scope of this thesis. Therefore, the simplified definitions for the adjacency tensors $T_{NH}$ and $T_{OH}$ are used in the examples below.\\
\end{remark}

\begin{example}[\textbf{Adjacency tensors}] \ \\
	Given the not oriented hypergraph $NH = \left(\mathcal{V}, \mathcal{E}_H\right)$ with the vertex set $\mathcal{V} =$\\ $\left\{v_1, v_2, v_3, v_4, v_5, v_6\right\}$ and hyperedges $\mathcal{E}_H = \left\{\left\{v_1, v_2, v_4\right\}, \left\{v_2, v_3, v_6\right\}, \left\{v_5, v_6\right\}\right\}$, then the maximum hyperedge cardinality is $max_e = 3$ and hence the hyperedges are described by the following indices:
	
	\begin{minipage}{.4\textwidth}\begin{tikzpicture}
				\tikzstyle{vertex} = [fill,shape=circle,node distance=80pt]
				\tikzstyle{edge} = [fill,opacity=.5,fill opacity=.5,line cap=round, line join=round, line width=30pt]
				\tikzstyle{elabel} =  [fill,shape=circle,node distance=40pt]
				
				\pgfdeclarelayer{background}
				\pgfsetlayers{background,main}
				
				\begin{scope}[every node/.style={circle,thick,draw}]
					\node (v1) at (0,0) {$v_1$};
					\node (v2) at (2,0) {$v_2$};
					\node (v3) at (4,0) {$v_3$};
					\node (v4) at (0,-2) {$v_4$};
					\node (v5) at (2,-2) {$v_5$};
					\node (v6) at (4,-2) {$v_6$};
				\end{scope}
				
				\begin{pgfonlayer}{background}
					\begin{scope}[transparency group,opacity=.5]
						\draw[edge,opacity=1,color=blue] (v1.center) -- (v2.center) -- (v4.center) -- (v1.center);
						\fill[edge,opacity=1,color=blue] (v1.center) -- (v2.center) -- (v4.center) -- (v1.center);
					\end{scope}
					\begin{scope}[transparency group,opacity=.5]
						\draw[edge,opacity=1,color=green] (v2.center) -- (v3.center) -- (v6.center) -- (v2.center);
						\fill[edge,opacity=1,color=green] (v2.center) -- (v3.center) -- (v6.center) -- (v2.center);
					\end{scope}
					\begin{scope}[transparency group,opacity=.5]
						\draw[edge,opacity=1,color=lightgray] (v5.center) -- (v6.center) -- (v5.center);
						\fill[edge,opacity=1,color=lightgray] (v5.center) -- (v6.center) -- (v5.center);
					\end{scope}
				\end{pgfonlayer}
				
				\node[text=blue] (a1) at (1,0.75) {$e_1$};
				\node[text=green] (a2) at (4.75,-1) {$e_2$};
				\node[text=lightgray] (a3) at (3,-2.75) {$e_3$};
	\end{tikzpicture}\end{minipage}
	\begin{minipage}{.6\textwidth} \begin{flushright}
		$\begin{aligned}[t]	
			e_1 = \left\{v_1, v_2, v_4\right\}: \quad & \left(1, 2, 4\right), \left(1, 4, 2\right), \left(2, 1, 4\right), &\\
			& \left(4, 1, 2\right), \left(2, 4, 1\right), \left(4, 2, 1\right) &\\
			e_2 = \left\{v_2, v_3, v_6\right\}: \quad & \left(2, 3, 6\right), \left(2, 6, 3\right), \left(3, 2, 6\right), &\\
			& \left(6, 2, 3\right), \left(3, 6, 2\right), \left(6, 3, 2\right) &\\
			e_3 = \left\{v_5, v_6\right\}: \quad & \left(5, 5, 6\right), \left(5, 6, 5\right), \left(6, 5, 5\right), &\\
			& \left(5, 6, 6\right), \left(6, 5, 6\right), \left(6, 6, 5\right) &\\
		\end{aligned}$
	\end{flushright}\end{minipage}
	
	And the adjacency tensor $T_{NH}$ is a $3^{\text{rd}}$-order $6$-dimensional tensor:
	
	$\left(T_{NH}\right)_{\cdot, \cdot, 1} = \begin{pmatrix}
		0 & 0 & 0 & 0 & 0 & 0\\
		0 & 0 & 0 & \frac{3}{2} & 0 & 0\\
		0 & 0 & 0 & 0 & 0 & 0\\
		0 & \frac{3}{2} & 0 & 0 & 0 & 0\\
		0 & 0 & 0 & 0 & 0 & 0\\
		0 & 0 & 0 & 0 & 0 & 0\\
		0 & 0 & 0 & 0 & 0 & 0\\
	\end{pmatrix}
	\qquad \qquad \quad \left(T_{NH}\right)_{\cdot, \cdot, 2} = \begin{pmatrix}
		0 & 0 & 0 & \frac{3}{2} & 0 & 0\\
		0 & 0 & 0 & 0 & 0 & 0\\
		0 & 0 & 0 & 0 & 0 & \frac{3}{2}\\
		\frac{3}{2} & 0 & 0 & 0 & 0 & 0\\
		0 & 0 & 0 & 0 & 0 & 0\\
		0 & 0 & 0 & 0 & 0 & 0\\
		0 & 0 & \frac{3}{2} & 0 & 0 & 0\\
	\end{pmatrix}$
	
	$\left(T_{NH}\right)_{\cdot, \cdot, 3} = \begin{pmatrix}
		0 & 0 & 0 & 0 & 0 & 0\\
		0 & 0 & 0 & 0 & 0 & \frac{3}{2}\\
		0 & 0 & 0 & 0 & 0 & 0\\
		0 & 0 & 0 & 0 & 0 & 0\\
		0 & 0 & 0 & 0 & 0 & 0\\
		0 & 0 & 0 & 0 & 0 & 0\\
		0 & \frac{3}{2} & 0 & 0 & 0 & 0\\
	\end{pmatrix}
	\qquad \qquad \quad \left(T_{NH}\right)_{\cdot, \cdot, 4} = \begin{pmatrix}
		0 & \frac{3}{2} & 0 & 0 & 0 & 0\\
		\frac{3}{2} & 0 & 0 & 0 & 0 & 0\\
		0 & 0 & 0 & 0 & 0 & 0\\
		0 & 0 & 0 & 0 & 0 & 0\\
		0 & 0 & 0 & 0 & 0 & 0\\
		0 & 0 & 0 & 0 & 0 & 0\\
		0 & 0 & 0 & 0 & 0 & 0\\
	\end{pmatrix}$

	$\left(T_{NH}\right)_{\cdot, \cdot, 5} = \begin{pmatrix}
		0 & 0 & 0 & 0 & 0 & 0\\
		0 & 0 & 0 & 0 & 0 & 0\\
		0 & 0 & 0 & 0 & 0 & 0\\
		0 & 0 & 0 & 0 & 0 & 0\\
		0 & 0 & 0 & 0 & 0 & 0\\
		0 & 0 & 0 & 0 & 0 & 1\\
		0 & 0 & 0 & 0 & 1 & 1\\
	\end{pmatrix}
	\qquad \qquad \quad \left(T_{NH}\right)_{\cdot, \cdot, 6} = \begin{pmatrix}
		0 & 0 & 0 & 0 & 0 & 0\\
		0 & 0 & \frac{3}{2} & 0 & 0 & 0\\
		0 & \frac{3}{2} & 0 & 0 & 0 & 0\\
		0 & 0 & 0 & 0 & 0 & 0\\
		0 & 0 & 0 & 0 & 0 & 0\\
		0 & 0 & 0 & 0 & 1 & 1\\
		0 & 0 & 0 & 0 & 1 & 0\\
	\end{pmatrix}$\\
	
	\clearpage
	Similarly, in the case of an oriented hypergraph $OH = \left(\mathcal{V}, \mathcal{A}_H\right)$ with vertex set $\mathcal{V} = \left\{v_1, v_2, v_3, v_4, v_5, v_6\right\}$ and the set of hyperarcs $\mathcal{A}_H =$ \\ $\left\{\left(\left\{v_1, v_2\right\}, \left\{v_4\right\}\right), \left(\left\{v_3\right\}, \left\{v_2, v_6\right\}\right), \left(\left\{v_5\right\}, \left\{v_6\right\}\right)\right\}$, the maximum hyperarc cardinality is given by $max_a = 3$ and therefore the hyperarcs are described by the following indices \big(Note: All other indices consisting of the vertices of the hyperarcs are not valid because it is not possible to separate them into output and input vertices according to the hyperarc without changing their order.\big):
	
	\begin{minipage}{.4\textwidth}\begin{tikzpicture}
				\tikzstyle{vertex} = [fill,shape=circle,node distance=80pt]
				\tikzstyle{edge} = [fill,opacity=.5,fill opacity=.5,line cap=round, line join=round, line width=30pt]
				\tikzstyle{elabel} =  [fill,shape=circle,node distance=40pt]
				
				\pgfdeclarelayer{background}
				\pgfsetlayers{background,main}
				
				\begin{scope}[every node/.style={circle,thick,draw}]
					\node (v1) at (0,0) {$v_1$};
					\node (v2) at (2,0) {$v_2$};
					\node (v3) at (4,0) {$v_3$};
					\node (v4) at (0,-2) {$v_4$};
					\node (v5) at (2,-2) {$v_5$};
					\node (v6) at (4,-2) {$v_6$};
				\end{scope}
				
				\begin{pgfonlayer}{background}
					\begin{scope}[transparency group,opacity=.2]
						\draw[edge,opacity=1,color=blue] (v1.center) -- (v2.center) -- (v4.center) -- (v1.center);
						\fill[edge,opacity=1,color=blue] (v1.center) -- (v2.center) -- (v4.center) -- (v1.center);
					\end{scope}
					\begin{scope}[transparency group,opacity=.2]
						\draw[edge,opacity=1,color=green] (v2.center) -- (v3.center) -- (v6.center) -- (v2.center);
						\fill[edge,opacity=1,color=green] (v2.center) -- (v3.center) -- (v6.center) -- (v2.center);
					\end{scope}
					\begin{scope}[transparency group,opacity=.2]
						\draw[edge,opacity=1,color=lightgray] (v5.center) -- (v6.center) -- (v5.center);
						\fill[edge,opacity=1,color=lightgray] (v5.center) -- (v6.center) -- (v5.center);
					\end{scope}
				\end{pgfonlayer}
				
				\node[text=blue] (a1) at (1,0.75) {$a_1$};
				\node[text=green] (a2) at (4.75,-1) {$a_2$};
				\node[text=lightgray] (a3) at (3,-2.75) {$a_3$};
				\node[text=blue] (v1a1) at (0.5, -0.5) {\footnotesize{$out$}};
				\node[text=blue] (v2a1) at (1.5, -0.5) {\footnotesize{$out$}};
				\node[text=blue] (v5a1) at (0.5, -1.5) {\footnotesize{$in$}};
				\node[text=green] (v3a2) at (4, -0.7) {\footnotesize{$out$}};
				\node[text=green] (v2a2) at (2.5, -0.5) {\footnotesize{$in$}};
				\node[text=green] (v7a2) at (4, -1.3) {\footnotesize{$in$}};
				\node[text=lightgray] (v6a3) at (2.75, -2) {\footnotesize{$out$}};
				\node[text=lightgray] (v7a3) at (3.35, -2) {\footnotesize{$in$}};
	\end{tikzpicture}\end{minipage}
	\begin{minipage}{.6\textwidth}\begin{flushright}
		$\begin{aligned}[t]	
			a_1 = \left(\left\{v_1, v_2\right\}, \left\{v_4\right\}\right): \quad & \left(1, 2, 4\right), \left(2, 1, 4\right) &\\
			a_2 = \left(\left\{v_3\right\}, \left\{v_2, v_6\right\}\right): \quad & \left(3, 2, 6\right), \left(3, 6, 2\right) &\\
			a_3 = \left(\left\{v_5\right\}, \left\{v_6\right\}\right): \quad & \left(5, 5, 6\right), \left(5, 6, 6\right) &\\
		\end{aligned}$
	\end{flushright}\end{minipage}
	
	The adjacency tensor $T_{OH}$ is again a $3^{\text{rd}}$-order $6$-dimensional tensor:
	
	$\left(T_{NH}\right)_{\cdot, \cdot, 1} = \begin{pmatrix}
		0 & 0 & 0 & 0 & 0 & 0\\
		0 & 0 & 0 & 0 & 0 & 0\\
		0 & 0 & 0 & 0 & 0 & 0\\
		0 & 0 & 0 & 0 & 0 & 0\\
		0 & 0 & 0 & 0 & 0 & 0\\
		0 & 0 & 0 & 0 & 0 & 0\\
		0 & 0 & 0 & 0 & 0 & 0\\
	\end{pmatrix}
	\qquad \qquad \quad \left(T_{NH}\right)_{\cdot, \cdot, 2} = \begin{pmatrix}
		0 & 0 & 0 & 0 & 0 & 0\\
		0 & 0 & 0 & 0 & 0 & 0\\
		0 & 0 & 0 & 0 & 0 & \frac{3}{2}\\
		0 & 0 & 0 & 0 & 0 & 0\\
		0 & 0 & 0 & 0 & 0 & 0\\
		0 & 0 & 0 & 0 & 0 & 0\\
		0 & 0 & 0 & 0 & 0 & 0\\
	\end{pmatrix}$
	
	$\left(T_{NH}\right)_{\cdot, \cdot, 3} = \begin{pmatrix}
		0 & 0 & 0 & 0 & 0 & 0\\
		0 & 0 & 0 & 0 & 0 & 0\\
		0 & 0 & 0 & 0 & 0 & 0\\
		0 & 0 & 0 & 0 & 0 & 0\\
		0 & 0 & 0 & 0 & 0 & 0\\
		0 & 0 & 0 & 0 & 0 & 0\\
		0 & 0 & 0 & 0 & 0 & 0\\
	\end{pmatrix}
	\qquad \qquad \quad \left(T_{NH}\right)_{\cdot, \cdot, 4} = \begin{pmatrix}
		0 & \frac{3}{2} & 0 & 0 & 0 & 0\\
		\frac{3}{2} & 0 & 0 & 0 & 0 & 0\\
		0 & 0 & 0 & 0 & 0 & 0\\
		0 & 0 & 0 & 0 & 0 & 0\\
		0 & 0 & 0 & 0 & 0 & 0\\
		0 & 0 & 0 & 0 & 0 & 0\\
		0 & 0 & 0 & 0 & 0 & 0\\
	\end{pmatrix}$
	
	$\left(T_{NH}\right)_{\cdot, \cdot, 5} = \begin{pmatrix}
		0 & 0 & 0 & 0 & 0 & 0\\
		0 & 0 & 0 & 0 & 0 & 0\\
		0 & 0 & 0 & 0 & 0 & 0\\
		0 & 0 & 0 & 0 & 0 & 0\\
		0 & 0 & 0 & 0 & 0 & 0\\
		0 & 0 & 0 & 0 & 0 & 0\\
		0 & 0 & 0 & 0 & 0 & 0\\
	\end{pmatrix}
	\qquad \qquad \quad \left(T_{NH}\right)_{\cdot, \cdot, 6} = \begin{pmatrix}
		0 & 0 & 0 & 0 & 0 & 0\\
		0 & 0 & 0 & 0 & 0 & 0\\
		0 & \frac{3}{2} & 0 & 0 & 0 & 0\\
		0 & 0 & 0 & 0 & 0 & 0\\
		0 & 0 & 0 & 0 & 0 & 0\\
		0 & 0 & 0 & 0 & 1 & 1\\
		0 & 0 & 0 & 0 & 0 & 0\\
	\end{pmatrix}$\\	
\end{example}

\clearpage
The next pages discuss the problem of how to deduce the unknown underlying oriented hypergraph from a given adjacency tensor $T_{OH}$. In order to accomplish this, some preliminary results are necessary, such as the number of indices $m_1, m_2, \dots m_{max_a}$ of the adjacency tensor $T_{OH}$, which encode a hyperarc $a_q \in \mathcal{A}_H$ with $\left\lvert a_q^{out} \cup a_q^{in}\right\rvert = max_a$ or $\left\lvert a_q^{out} \cup a_q^{in}\right\rvert = max_a - 1$.\\

\begin{lemma}[\textbf{Number of indices in $T_{OH}$ describing hyperarcs}]\label{TOHindi} \ \\
	In the adjacency tensor $T_{OH}$ of an oriented hypergraph $OH = \left(\mathcal{V}, \mathcal{A}_H\right)$, the number of indices describing hyperarcs $a_q \in \mathcal{A}_H$ with $\left\lvert a_q^{out} \cup a_q^{in}\right\rvert = max_a$ is given by
	\begin{equation}
		\left\lvert a_q^{out}\right\rvert! ~ \left\lvert a_q^{in}\right\rvert!
	\end{equation}
	and similarly the number of indices describing hyperarcs $a_q \in \mathcal{A}_H$ with $\left\lvert a_q^{out} \cup a_q^{in}\right\rvert = max_a - 1$ is given by
	\begin{equation}
		\frac{\left(\left\lvert a_q^{out}\right\rvert + 1\right)!}{2!} \left\lvert a_q^{out}\right\rvert \left\lvert a_q^{in}\right\rvert! ~ + ~ \left\lvert a_q^{out}\right\rvert! \frac{\left(\left\lvert a_q^{in}\right\rvert + 1\right)!}{2!} \left\lvert a_q^{in}\right\rvert .
	\end{equation}\vspace{0em}
\end{lemma}

\begin{proof}\ \\
	\textbf{First case:} $\left\lvert a_q^{out} \cup a_q^{in}\right\rvert = max_a$
	
	Given the adjacency tensor $T_{OH}$ of an oriented hypergraph $OH = \left(\mathcal{V}, \mathcal{A}_H\right)$ and a hyperarc $a_q \in \mathcal{A}_H$ with $\left\lvert a_q^{out} \cup a_q^{in}\right\rvert = max_a$, then the hyperarc consists of $n = \left\lvert a_q^{out}\right\rvert + \left\lvert a_q^{in}\right\rvert = max_a$ different vertices. The $max_a$ indices of every adjacency tensor entry describing $a_q$ are split into $\left\lvert a_q^{out}\right\rvert$ indices describing all $\left\lvert a_q^{out}\right\rvert$ output vertices and similarly $\left\lvert a_q^{in}\right\rvert$ indices describing all $\left\lvert a_q^{in}\right\rvert$ output vertices.\\
	
	There are $\left\lvert a_q^{out}\right\rvert!$ possible sequences describing the output vertices, since the order within does matter, and based on the same reason there are $\left\lvert a_q^{in}\right\rvert!$ possible sequences describing the input vertices. Since the indices are combined, overall there are
	\begin{equation*}
		\left\lvert a_q^{out}\right\rvert! \cdot \left\lvert a_q^{in}\right\rvert!
	\end{equation*}
	indices encoding hyperarc $a_q = \left(a_q^{out}, a_q^{in}\right)$.\\
	
	\clearpage
	\textbf{Second case:} $\left\lvert a_q^{out} \cup a_q^{in}\right\rvert = max_a - 1$
	
	Given a hyperarc $a_q \in \mathcal{A}_H$ with $\left\lvert a_q^{out} \cup a_q^{in}\right\rvert = max_a - 1$, then the hyperarc consists of $n = \left\lvert a_q^{out}\right\rvert + \left\lvert a_q^{in}\right\rvert = max_a - 1$ different vertices. This means, that any indices describing hyperarc $a_q$ have exactly one redundant vertex, either a duplicate of an output vertex or a duplicate of an input vertex.\\
	
	If the redundant vertex is an output vertex, then the number of valid indices describing the output vertices is given by
	\begin{equation*}
		\frac{\left(\left\lvert a_q^{out}\right\rvert + 1\right)!}{2!} \left\lvert a_q^{out}\right\rvert
	\end{equation*}
	where $\left(\left\lvert a_q^{out}\right\rvert + 1\right)!$ is the number of indices describing the output vertices, $\frac{1}{2!}$ is necessary because one vertex appears twice \big(with the duplicates of the index appearing twice being indistinguishable\big) and the factor $\left\lvert a_q^{out}\right\rvert$ appears, because there are $\left\lvert a_q^{out}\right\rvert$ remaining spots where the duplicate index can be.\\
	
	Since there is no redundant index for the input vertices, the number of valid indices describing the input vertices is given by $\left\lvert a_q^{in}\right\rvert!$.
	
	The indices are once again put together, and hence the number of indices describing the hyperarc $a_q$ with the redundancy in the output vertices indices can be calculated with:
	\begin{equation*}
		\frac{\left(\left\lvert a_q^{out}\right\rvert + 1\right)!}{2!} \cdot \left\lvert a_q^{out}\right\rvert \cdot \left\lvert a_q^{in}\right\rvert! .
	\end{equation*}
	Based on the same arguments, the number of indices describing the hyperarc $a_q$ with the redundancy being in the input vertices indices can be calculated with:
	\begin{equation*}
		\frac{\left(\left\lvert a_q^{in}\right\rvert + 1\right)!}{2!} \cdot \left\lvert a_q^{in}\right\rvert \cdot \left\lvert a_q^{out}\right\rvert! .
	\end{equation*}
	Therefore, the overall number of indices encoding hyperarc $a_q$ with $\left\lvert a_q^{out} \cup a_q^{in}\right\rvert = max_a - 1$ is given by:
	\begin{equation*}
		\frac{\left(\left\lvert a_q^{out}\right\rvert + 1\right)!}{2!} \left\lvert a_q^{out}\right\rvert \left\lvert a_q^{in}\right\rvert! ~ + ~ \left\lvert a_q^{out}\right\rvert! \frac{\left(\left\lvert a_q^{in}\right\rvert + 1\right)!}{2!} \left\lvert a_q^{in}\right\rvert .
	\end{equation*}\vspace{0em}
\end{proof}

\clearpage
Besides the number of indices in $T_{OH}$ describing a hyperarc $a_q \in \mathcal{A}_H$ with $\left\lvert a_q^{out} \cup a_q^{in}\right\rvert \in \left\{max_a - 1, max_a\right\}$, other properties of the adjacency tensor are necessary for the pseudo code of remark (\ref{TOHinter}) as well.\\

\begin{theorem}[\textbf{Properties of the adjacency tensor $T_{OH}$}]\label{TOHprop} \ \\
	All entries of the adjacency tensor $\left(T_{OH}\right)_{m_1, m_2, \dots m_{max_a}}$ of an oriented hypergraph $OH = \left(\mathcal{V}, \mathcal{A}_H\right)$ with maximum hyperarc cardinality $max_a$ have the following properties:
	\begin{itemize}
		\item[1)] $\left(T_{OH}\right)_{m_1, m_2, \dots m_{max_a}}$ can represent at most $max_a - 1$ different hyperarcs.
		\item[2)] $\left(T_{OH}\right)_{m_1, m_2, \dots m_1, m_{max_a}, \dots m_{max_a}}$ can represent at most one hyperarc.
		\item[3)] $\left(T_{OH}\right)_{m_1, m_2, \dots m_{max_a}} = \frac{\left\lvert\left\{a_q \in \mathcal{A}_H ~ \text{represented by} ~ m_1, m_2, \dots m_{max_a}\right\}\right\rvert \cdot n}{\sum_{\left\{l_1, l_2, \dots l_n \in \mathbb{N} \backslash \left\{0\right\}: ~ \sum_{k = 1}^{n} l_k = max_a\right\}} \frac{max_a!}{l_1! ~ l_2! ~ \cdots ~ l_n!}}$ .
	\end{itemize}	
\end{theorem}

\begin{proof}\ \\
	Given an oriented hypergraph $OH = \left(\mathcal{V}, \mathcal{A}_H\right)$ with maximum hyperarc cardinality $max_a$ and indices $m_1, m_2, \dots m_{max_a}$ for the adjacency tensor $T_{OH}$, then it holds true that:
	\begin{itemize}
		\item[1)] For all $s \in \left\{1, 2, \dots m_{max_a} - 1\right\}$ the partition of the indices $m_1, m_2, \dots m_{max_a}$ into two sets $\left(\left\{m_1, m_2, \dots m_s\right\}, \left\{m_{s + 1}, m_{s + 2}, \dots m_{max_a}\right\}\right)$ corresponds to a possible hyperarc $a_q \in \mathcal{A}_H$.
		
		Note: The indices $m_1, m_2, \dots m_{max_a}$ are not necessarily pairwise distinct and therefore some partitions do not constitute a feasible hyperarc. A partition $\left(\left\{m_1, m_2, \dots m_s\right\}, \left\{m_{s + 1}, m_{s + 2}, \dots m_{max_a}\right\}\right)$ of the indices with $m_t \in \left\{m_1, m_2, \dots m_s\right\}$, $m_u \in \left\{m_{s + 1}, m_{s + 2}, \dots m_{max_a}\right\}$ and $m_t = m_u$ does not represent a feasible hyperarc $a_q \in \mathcal{A}_H$ since $a_q^{out} \cap a_q^{in} \neq \emptyset$.
		
		Hence, $\left(T_{OH}\right)_{m_1, m_2, \dots m_{max_a}}$ encodes at most $max_a - 1$ different hyperarcs.\\
		
		\item[2)] Since $a_q^{out}, a_q^{in} \neq \emptyset$ and $a_q^{out} \cap a_q^{in} = \emptyset$ for all hyperarcs $a_q \in \mathcal{A}_H$, it is certain that:
		\begin{equation*}
			v_{m_1} \in a_q^{out} \quad \text{and} \quad v_{m_{max_a}} \in a_q^{in}
		\end{equation*}
		Thus, if $m_1$ or $m_{max_a}$ are also present in $m_2, m_3, \dots, m_{max_a - 1}$ then all the indices in between correspond to output or input vertices as well. In the case of $m_1, m_2, \dots m_1, m_{max_a}, \dots m_{max_a}$ the indices only allow one feasible partition $\left(\left\{m_1, m_2, \dots m_1\right\}, \left\{m_{max_a}, \dots m_{max_a - 1}, m_{max_a}\right\}\right)$ and hence encode at most one hyperarc.
		
		The partition $\left(\left\{m_1, m_2, \dots m_1\right\}, \left\{m_{max_a}, \dots m_{max_a - 1}, m_{max_a}\right\}\right)$ only corresponds to a feasible hyperarc $a_q \in \mathcal{A}_H$, if $\left\{m_1, m_2, \dots m_1\right\} \cap  \left\{m_{max_a}, \dots m_{max_a - 1}, m_{max_a}\right\} = \emptyset$.\\
		
		\item[3)] If the indices $m_1, m_2, \dots m_{max_a}$ do not correspond to any hyperarc then the set $\left\{a_q \in \mathcal{A}_H ~ \text{represented by} ~ m_1, m_2, \dots m_{max_a}\right\}$ is empty and it holds also true that $\left\lvert\left\{a_q \in \mathcal{A}_H ~ \text{represented by} ~ m_1, m_2, \dots m_{max_a}\right\}\right\rvert = 0$ and thus $\left(T_{OH}\right)_{m_1, m_2, \dots m_{max_a}} = \frac{\left\lvert\left\{a_q \in \mathcal{A}_H ~ \text{represented by} ~ m_1, m_2, \dots m_{max_a}\right\}\right\rvert \cdot n}{\sum_{\left\{l_1, l_2, \dots l_n \in \mathbb{N} \backslash \left\{0\right\}: ~ \sum_{k = 1}^{n} l_k = max_a\right\}} \frac{max_a!}{l_1! ~ l_2! ~ \cdots ~ l_n!}} = 0$ since $n \geq 2$. This aligns with the definition of the adjacency tensor $T_{OH}$ for entries which do not correspond to any hyperarcs $a_q \in \mathcal{A}_H$.\\
		
		If the indices $m_1, m_2, \dots m_{max_a}$ describe at least one hyperarc then all hyperarcs of $\left\{a_q \in \mathcal{A}_H ~ \text{represented by} ~ m_1, m_2, \dots m_{max_a}\right\}$ consist of the same $n$ vertices $\left\{v_{i_1}, v_{i_2}, \dots v_{i_n}\right\} = \left\{v_{m_1}, v_{m_2}, \dots v_{m_{max_a}}\right\}$ with the partition in output and input vertices varying. Since $n$ and $max_a$ are the same for all hyperarcs of the set $\left\{a_q \in \mathcal{A}_H ~ \text{represented by} ~ m_1, m_2, \dots m_{max_a}\right\}$, it holds true that:
		\begin{align*}
			\left(T_{OH}\right)_{m_1, m_2, \dots m_{max_a}} & = \sum_{\substack{a_q \in \mathcal{A}_H ~ \text{represented by} \\ m_1, m_2, \dots m_{max_a}}} \left(T_{OH}\right)_{m_1, m_2, \dots m_{max_a}}^{a_q}\\
			& = \sum_{\substack{a_q \in \mathcal{A}_H ~ \text{represented by} \\ m_1, m_2, \dots m_{max_a}}} \frac{n}{\sum_{\left\{l_1, l_2, \dots l_n \in \mathbb{N} \backslash \left\{0\right\}: ~ \sum_{k = 1}^{n} l_k = max_a\right\}} \frac{max_a!}{l_1! ~ l_2! ~ \cdots ~ l_n!}}\\
			& = \frac{\left\lvert\left\{a_q \in \mathcal{A}_H ~ \text{represented by} ~ m_1, m_2, \dots m_{max_a}\right\}\right\rvert \cdot n}{\sum_{\left\{l_1, l_2, \dots l_n \in \mathbb{N} \backslash \left\{0\right\}: ~ \sum_{k = 1}^{n} l_k = max_a\right\}} \frac{max_a!}{l_1! ~ l_2! ~ \cdots ~ l_n!}}
		\end{align*}
	
		Therefore, for all indices $m_1, m_2, \dots m_{max_a} \in \left\{1, 2, \dots N\right\}$ the equality
		\begin{flalign*}
			& \left(T_{OH}\right)_{m_1, m_2, \dots m_{max_a}} = \frac{\left\lvert\left\{a_q \in \mathcal{A}_H ~ \text{represented by} ~ m_1, m_2, \dots m_{max_a}\right\}\right\rvert \cdot n}{\sum_{\left\{l_1, l_2, \dots l_n \in \mathbb{N} \backslash \left\{0\right\}: ~ \sum_{k = 1}^{n} l_k = max_a\right\}} \frac{max_a!}{l_1! ~ l_2! ~ \cdots ~ l_n!}} &
		\end{flalign*}
	 	is fulfilled.
	\end{itemize}
\end{proof}

The following example shows why reading out an oriented hypergraph $OH = \left(\mathcal{V}, \mathcal{A}_H\right)$ from its adjacency tensor $T_{OH}$ is not as straight forward as in the not oriented case.\\

\begin{example}[\textbf{Difficulty with the interpretation of the adjacency tensor $T_{OH}$}]\label{TOHexam} \ \\
	Given the oriented hypergraph $OH = \left(\mathcal{V}, \mathcal{A}_H\right)$ with vertex set $\mathcal{V} = \left\{v_1, v_2, v_3\right\}$ and the set of hyperarcs $\mathcal{A}_H = \left\{\left(\left\{v_1\right\}, \left\{v_2, v_3\right\}\right), \left(\left\{v_2\right\}, \left\{v_1, v_3\right\}\right)\right\}$, then the two hyperarcs are encoded with the following indices in the adjacency tensor $T_{OH}$:
	
	$\begin{aligned}[t]	
		\qquad \qquad \qquad \qquad a_1 = \left(\left\{v_1\right\}, \left\{v_2, v_3\right\}\right): \quad & \left(1, 2, 3\right), \left(1, 3, 2\right) &\\
		\qquad \qquad \qquad \qquad a_2 = \left(\left\{v_2\right\}, \left\{v_1, v_3\right\}\right): \quad & \left(2, 1, 3\right), \left(2, 3, 1\right) &\\
	\end{aligned}$\\
	
	When only considering the nonzero indices of the adjacency tensor, it is not obvious that hyperarc $\left(\left\{v_1, v_2\right\}, \left\{v_3\right\}\right)$ is not a part of the underlying oriented hypergraph since both indices $\left(1, 2, 3\right)$ and $\left(2, 1, 3\right)$, encoding this hyperarc, have nonzero entries in the adjacency tensor.\\
\end{example}

Therefore, it is not enough to check for each possible hyperarc \big(based on the vertices $\mathcal{V}$\big) if all corresponding indices lead to a nonzero entry in the adjacency tensor $T_{OH}$ to find all hyperarcs of the underlying oriented hypergraph. Instead, an interative approach is necessary.\\

\begin{remark}[\textbf{Interpreting the adjacency tensor $T_{OH}$}]\label{TOHinter} \ \\
	Given an adjacency tensor entry $\left(T_{OH}\right)_{m_1, m_2, \dots m_{max_a}} > 0$, then the following steps should be completed in order to find the encoded hyperarcs $a_q = \left(a_q^{out}, a_q^{in}\right) \in \mathcal{A}_H$ of the oriented hypergraph $OH = \left(\mathcal{V}, \mathcal{A}_H\right)$.\\
	
	Since all hyperarcs represented by the indices $m_1, m_2, \dots m_{max_a}$ consist of the same $n = \left\lvert\left\{m_1, m_2, \dots m_{max_a}\right\}\right\rvert$ vertices given by $\left\{v_{i_1}, v_{i_2}, \dots v_{i_n}\right\} = \left\{v_{m_1}, v_{m_2}, \dots v_{m_{max_a}}\right\}$, the approach differentiates between two cases:
	\begin{equation*}
		n \leq max_a - 2 \qquad \text{and} \qquad n \in \left\{max_a - 1, max_a\right\}.
	\end{equation*}
	
	\textbf{First case:} $n \leq max_a - 2$
	\begin{itemize}
		\item[1)] As argued in theorem (\ref{TOHprop}), for all indices $m_1, m_2, \dots m_{max_a}$ with adjacency tensor entry $\left(T_{OH}\right)_{m_1, m_2, \dots m_{max_a}} > 0$ it is certain that $v_{m_1} \in a_q^{out}$ and $v_{m_{max_a}} \in a_q^{in}$ for all hyperarcs $a_q \in \mathcal{A}_H$ represented by $m_1, m_2, \dots m_{max_a}$. Furthermore, each encoded hyperarc can be described by a partition
		\begin{equation*}
			\left(\left\{v_{m_1}, v_{m_2}, \dots v_{m_s}\right\}, \left\{v_{m_{s + 1}}, v_{m_{s + 2}}, \dots v_{m_{max_a}}\right\}\right) ~ \text{with} ~ s \in \left\{1, 2, \dots max_a - 1\right\}.
		\end{equation*}
		
		In the case of $n \leq max_a - 2$, at least two of the indices in $m_1, m_2, \dots m_{max_a}$ are redundant and can therefore be substituted by $m_1$ and $m_{max_a}$, with the new indices still encoding the same hyperarc $a_q \in \mathcal{A}_H$. This is especially true for indices in the form of $m_1, m_2, \dots m_s, m_1, m_{max_a}, m_{s + 1}, \dots m_{max_a}$ uniquely describing the previous partition as shown in theorem (\ref{TOHprop}). Hence, it suffices to check whether the adjacency entries for such indices are positive, to find all hyperarcs $a_q \in \mathcal{A}_H$ represented by $m_1, m_2, \dots m_{max_a}$.\\
		
		\clearpage
		\begin{algorithmic}
			\State $\text{hyperarcs} = \left\{\right\}$
			\State $\text{remove two redundant indices from} ~ m_1, m_2, \dots m_{max_a}$
			\For{$s \in \left\{2, 3, \dots max_a - 1\right\}$}
				\If{$\left(T_{OH}\right)_{m_1, m_2, \dots m_s, m_1, m_{max_a}, m_{s + 1}, \dots m_{max_a}} > 0$}
					\State $\text{hyperarcs} \longleftarrow \text{hyperarcs} \cup \left(\left\{v_{m_1}, v_{m_2}, \dots v_{m_s}\right\}, \left\{v_{m_{s + 1}}, v_{m_{s + 2}}, \dots v_{m_{max_a}}\right\}\right)$
				\EndIf
			\EndFor
			\State $\textbf{return} ~ \text{hyperarcs}$
		\end{algorithmic}
		
		Note: $\left(T_{OH}\right)_{m_1, m_2, \dots m_s, m_1, m_{max_a}, m_{s + 1}, \dots m_{max_a}} = 0$ either indicates that the indices are not feasible for describing a hyperarc because
		\begin{equation*}
			\left\{m_1, m_2, \dots m_s\right\} \cap \left\{m_{s + 1}, m_{s + 2}, \dots m_{max_a}\right\} \neq \emptyset
		\end{equation*}
		or that the hyperarc $\left(\left\{v_{m_1}, v_{m_2}, \dots v_{m_s}\right\}, \left\{v_{m_{s + 1}}, v_{m_{s + 2}}, \dots v_{m_{max_a}}\right\}\right)$ is not part of the set of hyperarcs $\mathcal{A}_H$ of the underlying oriented hypergraph $OH = \left(\mathcal{V}, \mathcal{A}_H\right)$.\\
	\end{itemize}
	
	\textbf{Second case:} $n \in \left\{max_a - 1, max_a\right\}$
	
	In contrast to the first case, it is not possible to determine the hyperarcs encoded in an adjacency tensor entry $\left(T_{OH}\right)_{m_1, m_2, \dots m_{max_a}} > 0$ with $\left\lvert\left\{m_1, m_2, \dots m_{max_a}\right\}\right\rvert \geq max_a - 1$ without considering all other entries with indices made up of the same set $\left\{m_1, m_2, \dots m_{max_a}\right\}$. Therefore, the following steps should be applied to all nonzero entries of the adjacency tensor $T_{OH}$.
	\begin{itemize}
		\item[1)] It holds true that $v_{m_1} \in a_q^{out}$ and $v_{m_{max_a}} \in a_q^{in}$ for all hyperarcs $a_q \in \mathcal{A}_H$ represented by the indices $m_1, m_2, \dots m_{max_a}$ and feasible hyperarcs also fulfill the condition $a_q^{out} \cap a_q^{in} = \emptyset$. Thus, other certain output and input vertices can be found in the following way:\\
		
		\begin{algorithmic}
			\State $\text{outputvertices} = \left\{v_{m_1}\right\}$
			\State $\text{inputvertices} = \left\{v_{m_{max_a}}\right\}$
			\For{$t \in \left\{2, 3, \dots max_a - 1\right\}$}
				\If{$v_{m_t} \in \text{outputvertices}$}
					\State $\text{outputvertices} \longleftarrow \text{outputvertices} \cup \left\{v_{m_2}, v_{m_3}, \dots v_{m_t}\right\}$
				\EndIf
				\If{$v_{m_t} \in \text{inputvertices}$}
					\State $\text{inputvertices} \longleftarrow \text{inputvertices} \cup \left\{v_{m_t}, v_{m_{t + 1}}, \dots v_{m_{max_a - 1}}\right\}$
				\EndIf
			\EndFor
			\State $\textbf{return} ~ \text{outputvertices}, \text{inputvertices}$
		\end{algorithmic}
		
		\clearpage
		\item[2)] Given the vertex sets after the previous step as $\text{outputvertices} = \left\{v_{m_1}, v_{m_2}, \dots v_{m_t}\right\}$ and $\text{inputvertices} = \left\{v_{m_u}, v_{m_{u + 1}}, \dots v_{m_{max_a}}\right\}$ with $t \in \left\{1, 2, \dots max_a - 1\right\}$, $t < u$, and $u \in \left\{2, 3, \dots max_a\right\}$, then the feasible partitions $s \in \left\{t, t+1, \dots, u\right\}$ can be determined with the following steps:\\
		
		\begin{algorithmic}
			\State $\text{possiblepartitions} = \left\{t, t+1, \dots u - 1\right\}$
			\For{$s \in \left\{t, t+1, \dots u\right\}$}
				\If{$\left\{v_1, \dots v_{m_t}, v_{m_{t + 1}}, \dots v_{m_s}\right\} \cap \left\{v_{m_{s + 1}}, \dots v_{m_{u - 1}}, v_{m_u}, \dots v_{m_{max_a}}\right\} \neq \emptyset$}
					\State $\text{delete} ~ s ~ \text{from} ~ \text{possiblepartitions}$
				\EndIf
			\EndFor
			\State $\textbf{return} ~ \text{possiblepartitions}$
		\end{algorithmic}
	
		\item[3)] Based on the set ``possiblepartions'' for every nonzero entry in the adjacency tensor, an overall set of possible hyperarcs is created. In this set a ``possible hyperarc'' is added without checking if it is already part of the list \big(it can appear more than once\big) in order to record how many indices encode a specific hyperarc:\\

		\begin{algorithmic}
			\State $\text{possiblearcs} = \emptyset$
			\For{$\text{nonzeroindex} \in \text{nonzeroindices}$}
			\For{$s \in \text{possiblepartitions}$}
			\State $\text{add} ~ \left(\left\{v_1, v_2, \dots v_{m_s}\right\}, \left\{v_{m_{s + 1}}, v_{m_{s + 2}}, \dots v_{m_{max_a}}\right\}\right)~ \text{to possiblearcs}$
			\EndFor
			\EndFor
			\State $\textbf{return} ~ \text{possiblearcs}$
		\end{algorithmic}
	
		\item[4)] Since for every hyperarc $a_q \in \mathcal{A}_H$, all corresponding entries in the adjacency tensor $T_{OH}$ must be positive, it is possible to determine hyperarcs which are not part of the original oriented hypergraph by counting the nonzero entries in the adjacency tensor that should encode the ``possible hyperarc'':\\
		
		\clearpage
		\begin{algorithmic}
			\For{$h \in \text{possiblearcs}$}
				\If{$\text{number of vertices in} ~ h = max_a$}
					\State $\text{numberindices} = \left\lvert h^{out}\right\rvert! \cdot \left\lvert h^{in}\right\rvert!$
					\If{$\text{count of} ~ h \in \text{possiblearcs} < \text{numberindices}$}
						\State $\text{remove all} ~ h ~ \text{from possiblearcs}$
					\EndIf
				\EndIf
				\If{$\text{number of vertices in} ~ h = max_a - 1$}
					\State $\text{numberindices} = \frac{\left(\left\lvert h^{out}\right\rvert + 1\right)!}{2!} \cdot \left\lvert h^{out}\right\rvert \cdot \left\lvert h^{in}\right\rvert! + \left\lvert h^{out}\right\rvert! \cdot \frac{\left(\left\lvert h^{in}\right\rvert + 1\right)!}{2!} \cdot \left\lvert h^{in}\right\rvert$
					\If{$\text{count of} ~ h \in \text{possiblearcs} < \text{numberindices}$}
						\State $\text{remove all} ~ h ~ \text{from possiblearcs}$
					\EndIf
				\EndIf
			\EndFor
			\State $\textbf{return} ~ \text{possiblearcs}$
		\end{algorithmic}
	
		Note: ``$\text{count of} ~ h \in \text{possiblearcs}$'' is the number of times ``possible hyperarc'' $h$ appears in the set ``possiblearcs'', and thus this is also the number of indices possibly encoding $h$.\\
	
		\item[5)] The last step of finding the original hyperarcs $a_q \in \mathcal{A}_H$ is to calculate for each nonzero entry in the adjacency tensor how many hyperarcs are encoded in the specific entry. This information is then used to iteratively differentiate between possible index partitions definitely describing a hyperarc or based on step 4) definitely not describing a hyperarc of the underlying oriented hypergraph:\\
		
		\clearpage
		\begin{algorithmic}
			\State $\text{hyperarcs} = \left\{\right\}$
			\While{$\text{nonzeroindices} \neq \emptyset$}
				\For{$\text{nonzeroindex} \in \text{nonzeroindices}$}
					\State $\text{numberarcs} = \left(T_{OH}\right)_{\text{nonzeroindex}} \cdot \frac{2}{n}$
					\State $\text{determine possiplepartitions as hyperarcs based on step 2)}$
					\For{$\text{arc} \in \text{possiblepartitions}$}
						\If{$\text{arc} \notin \text{possiblearcs}$}
							\State $\text{remove arc from} ~ \text{possiblepartitions}$
						\EndIf
						\If{$\text{arc} \in \text{hyperarcs}$}
							\State $\text{remove arc from} ~ \text{possiblepartitions}$
							\State $\text{numberarcs} \longleftarrow \text{numberarcs} - 1$
						\EndIf
					\EndFor
					\If{$\text{numberarcs} = \left\lvert\text{possiblepartitions}\right\rvert$}
						\State{$\text{hyperarcs} \longleftarrow \text{hyperarcs} \cup \text{possiblepartitions}$}
						\State{$\text{remove nonzeroindex from nonzeroindices}$}
					\EndIf
					\If{$\text{numberarcs = 0}$}
						\State{$\text{remove nonzeroindex from nonzeroindices}$}
					\EndIf
				\EndFor
			\EndWhile
			\State $\textbf{return} ~ \text{hyperarcs}$
		\end{algorithmic}
	
		Note: Based on theorem (\ref{TOHprop}), the following equality holds true:
		\begin{equation*}
			\left(T_{OH}\right)_{\text{nonzeroindex}} = \frac{\left\lvert\left\{a_q \in \mathcal{A}_H ~ \text{represented by nonzeroindex}\right\}\right\rvert \cdot n}{2} \quad \Longleftrightarrow
		\end{equation*}
		\begin{equation*}
			\text{numberarcs} = \left\lvert\left\{a_q \in \mathcal{A}_H ~ \text{represented by nonzeroindex}\right\}\right\rvert = \left(T_{OH}\right)_{\text{nonzeroindex}} \cdot \frac{2}{n}
		\end{equation*}
	
		The algorithm iterates over all possible partitions for each nonzero index, until each partition of each nonzero index can either be declared a hyperarc $a_q \in \mathcal{A}_H$ or it is certain that the partition can not describe a hyperarc of the underlying not oriented hypergraph $OH = \left(\mathcal{V}, \mathcal{A}_H\right)$.\\
	\end{itemize}
\end{remark}

Several numerical tests suggest that the presented algorithm returns the correct oriented hypergraph in all cases, however the next example shows that it is possible to construct two hypergraphs, which are represented by the same adjacency tensor $T_{OH}$. Hence, the algorithm cannot differentiate between the two hypergraphs and would not necessarily return the correct underlying hypergraph.\\

\clearpage
\begin{example}[\textbf{Counterexample for uniqueness of the adjacency tensor $T_{OH}$}] \ \\
	Given the oriented hypergraph $OH_1 = \left(\mathcal{V}_1, \mathcal{A}_{H_1}\right)$ with vertex set $\mathcal{V}_1 = \left\{v_1, v_2, v_3, v_4\right\}$ and the set of hyperarcs $\mathcal{A}_{H_1} = \left\{\left(\left\{v_1, v_2\right\}, \left\{v_3, v_4\right\}\right), \left(\left\{v_1, v_3\right\}, \left\{v_2, v_4\right\}\right), \left(\left\{v_1, v_4\right\}, \left\{v_2, v_3\right\}\right),\right.$\\
	$\left.\left(\left\{v_2, v_3\right\}, \left\{v_1, v_4\right\}\right), \left(\left\{v_2, v_4\right\}, \left\{v_1, v_3\right\}\right), \left(\left\{v_3, v_4\right\}, \left\{v_1, v_2\right\}\right)\right\}$, then the hyperarcs are encoded with the following indices in the adjacency tensor $T_{OH}$:
	
	$\begin{aligned}[t]	
		\qquad \qquad a_1^1 = \left(\left\{v_1, v_2\right\}, \left\{v_3, v_4\right\}\right): \quad & \left(1, 2, 3, 4\right), \left(1, 2, 4, 3\right), \left(2, 1, 3, 4\right), \left(2, 1, 4, 3\right) &\\
		\qquad \qquad a_2^1 = \left(\left\{v_1, v_3\right\}, \left\{v_2, v_4\right\}\right): \quad & \left(1, 3, 2, 4\right), \left(1, 3, 4, 2\right), \left(3, 1, 2, 4\right), \left(3, 1, 4, 2\right) &\\
		\qquad \qquad a_3^1 = \left(\left\{v_1, v_4\right\}, \left\{v_2, v_3\right\}\right): \quad & \left(1, 4, 2, 3\right), \left(1, 4, 3, 2\right), \left(4, 1, 2, 3\right), \left(4, 1, 3, 2\right) &\\
		\qquad \qquad a_4^1 = \left(\left\{v_2, v_3\right\}, \left\{v_1, v_4\right\}\right): \quad & \left(2, 3, 1, 4\right), \left(2, 3, 4, 1\right), \left(3, 2, 1, 4\right), \left(3, 2, 4, 1\right) &\\
		\qquad \qquad a_5^1 = \left(\left\{v_2, v_4\right\}, \left\{v_1, v_3\right\}\right): \quad & \left(2, 4, 1, 3\right), \left(2, 4, 3, 1\right), \left(4, 2, 1, 3\right), \left(4, 2, 3, 1\right) &\\
		\qquad \qquad a_6^1 = \left(\left\{v_3, v_4\right\}, \left\{v_1, v_2\right\}\right): \quad & \left(3, 4, 1, 2\right), \left(3, 4, 2, 1\right), \left(4, 3, 1, 2\right), \left(4, 3, 2, 1\right) &\\
	\end{aligned}$\\

	Given a second oriented hypergraph $OH_2 = \left(\mathcal{V}_2, \mathcal{A}_{H_2}\right)$ with vertex set $\mathcal{V}_2 = \left\{v_1, v_2, v_3, v_4\right\}$ and the set of hyperarcs $\mathcal{A}_{H_2} = \left\{\left(\left\{v_1\right\}, \left\{v_2, v_3, v_4\right\}\right), \left(\left\{v_2\right\}, \left\{v_1, v_3, v_4\right\}\right),\right.$\\
	$\left.\left(\left\{v_3\right\}, \left\{v_1, v_2, v_4\right\}\right), \left(\left\{v_4\right\}, \left\{v_1, v_2, v_3\right\}\right)\right\}$, then the hyperarcs are encoded by the same indices as the first hypergraph:
	
	$\begin{aligned}[t]	
		\qquad \qquad a_1^2 = \left(\left\{v_1\right\}, \left\{v_2, v_3, v_4\right\}\right): \quad & \left(1, 2, 3, 4\right), \left(1, 2, 4, 3\right), \left(1, 3, 2, 4\right), \left(1, 3, 4, 2\right), &\\
		& \left(1, 4, 2, 3\right), \left(1, 4, 3, 2\right) &\\
		\qquad \qquad a_2^2 = \left(\left\{v_2\right\}, \left\{v_1, v_3, v_4\right\}\right): \quad & \left(2, 1, 3, 4\right), \left(2, 1, 4, 3\right), \left(2, 3, 1, 4\right), \left(2, 3, 4, 1\right), &\\
		& \left(2, 4, 1, 3\right), \left(2, 4, 3, 1\right) &\\
		\qquad \qquad a_3^2 = \left(\left\{v_3\right\}, \left\{v_1, v_2, v_4\right\}\right): \quad & \left(3, 1, 2, 4\right), \left(3, 1, 4, 2\right), \left(3, 2, 1, 4\right), \left(3, 2, 4, 1\right), &\\
		& \left(3, 4, 1, 2\right), \left(3, 4, 2, 1\right) &\\
		\qquad \qquad a_4^2 = \left(\left\{v_4\right\}, \left\{v_1, v_2, v_3\right\}\right): \quad & \left(4, 1, 2, 3\right), \left(4, 1, 3, 2\right), \left(4, 2, 1, 3\right), \left(4, 2, 3, 1\right), &\\
		& \left(4, 3, 1, 2\right), \left(4, 3, 2, 1\right) &\\
	\end{aligned}$\\

	Both hypergraphs have the same maximum hyperarc cardinality $max_a = 4$ and the same number of vertices $\left\lvert \mathcal{V}\right\rvert = N = 4$. Thus, both hypergraphs are represented by a $4^{\text{th}}$-order $4$-dimensional tensor with the same nonzero indices. Since the value of all nonzero entries only depends on the number of vertices in the encoded hyperarc \big(which is equal to $4$ for all hyperarcs in both hypergraphs\big), it is also not possible to distinguish between the hypergraphs based on the entries of the adjacency tensor.\\
	
	Overall, this constitutes an example of an adjacency tensor $T_{OH}$, which cannot be uniquely matched to the underlying oriented hypergraph and therefore the presented algorithm could return any of the two feasible hypergraphs \big(depending on the order of the set ``nonzeroindices'' in step 5) of remark (\ref{TOHinter})\big).\\
\end{example}

The issues described in the last example can be avoided by increasing the order of the adjacency tensor for oriented hypergraphs by two.\\

\begin{definition}[\textbf{Adjacency tensor with increased order $T_{OH}^+$ }]\label{T_OHincr} \ \\
	For an oriented hypergraph $OH = \left(\mathcal{V}, \mathcal{A}_H\right)$, the adjacency tensor with increased order $T_{OH}^+$ is defined as a $\left(max_a + 2\right)^{\text{th}}$-order $N$-dimensional tensor
	\begin{equation*}
		T_{OH}^+ \in \mathbb{R}^{\tiny{\underbrace{N \times N \times \dots \times N}_{\left(max_a + 2\right) -\text{times}}}} ,
	\end{equation*}
	where every hyperarc $a_q = \left(a_q^{out}, a_q^{in}\right) = \left(\left\{v_{i_1}, v_{i_2}, \dots v_{i_r}\right\}, \left\{v_{i_{r + 1}}, v_{i_{r + 2}}, \dots v_{i_n}\right\}\right) \in \mathcal{A}_H$ consisting of $2 \leq n \leq max_a$ vertices with $1 \leq r \leq n - 1$ is represented by all entries of the tensor $\left(T_{OH}^+\right)_{m_1, m_2, \dots m_s, m_{s + 1}, \dots m_{max_a + 2}}$ with $r \leq s \leq r + max_a + 2 - n$, for which it holds true that:
	\begin{equation*}
		\left\{i_1, i_2, \dots i_r\right\} = \left\{m_1, m_2, \dots m_s\right\}
	\end{equation*}
	\begin{equation*}
		\left\{i_{r + 1}, i_{r + 2}, \dots i_n\right\} = \left\{m_{s + 1}, m_{s + 2}, \dots m_{max_a + 2}\right\}.
	\end{equation*}
	This means that the first $s$ indices of $\left(T_{OH}^+\right)_{m_1, m_2, \dots m_{max_a + 2}}$ only consist of $\left\{i_1, i_2, \dots i_r\right\}$, with each of the vertex indices $\left\{i_1, i_2, \dots i_r\right\}$ being represented in $m_1, m_2, \dots m_s$ at least once. The remaining $max_a + 2 - s$ indices $m_{s + 1}, m_{s + 2}, \dots m_{max_a + 2}$ similarly only consist of the vertex indices $\left\{i_{r + 1}, i_{r + 2}, \dots i_n\right\}$.For all such indices $m_1, m_2, \dots m_{max_a + 2}$ encoding a hyperarc $a_q \in \mathcal{A}_H$, define:
	\begin{equation}
		\left(T_{OH}^+ \right)_{m_1, m_2, \dots m_{max_a + 2}}^{a_q} = \frac{n}{\sum_{\left\{l_1, l_2, \dots l_n \in \mathbb{N} \backslash \left\{0\right\}: ~ \sum_{k = 1}^{n} l_k = max_a + 2\right\}} \frac{\left(max_a + 2\right)!}{l_1! ~ l_2! ~ \cdots ~ l_n!}} > 0.
	\end{equation}
	Since indices $m_1, m_2, \dots m_{max_a + 2}$ can encode several hyperarcs at the same time, the entry of the adjacency tensor $T_{OH}^+$ is set to:
	\begin{equation}
		\left(T_{OH}^+\right)_{m_1, m_2, \dots m_{max_a + 2}} = \sum_{a_q \in \mathcal{A}_H ~ \text{represented by} ~ m_1, m_2, \dots m_{max_a + 2}} \left(T_{OH}^+\right)_{m_1, m_2, \dots m_{max_a} + 2}^{a_q}.
	\end{equation}
	All entries of $T_{OH}^+$, which do not correspond to any hyperarc $a_q \in \mathcal{A}_H$, are set to $0$.\\
\end{definition}

\begin{example}[\textbf{Adjacency tensor with increased order}] \ \\
	Using the before defined oriented hypergraph $OH = \left(\mathcal{V}, \mathcal{A}_H\right)$ with the set of vertices $\mathcal{V} = \left\{v_1, v_2, v_3, v_4, v_5, v_6\right\}$, hyperarcs $\mathcal{A}_H = \left\{\left(\left\{v_1, v_2\right\}, \left\{v_4\right\}\right), \left(\left\{v_3\right\}, \left\{v_2, v_6\right\}\right), \left(\left\{v_5\right\}, \left\{v_6\right\}\right)\right\}$ and maximum hyperarc cardinality $max_a = 3$, then the nonzero indices in the adjacency tensor with increased order $T_{OH}^+$ are given by:
	
	$\begin{aligned}[t]	
		\qquad a_1 = \left(\left\{v_1, v_2\right\}, \left\{v_4\right\}\right): \quad & \left(1, 1, 1, 2, 4\right), \left(1, 1, 2, 1, 4\right), \left(1, 2, 1, 1, 4\right), \left(2, 1, 1, 1, 4\right), &\\
		& \left(1, 1, 2, 2, 4\right), \left(1, 2, 1, 2, 4\right), \left(1, 2, 2, 1, 4\right), \left(2, 1, 1, 2, 4\right), &\\
		& \left(2, 1, 2, 1, 4\right), \left(2, 2, 1, 1, 4\right), \left(1, 2, 2, 2, 4\right), \left(2, 1, 2, 2, 4\right), &\\
		& \left(2, 2, 1, 2, 4\right), \left(2, 2, 2, 1, 4\right), \left(1, 1, 2, 4, 4\right), \left(1, 2, 1, 4, 4\right), &\\
		& \left(2, 1, 1, 4, 4\right), \left(1, 2, 2, 4, 4\right), \left(2, 1, 2, 4, 4\right), \left(2, 2, 1, 4, 4\right), &\\
		& \left(1, 2, 4, 4, 4\right), \left(2, 1, 4, 4, 4\right) &\\
	\end{aligned}$\\
	$\begin{aligned}[t]	
		\qquad a_2 = \left(\left\{v_3\right\}, \left\{v_2, v_6\right\}\right): \quad & \left(3, 3, 3, 2, 6\right), \left(3, 3, 3, 6, 2\right), \left(3, 3, 2, 2, 6\right), \left(3, 3, 2, 6, 2\right), &\\
		& \left(3, 3, 6, 2, 2\right), \left(3, 3, 2, 6, 6\right), \left(3, 3, 6, 2, 6\right), \left(3, 3, 6, 6, 2\right), &\\
		& \left(3, 2, 2, 2, 6\right), \left(3, 2, 2, 6, 2\right), \left(3, 2, 6, 2, 2\right), \left(3, 6, 2, 2, 2\right), &\\
		& \left(3, 2, 2, 6, 6\right), \left(3, 2, 6, 2, 6\right), \left(3, 2, 6, 6, 2\right), \left(3, 6, 2, 2, 6\right), &\\
		& \left(3, 6, 2, 6, 2\right), \left(3, 6, 6, 2, 2\right), \left(3, 2, 6, 6, 6\right), \left(3, 6, 2, 6, 6\right), &\\
		& \left(3, 6, 6, 2, 6\right), \left(3, 6, 6, 6, 2\right) &\\
		\qquad a_3 = \left(\left\{v_5\right\}, \left\{v_6\right\}\right): \quad & \left(5, 5, 5, 5, 6\right), \left(5, 5, 5, 6, 6\right), \left(5, 5, 6, 6, 6\right), \left(5, 6, 6, 6, 6\right) &\\
	\end{aligned}$\\
\end{example}

The last example shows that for the adjacency tensor of oriented hypergraphs there is a trade-off between the uniqueness of the represented hypergraph and the size of the tensor. The following theorem proves that the adjacency tensor with increased order $T_{OH}^+$ uniquely describes oriented hypergraphs and that the underlying hypergraph can be identified by using the first case of the algorithm introduced in remark (\ref{TOHinter}).\\

\begin{theorem}[\textbf{Adjacency tensor $T_{OH}^+$ uniquely describes oriented hypergraphs}]\label{TOH+proof} \ \\
	Given an oriented hypergraph $OH = \left(\mathcal{V}, \mathcal{A}_H\right)$, then the corresponding adjacency tensor with increased order $T_{OH}^+$ uniquely describes this hypergraph and the algorithm described in remark (\ref{TOHinter}) returns the correct underlying hypergraph.\\
\end{theorem}

\begin{proof}\ \\
	Given an adjacency tensor with increased order $T_{OH}^+$, then the set of vertices $\mathcal{V}$ of the underlying oriented hypergraph is clearly identifiable based on the dimension of the tensor. Furthermore, the maximum hyperarc cardinality $max_a$ can also be read out from the order of the adjacency tensor. Assume that $\mathcal{V} = \left\{v_1, v_2, \dots v_N\right\}$ and $max_a \in \mathbb{N} \backslash \left\{0, 1\right\}$.\\
	
	By checking the sign of the entries of the adjacency tensor, all possible hyperarcs \big(based on the vertices $\mathcal{V}$\big) can be partitioned into a set of ``hyperarcs'' and a set of ``no hyperarcs'':\\
	
	\textbf{First case:} $m_1, m_2, \dots m_{max_a + 2}$ with $\left(T_{OH}^+\right)_{m_1, m_2, \dots m_{max_a + 2}} = 0$
	
	For such indices it holds true that partitions $\left(\left\{m_1, m_2, \dots m_s\right\}, \left\{m_{s + 1}, m_{s + 2}, \dots m_{max_a + 2}\right\}\right)$ with $s \in \left\{1, 2, \dots max_a + 1\right\}$ are definitely no hyperarc of the underlying hypergraph. In the case $\left\{m_1, m_2, \dots m_s\right\} \cap \left\{m_{s + 1}, m_{s + 2}, \dots m_{max_a + 2}\right\} = \emptyset$, the partition $\left(\left\{m_1, m_2, \dots m_s\right\},\right.$\\
	$\left.\left\{m_{s + 1}, m_{s + 2}, \dots m_{max_a + 2}\right\}\right)$ can thus be added to the set of ``no hyperarcs''.\\
	
	\clearpage
	\textbf{Second case:} $m_1, m_2, \dots m_{max_a + 2}$ with $\left(T_{OH}^+\right)_{m_1, m_2, \dots m_{max_a + 2}} > 0$
	
	Since every hyperarc $a_q \in \mathcal{A}_H$ consists of at most $max_a$ different vertices, the condition $\left\lvert\left\{m_1, m_2, \dots m_{max_a + 2}\right\}\right\rvert \leq max_a$ is fulfilled for all indices $m_1, m_2, \dots m_{max_a + 2}$ with $\left(T_{OH}^+\right)_{m_1, m_2, \dots m_{max_a + 2}} > 0$. Otherwise, $\left\lvert\left\{m_1, m_2, \dots m_{max_a + 2}\right\}\right\rvert > max_a$ would imply that there is a hyperarc $a_q$ with more than $max_a$ different vertices, since the different indices each correspond to one vertex. Therefore, there are always at least two redundant indices in $m_1, m_2, \dots m_{max_a + 2}$ \big(which can be substituted by $m_1$ and $m_{max_a + 2}$\big) if $\left(T_{OH}^+\right)_{m_1, m_2, \dots m_{max_a + 2}} > 0$ holds true.
	
	As described in the first case of remark (\ref{TOHinter}), hyperarcs with at most two less vertices than there are positions for the indices \big(in this case there are $max_a + 2$ positions\big) can be clearly identified using the second property of theorem (\ref{TOHprop}). Hence, any hyperarc $a_q = \left(\left\{i_1, i_2, \dots i_r\right\}, \left\{i_{r + 1}, i_{r + 2}, \dots i_n\right\}\right) \in \mathcal{A}_H$ with $1 \leq r \leq n - 1$ can be uniquely described by indices $m_1, m_2, \dots m_s, m_1, m_{max_a + 2}, m_{s + 1}, m_{s + 2}, \dots m_{max_a + 2}$ with $\left\{i_1, i_2, \dots i_r\right\} = \left\{m_1, m_2, \dots m_s\right\}$ and $\left\{i_{r + 1}, i_{r + 2}, \dots i_n\right\} = \left\{m_{s + 1}, m_{s + 2}, \dots m_{max_a + 2}\right\}$.
	
	Thus, it suffices to check $\left(T_{OH}^+\right)_{m_1, m_2, \dots m_s, m_1, m_{max_a + 2}, m_{s + 1}, m_{s + 2}, \dots m_{max_a + 2}} > 0$ for all partitions $\left(\left\{m_1, m_2, \dots m_s\right\}, \left\{m_{s + 1}, m_{s + 2}, \dots m_{max_a + 2}\right\}\right)$ of the indices $m_1, m_2, \dots m_{max_a + 2}$ to differentiate between sure hyperarcs and no hyperarcs and find all hyperarcs of the underlying oriented hypergraph.\\
	
	The first case of the algorithm introduced in remark (\ref{TOHinter}) completes exactly the two necessary steps:
	\begin{itemize}
		\item Remove two of the redundant indices of $m_1, m_2, \dots m_{max_a + 2}$ with\\ $\left(T_{OH}^+\right)_{m_1, m_2, \dots m_{max_a + 2}} > 0$
		\item Check for all partitions $\left(\left\{m_1, m_2, \dots m_s\right\}, \left\{m_{s + 1}, m_{s + 2}, \dots m_{max_a + 2}\right\}\right)$ if\\ $\left(T_{OH}^+\right)_{m_1, m_2, \dots m_s, m_1, m_{max_a + 2}, m_{s + 1}, m_{s + 2}, \dots m_{max_a + 2}} > 0$
	\end{itemize}
	The second case of the algorithm is not needed, because the cases $\left(T_{OH}^+\right)_{m_1, m_2, \dots m_{max_a + 2}} > 0$ with $\left\lvert\left\{m_1, m_2, \dots m_{max_a + 2}\right\}\right\rvert = max_a + 1$ or $\left\lvert\left\{m_1, m_2, \dots m_{max_a + 2}\right\}\right\rvert = max_a + 2$ cannot occur in the adjacency tensor with increased order. Therefore, the algorithm can be applied to $T_{OH}^+$ and it will always find the correct underlying oriented hypergraph.\\
\end{proof}
\clearpage
\section{Additional definitions on normal graphs and hypergraphs}\label{3}

In this section additional definitions on normal graphs and hypergraphs are introduced, which are necessary for later definitions of gradients, adjoints, divergences, Laplacians and $p$-Laplacians.\\

Subsection (\ref{3.1}) introduces bipartite normal graphs, complete not oriented normal graphs, vertex degree functions \big($\deg$, $\deg_{out}$, and $\deg_{in}$\big) on normal graphs, the arc indicator function \big($\delta$\big), and vertex-arc indicator functions \big($\delta_{out}$ and $\delta_{in}$\big). Moreover, subsection (\ref{3.2}) gives definitions for the maximum hyperedge and hyperarc cardinality \big($max_e$ and $max_a$\big) as well as the minimum hyperedge and hyperarc cardinality \big($min_e$ and $min_a$\big), the vertex degree functions \big($\deg$, $\deg_{out}$, and $\deg_{in}$\big) on hypergraphs, and the vertex-hyperarc indicator functions \big($\delta_{out}$ and $\delta_{in}$\big).\\

\subsection{Additional definitions on normal graphs}\label{3.1}

The definitions of this subsection are based on the work of \cite{asratian1998bipartite}, \cite{augustson1970analysis} and \cite{bang2018basic}.\\

\begin{definition}[\textbf{Bipartite normal graph}]\label{bipartiteG} (Page 7 in \cite{asratian1998bipartite}: Bipartite graph) \ \\
	A not oriented normal graph $NG = \left(\mathcal{V}, \mathcal{E}_G\right)$ is called bipartite if the vertex set $\mathcal{V}$ can be partitioned into two subsets $X$ and $Y$, such that each edge connects a vertex from the first subset $X$ with a vertex from the second subset $Y$:
	\begin{equation}
		\mathcal{V} = X \cup Y ~ \text{with} ~ X \cap Y = \emptyset
	\end{equation}
	\begin{equation}
		\forall e_q \in \mathcal{E}_G: \quad e_q = \left\{v_i, v_j\right\} ~ \text{with} ~ v_i \in X ~ \text{and} ~ v_j \in Y .
	\end{equation}
	Similarly, an oriented normal graph $OG = \left(\mathcal{V}, \mathcal{A}_G\right)$ is called bipartite if the vertex set $\mathcal{V}$ can be partitioned as explained before, such that each arc connects a vertex from the first subset $X$ with a vertex from the second subset $Y$:
	\begin{equation}
		\mathcal{V} = X \cup Y ~ \text{with} ~ X \cap Y = \emptyset
	\end{equation}
	\begin{equation}
		\forall a_q \in \mathcal{A}_G: \quad a_q = \left(v_i, v_j\right) ~ \text{with} ~ v_i \in X, v_j \in Y ~ \text{or with} ~ v_i \in Y, v_j \in X .
	\end{equation}\vspace{0em}
\end{definition}

\begin{definition}[\textbf{Complete not oriented normal graph}]\label{completeNG} (Page 21 in \cite{bang2018basic}: Complete multigraph) \ \\
	A not oriented normal graph $NG = \left(\mathcal{V}, \mathcal{E}_G\right)$ is called complete if each vertex $v_i \in \mathcal{V}$ is connected to each other vertex $v_j \in \mathcal{V} \backslash \left\{v_i\right\}$ through exactly one edge $e_q \in \mathcal{E}_G$:
	\begin{equation}
		\mathcal{E}_G = \left\{\left\{v_i, v_j\right\} ~ \middle| ~ v_i \in \mathcal{V} ~ \text{and} ~ v_j \in \mathcal{V} \backslash \left\{v_i\right\}\right\}.
	\end{equation}\vspace{0em}
\end{definition}

\begin{definition}[\textbf{Connected components of normal graphs}]\label{connectedG} (Page 572 in \cite{augustson1970analysis}: Connected components) \ \\
	In a not oriented normal graph $NG = \left(\mathcal{V}, \mathcal{E}_G\right)$, two vertices $v_i, v_j \in \mathcal{V}$ are called connected if there exists a subset of edges in $\mathcal{E}_G$ such that it is possible to move along those edges from vertex $v_i$ to vertex $v_j$ and vice versa. Similarly, in an oriented normal graph $OG = \left(\mathcal{V}, \mathcal{A}_G\right)$, two vertices $v_i, v_j \in \mathcal{V}$ are called connected, if there exists a subset of arcs in $\mathcal{A}_G$ such that it is possible to move along those arcs \big(while ignoring their direction and hence treating them as edges\big) from one vertex to the other.\\
	
	A maximal connected subgraph of $NG$ or $OG$ is then called a connected component of the normal graph. This indicates that each vertex $v_i \in \mathcal{V}$, each edge $e_q \in \mathcal{E}_G$ and each arc $a_q \in \mathcal{A}_G$ always belong to exactly one connected component.\\
\end{definition}

\begin{definition}[\textbf{Vertex degree functions $\deg$, $\deg_{out}$, $\deg_{in}$}]\label{degG} (Page 5 in \cite{bang2018basic}: Out-degree, in-degree and degree) \ \\
	For a not oriented normal graph $NG = \left(\mathcal{V}, \mathcal{E}_G\right)$, only the general degree function $\deg$ is defined as
	\begin{equation*}
		\deg: ~ \mathcal{V} \longrightarrow \mathbb{Z}_{\geq 0} \qquad v_i \longmapsto \deg\left(v_i\right)
	\end{equation*}
	with the degree of vertex $v_i \in \mathcal{V}$ being set to
	\begin{equation}
		\deg\left(v_i\right) = \left\lvert \left\{e_q \in \mathcal{E}_G ~ \middle| ~ v_i \in e_q\right\}\right\rvert .
	\end{equation}
	Thus, the general degree function $\deg$ gives for every vertex $v_i \in \mathcal{V}$ the number of edges $e_q \in \mathcal{E}_G$ containing vertex $v_i$.\\
	
	For an oriented normal graph $OG = \left(\mathcal{V}, \mathcal{A}_G\right)$, the output degree function $\deg_{out}$ and the input degree function $\deg_{in}$ are defined as
	\begin{equation*}
		\deg_{out}: ~ \mathcal{V} \longrightarrow \mathbb{Z}_{\geq 0}
	\end{equation*}
	\begin{equation}
		v_i \longmapsto \deg_{out}\left(v_i\right) = \left\lvert \left\{a_q \in \mathcal{A}_G ~ \middle| ~ a_q = \left(v_i, v_j\right) ~ \text{for some} ~ v_j \in \mathcal{V}\backslash\left\{v_i\right\} \right\} \right\rvert
	\end{equation}
	and respectively
	\begin{equation*}
		\deg_{in}: ~ \mathcal{V} \longrightarrow \mathbb{Z}_{\geq 0}
	\end{equation*}
	\begin{equation}
		v_i \longmapsto \deg_{in}\left(v_i\right) = \left\lvert \left\{a_q \in \mathcal{A}_G ~ \middle| ~ a_q = \left(v_j, v_i\right) ~ \text{for some} ~ v_j \in \mathcal{V}\backslash\left\{v_i\right\} \right\} \right\rvert .
	\end{equation}

	Therefore, the output degree function $\deg_{out}$ returns the number of arcs $a_q \in \mathcal{A}_G$, which contain $v_i$  as an output vertex, and similarly the input degree function $\deg_{in}$ returns the number of arcs $a_q \in \mathcal{A}_G$, for which $v_i$ is an input vertex.
	
	\clearpage
	The general degree function $\deg$ for oriented normal graphs $OG = \left(\mathcal{V}, \mathcal{A}_G\right)$ can be defined using the output and input degree functions, which results in
	\begin{equation*}
		\deg: ~ \mathcal{V} \longrightarrow \mathbb{Z}_{\geq 0} \qquad v_i \longmapsto \deg\left(v_i\right)
	\end{equation*}
	with the degree of vertex $v_i \in \mathcal{V}$ being set to
	\begin{equation}
		\deg\left(v_i\right) = \deg_{out}\left(v_i\right) + \deg_{in}\left(v_i\right).
	\end{equation}\vspace{0em}
\end{definition}

\begin{remark}[\textbf{General degree function $\deg$ for oriented normal graphs counts connected arcs}]\label{degOG} \ \\
	The set of arcs $\mathcal{A}_G$ in an oriented normal graph $OG = \left(\mathcal{V}, \mathcal{A}_G\right)$ is defined as
	\begin{equation*}
		\mathcal{A}_G \subseteq \left\{\left(v_i, v_j\right) ~ \middle| ~ v_i, v_j \in \mathcal{V}, v_i \neq v_j\right\}
	\end{equation*}
	and hence it is impossible for one arc $a_q \in \mathcal{A}_G$ to be counted twice in $\deg\left(v_i\right)$ \big(once in $\deg_{out}\left(v_i\right)$ and once in $\deg_{in}\left(v_i\right)$\big) because any vertex $v_i \in \mathcal{V}$ cannot be the output vertex and the input vertex of $a_q$ at the same time. Therefore, it holds true that
	\begin{equation}
		\deg\left(v_i\right) = \left\lvert\left\{a_q \in \mathcal{A}_G ~ \middle| ~ a_q = \left(v_i, v_j\right) ~ \text{or} ~ a_q = \left(v_j, v_i\right) ~ \text{for} ~ v_j \in  \mathcal{V}\right\}\right\rvert ,
	\end{equation}
	where $v_j \neq v_i$ is automatically fulfilled. This means that $\deg\left(v_i\right)$ counts the arcs $a_q \in \mathcal{A}_G$, which contain vertex $v_i \in \mathcal{V}$.\\
\end{remark}

\begin{definition}[\textbf{Arc indicator function $\delta$}]\label{arcindiOG} \ \\
	For oriented normal graphs $OG = \left(\mathcal{V}, \mathcal{A}_G\right)$, the arc indicator function $\delta$ indicates for two vertices $v_i, v_j \in \mathcal{V}$ if it holds true that $a_q = \left(v_i, v_j\right) \in \mathcal{A}_G$ or if $a_q = \left(v_i, v_j\right) \notin \mathcal{A}_G$. Thus, the arc indicator function is defined as:
	\begin{equation}
		\delta: ~ \mathcal{V} \times \mathcal{V} \longrightarrow \left\{0, 1\right\} \qquad \left(v_i, v_j\right) \longmapsto \delta\left(v_i, v_j\right) = \left\{\begin{array}{ll}
			1 & \quad a_q = \left(v_i, v_j\right) \in \mathcal{A}_G\\
			0 & \quad \text{otherwise}
		\end{array}\right.
		.
	\end{equation}\vspace{0em}
\end{definition}

\begin{definition}[\textbf{Vertex-arc indicator functions $\delta_{out}$, $\delta_{in}$}]\label{vertexarcindiOG} \ \\
	For oriented normal graphs $OG = \left(\mathcal{V}, \mathcal{A}_G\right)$, functions indicating whether a vertex $v_i \in \mathcal{V}$ is an output or an input vertex of an arc $a_q \in \mathcal{A}_G$, can be defined. The output vertex-arc indicator function $\delta_{out}$ is given by
	\begin{equation*}
		\delta_{out}: ~ \mathcal{V} \times \mathcal{A}_G \longrightarrow \left\{0, 1\right\}
	\end{equation*}
	\begin{equation}
		\left(v_i, a_q\right) \longmapsto \delta_{out}\left(v_i, a_q\right) = \left\{\begin{array}{ll}
			1 & \quad \left(v_i, v_j\right) = a_q ~ \text{for some} ~ v_j \in \mathcal{V}\\
			0 & \quad \text{otherwise}
		\end{array}\right.
	\end{equation}
	and the input vertex-arc indicator function $\delta_{in}$ is given by:
	\begin{equation*}
		\delta_{in}: ~ \mathcal{V} \times \mathcal{A}_G \longrightarrow \left\{0, 1\right\}
	\end{equation*}
	\begin{equation}
		\left(v_i, a_q\right) \longmapsto \delta_{in}\left(v_i, a_q\right) = \left\{\begin{array}{ll}
			1 & \quad \left(v_j, v_i\right) = a_q ~ \text{for some} ~ v_j \in \mathcal{V}\\
			0 & \quad \text{otherwise}
		\end{array}\right..
	\end{equation}\vspace{0em}
\end{definition}
\subsection{Additional definitions on hypergraphs}\label{3.2}

Part of the content of this subsection is based on \cite{zhang2019introducing} and the remaining definitions generalize the normal graph case from the previous subsection.\\

\begin{definition}[\textbf{Maximum hyperedge cardinality $max_e$ and maximum hyperarc cardinality $max_a$}] \label{maxea} (Page 641 in \cite{zhang2019introducing}: Maximum cardinality of hyperedges $m.c.e\left(\mathcal{H}\right)$) \ \\
	For a not oriented hypergraph $NH = \left(\mathcal{V}, \mathcal{E}_H\right)$, the maximum hyperedge cardinality $max_e$ is defined as
	\begin{equation}
		max_e = \max \left\{\lvert e_q\rvert ~ \middle| ~ e_q \in \mathcal{E}_H\right\}.
	\end{equation}

	In case of an oriented hypergraph $OH = \left(\mathcal{V}, \mathcal{A}_H\right)$, the maximum hyperedge cardinality $max_e$ is then called maximum hyperarc cardinality $max_a$ with
	\begin{equation}
		max_a = \max \left\{\lvert a_q^{out} \rvert + \lvert a_q^{in} \rvert ~ \middle| ~ a_q \in \mathcal{A}_H\right\}.
	\end{equation}\vspace{0em}	
\end{definition}

\begin{definition}[\textbf{Minimum hyperedge cardinality $min_e$ and minimum hyperarc cardinality $min_a$}]\label{minea} \ \\
	Analogously to the maximum hyperedge cardinality $max_e$, for a not oriented hypergraph $NH = \left(\mathcal{V}, \mathcal{E}_H\right)$, the minimum hyperedge cardinality $min_e$ is defined as
	\begin{equation}
		min_e = \min \left\{\lvert e_q\rvert ~ \middle| ~ e_q \in \mathcal{E}_H\right\}.
	\end{equation}
	In case of an oriented hypergraph $OH = \left(\mathcal{V}, \mathcal{A}_H\right)$, the minimum hyperarc cardinality $min_a$ is defined as
	\begin{equation}
		min_a = \min \left\{\lvert a_q^{out} \rvert + \lvert a_q^{in} \rvert ~ \middle| ~ a_q \in \mathcal{A}_H\right\}.
	\end{equation}\vspace{0em}
\end{definition}

\begin{remark}[\textbf{Properties of the minimum and maximum hyperedge and hyperarc cardinality}]\label{maxminea} \ \\
	Given a not oriented hypergraph $NH = \left(\mathcal{V}, \mathcal{E}_H\right)$, then for the minimum hyperedge cardinality $min_e$ and the maximum hyperedge cardinality $max_e$, it holds true that
	\begin{equation}
		2 \leq min_e \leq max_e \leq \left\lvert\mathcal{V}\right\rvert = N,
	\end{equation}
	since all hyperedges $e_q \in \mathcal{E}_H$ fulfill the condition $\left\lvert e_q\right\rvert \geq 2$ and $e_q \subseteq \mathcal{V}$.\\
	
	Given an oriented hypergraph $OH = \left(\mathcal{V}, \mathcal{A}_H\right)$, then the minimum hyperarc cardinality $min_a$ and the maximum hyperarc cardinality $max_a$ fulfill the inequalities
	\begin{equation}
		2 \leq min_a \leq max_a \leq \left\lvert\mathcal{V}\right\rvert = N,
	\end{equation}
	since for all feasible hyperarcs $a_q \in \mathcal{A}_H$ it holds true that $\emptyset \subset a_q^{out}, a_q^{in} \subset \mathcal{V}$ and $a_q^{out} \cap a_q^{in} = \emptyset$.\\
\end{remark}

\begin{definition}[\textbf{Vertex degree functions $\deg$, $\deg_{out}$, $\deg_{in}$}]\label{degH} \ \\
	For a not oriented hypergraph $NH = \left(\mathcal{V}, \mathcal{E}_H\right)$, only the general degree function $\deg$ is defined
	\begin{equation*}
		\deg: ~ \mathcal{V} \longrightarrow \mathbb{Z}_{\geq 0} \qquad v_i \longmapsto \deg\left(v_i\right)
	\end{equation*}
	with the degree of vertex $v_i \in \mathcal{V}$ being set to
	\begin{equation}
		\deg\left(v_i\right) = \left\lvert \left\{e_q \in \mathcal{E}_H ~ \middle| ~ v_i \in e_q\right\}\right\rvert .
	\end{equation}
	As seen before, the general degree function $\deg$ returns, for every vertex $v_i \in \mathcal{V}$, the number of hyperedges $e_q \in \mathcal{E}_H$ containing vertex $v_i$.\\
	
	For an oriented hypergraph $OH = \left(\mathcal{V}, \mathcal{A}_H\right)$, the output degree function $\deg_{out}$ and the input degree function $\deg_{in}$ are defined as
	\begin{equation*}
		\deg_{out}: ~ \mathcal{V} \longrightarrow \mathbb{Z}_{\geq 0}
	\end{equation*}
	\begin{equation}
		v_i \longmapsto \deg_{out}\left(v_i\right) = \left\lvert \left\{a_q \in \mathcal{A}_H ~ \middle| ~ v_i \in a_q^{out} \right\} \right\rvert
	\end{equation}
	and respectively
	\begin{equation*}
		\deg_{in}: ~ \mathcal{V} \longrightarrow \mathbb{Z}_{\geq 0}
	\end{equation*}
	\begin{equation}
		v_i \longmapsto \deg_{in}\left(v_i\right) = \left\lvert \left\{a_q \in \mathcal{A}_H ~ \middle| ~ v_i \in a_q^{in} \right\} \right\rvert .
	\end{equation}
	Hence, the output degree function $\deg_{out}$ counts the number of hyperarcs $a_q \in \mathcal{A}_H$, which have output vertex $v_i$, and similarly the input degree function $\deg_{in}$ counts the number of hyperarcs $a_q \in \mathcal{A}_H$, which have input vertex $v_i$.
	
	The general degree function $\deg$ for oriented hypergraphs $OH = \left(\mathcal{V}, \mathcal{A}_H\right)$ can be defined using the output and input degree functions, which leads to
	\begin{equation*}
		\deg: ~ \mathcal{V} \longrightarrow \mathbb{Z}_{\geq 0} \qquad v_i \longmapsto \deg\left(v_i\right)
	\end{equation*}
	with the degree of vertex $v_i \in \mathcal{V}$ being set to
	\begin{equation}
		\deg\left(v_i\right) = \deg_{out}\left(v_i\right) + \deg_{in}\left(v_i\right).
	\end{equation}\vspace{0em}
\end{definition}

\clearpage
\begin{remark}[\textbf{General degree function $\deg$ for oriented hypergraphs counts connected hyperarcs}]\label{degOH} \ \\
	Since all hyperarcs $a_q \in \mathcal{A}_H$ in an oriented hypergraph $OH = \left(\mathcal{V}, \mathcal{A}_H\right)$ fulfill the condition
	\begin{equation*}
		a_q^{out} \cap a_q^{in} = \emptyset,
	\end{equation*}
	it is impossible for one hyperarc $a_q \in \mathcal{A}_H$ to be counted twice in $\deg\left(v_i\right)$ \big(once in $\deg_{out}\left(v_i\right)$ and once in $\deg_{in}\left(v_i\right)$\big), because any vertex $v_i \in \mathcal{V}$ can not be the output vertex and the input vertex of $a_q$ at the same time. Thus, it holds true that
	\begin{equation}
		\deg\left(v_i\right) = \left\lvert\left\{a_q \in \mathcal{A}_H ~ \middle| ~ v_i \in a_q^{out} ~ \text{or} ~ v_i \in a_q^{in}\right\}\right\rvert ,
	\end{equation}
	which means that the general degree function $\deg\left(v_i\right)$ counts the hyperarcs $a_q \in \mathcal{A}_H$, which contain vertex $v_i \in \mathcal{V}$.\\
\end{remark}

\begin{definition}[\textbf{Vertex-hyperarc indicator functions $\delta_{out}$, $\delta_{in}$}]\label{vertexarcindiOH} \ \\
	For oriented hypergraphs $OH = \left(\mathcal{V}, \mathcal{A}_H\right)$, functions indicating whether a vertex $v_i \in \mathcal{V}$ is an output or input vertex of a hyperarc $a_q \in \mathcal{A}_H$, can be defined. The output vertex-hyperarc indicator function $\delta_{out}$ is given by
	\begin{equation*}
		\delta_{out}: ~ \mathcal{V} \times \mathcal{A}_H \longrightarrow \left\{0, 1\right\}
	\end{equation*}	
	\begin{equation}
		\left(v_i, a_q\right) \longmapsto \delta_{out}\left(v_i, a_q\right) = \left\{\begin{array}{ll}
			1 & \quad v_i \in a_q^{out}\\
			0 & \quad \text{otherwise}
		\end{array}\right.
	\end{equation}
	and the input vertex-hyperarc indicator function $\delta_{in}$ is given by
	\begin{equation*}
		\delta_{in}: ~ \mathcal{V} \times \mathcal{A}_H \longrightarrow \left\{0, 1\right\}
	\end{equation*}
	\begin{equation}
		\left(v_i, a_q\right) \longmapsto \delta_{in}\left(v_i, a_q\right) = \left\{\begin{array}{ll}
			1 & \quad v_i \in a_q^{in}\\
			0 & \quad \text{otherwise}
		\end{array}\right..
	\end{equation}\vspace{0em}
\end{definition}
\clearpage
\section{Connections between normal graphs and hypergraphs}\label{4}

In order to further highlight why hypergraphs are a natural extension of normal graphs, subsection (\ref{4.1}) shows that every normal graph is a hypergraph and that applying the previously defined adjacency tensors $T_{NH}$ or $T_{OH}$ to normal graphs results in the adjacency matrices $A_{NG}$ or $A_{OG}$.\\

Moreover, subsection (\ref{4.2}) introduces two different normal graph representations of hypergraphs and it also analyzes the benefits and drawbacks of each representation. This section also includes examples for the two normal graph representations, which show that mappings of different hypergraphs to their normal graph representations can result in the same normal graph. Therefore, it is generally not possible to reconstruct the underlying hypergraph from its normal graph representation in a unique way without additional information.\\

\subsection{Every normal graph is a hypergraph}\label{4.1}

This subsection further explores the work of \cite{zhang2019introducing} and extends it to the oriented hypergraph case.\\

The following theorem proofs that any normal graph is ``just'' a special case of a hypergraph.\\

\begin{theorem}[\textbf{Every normal graph is a hypergraph}]\label{GaH} \ \\
	Normal graphs and hypergraphs are linked by the following connections:
	\begin{itemize}
		\item[1)] Every not oriented normal graph $NG = \left(\mathcal{V}, \mathcal{E}_G\right)$ is a not oriented hypergraph $NH = \left(\mathcal{V}, \mathcal{E}_H\right)$, where for every hyperedge $e_q \in \mathcal{E}_H$ it holds true that $\lvert e_q\rvert = 2 = max_e$.
		\item[2)] Every oriented normal graph $OG = \left(\mathcal{V}, \mathcal{A}_G\right)$ is an oriented hypergraph $OH = \left(\mathcal{V}, \mathcal{A}_H\right)$, where fore every hyperarc $a_q \in \mathcal{A}_H$ it holds true that $\lvert a_q\rvert = 2 = max_a$.\\
	\end{itemize}
\end{theorem}

\begin{proof}\ \\
	In order to show that every normal graph is a hypergraph, the vertex set $\mathcal{V}$, the hyperedge set $\mathcal{E}_H$ and the hyperarc set $\mathcal{A}_H$ are defined:
	\begin{itemize}
		\item[1)] The vertex set $\mathcal{V}$ of the not oriented hypergraph can be chosen the same as the vertex set of the not oriented normal graph. The hyperedges $e_q$ of $\mathcal{E}_H$ are then defined as:
		\begin{equation*}
			\mathcal{E}_H = \left\{e_q = \left\{v_i, v_j\right\} ~ \middle| ~ \left\{v_i, v_j\right\} \in \mathcal{E}_G\right\}.
		\end{equation*}
		Due to the definition of the set of hyperedges $\mathcal{E}_H$, for all hyperedges $e_q \in \mathcal{E}_H$ it holds true that $\lvert e_q\rvert = 2$, which automatically implies that the maximum hyperedge cardinality of the not oriented hypergraph satisfies $max_e = 2$.
		
		\item[2)] The vertex set $\mathcal{V}$ of the oriented hypergraph can be chosen to be the same as the vertex set of the oriented normal graph. The hyperarcs $a_q = a_q^{out} \cup a_q^{in}$ of $\mathcal{A}_H$ are then defined as:
		\begin{equation*}
			\mathcal{A}_H = \left\{a_q =  \left(a_q^{out}, a_q^{in}\right) = \left(\left\{v_i\right\}, \left\{v_j\right\}\right) ~ \middle| ~ \left(v_i, v_j\right) \in \mathcal{A}_G\right\}.
		\end{equation*}
		Based on this definition of the set of hyperarcs, all hyperarcs $a_q \in \mathcal{A}_H$ fulfill the equality $\lvert a_q^{out} \rvert + \lvert a_q^{in} \rvert = 1 + 1 = 2$ and thus the maximum hyperarc cardinality of the oriented hypergraph is given as $max_a = 2$.
	\end{itemize}
\end{proof}

Due to the previous theorem, the algebraic representations of hypergraphs as adjacency tensors can now also be applied to normal graphs.\\

\begin{theorem}[\textbf{Equivalence of adjacency tensor $T_{NG}$ and adjacency matrix $A_{NG}$ for not oriented normal graphs}]\label{TNG} (Page 641 in \cite{zhang2019introducing})\ \\
	Applying the adjacency tensor for not oriented hypergraphs $T_{NG}$ from definition (\ref{TNH}) to a not oriented normal graph $NG = \left(\mathcal{V}, \mathcal{E}_G\right)$ results in the adjacency matrix $A_{NG}$ for not oriented normal graphs of definition (\ref{ANG}).\\
\end{theorem}

\begin{proof}\ \\
	Since every not oriented normal graph $NG = \left(\mathcal{V}, \mathcal{E}_G\right)$ is a not oriented hypergraph with the property $\lvert e_q\rvert = 2 = max_e$ for all hyperedges $e_q \in \mathcal{E}_H$, the corresponding adjacency tensor $T_{NG}$ is $N$-dimensional and of $2^{\text{nd}}$-order. Hence, $T_{NG}$ can be written as an $N \times N$ matrix.
	
	Furthermore, for every edge $\left\{v_i, v_j\right\} \in \mathcal{E}_G$ it holds true that it corresponds to the hyperedge $e_q = \left\{v_i, v_j\right\} \in \mathcal{E}_H$ with $\lvert e_q\rvert = n = 2$ in the corresponding hypergraph of $NG$. Using the definition of the adjacency tensor yields the following:
	\begin{align*}
		\left(T_{NG}\right)_{i j} = \left(T_{NG}\right)_{j i} & = \frac{2}{\sum_{l_1, l_2 \in \mathbb{N} \backslash \left\{0\right\}: ~ \sum_{k = 1}^{2} l_k = 2} \frac{2!}{l_1! ~ l_2!}} = \frac{2}{\sum_{l_1 = 1, l_2 = 1} \frac{2!}{l_1! ~ l_2!}} = \frac{2}{\frac{2!}{1! ~ 1!}} = 1\\
		& = \left(A_{NG}\right)_{i j} = \left(A_{NG}\right)_{j i}
	\end{align*}
	where the last equality is given by the symmetry of $A_{NG}$. \big(Note: The equalities $\left(T_{NG}\right)_{i j} = \left(T_{NG}\right)_{j i} = 1 = \left(A_{NG}\right)_{i j} = \left(A_{NG}\right)_{j i}$ also hold true if the simplified adjacency tensor definition $\left(T_{NG}\right)_{i j} = \left(T_{NG}\right)_{j i} = \frac{n}{2}$ is used.\big)
	
	\clearpage
	Since all other entries \big(entries not corresponding to edges or hyperedges\big) of the adjacency matrix $A_{NG}$ and the adjacency tensor $T_{NG}$ are set to zero by definition, the equality $T_{NG} = A_{NG}$ holds true.\\
\end{proof}

\begin{theorem}[\textbf{Equivalence of adjacency tensor $T_{OG}$ and adjacency matrix $A_{OG}$ for oriented normal graphs}]\label{TOG} \ \\	
	Applying the adjacency tensor for oriented hypergraphs $T_{OG}$ from definition (\ref{TOH}) to an oriented normal graph $OG = \left(\mathcal{V}, \mathcal{A}_G\right)$ results in the adjacency matrix $A_{OG}$ for oriented normal graphs of definition (\ref{AOG}).\\
\end{theorem}

\begin{proof}\ \\
	As argued before, every oriented normal graph $OG = \left(\mathcal{V}, \mathcal{A}_G\right)$ is an oriented hypergraph with the property $\lvert a_q\rvert = 2 = max_a$ for all hyperarcs $a_q \in \mathcal{A}_H$. Thus, the corresponding adjacency tensor $T_{OG}$ is $N$-dimensional and of $2^{\text{nd}}$-order and can be written as an $N \times N$ matrix.
	
	Moreover, for every arc $\left(v_i, v_j\right) \in \mathcal{A}_G$ it holds true that it corresponds to the hyperarc $a_q = \left(v_i, v_j\right) \in \mathcal{A}_H$ with $\lvert a_q^{out} \cup a_q^{in}\rvert = n = 2$ in the corresponding hypergraph of $OG$. Using the definition of the adjacency tensor yields the following:
	\begin{align*}
		\left(T_{OG}\right)_{i j} & = \frac{2}{\sum_{l_1, l_2 \in \mathbb{N} \backslash \left\{0\right\}: ~ \sum_{k = 1}^{2} l_k = 2} \frac{2!}{l_1! ~ l_2!}} = \frac{2}{\sum_{l_1 = 1, l_2 = 1} \frac{2!}{l_1! ~ l_2!}} = \frac{2}{\frac{2!}{1! ~ 1!}} = 1\\
		& = \left(A_{OG}\right)_{i j} .
	\end{align*}
	\big(Note: The equalities $\left(T_{OG}\right)_{i j} = 1 = \left(A_{OG}\right)_{i j}$ also hold true if the simplified adjacency tensor definition $\left(T_{OG}\right)_{i j} = \frac{n}{2}$ is used.\big)
	
	Since all other entries \big(entries not corresponding to arcs or hyperarcs\big) of the adjacency matrix $A_{OG}$ and the adjacency tensor $T_{OG}$ are set to zero by definition, the equality $T_{OG} = A_{OG}$ holds true.\\
\end{proof}

Note: Applying the adjacency tensor with increased order $T_{OG}^+$ from definition (\ref{T_OHincr}) to an oriented normal graph $OG = \left(\mathcal{V}, \mathcal{A}_G\right)$ results in an $N \times N \times N \times N$ matrix and hence does not result in the adjacency matrix $A_{OG}$.\\
\clearpage
\subsection{Hypergraphs represented as normal graphs}\label{4.2}

The following two normal graph representations of hypergraphs underline why hypergraphs are an important generalization of normal graphs and why it is generally impossible to capture all properties of hypergraphs effectively using normal graph representations. The first normal graph representation is briefly discussed in \cite{mori2015peeling} and the second normal graph representation was mentioned in the introduction of \cite{zhou2006learning}.\\

\begin{definition}[\textbf{First normal graph representation of not oriented hypergraphs}]\label{NNGR1} (Page 3 in \cite{mori2015peeling}: Bipartite graph representation of hypergraphs)\ \\
	Given a not oriented hypergraph $NH = \left(\mathcal{V}_H, \mathcal{E}_H\right)$ with the vertex set $\mathcal{V}_H = \left\{v_1, v_2, \dots v_N\right\}$, then a normal graph $NG = \left(\mathcal{V}_G, \mathcal{E}_G\right)$ encoding the hypergraph can be defined with
	\begin{equation}
		\mathcal{V}_G = \mathcal{V}_H \cup \left\{v_{N + q} ~ \middle| ~ e_q \in \mathcal{E}_H\right\},
	\end{equation}
	\begin{equation}
		\mathcal{E}_G = \left\{\left\{v_i, v_{N + q}\right\} ~ \middle| ~ v_i \in e_q ~ \text{for} ~ e_q \in \mathcal{E}_H\right\},
	\end{equation}
	where for each hyperedge $e_q \in \mathcal{E}_H$ a new vertex $v_{N + q}$ is introduced. The edges then connect each vertex $v_{N + q}$ representing a hyperedge $e_q$ with the original vertices $v_i \in \mathcal{V}_H$, which were part of the hyperedge.\\
\end{definition}

\begin{example}[\textbf{First normal graph representation of not oriented hypergraphs}]\ \\
	Given the not oriented hypergraph $NH = \left(\mathcal{V}_H, \mathcal{E}_H\right)$ with $\mathcal{V}_H = \left\{v_1, v_2, v_3, v_4, v_5, v_6, v_7, v_8\right\}$ and $\mathcal{E}_H = \left\{\left\{v_1, v_2, v_5\right\}, \left\{v_2, v_3, v_7, v_8\right\}, \left\{v_6, v_7\right\}\right\}$, then the not oriented normal graph $NG = \left(\mathcal{V}_G, \mathcal{E}_G\right)$ representing this hypergraph is given by:
	\begin{equation*}
		\mathcal{V}_G = \left\{v_1, v_2, v_3, v_4, v_5, v_6, v_7, v_8\right\} \cup \left\{v_9, v_{10}, v_{11}\right\}
	\end{equation*}
	\begin{equation*}
		\mathcal{E}_G = \left\{\left\{v_1, v_9\right\}, \left\{v_2, v_9\right\}, \left\{v_5, v_9\right\}, \left\{v_2, v_{10}\right\}, \left\{v_3, v_{10}\right\},\right.
	\end{equation*}
	\begin{equation*}
		 \left.\left\{v_7, v_{10}\right\}, \left\{v_8, v_{10}\right\}, \left\{v_6, v_{11}\right\}, \left\{v_7, v_{11}\right\}\right\}
	\end{equation*}
	and can be visualized in the following way:
	
	\begin{minipage}{.5\textwidth}\begin{tikzpicture}
				\tikzstyle{vertex} = [fill,shape=circle,node distance=80pt]
				\tikzstyle{edge} = [fill,opacity=.5,fill opacity=.5,line cap=round, line join=round, line width=30pt]
				\tikzstyle{elabel} =  [fill,shape=circle,node distance=40pt]
				
				\pgfdeclarelayer{background}
				\pgfsetlayers{background,main}
				
				\begin{scope}[every node/.style={circle,thick,draw}]
					\node (v1) at (0,0) {$v_1$};
					\node (v2) at (2,0) {$v_2$};
					\node (v3) at (4,0) {$v_3$};
					\node (v4) at (6, 0) {$v_4$};
					\node (v5) at (0,-2) {$v_5$};
					\node (v6) at (2,-2) {$v_6$};
					\node (v7) at (4,-2) {$v_7$};
					\node (v8) at (6,-2) {$v_8$};
				\end{scope}
				
				\begin{pgfonlayer}{background}
					\begin{scope}[transparency group,opacity=.5]
						\draw[edge,opacity=1,color=blue] (v1.center) -- (v2.center) -- (v5.center) -- (v1.center);
						\fill[edge,opacity=1,color=blue] (v1.center) -- (v2.center) -- (v5.center) -- (v1.center);
					\end{scope}
					\begin{scope}[transparency group,opacity=.5]
						\draw[edge,opacity=1,color=green] (v2.center) -- (v3.center) -- (v8.center) -- (v7.center) -- (v2.center);
						\fill[edge,opacity=1,color=green] (v2.center) -- (v3.center) -- (v8.center) -- (v7.center) -- (v2.center);
					\end{scope}
					\begin{scope}[transparency group,opacity=.5]
						\draw[edge,opacity=1,color=lightgray] (v6.center) -- (v7.center) -- (v6.center);
						\fill[edge,opacity=1,color=lightgray] (v6.center) -- (v7.center) -- (v6.center);
					\end{scope}
				\end{pgfonlayer}
				
				\node[text=blue] (a1) at (1,0.75) {$e_1$};
				\node[text=green] (a2) at (6,-1) {$e_2$};
				\node[text=lightgray] (a3) at (3,-2.75) {$e_3$};
	\end{tikzpicture}\end{minipage}
	\begin{minipage}{.5\textwidth}\begin{flushright}\begin{tikzpicture}
				\tikzstyle{vertex} = [fill,shape=circle,node distance=80pt]
				\tikzstyle{edge} = [fill,opacity=.5,fill opacity=.5,line cap=round, line join=round, line width=30pt]
				\tikzstyle{elabel} =  [fill,shape=circle,node distance=40pt]
				
				\pgfdeclarelayer{background}
				\pgfsetlayers{background,main}
				
				\begin{scope}[every node/.style={circle,thick,draw}]
					\node (v1) at (2,0) {$v_1$};
					\node (v2) at (2,-1) {$v_2$};
					\node (v3) at (2,-2) {$v_3$};
					\node (v4) at (2, -3) {$v_4$};
					\node (v5) at (2,-4) {$v_5$};
					\node (v6) at (2,-5) {$v_6$};
					\node (v7) at (2,-6) {$v_7$};
					\node (v8) at (2,-7) {$v_8$};
					\node (e1) at (7, -1.5) {$v_9$};
					\node (e2) at (7, -3.5) {$v_{10}$};
					\node (e3) at (7, -5.5) {$v_{11}$};
				\end{scope}
				\begin{scope}[>={Stealth[black]},
					every node/.style={circle},
					every edge/.style={draw = blue,very thick}]
					\path [-] (v1) edge node[] {} (e1);
					\path [-] (v2) edge node[] {} (e1);
					\path [-] (v5) edge node[] {} (e1);
				\end{scope}
				\begin{scope}[>={Stealth[black]},
					every node/.style={circle},
					every edge/.style={draw = green,very thick}]
					\path [-] (v2) edge node[] {} (e2);
					\path [-] (v3) edge node[] {} (e2);
					\path [-] (v7) edge node[] {} (e2);
					\path [-] (v8) edge node[] {} (e2);
				\end{scope}
				\begin{scope}[>={Stealth[black]},
					every node/.style={circle},
					every edge/.style={draw = lightgray,very thick}]
					\path [-] (v6) edge node[] {} (e3);
					\path [-] (v7) edge node[] {} (e3);
				\end{scope}
	\end{tikzpicture}\end{flushright}\end{minipage}\vspace{0em}
\end{example}

This normal graph representation can also be defined for oriented hypergraphs.\\

\begin{definition}[\textbf{First normal graph representation of oriented hypergraphs}]\label{ONGR1} \ \\
	For an oriented hypergraph $OH = \left(\mathcal{V}_H, \mathcal{A}_H\right)$ with $\mathcal{V}_H = \left\{v_1, v_2, \dots v_N\right\}$, a normal graph $OG = \left(\mathcal{V}_G, \mathcal{A}_G\right)$ encoding the hypergraph can be defined with
	\begin{equation}
		\mathcal{V}_G = \mathcal{V}_H \cup \left\{v_{N + q} ~ \middle| ~ a_q \in \mathcal{A}_H\right\},
	\end{equation}
	\begin{equation}
		\mathcal{A}_G = \left\{\left(v_i, v_{N + q}\right) ~ \middle| ~ v_i \in a_q^{out} ~ \text{for} ~ a_q \in \mathcal{A}_H\right\} \cup \left\{\left(v_{N + q}, v_i\right) ~ \middle| ~ v_i \in a_q^{in} ~ \text{for} ~ a_q \in \mathcal{A}_H\right\},
	\end{equation}
	where the orientation of the arcs $\left(v_i, v_{N + q}\right)$ and $\left(v_{N + q}, v_i\right)$ indicates whether the original vertex $v_i \in \mathcal{V}_H$ is an output vertex or an input vertex of the hyperarc $a_q \in \mathcal{A}_H$.\\
\end{definition}

\begin{example}[\textbf{First normal graph representation of oriented hypergraphs}]\ \\
	Given the oriented hypergraph $OH = \left(\mathcal{V}_H, \mathcal{A}_H\right)$ with $\mathcal{V}_H = \left\{v_1, v_2, v_3, v_4, v_5, v_6, v_7, v_8\right\}$ and $\mathcal{A}_H = \left\{\left(\left\{v_1, v_2\right\}, \left\{v_5\right\}\right), \left(\left\{v_3\right\}, \left\{v_2, v_7, v_8\right\}\right), \left(\left\{v_6\right\}, \left\{v_7\right\}\right)\right\}$, then the oriented normal graph $OG = \left(\mathcal{V}_G, \mathcal{A}_G\right)$ encoding this hypergraph is given by:
	\begin{equation*}
		\mathcal{V}_G = \left\{v_1, v_2, v_3, v_4, v_5, v_6, v_7, v_8\right\} \cup \left\{v_9, v_{10}, v_{11}\right\}
	\end{equation*}
	\begin{equation*}
		\mathcal{A}_G = \left\{\left(v_1, v_9\right), \left(v_2, v_9\right), \left(v_9, v_5\right), \left(v_3, v_{10}\right), \left(v_{10}, v_2\right),\right.
	\end{equation*}
	\begin{equation*}
		 \left.\left(v_{10}, v_7\right), \left(v_{10}, v_8\right), \left(v_6, v_{11}\right), \left(v_{11}, v_7\right)\right\}
	\end{equation*}
	and can be visualized in the following way:
	
	\begin{minipage}{.5\textwidth}\begin{tikzpicture}
				\tikzstyle{vertex} = [fill,shape=circle,node distance=80pt]
				\tikzstyle{edge} = [fill,opacity=.5,fill opacity=.5,line cap=round, line join=round, line width=30pt]
				\tikzstyle{elabel} =  [fill,shape=circle,node distance=40pt]
				
				\pgfdeclarelayer{background}
				\pgfsetlayers{background,main}
				
				\begin{scope}[every node/.style={circle,thick,draw}]
					\node (v1) at (0,0) {$v_1$};
					\node (v2) at (2,0) {$v_2$};
					\node (v3) at (4,0) {$v_3$};
					\node (v4) at (6, 0) {$v_4$};
					\node (v5) at (0,-2) {$v_5$};
					\node (v6) at (2,-2) {$v_6$};
					\node (v7) at (4,-2) {$v_7$};
					\node (v8) at (6,-2) {$v_8$};
				\end{scope}
				
				\begin{pgfonlayer}{background}
					\begin{scope}[transparency group,opacity=.2]
						\draw[edge,opacity=1,color=blue] (v1.center) -- (v2.center) -- (v5.center) -- (v1.center);
						\fill[edge,opacity=1,color=blue] (v1.center) -- (v2.center) -- (v5.center) -- (v1.center);
					\end{scope}
					\begin{scope}[transparency group,opacity=.2]
						\draw[edge,opacity=1,color=green] (v2.center) -- (v3.center) -- (v8.center) -- (v7.center) -- (v2.center);
						\fill[edge,opacity=1,color=green] (v2.center) -- (v3.center) -- (v8.center) -- (v7.center) -- (v2.center);
					\end{scope}
					\begin{scope}[transparency group,opacity=.2]
						\draw[edge,opacity=1,color=lightgray] (v6.center) -- (v7.center) -- (v6.center);
						\fill[edge,opacity=1,color=lightgray] (v6.center) -- (v7.center) -- (v6.center);
					\end{scope}
				\end{pgfonlayer}
				
				\node[text=blue] (v1a1) at (0.5, -0.5) {\footnotesize{$out$}};
				\node[text=blue] (v2a1) at (1.5, -0.5) {\footnotesize{$out$}};
				\node[text=blue] (v5a1) at (0.5, -1.5) {\footnotesize{$in$}};
				\node[text=green] (v3a2) at (4, -0.7) {\footnotesize{$out$}};
				\node[text=green] (v2a2) at (2.5, -0.5) {\footnotesize{$in$}};
				\node[text=green] (v7a2) at (4, -1.3) {\footnotesize{$in$}};
				\node[text=green] (v8a2) at (5.5, -1.5) {\footnotesize{$in$}};
				\node[text=lightgray] (v6a3) at (2.75, -2) {\footnotesize{$out$}};
				\node[text=lightgray] (v7a3) at (3.35, -2) {\footnotesize{$in$}};
				\node[text=blue] (a1) at (1,0.75) {$a_1$};
				\node[text=green] (a2) at (6,-1) {$a_2$};
				\node[text=lightgray] (a3) at (3,-2.75) {$a_3$};
	\end{tikzpicture}\end{minipage}
	\begin{minipage}{.5\textwidth}\begin{flushright}\begin{tikzpicture}
				\tikzstyle{vertex} = [fill,shape=circle,node distance=80pt]
				\tikzstyle{edge} = [fill,opacity=.5,fill opacity=.5,line cap=round, line join=round, line width=30pt]
				\tikzstyle{elabel} =  [fill,shape=circle,node distance=40pt]
				
				\pgfdeclarelayer{background}
				\pgfsetlayers{background,main}
				
				\begin{scope}[every node/.style={circle,thick,draw}]
					\node (v1) at (2,0) {$v_1$};
					\node (v2) at (2,-1) {$v_2$};
					\node (v3) at (2,-2) {$v_3$};
					\node (v4) at (2, -3) {$v_4$};
					\node (v5) at (2,-4) {$v_5$};
					\node (v6) at (2,-5) {$v_6$};
					\node (v7) at (2,-6) {$v_7$};
					\node (v8) at (2,-7) {$v_8$};
					\node (e1) at (7, -1.5) {$v_9$};
					\node (e2) at (7, -3.5) {$v_{10}$};
					\node (e3) at (7, -5.5) {$v_{11}$};
				\end{scope}
				\begin{scope}[>={Stealth[blue]},
					every node/.style={circle},
					every edge/.style={draw = blue,very thick}]
					\path [->] (v1) edge node[] {} (e1);
					\path [->] (v2) edge node[] {} (e1);
					\path [<-] (v5) edge node[] {} (e1);
				\end{scope}
				\begin{scope}[>={Stealth[green]},
					every node/.style={circle},
					every edge/.style={draw = green,very thick}]
					\path [<-] (v2) edge node[] {} (e2);
					\path [->] (v3) edge node[] {} (e2);
					\path [<-] (v7) edge node[] {} (e2);
					\path [<-] (v8) edge node[] {} (e2);
				\end{scope}
				\begin{scope}[>={Stealth[lightgray]},
					every node/.style={circle},
					every edge/.style={draw = lightgray,very thick}]
					\path [->] (v6) edge node[] {} (e3);
					\path [<-] (v7) edge node[] {} (e3);
				\end{scope}
	\end{tikzpicture}\end{flushright}\end{minipage}\vspace{0em}
\end{example}

In the case of the first normal graph representation, the vertex degrees of the resulting normal graph have different interpretations, depending on whether it is the vertex degree of an original vertex in the vertex set $\mathcal{V}_H$, or the vertex degree of a new vertex in the vertex set $\mathcal{V}_G \backslash \mathcal{V}_H$.\\

\begin{remark}[\textbf{Vertex degrees in the first normal graph representation}]\label{Vdeg} \ \\
	Given a not oriented hypergraph $NH = \left(\mathcal{V}_H, \mathcal{E}_H\right)$ and its not oriented normal graph representation $NG = \left(\mathcal{V}_G, \mathcal{E}_G\right)$, then for the vertex degrees of vertices $v_i \in \mathcal{V}_H$ it holds true that
	\begin{equation}
		\deg\left(v_i\right) = \left\lvert\left\{e_q \in \mathcal{E}_H ~ \middle| ~ v_i \in e_q\right\}\right\rvert
	\end{equation}
	and for vertices $v_{N + q} \in \mathcal{V}_G \backslash \mathcal{V}_H$ it holds true that
	\begin{equation}
		\deg\left(v_{N + q}\right) = \left\lvert e_q\right\rvert,
	\end{equation}
	since every original vertex $v_i \in \mathcal{V}_H$ is only connected to vertices $v_{N + q} \in \mathcal{V}_G \backslash \mathcal{V}_H$, which encode a hyperedge $e_q \in \mathcal{E}_H$ containing vertex $v_i$, and vice versa.\\
	
	Analogously, given an oriented hypergraph $OH = \left(\mathcal{V}_H, \mathcal{A}_H\right)$ with its respective normal graph representations $OG = \left(\mathcal{V}_G, \mathcal{A}_G\right)$, for all original vertices $v_i \in \mathcal{V}_H$ it holds true that 
	\begin{equation}
		\deg\left(v_i\right) = \left\lvert\left\{a_q \in \mathcal{A}_H ~ \middle| ~ v_i \in a_q^{out} ~ \text{or} ~ v_i \in a_q^{in}\right\}\right\rvert
	\end{equation}
	and for the vertex degrees of vertices $v_{N + q} \in \mathcal{V}_G \backslash \mathcal{V}_H$ it holds true that
	\begin{equation}
		\deg\left(v_{N + q}\right) = \left\lvert a_q^{out} \cup a_q^{in}\right\rvert,
	\end{equation}
	since every original vertex $v_i \in \mathcal{V}_H$ is only connected to vertices $v_{N + q} \in \mathcal{V}_G \backslash \mathcal{V}_H$ encoding a hyperarc $a_q \in \mathcal{A}_H$, which contains $v_i$ either as an output or as an input vertex.\\
\end{remark}

The following theorem shows that the resulting normal graph of the first normal graph representation is a bipartite normal graph.\\

\begin{theorem}[\textbf{First normal graph representation results in bipartite graph}]\label{NGR1bipartite} (Page 3 in \cite{mori2015peeling}: Bipartite graph representation of hypergraphs) \ \\
	Both, in the case of a not oriented hypergraph $NH = \left(\mathcal{V}_H, \mathcal{E}_H\right)$ and in the case of an oriented hypergraph $OH = \left(\mathcal{V}_H, \mathcal{A}_H\right)$, the resulting first normal graph representation is either a bipartite not oriented normal graph $NG = \left(\mathcal{V}_G, \mathcal{E}_G\right)$ or a bipartite oriented normal graph $OG = \left(\mathcal{V}_G, \mathcal{A}_G\right)$.\\
\end{theorem}

\begin{proof}\ \\
	\textbf{First case:} Not oriented hypergraph
	
	Given a not oriented hypergraph $NH = \left(\mathcal{V}_H, \mathcal{E}_H\right)$ with $\mathcal{V}_H = \left\{v_1, v_2, \dots v_N\right\}$, $\mathcal{E}_H = \left\{e_1, e_2, \dots e_M\right\}$ and the corresponding normal graph representation $NG = \left(\mathcal{V}_G, \mathcal{E}_G\right)$, then the vertex set $\mathcal{V}_G$ can be partitioned in the following way:
	\begin{align*}
		\mathcal{V}_G & = \left\{v_1, v_2, \dots v_N, v_{N + 1}, v_{N + 2}, \dots v_{N + M}\right\}\\
		& = \left\{v_1, v_2, \dots, v_N\right\} \cup \left\{v_{N + 1}, v_{N + 2}, \dots v_{N + M}\right\}\\
		& = \mathcal{V}_H \cup \left\{v_{N + 1}, v_{N + 2}, \dots v_{N + M}\right\},
	\end{align*}
	which automatically fulfills the condition $\mathcal{V}_H \cap \left\{v_{N + 1}, v_{N + 2}, \dots v_{N + M}\right\} = \emptyset$. Hence, the partition of the vertex set $\mathcal{V}_G$ is valid based on the first condition for bipartite not oriented normal graphs from definition (\ref{bipartiteG}).\\
	
	Since the definition of the edge set $\mathcal{E}_G$ is based on the first normal graph representation, it holds true that:
	\begin{align*}
		\mathcal{E}_G & = \left\{\left\{v_i, v_{N + q}\right\} ~ \middle| ~ v_i \in e_q ~ \text{for} ~ e_q \in \mathcal{E}_H\right\}\\
		& = \left\{\left\{v_i, v_j\right\} ~ \middle| ~ v_i \in e_q \subseteq \mathcal{V}_H ~ \text{for} ~ e_q \in \mathcal{E}_H ~ \text{and} ~ v_{j} \in \left\{v_{N + 1}, v_{N + 2}, \dots v_{N + M}\right\} ~ \text{with} ~ j = N + q\right\}
	\end{align*}
	Thus, the vertex set partition of $\mathcal{V}_G$ together with the edge set $\mathcal{E}_G$ also fulfills the second condition for bipartite not oriented normal graphs
	\begin{equation*}
		\forall e_q \in \mathcal{E}_G: \quad e_q = \left\{v_i, v_j\right\} ~ \text{with} ~ v_i \in \mathcal{V}_H ~ \text{and} ~ v_j \in \left\{v_{N + 1}, v_{N + 2}, \dots v_{N + M}\right\}.
	\end{equation*}

	This proves that the normal graph representation $NG = \left(\mathcal{V}_G, \mathcal{E}_G\right)$ of a not oriented hypergraph $NH = \left(\mathcal{V}_H, \mathcal{E}_H\right)$ results in a bipartite not oriented normal graph.\\
	
	\textbf{Second case:} Oriented hypergraph
	
	Similarly, for an oriented hypergraph $OH = \left(\mathcal{V}_H, \mathcal{A}_H\right)$ with $\mathcal{V}_H = \left\{v_1, v_2, \dots v_N\right\}$, $\mathcal{A}_H = \left\{a_1, a_2, \dots a_M\right\}$ and its corresponding normal graph representation $OG = \left(\mathcal{V}_G, \mathcal{A}_G\right)$, the vertex set $\mathcal{V}_G$ can be partitioned to
	\begin{align*}
		\mathcal{V}_G & = \left\{v_1, v_2, \dots v_N, v_{N + 1}, v_{N + 2}, \dots v_{N + M}\right\}\\
		& = \left\{v_1, v_2, \dots, v_N\right\} \cup \left\{v_{N + 1}, v_{N + 2}, \dots v_{N + M}\right\}\\
		& = \mathcal{V}_H \cup \left\{v_{N + 1}, v_{N + 2}, \dots v_{N + M}\right\},
	\end{align*}
	with $\mathcal{V}_H \cap \left\{v_{n + 1}, v_{n + 2}, \dots v_{n + m}\right\} = \emptyset$. Therefore, the vertex set $\mathcal{V}_G$ fulfills the first condition of a bipartite oriented normal graph introduced in definition (\ref{bipartiteG}).\\
	
	Moreover, for the arc set $\mathcal{A}_G$ it holds true that:
	\begin{align*}
		\mathcal{A}_G = & \left\{\left(v_i, v_{N + q}\right) ~ \middle| ~ v_i \in a_q^{out} ~ \text{for} ~ a_q \in \mathcal{A}_H\right\}\\
		& \cup \left\{\left(v_{N + q}, v_i\right) ~ \middle| ~ v_i \in a_q^{in} ~ \text{for} ~ a_q \in \mathcal{A}_H\right\}\\
		= & \left\{\left(v_i, v_j\right) ~ \middle| ~ v_i \in a_q^{out} \subseteq \mathcal{V}_H ~ \text{for} ~ a_q \in \mathcal{A}_H ~ \text{and} ~ v_{j} \in \left\{v_{N + 1}, v_{N + 2}, \dots v_{N + M}\right\}, j = N + q\right\}\\
		& \cup \left\{\left(v_j, v_i\right) ~ \middle| ~ v_i \in a_q^{in} \subseteq \mathcal{V}_H ~ \text{for} ~ a_q \in \mathcal{A}_H ~ \text{and} ~ v_{j} \in \left\{v_{N + 1}, v_{N + 2}, \dots v_{N + M}\right\}, j = N + q\right\}
	\end{align*}
	
	\clearpage
	This shows that the second condition for bipartite oriented normal graphs
	\begin{equation*}
		\forall a_q \in \mathcal{A}_G: \quad a_q = \left(v_k, v_l\right) ~ \text{with} ~ v_k \in \mathcal{V}_H, v_l \in \left\{v_{N + 1}, v_{N + 2}, \dots v_{N + M}\right\}
	\end{equation*}
	\begin{equation*}
		\text{or with} ~ v_k \in \left\{v_{N + 1}, v_{N + 2}, \dots v_{N + M}\right\}, v_l \in \mathcal{V}_H
	\end{equation*}
 	is also fulfilled.\\

	Hence, the normal graph representation $OG = \left(\mathcal{V}_G, \mathcal{A}_G\right)$ of an oriented hypergraph $OH = \left(\mathcal{V}_H, \mathcal{A}_H\right)$ results in a bipartite oriented normal graph as well.\\
\end{proof}

The following theorem describes necessary additional information with which it is possible to reconstruct the underlying hypergraph from its normal graph representation uniquely.\\

\begin{theorem}[\textbf{Unique first normal graph representation}]\label{NGR1unique} \ \\
	Given a not oriented hypergraph $NH = \left(\mathcal{V}_H, \mathcal{E}_H\right)$ or an oriented hypergraph $OH = \left(\mathcal{V}_H, \mathcal{A}_H\right)$ and the corresponding not oriented normal graph $NG = \left(\mathcal{V}_G, \mathcal{E}_G\right)$ or oriented normal graph $OG = \left(\mathcal{V}_G, \mathcal{A}_G\right)$ based on the first normal graph representation. Then the hypergraph $NH$ or $OH$ can be reconstructed uniquely from $NG$ or $OG$ if for each connected component of $NG$ or $OG$, there exists at least one vertex $v \in \mathcal{V}_G$, where it is known, whether $v \in \mathcal{V}_H$ or $v \notin \mathcal{V}_H$.
\end{theorem}

\begin{proof}\ \\
	For a not oriented hypergraph $NH = \left(\mathcal{V}_H, \mathcal{E}_H\right)$ or an oriented hypergraph $OH = \left(\mathcal{V}_H, \mathcal{A}_H\right)$, the normal graph representations $NG = \left(\mathcal{V}_G, \mathcal{E}_G\right)$ or $OG = \left(\mathcal{V}_G, \mathcal{A}_G\right)$ are always bipartite normal graphs as shown in theorem (\ref{NGR1bipartite}). This also indicates that all connected components are bipartite subgraphs and for each connected components $G$ it holds true that:\\
	
	If vertex $v$ is the only vertex in $G$ then it is already clear that $v \in \mathcal{V}_H$, because for $v \notin \mathcal{V}_H$ \big(automatically implying $v \in \mathcal{V}_G \backslash \mathcal{V}_H$\big), $v$ would need to be connected to at least $min_e \geq 2$ or $min_a \geq 2$ original vertices and thus $G$ would need to consist of at least three vertices.\\
	
	Based on the last argument a connected component with only two vertices is also impossible.\\
	
	Otherwise, if $G$ consists of at least three vertices, then the set of all vertices of $G$ can be partitioned, such that the two vertex sets are valid with regard to the bipartite property of the normal graph representation. For each connected component the partition of vertices necessary for the bipartite property is unique and for one vertex $v$ in $G$ it is already known whether $v \in \mathcal{V}_H$ or $v \in \mathcal{V}_G \backslash \mathcal{V}_H$ based on the assumption. Thus, the two vertex sets of $G$ can be clearly identified as $\mathcal{V}_H$ and $\mathcal{V}_G \backslash \mathcal{V}_H$ based on vertex $v$. Hence, all hyperedges or hyperarcs encoded as vertices in $\mathcal{V}_G \backslash \mathcal{V}_H$ together with their respective vertices in $\mathcal{V}_H$ can be decoded resulting in a unique not oriented or oriented hypergraph.\\
	
	This reasoning holds true for all connected components $G$ of the normal graph representation and therefore for the whole normal graph. Thus, for any given normal graph representation with a vertex $v \in \mathcal{V}_H$ or $v \notin \mathcal{V}_H$ for each connected component, it is possible to uniquely reconstruct the underlying hypergraph, which indicates a form of duality between hypergraphs and their normal graph representations.\\
\end{proof}

\begin{corollary}[\textbf{Unique first normal graph representation property}]\label{NGR1uniqueprop} \ \\
	Given a not oriented hypergraph $NH = \left(\mathcal{V}_H, \mathcal{E}_H\right)$ or an oriented hypergraph $OH = \left(\mathcal{V}_H, \mathcal{A}_H\right)$ with its normal graph representation $NG = \left(\mathcal{V}_G, \mathcal{E}_G\right)$ or $OG = \left(\mathcal{V}_G, \mathcal{A}_G\right)$, then the mapping from the normal graph representation to the underlying hypergraph is unique if for each connected component at least one of the two properties is fulfilled:
	\begin{itemize}
		\item[1)] There exists at least one vertex $v \in \mathcal{V}_G$ with $\deg\left(v\right) < min_e$ or $\deg\left(v\right) < min_a$.
		\item[2)] There exists at least one vertex $v \in \mathcal{V}_G$ with $\deg\left(v\right) > max_e$ or $\deg\left(v\right) > max_a$.
	\end{itemize}
\end{corollary}

\begin{proof}\ \\
	For the normal graph representations $NG = \left(\mathcal{V}_G, \mathcal{E}_G\right)$ or $OG = \left(\mathcal{V}_G, \mathcal{A}_G\right)$ of hypergraphs $NH = \left(\mathcal{V}_H, \mathcal{E}_H\right)$ or $OH = \left(\mathcal{V}_H, \mathcal{A}_H\right)$, it holds true that:
	
	\begin{itemize}
		\item[1)] If for any vertex $v \in \mathcal{V}_G$ it is the case that $\deg\left(v\right) < min_e$ or $\deg\left(v\right) < min_a$, then this implies $v \in \mathcal{V}_H$.
		
		Assume that $v \in \mathcal{V}_G \backslash \mathcal{V}_H$ would be the case. This means vertex $v$ represents a hyperedge $e_q \in  \mathcal{E}_H$ or a hyperarc $a_q \in \mathcal{A}_H$, while only being connected to $\deg\left(v\right) < min_e$ or $\deg\left(v\right) < min_a$ original vertices, which is a contradiction to the minimum hyperedge or hyperarc cardinality in definition (\ref{minea}).
		
		\item[2)] Similarly, if for any vertex $v \in \mathcal{V}_G$ it is the case that $\deg\left(v\right) > max_e$ or $\deg\left(v\right) > max_a$, then this implies $v \in \mathcal{V}_H$.
		
		Assume that $v \in \mathcal{V}_G \backslash \mathcal{V}_H$ would hold true, then this means vertex $v$ represents a hyperedge $e_q \in  \mathcal{E}_H$ or a hyperarc $a_q \in \mathcal{A}_H$ consisting of $\deg\left(v\right) > max_e$ or $\deg\left(v\right) > max_a$ original vertices, which is a contradiction to the maximum hyperedge or hyperarc cardinality in definition (\ref{maxea}).
	\end{itemize}
	
	Therefore, if at least one of the two conditions is fulfilled for each connected component of the normal graph representation, then for each connected component there exists at least one vertex $v$, where it is known that $v \in \mathcal{V}_H$. Hence, theorem (\ref{NGR1unique}) can be applied and based on this, the mapping from the normal graph representation back to the underlying hypergraph is unique.\\
\end{proof}

As theorem (\ref{NGR1unique}) indicates, for a general hypergraph it does not hold true that the resulting first normal graph representation can be uniquely traced back to the underlying hypergraph, rather several hypergraphs can have the same normal graph representation.\\

\begin{example}[\textbf{Counterexamples for the uniqueness of the first normal graph representation}]\ \\
	For the two not oriented hypergraphs $NH_1 = \left(\mathcal{V}_{H_1}, \mathcal{E}_{H_1}\right)$ with
	\begin{equation*}
		\mathcal{V}_{H_1} = \left\{v_1, v_2, v_3, v_4, v_5, v_6\right\}
	\end{equation*}
	\begin{equation*}
		\mathcal{E}_{H_1} = \left\{\left\{v_1, v_2, v_3\right\}, \left\{v_1, v_4, v_5\right\}, \left\{v_2, v_4, v_6\right\}, \left\{v_3, v_5, v_6\right\}\right\}
	\end{equation*}
	and $NH_2 = \left(\mathcal{V}_{H_2}, \mathcal{E}_{H_2}\right)$ with
	\begin{equation*}
		\mathcal{V}_{H_2} = \left\{v_A, v_B, v_C, v_D\right\}
	\end{equation*}
	\begin{equation*}
		\mathcal{E}_{H_2} = \left\{\left\{v_A, v_B\right\}, \left\{v_A, v_C\right\}, \left\{v_A, v_D\right\}, \left\{v_B, v_C\right\}, \left\{v_B, v_D\right\}, \left\{v_C, v_D\right\}\right\},
	\end{equation*}
	the resulting first normal graph representation is equivalent.
	
	\begin{minipage}{.5\textwidth}\begin{tikzpicture}
				\tikzstyle{vertex} = [fill,shape=circle,node distance=80pt]
				\tikzstyle{edge} = [fill,opacity=.5,fill opacity=.5,line cap=round, line join=round, line width=30pt]
				\tikzstyle{elabel} =  [fill,shape=circle,node distance=40pt]
				
				\pgfdeclarelayer{background}
				\pgfsetlayers{background,main}
				
				\begin{scope}[every node/.style={circle,thick,draw}]
					\node (v1) at (0,0) {$v_1$};
					\node (v2) at (3,0) {$v_2$};
					\node (v3) at (6,0) {$v_3$};
					\node (v4) at (0, -4) {$v_4$};
					\node (v5) at (3,-4) {$v_5$};
					\node (v6) at (6,-4) {$v_6$};
				\end{scope}
				
				\begin{pgfonlayer}{background}
					\begin{scope}[transparency group,opacity=.4]
						\draw[edge,opacity=1,color=blue] (v1.center) -- (v2.center) -- (v3.center) -- (v1.center);
						\fill[edge,opacity=1,color=blue] (v1.center) -- (v2.center) -- (v3.center) -- (v1.center);
					\end{scope}
					\begin{scope}[transparency group,opacity=.4]
						\draw[edge,opacity=1,color=green] (v1.center) -- (v4.center) -- (v5.center) -- (v1.center);
						\fill[edge,opacity=1,color=green] (v1.center) -- (v4.center) -- (v5.center) -- (v1.center);
					\end{scope}
					\begin{scope}[transparency group,opacity=.4]
						\draw[edge,opacity=1,color=lightgray] (v2.center) -- (v4.center) -- (v2.center) -- (v6.center) -- (v2.center);
						\fill[edge,opacity=1,color=lightgray] (v2.center) -- (v4.center) -- (v2.center) -- (v6.center) -- (v2.center);
					\end{scope}
					\begin{scope}[transparency group,opacity=.4]
						\draw[edge,opacity=1,color=cyan] (v3.center) -- (v5.center) -- (v6.center) -- (v3.center);
						\fill[edge,opacity=1,color=cyan] (v3.center) -- (v5.center) -- (v6.center) -- (v3.center);
					\end{scope}
				\end{pgfonlayer}
				
				\node[text=blue] (a1) at (-0.75,0) {$e_A$};
				\node[text=green] (a2) at (-0.75,-2) {$e_B$};
				\node[text=lightgray] (a3) at (-0.75,-4) {$e_C$};
				\node[text=cyan] (a4) at (6.75,-2) {$e_D$};
	\end{tikzpicture}\end{minipage}
	\begin{minipage}{.5\textwidth}\begin{flushright}\begin{tikzpicture}
				\tikzstyle{vertex} = [fill,shape=circle,node distance=80pt]
				\tikzstyle{edge} = [fill,opacity=.5,fill opacity=.5,line cap=round, line join=round, line width=30pt]
				\tikzstyle{elabel} =  [fill,shape=circle,node distance=40pt]
				
				\pgfdeclarelayer{background}
				\pgfsetlayers{background,main}
				
				\begin{scope}[every node/.style={circle,thick,draw}]
					\node (v1) at (2,-1) {$v_1$};
					\node (v2) at (2,-2) {$v_2$};
					\node (v3) at (2,-3) {$v_3$};
					\node (v4) at (2,-4) {$v_4$};
					\node (v5) at (2,-5) {$v_5$};
					\node (v6) at (2,-6) {$v_6$};
					\node (eA) at (7,-1) {$v_A$};
					\node (eB) at (7,-2.67) {$v_B$};
					\node (eC) at (7,-4.33) {$v_C$};
					\node (eD) at (7,-6) {$v_D$};
				\end{scope}
				\begin{scope}[>={Stealth[black]},
					every node/.style={circle},
					every edge/.style={draw = blue,very thick}]
					\path [-] (v1) edge node[] {} (eA);
					\path [-] (v2) edge node[] {} (eA);
					\path [-] (v3) edge node[] {} (eA);
				\end{scope}
				\begin{scope}[>={Stealth[black]},
					every node/.style={circle},
					every edge/.style={draw = green,very thick}]
					\path [-] (v1) edge node[] {} (eB);
					\path [-] (v4) edge node[] {} (eB);
					\path [-] (v5) edge node[] {} (eB);
				\end{scope}
				\begin{scope}[>={Stealth[black]},
					every node/.style={circle},
					every edge/.style={draw = lightgray,very thick}]
					\path [-] (v2) edge node[] {} (eC);
					\path [-] (v4) edge node[] {} (eC);
					\path [-] (v6) edge node[] {} (eC);
				\end{scope}
				\begin{scope}[>={Stealth[black]},
					every node/.style={circle},
					every edge/.style={draw = cyan,very thick}]
					\path [-] (v3) edge node[] {} (eD);
					\path [-] (v5) edge node[] {} (eD);
					\path [-] (v6) edge node[] {} (eD);
				\end{scope}
	\end{tikzpicture}\end{flushright}\end{minipage}

	\begin{minipage}{.5\textwidth}\begin{tikzpicture}
				\tikzstyle{vertex} = [fill,shape=circle,node distance=80pt]
				\tikzstyle{edge} = [fill,opacity=.5,fill opacity=.5,line cap=round, line join=round, line width=30pt]
				\tikzstyle{elabel} =  [fill,shape=circle,node distance=40pt]
				
				\pgfdeclarelayer{background}
				\pgfsetlayers{background,main}
				
				\begin{scope}[every node/.style={circle,thick,draw}]
					\node (vA) at (0,0) {$v_A$};
					\node (vB) at (4,0) {$v_B$};
					\node (vC) at (0, -4) {$v_C$};
					\node (vD) at (4,-4) {$v_D$};
				\end{scope}
				
				\begin{pgfonlayer}{background}
					\begin{scope}[transparency group,opacity=.4]
						\draw[edge,opacity=1,color=blue] (vA.center) -- (vB.center) -- (vA.center);
						\fill[edge,opacity=1,color=blue] (vA.center) -- (vB.center) -- (vA.center);
					\end{scope}
					\begin{scope}[transparency group,opacity=.4]
						\draw[edge,opacity=1,color=green] (vA.center) -- (vC.center) -- (vA.center);
						\fill[edge,opacity=1,color=green] (vA.center) -- (vC.center) -- (vA.center);
					\end{scope}
					\begin{scope}[transparency group,opacity=.4]
						\draw[edge,opacity=1,color=lightgray] (vA.center) -- (vD.center) -- (vA.center);
						\fill[edge,opacity=1,color=lightgray] (vA.center) -- (vD.center) -- (vA.center);
					\end{scope}
					\begin{scope}[transparency group,opacity=.4]
						\draw[edge,opacity=1,color=cyan] (vB.center) -- (vC.center) -- (vB.center);
						\fill[edge,opacity=1,color=cyan] (vB.center) -- (vC.center) -- (vB.center);
					\end{scope}
					\begin{scope}[transparency group,opacity=.4]
						\draw[edge,opacity=1,color=lime] (vB.center) -- (vD.center) -- (vB.center);
						\fill[edge,opacity=1,color=lime] (vB.center) -- (vD.center) -- (vB.center);
					\end{scope}
					\begin{scope}[transparency group,opacity=.4]
						\draw[edge,opacity=1,color=teal] (vC.center) -- (vD.center) -- (vC.center);
						\fill[edge,opacity=1,color=teal] (vC.center) -- (vD.center) -- (vC.center);
					\end{scope}
				\end{pgfonlayer}
				
				\node[text=blue] (a1) at (2,0.75) {$e_1$};
				\node[text=green] (a2) at (-0.75,-2) {$e_2$};
				\node[text=lightgray] (a3) at (-0.75,0) {$e_3$};
				\node[text=cyan] (a4) at (-0.75,-4) {$e_4$};
				\node[text=lime] (a2) at (4.75,-2) {$e_5$};
				\node[text=teal] (a6) at (2,-4.75) {$e_6$};
	\end{tikzpicture}\end{minipage}
	\begin{minipage}{.5\textwidth}\begin{flushright}\begin{tikzpicture}
				\tikzstyle{vertex} = [fill,shape=circle,node distance=80pt]
				\tikzstyle{edge} = [fill,opacity=.5,fill opacity=.5,line cap=round, line join=round, line width=30pt]
				\tikzstyle{elabel} =  [fill,shape=circle,node distance=40pt]
				
				\pgfdeclarelayer{background}
				\pgfsetlayers{background,main}
				
				\begin{scope}[every node/.style={circle,thick,draw}]
					\node (v1) at (2,-1) {$v_1$};
					\node (v2) at (2,-2) {$v_2$};
					\node (v3) at (2,-3) {$v_3$};
					\node (v4) at (2,-4) {$v_4$};
					\node (v5) at (2,-5) {$v_5$};
					\node (v6) at (2,-6) {$v_6$};
					\node (eA) at (7,-1) {$v_A$};
					\node (eB) at (7,-2.67) {$v_B$};
					\node (eC) at (7,-4.33) {$v_C$};
					\node (eD) at (7,-6) {$v_D$};
				\end{scope}
				\begin{scope}[>={Stealth[black]},
					every node/.style={circle},
					every edge/.style={draw = blue,very thick}]
					\path [-] (v1) edge node[] {} (eA);
					\path [-] (v1) edge node[] {} (eB);
				\end{scope}
				\begin{scope}[>={Stealth[black]},
					every node/.style={circle},
					every edge/.style={draw = green,very thick}]
					\path [-] (v2) edge node[] {} (eA);
					\path [-] (v2) edge node[] {} (eC);
				\end{scope}
				\begin{scope}[>={Stealth[black]},
					every node/.style={circle},
					every edge/.style={draw = lightgray,very thick}]
					\path [-] (v3) edge node[] {} (eA);
					\path [-] (v3) edge node[] {} (eD);
				\end{scope}
				\begin{scope}[>={Stealth[black]},
					every node/.style={circle},
					every edge/.style={draw = cyan,very thick}]
					\path [-] (v4) edge node[] {} (eB);
					\path [-] (v4) edge node[] {} (eC);
				\end{scope}
				\begin{scope}[>={Stealth[black]},
					every node/.style={circle},
					every edge/.style={draw = lime,very thick}]
					\path [-] (v5) edge node[] {} (eB);
					\path [-] (v5) edge node[] {} (eD);
				\end{scope}
				\begin{scope}[>={Stealth[black]},
					every node/.style={circle},
					every edge/.style={draw = teal,very thick}]
					\path [-] (v6) edge node[] {} (eC);
					\path [-] (v6) edge node[] {} (eD);
				\end{scope}
	\end{tikzpicture}\end{flushright}\end{minipage}
	
	\clearpage
	Thus, without further information, such as the number of original vertices $\left\lvert\mathcal{V}_H\right\rvert$, the number of hyperedges $\left\lvert\mathcal{E}_H\right\rvert$, the minimum or maximum hyperedge cardinality or the knowledge for one vertex $v \in \mathcal{V}_G$ whether $v \in \mathcal{V}_H$ or $v \in \mathcal{V}_G \backslash \mathcal{V}_H$, it is impossible to reconstruct the underlying hypergraph uniquely. Both hypergraphs have the same first normal graph representation $NG = \left(\mathcal{V}_G, \mathcal{E}_G\right)$ with:
	\begin{equation*}
		\mathcal{V}_G = \left\{v_1, v_2, v_3, v_4, v_5, v_6, v_A, v_B, v_C, v_D\right\}
	\end{equation*}
	\begin{equation*}
		\mathcal{E}_G = \left\{\left\{v_1, v_A\right\}, \left\{v_1, v_B\right\}, \left\{v_2, v_A\right\}, \left\{v_2, v_C\right\}, \left\{v_3, v_A\right\}, \left\{v_3, v_D\right\}, \right.
	\end{equation*}
	\begin{equation*}
		\left.\left\{v_4, v_B\right\}, \left\{v_4, v_C\right\}, \left\{v_5, v_B\right\}, \left\{v_5, v_D\right\}, \left\{v_6, v_C\right\}, \left\{v_6, v_D\right\}\right\}.
	\end{equation*}\vspace{0em}
\end{example}

Besides the issue of the non-uniqueness when going back from the normal graph to the hypergraph, the first normal graph representation can also lead to an exponentially growing number of new vertices and edges or arcs.\\

\begin{remark}[\textbf{Drawback of the first normal graph representation}]\label{NGR1drawback} \ \\
	Assuming that hyperedges and hyperarcs in a hypergraph have to be unique, as discussed in remark (\ref{Hfinite}), then the number of hyperedges in a not oriented hypergraph $NH = \left(\mathcal{V}_H, \mathcal{E}_H\right)$ is limited tightly by
	\begin{equation}
		\left\lvert\mathcal{E}_H\right\rvert \leq \left\lvert 2^\mathcal{V} \backslash \left\{\emptyset, \left\{v_1\right\}, \left\{v_2\right\}, \dots \left\{v_N\right\}\right\}\right\rvert = \left\lvert 2^\mathcal{V}\right\rvert - \left(N + 1\right) = 2^N - \left(N + 1\right)
	\end{equation}
	and the number of hyperarcs in a not oriented hypergraph $OH = \left(\mathcal{V}_H, \mathcal{A}_H\right)$ is limited tightly by
	\begin{equation}
		\left\lvert\mathcal{A}_H\right\rvert \leq C ~ \text{with} ~ 2^N - \left(N + 1\right) < C,
	\end{equation}
	because every possible hyperedge $\left\{v_{i_1}, v_{i_2}, \dots v_{i_n}\right\} \subseteq \mathcal{V}$ can be partitioned into at least two feasible and unique hyperarcs $\left(\left\{v_{i_1}\right\}, \left\{v_{i_2}, v_{i_3}\dots v_{i_n}\right\}\right)$ and $\left(\left\{v_{i_2}, v_{i_3}\dots v_{i_n}\right\}, \left\{v_{i_1}\right\}\right)$. Thus, the maximum number of possible hyperarcs is at least twice the number of possible hyperedges given a set of vertices $\mathcal{V} = \left\{v_1, v_2, \dots v_N\right\}$.\\
	
	This indicates, that the cardinality of the set of new vertices $\mathcal{V}_G \backslash \mathcal{V}_H$, representing the hyperedges or hyperarcs of the underlying hypergraph, is possibly exponential in the number of original vertices $N$. Moreover, since every new vertex $v \in \mathcal{V}_G \backslash \mathcal{V}_H$ is connected to at least $2$ and at most $N$ original vertices, the number of edges or arcs possibly also grows exponential in the number of original vertices $N$.\\
\end{remark}

In order to tackle the problem of a possibly exponential number of new vertices, a second normal graph representation is introduced.\\

\clearpage
\begin{definition}[\textbf{Second normal graph representation of not oriented hypergraphs}]\label{NNGR2} (Introduction in \cite{zhou2006learning})\ \\
	Given a not oriented hypergraph $NH = \left(\mathcal{V}_H, \mathcal{E}_H\right)$, then a not oriented normal graph $NG = \left(\mathcal{V}_G, \mathcal{E}_G\right)$ encoding the hypergraph can be defined with
	\begin{equation}
		\mathcal{V}_G = \mathcal{V}_H,
	\end{equation}
	\begin{equation}
		\mathcal{E}_G = \bigcup_{e_q \in \mathcal{E}_H} \left\{\left\{v_i, v_j\right\} ~ \middle| ~ v_i, v_j \in e_q ~ \text{and} ~ v_i \neq v_j\right\},
	\end{equation}
	where for each hyperedge $e_q \in \mathcal{E}_H$ and each pair of vertices $v_i, v_j \in e_q \subseteq \mathcal{V}_H$ a new edge is introduced. Therefore, each hyperedge $e_q \in \mathcal{E}_H$ of the underlying hypergraph results in a complete subgraph in the normal graph representation.\\
	
	Note: If two vertices $v_i, v_j \in \mathcal{V}_H$ are together part of more than one hyperedge $e_q \in \mathcal{E}_H$, then the new edge $\left\{v_i, v_j\right\}$ appears in $\mathcal{E}_G$ more than once \big(exactly $\left\lvert\left\{e_q \in \mathcal{E}_H ~ \middle| ~ v_i, v_j \in e_q\right\}\right\rvert$ times\big).\\
\end{definition}

\begin{example}[\textbf{Second normal graph representation of not oriented hypergraphs}] \ \\
	Using a slight variation of the previous not oriented hypergraph $NH = \left(\mathcal{V}_H, \mathcal{E}_H\right)$ with $\mathcal{V}_H = \left\{v_1, v_2, v_3, v_4, v_5, v_6, v_7, v_8\right\}$ and $\mathcal{E}_H = \left\{\left\{v_1, v_2, v_3\right\}, \left\{v_2, v_3, v_7, v_8\right\}, \left\{v_6, v_7\right\}\right\}$, then the not oriented normal graph $NG = \left(\mathcal{V}_G, \mathcal{E}_G\right)$ representing this hypergraph is given by
	\begin{equation*}
		\mathcal{V}_G = \mathcal{V}_H = \left\{v_1, v_2, v_3, v_4, v_5, v_6, v_7, v_8\right\},
	\end{equation*}
	\begin{equation*}
		\mathcal{E}_G = \left\{\left\{v_1, v_2\right\}, \left\{v_1, v_3\right\}, \left\{v_2, v_3\right\}, \left\{v_2, v_3\right\}, \left\{v_2, v_7\right\}, \left\{v_2, v_8\right\}, \right.
	\end{equation*}
	\begin{equation*}
		\left. \left\{v_3, v_7\right\}, \left\{v_3, v_8\right\}, \left\{v_7, v_8\right\}, \left\{v_6, v_7\right\}\right\}
	\end{equation*}
	and can be visualized in the following way:
	
	\begin{minipage}{.5\textwidth}\begin{tikzpicture}
				\tikzstyle{vertex} = [fill,shape=circle,node distance=80pt]
				\tikzstyle{edge} = [fill,opacity=.5,fill opacity=.5,line cap=round, line join=round, line width=30pt]
				\tikzstyle{elabel} =  [fill,shape=circle,node distance=40pt]
				
				\pgfdeclarelayer{background}
				\pgfsetlayers{background,main}
				
				\begin{scope}[every node/.style={circle,thick,draw}]
					\node (v1) at (0,0) {$v_1$};
					\node (v2) at (2,0) {$v_2$};
					\node (v3) at (4,0) {$v_3$};
					\node (v4) at (6, 0) {$v_4$};
					\node (v5) at (0,-2) {$v_5$};
					\node (v6) at (2,-2) {$v_6$};
					\node (v7) at (4,-2) {$v_7$};
					\node (v8) at (6,-2) {$v_8$};
				\end{scope}
				
				\begin{pgfonlayer}{background}
					\begin{scope}[transparency group,opacity=.5]
						\draw[edge,opacity=1,color=blue] (v1.center) -- (v2.center) -- (v3.center) -- (v1.center);
						\fill[edge,opacity=1,color=blue] (v1.center) -- (v2.center) -- (v3.center) -- (v1.center);
					\end{scope}
					\begin{scope}[transparency group,opacity=.5]
						\draw[edge,opacity=1,color=green] (v2.center) -- (v3.center) -- (v8.center) -- (v7.center) -- (v2.center);
						\fill[edge,opacity=1,color=green] (v2.center) -- (v3.center) -- (v8.center) -- (v7.center) -- (v2.center);
					\end{scope}
					\begin{scope}[transparency group,opacity=.5]
						\draw[edge,opacity=1,color=lightgray] (v6.center) -- (v7.center) -- (v6.center);
						\fill[edge,opacity=1,color=lightgray] (v6.center) -- (v7.center) -- (v6.center);
					\end{scope}
				\end{pgfonlayer}
				
				\node[text=blue] (a1) at (2,0.75) {$e_1$};
				\node[text=green] (a2) at (6,-1) {$e_2$};
				\node[text=lightgray] (a3) at (3,-2.75) {$e_3$};
	\end{tikzpicture}\end{minipage}
	\begin{minipage}{.5\textwidth}\begin{flushright}\begin{tikzpicture}
				\tikzstyle{vertex} = [fill,shape=circle,node distance=80pt]
				\tikzstyle{edge} = [fill,opacity=.5,fill opacity=.5,line cap=round, line join=round, line width=30pt]
				\tikzstyle{elabel} =  [fill,shape=circle,node distance=40pt]
				
				\pgfdeclarelayer{background}
				\pgfsetlayers{background,main}
				
				\begin{scope}[every node/.style={circle,thick,draw}]
					\node (v1) at (0,0) {$v_1$};
					\node (v2) at (1.5,1.5) {$v_2$};
					\node (v3) at (3.5,1.5) {$v_3$};
					\node (v4) at (5,0) {$v_4$};
					\node (v5) at (0,-2) {$v_5$};
					\node (v6) at (1.5,-3.5) {$v_6$};
					\node (v7) at (3.5,-3.5) {$v_7$};
					\node (v8) at (5,-2) {$v_8$};
				\end{scope}

				\begin{scope}[>={Stealth[black]},
					every node/.style={circle},
					every edge/.style={draw = blue,very thick}]
					\path [-] (v1) edge node[] {} (v2);
					\path [-] (v1) edge node[] {} (v3);
					\path [-] (v2) edge node[] {} (v3);
				\end{scope}
			
				\begin{scope}[>={Stealth[black]},
					every node/.style={circle},
					every edge/.style={draw = green,very thick}]
					\path [-] (v2) edge[bend left=20] node[] {} (v3);
					\path [-] (v2) edge node[] {} (v7);
					\path [-] (v2) edge node[] {} (v8);
					\path [-] (v3) edge node[] {} (v7);
					\path [-] (v3) edge node[] {} (v8);
					\path [-] (v7) edge node[] {} (v8);
				\end{scope}
			
				\begin{scope}[>={Stealth[black]},
					every node/.style={circle},
					every edge/.style={draw = lightgray,very thick}]
					\path [-] (v6) edge node[] {} (v7);
				\end{scope}
	\end{tikzpicture}\end{flushright}\end{minipage}\vspace{0em}
\end{example}

\clearpage
\begin{definition}[\textbf{Second normal graph representation of oriented hypergraphs}]\label{ONGR2} \ \\
	Given an oriented hypergraph $OH = \left(\mathcal{V}_H, \mathcal{A}_H\right)$, then an oriented normal graph $OG = \left(\mathcal{V}_G, \mathcal{A}_G\right)$ encoding the hypergraph can be defined with
	\begin{equation}
		\mathcal{V}_G = \mathcal{V}_H,
	\end{equation}
	\begin{equation}
		\mathcal{A}_G = \bigcup_{a_q \in \mathcal{A}_H} \left\{\left(v_i, v_j\right) ~ \middle| ~ v_i \in a_q^{out}, v_j \in a_q^{in} ~ \text{or} ~ v_i, v_j \in a_q^{out} ~ \text{or} ~ v_i, v_j \in a_q^{in}, v_i \neq v_j\right\},
	\end{equation}
	where a newly introduced arc $\left(v_i, v_j\right)$ either indicates that the vertices $v_i$ and $v_j$ are part of a hyperarc $a_q \in \mathcal{A}_H$ with $v_i$ being an output vertex and $v_j$ being an input vertex, or the vertices $v_i$ and $v_j$ are both part of the output or input vertex set of a hyperarc $a_q \in \mathcal{A}_H$.\\
	
	Note: If for two vertices $v_i, v_j \in \mathcal{V}_H$, there exists more than one hyperarc $a_q \in \mathcal{A}_H$ such that $v_i \in a_q^{out}$ and $v_j \in a_q^{in}$, then the new arc $\left(v_i, v_j\right)$ is part of the arc set $\mathcal{A}_G$ more than once and similarly, if two vertices $v_i, v_j \in \mathcal{V}_H$ are part of the same output or input vertex set of more than one hyperarc $a_q \in \mathcal{A}_H$, then the arcs $\left(v_i, v_j\right)$ and $\left(v_j, v_i\right)$ appear more than once in $\mathcal{A}_G$. Overall, an arc $\left(v_i, v_j\right)$ with $v_i \neq v_j$ is part of $\mathcal{A}_G$
	\begin{equation*}
		\left\lvert\left\{a_q \in \mathcal{A}_H ~ \middle| ~ v_i \in a_q^{out}, v_j \in a_q^{in} ~ \text{or} ~ v_i, v_j \in a_q^{out} ~ \text{or} ~ v_i, v_j \in a_q^{in}\right\}\right\rvert
	\end{equation*}
	times.\\
\end{definition}

\begin{example}[\textbf{Second normal graph representation of oriented hypergraphs}]\ \\
	Given the oriented hypergraph $OH = \left(\mathcal{V}_H, \mathcal{A}_H\right)$ with $\mathcal{V}_H = \left\{v_1, v_2, v_3, v_4, v_5, v_6, v_7, v_8\right\}$ and $\mathcal{A}_H = \left\{\left(\left\{v_1, v_2\right\}, \left\{v_3\right\}\right), \left(\left\{v_3\right\}, \left\{v_2, v_7, v_8\right\}\right), \left(\left\{v_3\right\}, \left\{v_4, v_8\right\}\right), \left(\left\{v_6\right\}, \left\{v_7\right\}\right)\right\}$, then the oriented normal graph $OG = \left(\mathcal{V}_G, \mathcal{A}_G\right)$ encoding this hypergraph is given by
	\begin{equation*}
		\mathcal{V}_G = \mathcal{V}_H = \left\{v_1, v_2, v_3, v_4, v_5, v_6, v_7, v_8\right\},
	\end{equation*}
	\begin{equation*}
		\mathcal{A}_G = \left\{\left(v_1, v_3\right), \left(v_2, v_3\right), \left(v_1, v_2\right), \left(v_2, v_1\right)\right\} \cup
	\end{equation*}
	\begin{equation*}
		\left\{\left(v_3, v_2\right), \left(v_3, v_7\right), \left(v_3, v_8\right), \left(v_2, v_7\right), \left(v_7, v_2\right), \left(v_2, v_8\right), \left(v_8, v_2\right), \left(v_7, v_8\right), \left(v_8, v_7\right)\right\} \cup
	\end{equation*}
	\begin{equation*}
		\left\{\left(v_3, v_4\right), \left(v_3, v_8\right), \left(v_4, v_8\right), \left(v_8, v_4\right)\right\} \cup
		\left\{\left(v_6, v_7\right)\right\}
	\end{equation*}
	and can be visualized in the following way:
	
	\begin{minipage}{.5\textwidth}\begin{tikzpicture}
			\tikzstyle{vertex} = [fill,shape=circle,node distance=80pt]
			\tikzstyle{edge} = [fill,opacity=.5,fill opacity=.5,line cap=round, line join=round, line width=30pt]
			\tikzstyle{elabel} =  [fill,shape=circle,node distance=40pt]
			
			\pgfdeclarelayer{background}
			\pgfsetlayers{background,main}
			
			\begin{scope}[every node/.style={circle,thick,draw}]
				\node (v1) at (0,0) {$v_1$};
				\node (v2) at (2,0) {$v_2$};
				\node (v3) at (4,0) {$v_3$};
				\node (v4) at (6, 0) {$v_4$};
				\node (v5) at (0,-2) {$v_5$};
				\node (v6) at (2,-2) {$v_6$};
				\node (v7) at (4,-2) {$v_7$};
				\node (v8) at (6,-2) {$v_8$};
			\end{scope}
			
			\begin{pgfonlayer}{background}
				\begin{scope}[transparency group,opacity=.2]
					\draw[edge,opacity=1,color=blue] (v1.center) -- (v2.center) -- (v3.center) -- (v1.center);
					\fill[edge,opacity=1,color=blue] (v1.center) -- (v2.center) -- (v3.center) -- (v1.center);
				\end{scope}
				\begin{scope}[transparency group,opacity=.2]
					\draw[edge,opacity=1,color=green] (v2.center) -- (v3.center) -- (v8.center) -- (v7.center) -- (v2.center);
					\fill[edge,opacity=1,color=green] (v2.center) -- (v3.center) -- (v8.center) -- (v7.center) -- (v2.center);
				\end{scope}
				\begin{scope}[transparency group,opacity=.2]
					\draw[edge,opacity=1,color=lightgray] (v6.center) -- (v7.center) -- (v6.center);
					\fill[edge,opacity=1,color=lightgray] (v6.center) -- (v7.center) -- (v6.center);
				\end{scope}
				\begin{scope}[transparency group,opacity=.2]
					\draw[edge,opacity=1,color=cyan] (v3.center) -- (v4.center) -- (v8.center) -- (v3.center);
					\fill[edge,opacity=1,color=cyan] (v3.center) -- (v4.center) -- (v8.center) -- (v3.center);
				\end{scope}
			\end{pgfonlayer}
			
			\node[text=blue] (v1a1) at (0.75, 0) {\footnotesize{$out$}};
			\node[text=blue] (v2a1) at (2.75, 0) {\footnotesize{$out$}};
			\node[text=blue] (v5a1) at (3.35, 0) {\footnotesize{$in$}};
			\node[text=green] (v3a2) at (4, -0.7) {\footnotesize{$out$}};
			\node[text=green] (v2a2) at (2.5, -0.5) {\footnotesize{$in$}};
			\node[text=green] (v7a2) at (4, -1.3) {\footnotesize{$in$}};
			\node[text=green] (v8a2) at (5.35, -2) {\footnotesize{$in$}};
			\node[text=cyan] (v3a3) at (4.75, 0) {\footnotesize{$out$}};
			\node[text=cyan] (v4a3) at (5.5, -0.5) {\footnotesize{$in$}};
			\node[text=cyan] (v8a3) at (6, -1.3) {\footnotesize{$in$}};
			\node[text=lightgray] (v6a4) at (2.75, -2) {\footnotesize{$out$}};
			\node[text=lightgray] (v7a4) at (3.35, -2) {\footnotesize{$in$}};
			\node[text=blue] (a1) at (2,0.75) {$a_1$};
			\node[text=green] (a2) at (5,-2.75) {$a_2$};
			\node[text=cyan] (a3) at (5,0.75) {$a_3$};
			\node[text=lightgray] (a4) at (3,-2.75) {$a_4$};
	\end{tikzpicture}\end{minipage}
	\begin{minipage}{.5\textwidth}\begin{flushright}\begin{tikzpicture}
				\tikzstyle{vertex} = [fill,shape=circle,node distance=80pt]
				\tikzstyle{edge} = [fill,opacity=.5,fill opacity=.5,line cap=round, line join=round, line width=30pt]
				\tikzstyle{elabel} =  [fill,shape=circle,node distance=40pt]
				
				\pgfdeclarelayer{background}
				\pgfsetlayers{background,main}
				
				\begin{scope}[every node/.style={circle,thick,draw}]
					\node (v1) at (0,0) {$v_1$};
					\node (v2) at (1.5,1.5) {$v_2$};
					\node (v3) at (3.5,1.5) {$v_3$};
					\node (v4) at (5,0) {$v_4$};
					\node (v5) at (0,-2) {$v_5$};
					\node (v6) at (1.5,-3.5) {$v_6$};
					\node (v7) at (3.5,-3.5) {$v_7$};
					\node (v8) at (5,-2) {$v_8$};
				\end{scope}

				\begin{scope}[>={Stealth[blue]},
					every node/.style={circle},
					every edge/.style={draw = blue,very thick}]
					\path [->] (v1) edge node[] {} (v3);
					\path [->] (v2) edge node[] {} (v3);
					\path [<->] (v1) edge node[] {} (v2);
				\end{scope}
				
				\begin{scope}[>={Stealth[green]},
					every node/.style={circle},
					every edge/.style={draw = green,very thick}]
					\path [<-] (v2) edge[bend left=20] node[] {} (v3);
					\path [->] (v3) edge node[] {} (v7);
					\path [->] (v3) edge[bend right=15] node[] {} (v8);
					\path [<->] (v2) edge node[] {} (v7);
					\path [<->] (v2) edge[bend right=10] node[] {} (v8);
					\path [<->] (v7) edge node[] {} (v8);
				\end{scope}
			
				\begin{scope}[>={Stealth[cyan]},
					every node/.style={circle},
					every edge/.style={draw = cyan,very thick}]
					\path [->] (v3) edge node[] {} (v4);
					\path [->] (v3) edge node[] {} (v8);
					\path [<->] (v4) edge node[] {} (v8);
				\end{scope}
				
				\begin{scope}[>={Stealth[lightgray]},
					every node/.style={circle},
					every edge/.style={draw = lightgray,very thick}]
					\path [->] (v6) edge node[] {} (v7);
				\end{scope}
	\end{tikzpicture}\end{flushright}\end{minipage}

	In order to keep clarity in the normal graph representation as best as possible, the case $\left\{\left(v_i, v_j\right), \left(v_j, v_i\right) ~ \middle| ~ v_i, v_j \in a_q^{out} ~ \text{or} ~ v_i, v_j \in a_q^{in}, v_i \neq v_j\right\}$ for any hyperarc $a_q \in \mathcal{A}_H$ is represented by a combined arc $\longleftrightarrow$ instead of two individual arcs $\longleftarrow$ and $\longrightarrow$.\\
\end{example}

Similar to the first normal graph representation, the mapping back from the normal graph to the underlying hypergraph is not necessarily unique.\\

\begin{example}[\textbf{Counterexamples for the uniqueness of the second normal graph representation}]\ \\
	For the two not oriented hypergraphs $NH_1 = \left(\mathcal{V}_{H_1}, \mathcal{E}_{H_1}\right)$ with
	\begin{equation*}
		\mathcal{V}_{H_1} = \left\{v_1, v_2, v_3, v_4\right\},
	\end{equation*}
	\begin{equation*}
		\mathcal{E}_{H_1} = \left\{\left\{v_1, v_2, v_3\right\}, \left\{v_1, v_2, v_4\right\}, \left\{v_1, v_3, v_4\right\}\right\}
	\end{equation*}
	and $NH_2 = \left(\mathcal{V}_{H_2}, \mathcal{E}_{H_2}\right)$ with
	\begin{equation*}
		\mathcal{V}_{H_2} = \left\{v_1, v_2, v_3, v_4\right\},
	\end{equation*}
	\begin{equation*}
		\mathcal{E}_{H_2} = \left\{\left\{v_1, v_2\right\}, \left\{v_1, v_3\right\}, \left\{v_1, v_4\right\}, \left\{v_1, v_2, v_3, v_4\right\}\right\},
	\end{equation*}
	the resulting second normal graph representation is equivalent.
	
	\begin{minipage}{.5\textwidth}\begin{tikzpicture}
			\tikzstyle{vertex} = [fill,shape=circle,node distance=80pt]
			\tikzstyle{edge} = [fill,opacity=.5,fill opacity=.5,line cap=round, line join=round, line width=30pt]
			\tikzstyle{elabel} =  [fill,shape=circle,node distance=40pt]
			
			\pgfdeclarelayer{background}
			\pgfsetlayers{background,main}
			
			\begin{scope}[every node/.style={circle,thick,draw}]
				\node (v1) at (0,0) {$v_1$};
				\node (v2) at (4,0) {$v_2$};
				\node (v3) at (0, -4) {$v_3$};
				\node (v4) at (4,-4) {$v_4$};
			\end{scope}
			
			\begin{pgfonlayer}{background}
				\begin{scope}[transparency group,opacity=.4]
					\draw[edge,opacity=1,color=blue] (v1.center) -- (v2.center) -- (v3.center) -- (v1.center);
					\fill[edge,opacity=1,color=blue] (v1.center) -- (v2.center) -- (v3.center) -- (v1.center);
				\end{scope}
				\begin{scope}[transparency group,opacity=.4]
					\draw[edge,opacity=1,color=green] (v1.center) -- (v2.center) -- (v4.center) -- (v1.center);
					\fill[edge,opacity=1,color=green] (v1.center) -- (v2.center) -- (v4.center) -- (v1.center);
				\end{scope}
				\begin{scope}[transparency group,opacity=.4]
					\draw[edge,opacity=1,color=lightgray] (v1.center) -- (v3.center) -- (v4.center) -- (v1.center);
					\fill[edge,opacity=1,color=lightgray] (v1.center) -- (v3.center) -- (v4.center) -- (v1.center);
				\end{scope}
			\end{pgfonlayer}
		
			\node[text=blue] (a1) at (0,0.75) {$e_1$};
			\node[text=green] (a2) at (4,0.75) {$e_2$};
			\node[text=lightgray] (a3) at (0,-4.75) {$e_3$};
	\end{tikzpicture}\end{minipage}
	\begin{minipage}{.5\textwidth}\begin{flushright}\begin{tikzpicture}
				\tikzstyle{vertex} = [fill,shape=circle,node distance=80pt]
				\tikzstyle{edge} = [fill,opacity=.5,fill opacity=.5,line cap=round, line join=round, line width=30pt]
				\tikzstyle{elabel} =  [fill,shape=circle,node distance=40pt]
				
				\pgfdeclarelayer{background}
				\pgfsetlayers{background,main}
				
				\begin{scope}[every node/.style={circle,thick,draw}]
					\node (v1) at (0,0) {$v_1$};
					\node (v2) at (4,0) {$v_2$};
					\node (v3) at (0, -4) {$v_3$};
					\node (v4) at (4,-4) {$v_4$};
				\end{scope}
				
				\begin{scope}[>={Stealth[black]},
					every node/.style={circle},
					every edge/.style={draw = blue,very thick}]
					\path [-] (v1) edge[bend left=10] node[] {} (v2);
					\path [-] (v1) edge[bend right=10] node[] {} (v3);
					\path [-] (v2) edge node[] {} (v3);
				\end{scope}
				\begin{scope}[>={Stealth[black]},
					every node/.style={circle},
					every edge/.style={draw = green,very thick}]
					\path [-] (v1) edge[bend right=10] node[] {} (v2);
					\path [-] (v1) edge[bend left=10] node[] {} (v4);
					\path [-] (v2) edge node[] {} (v4);
				\end{scope}
				\begin{scope}[>={Stealth[black]},
					every node/.style={circle},
					every edge/.style={draw = lightgray,very thick}]
					\path [-] (v1) edge[bend left=10] node[] {} (v3);
					\path [-] (v1) edge[bend right=10] node[] {} (v4);
					\path [-] (v3) edge node[] {} (v4);
				\end{scope}
	\end{tikzpicture}\end{flushright}\end{minipage}
	
	\begin{minipage}{.5\textwidth}\begin{tikzpicture}
			\tikzstyle{vertex} = [fill,shape=circle,node distance=80pt]
			\tikzstyle{edge} = [fill,opacity=.5,fill opacity=.5,line cap=round, line join=round, line width=30pt]
			\tikzstyle{elabel} =  [fill,shape=circle,node distance=40pt]
			
			\pgfdeclarelayer{background}
			\pgfsetlayers{background,main}
			
			\begin{scope}[every node/.style={circle,thick,draw}]
				\node (v1) at (0,0) {$v_1$};
				\node (v2) at (4,0) {$v_2$};
				\node (v3) at (0, -4) {$v_3$};
				\node (v4) at (4,-4) {$v_4$};
			\end{scope}
			
			\begin{pgfonlayer}{background}
				\begin{scope}[transparency group,opacity=.4]
					\draw[edge,opacity=1,color=cyan] (v1.center) -- (v2.center) -- (v4.center) -- (v3.center) -- (v1.center);
					\fill[edge,opacity=1,color=cyan] (v1.center) -- (v2.center) -- (v4.center) -- (v3.center) -- (v1.center);
				\end{scope}
				\begin{scope}[transparency group,opacity=.4]
					\draw[edge,opacity=1,color=blue] (v1.center) -- (v2.center) -- (v1.center);
					\fill[edge,opacity=1,color=blue] (v1.center) -- (v2.center) -- (v1.center);
				\end{scope}
				\begin{scope}[transparency group,opacity=.4]
					\draw[edge,opacity=1,color=green] (v1.center) -- (v3.center) -- (v1.center);
					\fill[edge,opacity=1,color=green] (v1.center) -- (v3.center) -- (v1.center);
				\end{scope}
				\begin{scope}[transparency group,opacity=.4]
					\draw[edge,opacity=1,color=lightgray] (v1.center) -- (v4.center) -- (v1.center);
					\fill[edge,opacity=1,color=lightgray] (v1.center) -- (v4.center) -- (v1.center);
				\end{scope}
			\end{pgfonlayer}
			
			\node[text=blue] (a1) at (4,0.75) {$e_1$};
			\node[text=green] (a2) at (0,-4.75) {$e_2$};
			\node[text=lightgray] (a3) at (0,0.75) {$e_3$};
			\node[text=cyan] (a4) at (4,-4.75) {$e_4$};
	\end{tikzpicture}\end{minipage}
	\begin{minipage}{.5\textwidth}\begin{flushright}\begin{tikzpicture}
				\tikzstyle{vertex} = [fill,shape=circle,node distance=80pt]
				\tikzstyle{edge} = [fill,opacity=.5,fill opacity=.5,line cap=round, line join=round, line width=30pt]
				\tikzstyle{elabel} =  [fill,shape=circle,node distance=40pt]
				
				\pgfdeclarelayer{background}
				\pgfsetlayers{background,main}
				
				\begin{scope}[every node/.style={circle,thick,draw}]
					\node (v1) at (0,0) {$v_1$};
					\node (v2) at (4,0) {$v_2$};
					\node (v3) at (0, -4) {$v_3$};
					\node (v4) at (4,-4) {$v_4$};
				\end{scope}
				
				\begin{scope}[>={Stealth[black]},
					every node/.style={circle},
					every edge/.style={draw = blue,very thick}]
					\path [-] (v1) edge[bend left=20] node[] {} (v2);
				\end{scope}
				\begin{scope}[>={Stealth[black]},
					every node/.style={circle},
					every edge/.style={draw = green,very thick}]
					\path [-] (v1) edge[bend right=20] node[] {} (v3);
				\end{scope}
				\begin{scope}[>={Stealth[black]},
					every node/.style={circle},
					every edge/.style={draw = lightgray,very thick}]
					\path [-] (v1) edge[bend right=20] node[] {} (v4);
				\end{scope}
				\begin{scope}[>={Stealth[black]},
					every node/.style={circle},
					every edge/.style={draw = cyan,very thick}]
					\path [-] (v1) edge node[] {} (v2);
					\path [-] (v1) edge node[] {} (v3);
					\path [-] (v1) edge node[] {} (v4);
					\path [-] (v2) edge node[] {} (v3);
					\path [-] (v2) edge node[] {} (v4);
					\path [-] (v3) edge node[] {} (v4);
				\end{scope}
	\end{tikzpicture}\end{flushright}\end{minipage}
	
	Therefore, without further information, such as the number of hyperedges $\left\lvert\mathcal{E}_H\right\rvert$ or the minimum or maximum hyperedge cardinality, it is impossible to reconstruct the underlying hypergraph uniquely. Both hypergraphs have the same second normal graph representation $NG = \left(\mathcal{V}_G, \mathcal{E}_G\right)$ with:
	\begin{equation*}
		\mathcal{V}_G = \left\{v_1, v_2, v_3, v_4\right\}
	\end{equation*}
	\begin{equation*}
		\mathcal{E}_G = \left\{\left\{v_1, v_2\right\}, \left\{v_1, v_2\right\}, \left\{v_1, v_3\right\}, \left\{v_1, v_3\right\}, \left\{v_1, v_4\right\}, \left\{v_1, v_4\right\}, \right.
	\end{equation*}
	\begin{equation*}
		\left.\left\{v_2, v_3\right\}, \left\{v_2, v_4\right\}, \left\{v_3, v_4\right\}\right\}.
	\end{equation*}\vspace{0em}
\end{example}

Even though the second normal graph representation does not introduce new vertices and hence has no issue with a possibly exponential number of new vertices, the number of edges or arcs necessary for describing a hyperedge or hyperarc still grows quadratic in the number of vertices in the hyperedge or hyperarc.\\

\begin{remark}[\textbf{Drawback of the second normal graph representation}]\ \\
	Given a not oriented hypergraph $NH = \left(\mathcal{V}_H, \mathcal{E}_H\right)$ together with its normal graph representation $NG = \left(\mathcal{V}_G, \mathcal{E}_G\right)$, then the number of new edges necessary to describe a hyperedge $e_q \in \mathcal{E}_H$ is given by
	\begin{equation}
		\left(\begin{array}{c}
			\left\lvert e_q\right\rvert\\
			2
		\end{array}\right) = \frac{\left\lvert e_q\right\rvert !}{2 ! \left(\left\lvert e_q\right\rvert - 2\right) !} = \frac{\left\lvert e_q\right\rvert \left(\left\lvert e_q\right\rvert - 1\right)}{2} \in \mathcal{O}\left(\left\lvert e_q\right\rvert^2\right)
	\end{equation}
	since every pair of vertices from the hyperedge is connected by exactly one new edge \big(the order of the vertices does not matter\big) in the normal graph representation.\\
	
	Similarly, in the case of an oriented hypergraph $OH = \left(\mathcal{V}_H, \mathcal{A}_H\right)$ together with the normal graph representations $OG = \left(\mathcal{V}_G, \mathcal{A}_G\right)$, the number of new arcs necessary to describe a hyperarc $a_q = \left(a_q^{out}, a_q^{in}\right) \in \mathcal{A}_H$ is given by
	\begin{equation*}
		2 \left(\begin{array}{c}
			\left\lvert a_q^{out}\right\rvert\\
			2
		\end{array}\right)
		+ 2 \left(\begin{array}{c}
			\left\lvert a_q^{in}\right\rvert\\
			2
		\end{array}\right)
		+ \left\lvert a_q^{out}\right\rvert\left\lvert a_q^{in}\right\rvert =
	\end{equation*}
	\begin{equation}	
		\left\lvert a_q^{out}\right\rvert \left(\left\lvert a_q^{out}\right\rvert - 1\right) + \left\lvert a_q^{in}\right\rvert \left(\left\lvert a_q^{in}\right\rvert - 1\right) + \left\lvert a_q^{out}\right\rvert\left\lvert a_q^{in}\right\rvert \in \mathcal{O}\left(\left(\left\lvert a_q^{out}\right\rvert + \left\lvert a_q^{in}\right\rvert\right)^2\right)
	\end{equation}
	since every pair of vertices within the set of output vertices and input vertices of the hyperarc is connected by two arcs. Moreover, every output vertex is connected to each input vertex of the hyperarc, which results in additional $\left\lvert a_q^{out}\right\rvert\left\lvert a_q^{in}\right\rvert$ new arcs.\\
	
	Thus, the number of new edges or arcs necessary for encoding a hyperedge or hyperarc with the second normal graph representation is quadratic in the number of vertices in the hyperedge or hyperarc.\\
\end{remark}

The issue with the non-uniqueness of the mapping from the normal graph representation back to the hypergraph and the increasing number of vertices and $\backslash$ or edges and arcs in the normal graph representation clearly showcase why the extension of normal graphs to hypergraphs is necessary.\\

Hypergraphs are extremely relevant for modeling connections between multiple objects consistently, such as interactions within a social network or chemical reactions. Furthermore, the presented normal graph representations both suggest computational advantages in modeling complex links between objects with hypergraphs instead of normal graphs. These considerations prompt the extension of concepts from the normal graph setting to hypergraphs and at the same time justifies the more complex definitions on hypergraphs compared to normal graphs in later sections.\\
\clearpage
\section{Real functions on normal graphs}\label{5}

On the vertex set $\mathcal{V}$ and the edge set $\mathcal{E}_G$ or arc set $\mathcal{A}_G$ of normal graphs, real functions can be defined in order to encode further information about the original data set $v_1, v_2, \dots v_N$. These functions can be used to define gradients, adjoints, divergences, Laplacians and $p$-Laplacians in later sections.\\

Subsection (\ref{5.1}) introduces vertex functions, which can further describe the individual data points $v_1, v_2, \dots v_N$, for example with a vertex weight function. This subsection also includes the definition of the space of all vertex functions as a Hilbert space with the corresponding inner product, absolute value, and $\mathcal{L}^p$-norm of vertex functions.\\

Moreover, subsection (\ref{5.2}) defines edge functions on not oriented normal graphs and arc functions on oriented normal graphs for encoding more information about the pairwise relationship between the data points $v_1, v_2, \dots v_N$. As an example for edge and arc functions the edge weight and arc weight functions are introduced together with their symmetry property and the possibility to normalize them. As in the previous subsection, the definition of the space of all edge or arc functions as a Hilbert space with the corresponding inner product, absolute value, and $\mathcal{L}^p$-norm is given as well.

\subsection{Real vertex functions on normal graphs}\label{5.1}

To further describe the data points, interpreted as vertices $v_1, v_2, \dots v_N$, real functions on the vertex set $\mathcal{V}$ can be used. Extensive work about this topic was already published in \cite{elmoataz2015p} and hence this subsection generalizes their definitions and concepts.\\

\begin{definition}[\textbf{Real vertex function $f$}]\label{f} (Page 2418 in \cite{elmoataz2015p}: Real valued function on vertices) \ \\	
	A real vertex function $f$ is defined on the domain of the vertices $\mathcal{V}$ of a not oriented normal graph $NG = \left(\mathcal{V}, \mathcal{E}_G\right)$ or of an oriented normal graph $OG = \left(\mathcal{V}, \mathcal{A}_G\right)$ and can be written as
	\begin{equation}
		f: ~ \mathcal{V} \longrightarrow \mathbb{R} \qquad v_i \longmapsto f\left(v_i\right).
	\end{equation}\vspace{0em}
\end{definition}

An example for a real vertex function is the vertex weight function $w$ assigning a positive weight $w$ to each vertex $v_i \in \mathcal{V}$.\\

\clearpage
\begin{definition}[\textbf{Vertex weight function $w$}]\label{w} \ \\
	The vertex weight function $w$ can be defined for both, a not oriented normal graph $NG = \left(\mathcal{V}, \mathcal{E}_G\right)$ and an oriented normal graph $OG = \left(\mathcal{V}, \mathcal{A}_G\right)$, as
	\begin{equation}
		w: ~ \mathcal{V} \longrightarrow \mathbb{R}_{> 0} \qquad v_i \longmapsto w\left(v_i\right).
	\end{equation}
	A not oriented normal graph and an oriented normal graph together with a vertex weight function $w$, can be written as triple $NG = \left(\mathcal{V}, \mathcal{E}_G, w\right)$ or $OG = \left(\mathcal{V}, \mathcal{A}_G, w\right)$.\\
\end{definition}

Generally it is possible to define multiple vertex functions on the same normal graph. Therefore, in this thesis, different vertex functions are indicated by their subscript.\\

\begin{remark}[\textbf{Vertex weight functions $w_I$ and $w_G$}]\label{w_Iw_G} \ \\
	Since different vertex weight functions $w$ of a weighted oriented normal graph $OG = \left(\mathcal{V}, \mathcal{A}_G, w\right)$ will be used to define the inner product of the space of real vertex functions and the vertex and arc gradients, different subscripts will indicate if the vertex weight function stems from the inner product $w_I$ or the gradient definitions $w_G$. Both vertex weight functions $w_I$ and $w_G$ are defined on the same vertex set $\mathcal{V}$ of the weighted oriented normal graph $OG$.\\
\end{remark}

Other important examples for real vertex functions on not oriented normal graphs $NG = \left(\mathcal{V}, \mathcal{E}_G\right)$ or on oriented normal graphs $OG = \left(\mathcal{V}, \mathcal{A}_G\right)$ are the before defined vertex degree functions $\deg$, $\deg_{out}$, and $\deg_{in}$ of definition (\ref{degG}).\\

Since the number of vertices $|\mathcal{V}| = N \in \mathbb{N}$ is assumed to be finite, the space of all real vertex functions $f$ can be associated with an $N-$dimensional Hilbert space.\\

\begin{definition}[\textbf{Space of real vertex functions $\mathcal{H}\left(\mathcal{V}\right)$}]\label{H(V)} \ \\
	The space of real vertex functions $f$ is given by
	\begin{equation}
		\mathcal{H}\left(\mathcal{V}\right) = \left\{f ~ \middle| ~ f: ~ \mathcal{V} \longrightarrow \mathbb{R}\right\}
	\end{equation}
	with the corresponding inner product
	\begin{equation}
		{\langle f, g \rangle}_{\mathcal{H}\left(\mathcal{V}\right)} = \sum_{v_i \in \mathcal{V}} w_I\left(v_i\right)^\alpha f\left(v_i\right) g\left(v_i\right)
	\end{equation}
	for any two real vertex functions $f, g \in \mathcal{H}\left(\mathcal{V}\right)$, vertex weight function $w_I$, and parameter $\alpha \in \mathbb{R}$. \big(Note: $w_I\left(v_i\right)^\alpha = \left(w_I\left(v_i\right)\right)^\alpha$\big)\\
\end{definition}

\clearpage
\begin{remark}[\textbf{Parameter choice for the inner product of real vertex functions}]\label{inproH(V)} \ \\
	By choosing $\alpha = 0$, the common definition of the inner product of real vertex functions (Page 2418 in \cite{elmoataz2015p}: Inner product on $\mathcal{H}\left(\mathcal{V}\right)$) can be obtained
	\begin{equation}
		{\langle f, g \rangle}_{\mathcal{H}\left(\mathcal{V}\right)} = \sum_{v_i \in \mathcal{V}} f\left(v_i\right) g\left(v_i\right)
	\end{equation}
	for any two real vertex functions $f, g \in \mathcal{H}\left(\mathcal{V}\right)$.\\
\end{remark}

Moreover, based on the definitions for real valued continuous functions, the absolute value and the $\mathcal{L}^p$-norm for real vertex functions can be defined.\\

\begin{definition}[\textbf{Absolute value and $\mathcal{L}^p$-norm on the space of real vertex functions}]\label{Lpf} (Page 2418 in \cite{elmoataz2015p}: $\mathcal{L}^p\left(\mathcal{V}\right)$) \ \\
	For a real vertex function $f \in \mathcal{H}\left(\mathcal{V}\right)$ at a vertex $v_i \in \mathcal{V}$, the absolute value is defined as:
	\begin{equation*}
		\left\lvert ~ \cdot ~ \right\rvert: ~ \mathbb{R} \longrightarrow \mathbb{R}_{\geq 0}
	\end{equation*}
	\begin{equation}
		f\left(v_i\right) \longmapsto \left\lvert f\left(v_i\right)\right\rvert = \left\{\begin{array}{ll}
			f\left(v_i\right) & \quad f\left(v_i\right) \geq 0\\
			-f\left(v_i\right) & \quad \text{otherwise}
		\end{array}\right..
	\end{equation}
	Furthermore, the $\mathcal{L}^p$-norm of a real vertex function $f \in \mathcal{H}\left(\mathcal{V}\right)$ is defined as:
	\begin{equation*}
		\left\lvert \left\lvert ~ \cdot ~ \right\rvert\right\rvert_p: ~ \mathcal{H}\left(\mathcal{V}\right) \longrightarrow \mathbb{R}_{\geq 0}
	\end{equation*}
	\begin{equation}
		f \longmapsto \left\lvert \left\lvert f \right\rvert\right\rvert_p = \left\{\begin{array}{ll}
			\left(\sum_{v_i \in \mathcal{V}} \left\lvert f\left(v_i\right)\right\rvert^p\right)^{\frac{1}{p}} & \quad 1 \leq p < \infty\\
			\max_{v_i \in \mathcal{V}} \left(\left\lvert f\left(v_i\right) \right\rvert\right) & \quad p = \infty
		\end{array}\right..
	\end{equation}\vspace{0em}
\end{definition}
\subsection{Real edge and arc functions on normal graphs}\label{5.2}

Further information about the pairwise relationship between any two data points $v_i, v_j \in \left\{v_1, v_2, \dots v_N\right\}$ with $v_i \neq v_j$ can be encoded by edges or arcs and so-called edge functions or arc functions defined on the edge set $\mathcal{E}_G$ or the arc set $\mathcal{A}_G$. This subsection again generalizes the work of \cite{elmoataz2015p}.\\

\begin{definition}[\textbf{Real edge or arc function $F$}]\label{F} (Page 2418 in \cite{elmoataz2015p}: Real valued function on edges) \ \\
	A real edge or arc function $F$ is defined on the domain of the edges $\mathcal{E}_G$ of a not oriented normal graph $NG = \left(\mathcal{V}, \mathcal{E}_G\right)$ or on the domain of the arcs $\mathcal{A}_G$ of an oriented normal graph $OG = \left(\mathcal{V}, \mathcal{A}_G\right)$ as
	\begin{equation}
		F: ~ \mathcal{E}_G \longrightarrow \mathbb{R} \qquad e_q = \left\{v_i, v_j\right\} \longmapsto F\left(e_q\right) = F\left(\left\{v_i, v_j\right\}\right)
	\end{equation}
	or alternatively in case of an oriented normal graph as
	\begin{equation}
		F: ~ \mathcal{A}_G \longrightarrow \mathbb{R} \qquad a_q = \left(v_i, v_j\right) \longmapsto F\left(a_q\right) = F\left(\left(v_i, v_j\right)\right).
	\end{equation}\vspace{0em}
\end{definition}

Because of the lack of order in edges $\left\{v_i, v_j\right\} \in \mathcal{E}_G$, real edge functions are always symmetric.\\

\begin{remark}[\textbf{Symmetry of a real edge function}]\label{Fsymm} \ \\
	For a not oriented normal graph $NG = \left(\mathcal{V}, \mathcal{E}_G\right)$, an edge function $F: ~ \mathcal{E}_G \longrightarrow \mathbb{R}$ is always symmetric, which means
	\begin{equation}
		\forall \left\{v_i, v_j\right\} \in \mathcal{E}_G: ~ F\left(\left\{v_i, v_j\right\}\right) = F\left(\left\{v_j, v_i\right\}\right).
	\end{equation}

	This property is generally not fulfilled by a real arc function $F: ~ \mathcal{A}_G \longrightarrow \mathbb{R}$, hence $F\left(\left(v_i, v_j\right)\right) = F\left(\left(v_j, v_i\right)\right)$ can not be assumed to hold true for all arcs $\left(v_i, v_j\right) \in \mathcal{A}_G$. \big(Note: $\left(v_i, v_j\right) \in \mathcal{A}_G$ does not imply $\left(v_j, v_i\right) \in \mathcal{A}_G$ and therefore $F\left(\left(v_j, v_i\right)\right)$ might not even be well-defined.\big)\\
\end{remark}

An example for a real edge or arc function $F$ is the edge or arc weight function $W$, which assigns positive weights to edges or arcs.\\

\begin{definition}[\textbf{Edge or arc weight function $W$}]\label{W} \ \\
	The edge or arc weight function $W$ can be defined for both a not oriented normal graph $NG = \left(\mathcal{V}, \mathcal{E}_G\right)$ and an oriented normal graph $OG = \left(\mathcal{V}, \mathcal{A}_G\right)$:
	\begin{itemize}
		\item[1)] In case of a not oriented normal graph $NG$, the edge weight function $W$ is given by
		\begin{equation}
			W: ~ \mathcal{E}_G \longrightarrow \mathbb{R}_{> 0} \qquad e_q = \left\{v_i, v_j\right\} \longmapsto W\left(e_q\right) = W\left(\left\{v_i, v_j\right\}\right)
		\end{equation}
		with the automatic symmetry property $W\left(\left\{v_i, v_j\right\}\right) = W\left(\left\{v_j, v_i\right\}\right)$ for all edges $\left\{v_i, v_j\right\} \in \mathcal{E}_G$.
		\item[2)] Similarly, in case of an oriented normal graph $OG$, the arc weight function $W$ is defined on the domain of all arcs $\mathcal{A}_G$ instead of all edges $\mathcal{E}_G$
		\begin{equation}
			W: ~ \mathcal{A}_G \longrightarrow \mathbb{R}_{> 0} \qquad a_q = \left(v_i, v_j\right) \longmapsto W\left(a_q\right) = W\left(\left(v_i, v_j\right)\right)
		\end{equation}
		with the symmetry property $W\left(\left(v_i, v_j\right)\right) = W\left(\left(v_j, v_i\right)\right)$ not necessarily holding true for all arcs $\left(v_i, v_j\right) \in \mathcal{A}_G$.
	\end{itemize}
	A not oriented normal graph $NG = \left(\mathcal{V}, \mathcal{E}_G\right)$ together with an edge weight function $W$ or an oriented normal graph $OG = \left(\mathcal{V}, \mathcal{A}_G\right)$ together with an arc weight function $W$ is called weighted \big(not\big) oriented normal graph and can be written as triple $NG = \left(\mathcal{V}, \mathcal{E}_G, W\right)$ or $OG = \left(\mathcal{V}, \mathcal{A}_G, W\right)$.\\
\end{definition}

It is also possible to define multiple edge or arc functions on the same normal graph. Therefore, different edge or arc functions are indicated by their subscript.\\

\begin{remark}[\textbf{Arc weight functions $W_I$ and $W_G$}]\label{W_IW_G} \ \\
	Different arc weight functions $W$ of a weighted oriented normal graph $OG = \left(\mathcal{V}, \mathcal{A}_G, W\right)$ will be used for the inner product of the space of real arc functions and also for the vertex and arc gradient definitions. Thus, different subscripts will indicate if the arc weight function stems from the inner product $W_I$ or the gradient definitions $W_G$, however both arc weight functions $W_I$ and $W_G$ are defined on the same set of arcs $\mathcal{A}_G$ of the weighted oriented normal graph $OG$.\\
\end{remark}

As discussed before, a real arc function is not necessarily symmetric. Thus, the symmetry of an arc weight function $W$ on an oriented normal graph must not be taken for granted.\\

\begin{definition}[\textbf{Symmetric arc weight function}]\label{Wsymm} \ \\
	An arc weight function $W$ on an oriented normal graph $OG = \left(\mathcal{V}, \mathcal{A}_G\right)$ is called symmetric if for all arcs $a_q = \left(v_i, v_j\right)\in \mathcal{A}_G$ it holds true that
	\begin{equation}
		W\left(a_q\right) = W\left(\widetilde{a_q}\right)
	\end{equation}
	for $\widetilde{a_q} = \left(v_j, v_i\right) \in \widetilde{\mathcal{A}_G}$, where $\widetilde{\mathcal{A}_G}$ can be found in definition (\ref{switchedOG}) of the oriented normal graph with switched orientation.
	
	Note: This property either needs an extension of the domain of the arc weight function $W$ from $\mathcal{A}_G$ to $\mathcal{A}_G \cup \widetilde{\mathcal{A}_G}$, or an oriented normal graph with the property $\mathcal{A}_G = \widetilde{\mathcal{A}_G}$ which then implies for all arcs $\left(v_i, v_j\right) \in \mathcal{A}_G$ that $\left(v_j, v_i\right) \in \mathcal{A}_G$ also holds true.\\
\end{definition}

Every edge or arc weight function can be rescaled such that they only map to positive values less or equal to $1$.\\

\begin{theorem}[\textbf{Normalization of the edge or arc weight function}]\label{Wnorm} \ \\
	Every edge or arc weight function $W$ can be normalized in order to map to values in $\left(0, 1\right]$ instead of $\mathbb{R}_{> 0}$.\\
\end{theorem}

\begin{proof}\ \\
	All arguments of the following proof can be applied to both, a not oriented normal graph $NG = \left(\mathcal{V}, \mathcal{E}_G\right)$ and an oriented normal graph $OG = \left(\mathcal{V}, \mathcal{A}_G\right)$, since the orientation of the edges or arcs and the symmetry of $W$ are not used in the argumentation. Therefore, without loss of generality, it is assumed that $W$ is defined on the domain of the edge set $\mathcal{E}_G$ of a not oriented normal graph $NG$.\\
	
	\clearpage
	Then for every edge $e_q \in \mathcal{E}_G$ it holds true that $W\left(e_q\right) \in \mathbb{R}_{> 0}$ and hence
	\begin{equation*}
		m := \max_{e_q \in \mathcal{E}_G} W\left(e_q\right) \in \mathbb{R}_{> 0}
	\end{equation*}
	is well-defined, due to the finiteness of the edge set $\mathcal{E}_G$. With the maximum weight $m$, a new edge weight function $\overline{W}$ can be defined
	\begin{equation}
		\overline{W}: ~ \mathcal{E}_G \longrightarrow \left(0, 1\right] \qquad e_q = \left\{v_i, v_j\right\} \longmapsto \frac{W\left(e_q\right)}{m} = \frac{W\left(\left\{v_i, v_j\right\}\right)}{m},
	\end{equation}
	where $\frac{W\left(e_q\right)}{m} \in \left(0, 1\right]$ for all edges $e_q \in \mathcal{E}_G$ since the following properties are fulfilled:
	\begin{itemize}
		\item $W\left(e_q\right), m > 0 \quad \Longrightarrow \quad \overline{W}\left(e_q\right) = \frac{W\left(e_q\right)}{m} > 0$
		\item $W\left(e_q\right) \leq  m \quad \Longrightarrow \quad \overline{W}\left(e_q\right) = \frac{W\left(e_q\right)}{m} \leq \frac{m}{m} = 1$
	\end{itemize}
	Thus, any edge or arc weight function $W$ can be normalized in order to map to values in $\left(0, 1\right]$.\\
\end{proof}

The following remark shows that the number of edges in a not oriented normal graph $NG = \left(\mathcal{V}, \mathcal{E}_G\right)$ and the number of arcs in an oriented normal graph $OG = \left(\mathcal{V}, \mathcal{A}_G\right)$ are limited by $N^2$, which enables the association of the space of real edge or arc functions with an at most $N^2$-dimensional Hilbert space.\\

\begin{remark}[\textbf{Number of edges or arcs in normal graphs}]\label{numbedgearc} \ \\
	As described in remark (\ref{Gfinite}), it is assumed that no edges or arcs exist more than once. Together with the property $\mathcal{E}_G \subseteq \left\{\left\{v_i, v_j\right\} ~ \middle| ~ v_i, v_j \in \mathcal{V}\right\}$, which means that each vertex $v_i \in \mathcal{V}$ is at most connected to all other vertices $v_j \in \mathcal{V} \backslash\left\{v_i\right\}$, this results in the number of edges in a not oriented normal graph $NG = \left(\mathcal{V}, \mathcal{E}_G\right)$ being limited by
	\begin{equation}
		\left\lvert \mathcal{E}_G \right\rvert \leq \left\lvert \mathcal{V} \right\rvert \cdot  \left\lvert\mathcal{V} \right\rvert = N^2.
	\end{equation}
	
	Similarly, the number of arcs in an oriented normal graph $OG = \left(\mathcal{V}, \mathcal{A}_G\right)$ is limited by
	\begin{equation}
		\left\lvert \mathcal{A}_G \right\rvert \leq \left\lvert \mathcal{V} \right\rvert \cdot \left\lvert\mathcal{V} \right\rvert = N^2
	\end{equation}
	since $\mathcal{A}_G \subseteq \left\{\left(v_i, v_j\right) ~ \middle| ~ v_i, v_j \in \mathcal{V}\right\} = \mathcal{V} \times \mathcal{V}$, which means that each vertex $v_i \in \mathcal{V}$ is at most connected to all other vertices $v_j \in \mathcal{V} \backslash\left\{v_i\right\}$, once as the output vertex in $\left(v_i, v_j\right)$ and once as the input vertex in $\left(v_j, v_i\right)$.\\
\end{remark}

With the extension of the domain of all edge functions from $\mathcal{E}_G$ to $\left\{\left\{v_i, v_j\right\} ~ \middle| ~ v_i, v_j \in \mathcal{V}\right\}$ and the extension of the domain of all arc functions from $\mathcal{A}_G$ to $\mathcal{V} \times \mathcal{V}$ by setting $F\left(\left\{v_i, v_j\right\}\right) = 0$ for any edge $\left\{v_i, v_j\right\} \notin \mathcal{E}_G$ and $F\left(\left(v_i, v_j\right)\right) = 0$ for any arc $\left(v_i, v_j\right) \notin \mathcal{A}_G$, the space of all real edge or arc functions can be identified as an $N^2$-dimensional Hilbert space.\\

\begin{definition}[\textbf{Space of real edge functions $\mathcal{H}\left(\mathcal{E}_G\right)$ and real arc functions $\mathcal{H}\left(\mathcal{A}_G\right)$}]\label{H(E)H(A)} \ \\
	The space of real edge functions $F$ for a not oriented normal graph $NG = \left(\mathcal{V}, \mathcal{E}_G, W\right)$ is given by
	\begin{equation}
		\mathcal{H}\left(\mathcal{E}_G\right) = \left\{F ~ \middle| ~ F: ~ \mathcal{E}_G \longrightarrow \mathbb{R}\right\}
	\end{equation}
	with the corresponding inner product
	\begin{equation}
		{\langle F, G \rangle}_{\mathcal{H}\left(\mathcal{E}_G\right)} = \sum_{e_q \in \mathcal{E}_G} W_I \left(e_q\right)^\beta F\left(e_q\right) G\left(e_q\right)
	\end{equation}
	for any two real edge functions $F, G \in \mathcal{H}\left(\mathcal{E}_G\right)$ and parameter $\beta \in \mathbb{R}$. Similarly, in the case of an oriented normal graph $OG = \left(\mathcal{V}, \mathcal{A}_G, W\right)$, the space of real arc functions is defined as
	\begin{equation}
		\mathcal{H}\left(\mathcal{A}_G\right) = \left\{F ~ \middle| ~ F: ~ \mathcal{A}_G \longrightarrow \mathbb{R}\right\}
	\end{equation}
	with the corresponding inner product
	\begin{equation}
		{\langle F, G \rangle}_{\mathcal{H}\left(\mathcal{A}_G\right)} = \sum_{a_q \in \mathcal{A}_G} W_I \left(a_q\right)^\beta F\left(a_q\right) G\left(a_q\right)
	\end{equation}
	for any two real arc functions $F, G \in \mathcal{H}\left(\mathcal{A}_G\right)$ and parameter $\beta \in \mathbb{R}$.\\
	
	\big(Note: $W_I \left(e_q\right)^\beta = \left(W_I \left(e_q\right)\right)^\beta$ and $W_I \left(a_q\right)^\beta = \left(W_I \left(a_q\right)\right)^\beta$\big)\\
\end{definition}

\begin{remark}[\textbf{Parameter choice for the inner product of real edge and arc functions}]\label{inproH(E)H(A)} \ \\
	Choosing $\beta = 1$, $W\left(e_q\right) \equiv \frac{1}{2}$ for all edges $e_q \in \mathcal{E}_G$ and $W\left(a_q\right) \equiv \frac{1}{2}$ for all arcs $a_q \in \mathcal{A}_G$ results in the common definition of the inner product of real edge functions (Page 2418 in \cite{elmoataz2015p}: Inner product on $\mathcal{H}\left(\mathcal{E}\right)$)
	\begin{equation}
		{\langle F, G \rangle}_{\mathcal{H}\left(\mathcal{E}_G\right)} = \frac{1}{2} \sum_{e_q \in \mathcal{E}_G} F\left(e_q\right) G\left(e_q\right)
	\end{equation}
	for any two real edge functions $F, G \in \mathcal{H}\left(\mathcal{E}_G\right)$. Similarly, the common definition of the inner product of real arc functions can be obtained
	\begin{equation}
		{\langle F, G \rangle}_{\mathcal{H}\left(\mathcal{A}_G\right)} = \frac{1}{2} \sum_{a_q \in \mathcal{A}_G} F\left(a_q\right) G\left(a_q\right)
	\end{equation}
	for any two real arc functions $F, G \in \mathcal{H}\left(\mathcal{A}_G\right)$.\\
\end{remark}

Furthermore, the absolute value and the $\mathcal{L}^p$-norm for real edge and arc functions can be defined, which are necessary for later sections.\\

\begin{definition}[\textbf{Absolute value and $\mathcal{L}^p$-norm on the space of real edge and arc functions}]\label{LpF} \ \\
	The absolute value and the $\mathcal{L}^p$-norm can be defined for both a real edge function $F \in \mathcal{H}\left(\mathcal{E}_G\right)$ and a real arc function $F \in \mathcal{H}\left(\mathcal{A}_G\right)$:
	\begin{itemize}
		\item[1)] For a real edge function $F \in \mathcal{H}\left(\mathcal{E}_G\right)$ at an edge $e_q \in \mathcal{E}_G$, the absolute value is defined as
		\begin{equation*}
			\left\lvert ~ \cdot ~ \right\rvert: ~ \mathbb{R} \longrightarrow \mathbb{R}_{\geq 0}
		\end{equation*}
		\begin{equation}
			F\left(e_q\right) \longmapsto \left\lvert F\left(e_q\right)\right\rvert = \left\{\begin{array}{ll}
				F\left(e_q\right) & \quad F\left(e_q\right) \geq 0\\
				-F\left(e_q\right) & \quad \text{otherwise}
			\end{array}\right.
		\end{equation}
	
		and the $\mathcal{L}^p$-norm of the edge function $F$ is given by
		\begin{equation*}
			\left\lvert \left\lvert ~ \cdot ~ \right\rvert\right\rvert_p: ~ \mathcal{H}\left(\mathcal{E}_G\right) \longrightarrow \mathbb{R}_{\geq 0}
		\end{equation*}
		\begin{equation}
			F \longmapsto \left\lvert \left\lvert F \right\rvert\right\rvert_p = \left\{\begin{array}{ll}
				\left(\frac{1}{2} \sum_{e_q \in \mathcal{E}_G} \left\lvert F\left(e_q\right)\right\rvert^p\right)^{\frac{1}{p}} & \quad 1 \leq p < \infty\\
				\max_{e_q \in \mathcal{E}_G} \left(\left\lvert F\left(e_q\right) \right\rvert\right) & \quad p = \infty
			\end{array}\right..
		\end{equation}\vspace{0em}
		
		\item[2)] Similarly, for a real arc function $F \in \mathcal{H}\left(\mathcal{A}_G\right)$ at an arc $a_q \in \mathcal{A}_G$, the absolute value is defined as
		\begin{equation*}
			\left\lvert ~ \cdot ~ \right\rvert: ~ \mathbb{R} \longrightarrow \mathbb{R}_{\geq 0}
		\end{equation*}
		\begin{equation}
			F\left(a_q\right) \longmapsto \left\lvert F\left(a_q\right)\right\rvert = \left\{\begin{array}{ll}
				F\left(a_q\right) & \quad F\left(a_q\right) \geq 0\\
				-F\left(a_q\right) & \quad \text{otherwise}
			\end{array}\right.
		\end{equation}
	
		and the $\mathcal{L}^p$-norm of the arc function $F$ is given by
		\begin{equation*}
			\left\lvert \left\lvert ~ \cdot ~ \right\rvert\right\rvert_p: ~ \mathcal{H}\left(\mathcal{A}_G\right) \longrightarrow \mathbb{R}_{\geq 0}
		\end{equation*}
		\begin{equation}
			F \longmapsto \left\lvert \left\lvert F \right\rvert\right\rvert_p = \left\{\begin{array}{ll}
				\left(\frac{1}{2} \sum_{a_q \in \mathcal{A}_G} \left\lvert F\left(a_q\right)\right\rvert^p\right)^{\frac{1}{p}} & \quad 1 \leq p < \infty\\
				\max_{a_q \in \mathcal{A}_G} \left(\left\lvert F\left(a_q\right) \right\rvert\right) & \quad p = \infty
			\end{array}\right..
		\end{equation}\vspace{0em}
	\end{itemize}

\end{definition}

In order to simplify the notation in more complex formulas, any double brackets are from now on replaced with single brackets.

\begin{remark}[\textbf{Simplification of the notation of arcs}]\label{arcnot} \ \\
	For simplification purposes, the double brackets in arc functions are from now on omitted
	\begin{equation*}
		F\left(\left(v_i, v_j\right)\right) \longleftarrow F\left(v_i, v_j\right)
	\end{equation*}
	for every arc $\left(v_i, v_j\right) \in \mathcal{A}_G$ of an oriented normal graph $OG = \left(\mathcal{V}, \mathcal{A}_G\right)$.\\
\end{remark}
\clearpage
\section{Gradient and adjoint operators on oriented normal graphs}\label{6}

Using the definitions of vertex functions and edge or arc functions, the standard differential operators from the continuum setting can be extended to the discrete setting of a weighted oriented normal graph $OG = \left(\mathcal{V}, \mathcal{A}_G, w, W\right)$ with vertex weight function $w$ and arc weight function $W$.\\

Thus, this section introduces the gradient and adjoint operators on weighted oriented normal graphs, which are needed for the divergence, Laplacian and $p$-Laplacian operators later on. Subsection (\ref{6.1}) defines the vertex gradient operator and the arc gradient operator as well as their specific properties, which match the common properties of gradients in the continuum setting. Next, subsection (\ref{6.2}) derives the corresponding vertex adjoint operator and arc adjoint operator based on the before defined vertex and arc gradients. In both subsections the definitions are chosen in a way such that they match the vertex gradient and vertex adjoint defined in \cite{elmoataz2015p}.\\

\subsection{Gradient operators on oriented normal graphs}\label{6.1}

In case of a weighted oriented normal graph $OG = \left(\mathcal{V}, \mathcal{A}_G, w, W\right)$, the gradient operator can be defined both for the set of all vertices $\mathcal{V}$ and for the set of all arcs $\mathcal{A}_G$.\\

\begin{definition}[\textbf{Vertex gradient operator $\nabla_v$}]\label{Gnabla_v} \ \\
	For a weighted oriented normal graph $OG = \left(\mathcal{V}, \mathcal{A}_G, w, W\right)$, the vertex gradient operator $\nabla_v$ is given by
	\begin{equation}
		\nabla_v: ~ \mathcal{H}\left(\mathcal{V}\right) \longrightarrow \mathcal{H}\left(\mathcal{A}_G\right) \qquad f \longmapsto \nabla_v f
	\end{equation}
	with the weighted difference of vertex function $f \in \mathcal{H}\left(\mathcal{V}\right)$ along an arc $a_q = \left(v_i, v_j\right) \in \mathcal{A}_G$ being defined as
	\begin{equation*}
		\nabla_v f: ~ \mathcal{A}_G \longrightarrow \mathbb{R} \qquad a_q = \left(v_i, v_j\right) \longmapsto \nabla_v f \left(a_q\right) = 
	\end{equation*}
	\begin{equation}
		W_G \left(v_i, v_j\right)^\gamma \left(w_I \left(v_j\right)^\alpha w_G \left(v_j\right)^\epsilon f\left(v_j\right) - w_I \left(v_i\right)^\alpha w_G \left(v_i\right)^\eta f\left(v_i\right)\right).
	\end{equation}\vspace{0em}
\end{definition}

\begin{remark}[\textbf{Parameter choice for the vertex gradient operator}]\label{Gnabla_vparam} \ \\
	By choosing $\alpha = 0$, $\gamma = \frac{1}{2}$, $\epsilon = 0$ and $\eta = 0$ the common definition of the vertex gradient (Page 2418 in \cite{elmoataz2015p}: Weighted partial difference) for all arcs $a_q = \left(v_i, v_j\right) \in \mathcal{A}_G$ can be obtained:
	\begin{equation}
		\nabla_v f \left(a_q\right) = \nabla_v f \left(v_i, v_j\right) = \sqrt{W_G \left(v_i, v_j\right)} \left(f\left(v_j\right) - f\left(v_i\right)\right).
	\end{equation}\vspace{0em}
\end{remark}

The weighted difference $\nabla_v f$ has similar properties compared to gradients in the continuum setting.\\

\begin{theorem}[\textbf{Vertex gradient operator properties}]\label{Gnabla_vprop} \ \\
	The weighted difference of a vertex function $f \in \mathcal{H}\left(\mathcal{V}\right)$ on a weighted oriented normal graph $OG = \left(\mathcal{V}, \mathcal{A}_G, w, W\right)$ fulfills the following properties:
	\begin{itemize}
		\item[1)] Gradient of a constant vertex function is equal to zero:
		
		$f\left(v_i\right) \equiv \overline{f} \in \mathbb{R}$ for all vertices $v_i \in \mathcal{V}$ and $w_I \left(v_k\right)^\alpha w_G \left(v_k\right)^\epsilon = w_I \left(v_j\right)^\alpha w_G \left(v_j\right)^\eta$ for all arcs $a_q = \left(v_j, v_k\right) \in \mathcal{A}_G$ \quad $\Longrightarrow$ \quad $\nabla_v f \left(a_q\right) = 0$ for all arcs $a_q = \left(v_j, v_k\right) \in \mathcal{A}_G$
		
		\item[2)] Antisymmetry:
		
		Symmetric arc weight function $W_G$ and $\epsilon = \eta$ \quad $\Longrightarrow$ \quad $\nabla_v f \left(a_q\right) = - \nabla_v f \left(\widetilde{a_q}\right)$ for all arcs $a_q \in \mathcal{A}_G$ \big(with the corresponding arc $\widetilde{a_q} \in \widetilde{\mathcal{A}_G}$ of the graph with switched orientation $\widetilde{OG}$ from definition (\ref{switchedOG})\big)\\
	\end{itemize}
\end{theorem}

\begin{proof}\ \\
	Given a vertex function $f \in \mathcal{H}\left(\mathcal{V}\right)$ on a weighted oriented normal graph $OG = \left(\mathcal{V}, \mathcal{A}_G, w, W\right)$, then it holds true that:
	\begin{itemize}
		\item[1)] $f\left(v_i\right) \equiv \overline{f} \in \mathbb{R}$ for all vertices $v_i \in \mathcal{V}$ and $w_I \left(v_k\right)^\alpha w_G \left(v_k\right)^\epsilon = w_I \left(v_j\right)^\alpha w_G \left(v_j\right)^\eta$ for all arcs $a_q = \left(v_j, v_k\right) \in \mathcal{A}_G$ implies that:
		
		$\begin{aligned}[t]
			\nabla_v f \left(a_q\right) & = W_G \left(v_j, v_k\right)^\gamma \left(w_I \left(v_k\right)^\alpha w_G \left(v_k\right)^\epsilon f\left(v_k\right) - w_I \left(v_j\right)^\alpha w_G \left(v_j\right)^\eta f\left(v_j\right)\right) &\\
			& = W_G \left(v_j, v_k\right)^\gamma \left(w_I \left(v_k\right)^\alpha w_G \left(v_k\right)^\epsilon \overline{f} - w_I \left(v_k\right)^\alpha w_G \left(v_k\right)^\epsilon \overline{f}\right) &\\
			& = W_G \left(v_j, v_k\right)^\gamma \cdot 0 &\\
			& = 0 &\\
		\end{aligned}$
	
		Where the last equality is feasible due to the arc weight function $W_G$ being a real function and therefore mapping to finite values.\\

		\item[2)] $W_G$ being a symmetric arc weight function and $\epsilon = \eta$ holding true, implies for every arc $a_q = \left(v_i, v_j\right) \in \mathcal{A}_G$:
		
		$\begin{aligned}[t]
			\nabla_v f \left(a_q\right) & = \nabla_v f \left(v_i, v_j\right) &\\
			& = W_G \left(v_i, v_j\right)^\gamma \left(w_I \left(v_j\right)^\alpha w_G \left(v_j\right)^\epsilon f\left(v_j\right) - w_I \left(v_i\right)^\alpha w_G \left(v_i\right)^\eta f\left(v_i\right)\right) &\\
			& = - W_G \left(v_i, v_j\right)^\gamma \left(w_I \left(v_i\right)^\alpha w_G \left(v_i\right)^\eta f\left(v_i\right) - w_I \left(v_j\right)^\alpha w_G \left(v_j\right)^\epsilon f\left(v_j\right)\right) &\\
			& = - W_G \left(v_i, v_j\right)^\gamma \left(w_I \left(v_i\right)^\alpha w_G \left(v_i\right)^\epsilon f\left(v_i\right) - w_I \left(v_j\right)^\alpha w_G \left(v_j\right)^\eta f\left(v_j\right)\right) &\\
			& = - W_G \left(v_j, v_i\right)^\gamma \left(w_I \left(v_i\right)^\alpha w_G \left(v_i\right)^\epsilon f\left(v_i\right) - w_I \left(v_j\right)^\alpha w_G \left(v_j\right)^\eta f\left(v_j\right)\right) &\\
			& = - \nabla_v f \left(v_j, v_i\right) = - \nabla_v f \left(\widetilde{a_q}\right) &\\
		\end{aligned}$
		
		Where the equality from the fourth to the fifth line is obtained by using the symmetry of the arc weight function $W_G$.
	\end{itemize}
\end{proof}

Similar to the vertex gradient operator $\nabla_v$, a linear gradient operator can be defined on the domain of arc functions $\mathcal{H}\left(\mathcal{A}_G\right)$.\\

\begin{definition}[\textbf{Arc gradient operator $\nabla_a$}]\label{Gnabla_a} \ \\
	For a weighted oriented normal graph $OG = \left(\mathcal{V}, \mathcal{A}_G, w, W\right)$, the arc gradient operator $\nabla_a$ is defined as
	\begin{equation}
		\nabla_a: ~ \mathcal{H}\left(\mathcal{A}_G\right) \longrightarrow \mathcal{H}\left(\mathcal{V}\right) \qquad F \longmapsto \nabla_a F
	\end{equation}
	with the weighted difference of arc function $F \in \mathcal{H}\left(\mathcal{A}_G\right)$ at a vertex $v_i \in \mathcal{V}$ being given as
	\begin{equation*}
		\nabla_a F: ~ \mathcal{V} \longrightarrow \mathbb{R} \qquad v_i \longmapsto \nabla_a F \left(v_i\right) = 
	\end{equation*}
	\begin{equation}
		w_G \left(v_i\right)^\zeta \sum_{a_q \in \mathcal{A}_G} \left(\frac{\delta_{in}\left(v_i, a_q\right)}{\deg_{in}\left(v_i\right)} - \frac{\delta_{out}\left(v_i, a_q\right)}{\deg_{out}\left(v_i\right)}\right) W_I \left(a_q\right)^\beta W_G \left(a_q\right)^\theta F\left(a_q\right).
	\end{equation}\vspace{0em}
\end{definition}

\begin{remark}[\textbf{Parameter choice for the arc gradient operator}]\label{Gnabla_aparam} \ \\
	Choosing $\beta = 1$, $W_I \left(a_q\right) \equiv \frac{1}{2}$ for all arcs $a_q \in \mathcal{A}_G$, $\zeta = \frac{1}{2}$ and $\theta = 0$ together with leaving out the degree functions $\deg_{in}\left(v_i\right)$ and $\deg_{out}\left(v_i\right)$ for all vertices $v_i \in \mathcal{V}$, results in a simplified definition of the arc gradient for all vertices $v_i \in \mathcal{V}$:
	\begin{equation}
		\nabla_a F \left(v_i\right) = \frac{\sqrt{w_G \left(v_i\right)}}{2} \sum_{a_q \in \mathcal{A}_G} \left(\delta_{in}\left(v_i, a_q\right) - \delta_{out}\left(v_i, a_q\right)\right) F\left(a_q\right).
	\end{equation}\vspace{1em}
\end{remark}

As the weighted difference $\nabla_v f$ of vertex functions $f \in \mathcal{H}\left(\mathcal{V}\right)$, the weighted difference $\nabla_a F$ of arc functions $F \in \mathcal{H}\left(\mathcal{A}_G\right)$ also has similar properties compared to gradients in the continuum setting.\\

\begin{theorem}[\textbf{Arc gradient operator properties}]\label{Gnabla_aprop} \ \\
	The weighted difference of an arc function $F \in \mathcal{H}\left(\mathcal{A}_G\right)$ on a weighted oriented normal graph $OG = \left(\mathcal{V}, \mathcal{A}_G, w, W\right)$ fulfills the following properties:
	\begin{itemize}
		\item[1)] Gradient of a constant arc function is equal to zero:
		
		$F\left(a_q\right) \equiv \overline{F} \in \mathbb{R}$ and $W_I \left(a_q\right)^\beta W_G \left(a_q\right)^\theta \equiv \overline{W} \in \mathbb{R}$ for all arcs $a_q \in \mathcal{A}_G$
		
		$\Longrightarrow$ \quad $\nabla_a F \left(v_i\right) = 0$ for all vertices $v_i \in \mathcal{V}$
		
		\item[2)] Antisymmetry:
		
		Symmetric arc weight functions $W_I$ and $W_G$ and $F\left(a_q\right) = F\left(\widetilde{a_q}\right)$ for all arcs $a_q \in \mathcal{A}_G$, $\widetilde{a_q} \in \widetilde{\mathcal{A}_G}$
		
		$\Longrightarrow$ \quad $\nabla_a F \left(v_i\right) = - \nabla_{\widetilde{a}} F \left(v_i\right)$ for all vertices $v_i \in \mathcal{V}$ \big(with $\nabla_{\widetilde{a}} F$ being the arc gradient applied to the graph with switched orientation $\widetilde{OG}$\big)\\
	\end{itemize}
\end{theorem}

\begin{proof}\ \\
	Given an arc function $F \in \mathcal{H}\left(\mathcal{A}_G\right)$ on a weighted oriented normal graph $OG = \left(\mathcal{V}, \mathcal{A}_G, w, W\right)$, then it holds true that:
	\begin{itemize}
		\item[1)] $F\left(a_q\right) \equiv \overline{F} \in \mathbb{R}$ and $W_I \left(a_q\right)^\beta W_G \left(a_q\right)^\theta \equiv \overline{W} \in \mathbb{R}$ for all arcs $a_q \in \mathcal{A}_G$ implies for all vertices $v_i \in \mathcal{V}$ that:
		
		$\begin{aligned}[t]
			\nabla_a F \left(v_i\right) & = w_G \left(v_i\right)^\zeta \sum_{a_q \in \mathcal{A}_G} \left(\frac{\delta_{in}\left(v_i, a_q\right)}{\deg_{in}\left(v_i\right)} - \frac{\delta_{out}\left(v_i, a_q\right)}{\deg_{out}\left(v_i\right)}\right) W_I \left(a_q\right)^\beta W_G \left(a_q\right)^\theta F\left(a_q\right) &\\
			& = w_G \left(v_i\right)^\zeta \sum_{a_q \in \mathcal{A}_G} \left(\frac{\delta_{in}\left(v_i, a_q\right)}{\deg_{in}\left(v_i\right)} - \frac{\delta_{out}\left(v_i, a_q\right)}{\deg_{out}\left(v_i\right)}\right) \overline{W} ~ \overline{F} &\\
			& = w_G \left(v_i\right)^\zeta \left(\sum_{a_q \in \mathcal{A}_G} \frac{\delta_{in}\left(v_i, a_q\right)}{\deg_{in}\left(v_i\right)} - \sum_{a_r \in \mathcal{A}_G} \frac{\delta_{out}\left(v_i, a_r\right)}{\deg_{out}\left(v_i\right)}\right) \overline{W} ~ \overline{F} &\\
			& = w_G \left(v_i\right)^\zeta \left( \frac{\sum_{a_q \in \mathcal{A}_G} \delta_{in}\left(v_i, a_q\right)}{\deg_{in}\left(v_i\right)} - \frac{\sum_{a_r \in \mathcal{A}_G} \delta_{out}\left(v_i, a_r\right)}{\deg_{out}\left(v_i\right)}\right) \overline{W} ~ \overline{F} &\\
			& = w_G \left(v_i\right)^\zeta \left( \frac{\deg_{in}\left(v_i\right)}{\deg_{in}\left(v_i\right)} - \frac{\deg_{out}\left(v_i\right)}{\deg_{out}\left(v_i\right)}\right) \overline{W} ~ \overline{F} &\\
			& = w_G \left(v_i\right)^\zeta \cdot 0 \cdot \overline{W} ~ \overline{F} &\\
			& = 0 &\\
		\end{aligned}$
	
		Where the last equality is feasible due to the vertex weight function $w_G$ being a real function and hence mapping to finite values.\\
		
		\item[2)] $W_I$ and $W_G$ being symmetric arc weight functions and $F\left(a_q\right) = F\left(\widetilde{a_q}\right)$ for all arcs $a_q \in \mathcal{A}_G$, $\widetilde{a_q} \in \widetilde{\mathcal{A}_G}$ implies for all vertices $v_i \in \mathcal{V}$:
		
		$\begin{aligned}[t]
			\nabla_a F \left(v_i\right) & = w_G \left(v_i\right)^\zeta \sum_{a_q \in \mathcal{A}_G} \left(\frac{\delta_{in}\left(v_i, a_q\right)}{\deg_{in}\left(v_i\right)} - \frac{\delta_{out}\left(v_i, a_q\right)}{\deg_{out}\left(v_i\right)}\right) W_I \left(a_q\right)^\beta W_G \left(a_q\right)^\theta F\left(a_q\right) &\\
			& = - w_G \left(v_i\right)^\zeta \sum_{a_q \in \mathcal{A}_G} \left(\frac{\delta_{out}\left(v_i, a_q\right)}{\deg_{out}\left(v_i\right)} - \frac{\delta_{in}\left(v_i, a_q\right)}{\deg_{in}\left(v_i\right)}\right) W_I \left(a_q\right)^\beta W_G \left(a_q\right)^\theta F\left(a_q\right) &\\
			& = - w_G \left(v_i\right)^\zeta \sum_{a_q \in \mathcal{A}_G} \left(\frac{\delta_{out}\left(v_i, a_q\right)}{\deg_{out}\left(v_i\right)} - \frac{\delta_{in}\left(v_i, a_q\right)}{\deg_{in}\left(v_i\right)}\right) W_I \left(\widetilde{a_q}\right)^\beta W_G \left(\widetilde{a_q}\right)^\theta F\left(\widetilde{a_q}\right) &\\
			& = - w_G \left(v_i\right)^\zeta \sum_{\widetilde{a_q} \in \widetilde{\mathcal{A}_G}} \left(\frac{\delta_{in}\left(v_i, \widetilde{a_q}\right)}{\widetilde{\deg_{in}}\left(v_i\right)} - \frac{\delta_{out}\left(v_i, \widetilde{a_q}\right)}{\widetilde{\deg_{out}}\left(v_i\right)}\right) W_I \left(\widetilde{a_q}\right)^\beta W_G \left(\widetilde{a_q}\right)^\theta F\left(\widetilde{a_q}\right) &\\
			& = - \nabla_{\widetilde{a}} F \left(v_i\right)\\
		\end{aligned}$
		
		Where the equality from the second to the third line is obtained by using the symmetry of the arc weight functions $W_I$ and $W_G$, while the reformulation from the third to the fourth line holds true due to fact that $\delta_{out}\left(v_i, a_q\right) = \delta_{in}\left(v_i, \widetilde{a_q}\right)$ and $\delta_{in}\left(v_i, a_q\right) = \delta_{out}\left(v_i, \widetilde{a_q}\right)$ for all vertices $v_i \in \mathcal{V}$, since the output and input vertices of any arc $a_q \in \mathcal{A}_G$ are exchanged to get $\widetilde{a_q} \in \widetilde{\mathcal{A}_G}$.
		
		Furthermore, $\deg_{out}\left(v_i\right)$ becomes $\widetilde{\deg_{in}}\left(v_i\right)$, when switching from $OG = \left(\mathcal{V}, \mathcal{A}_G\right)$ to $\widetilde{OG} = \left(\mathcal{V}, \widetilde{\mathcal{A}_G}\right)$ and respectively $\deg_{in}\left(v_i\right)$ becomes $\widetilde{\deg_{out}}\left(v_i\right)$ for all vertices $v_i \in \mathcal{V}$.\\
	\end{itemize}
\end{proof}
\subsection{Adjoint operators on oriented normal graphs}\label{6.2}

On a weighted oriented normal graph $OG = \left(\mathcal{V}, \mathcal{A}_G, w, W\right)$, the adjoint operators $\nabla^*_v$ and $\nabla^*_a$ are the dual counterpart to the weighted gradient operators $\nabla_v$ and $\nabla_a$ with regard to the inner products defined on $\mathcal{H}\left(\mathcal{V}\right)$ and $\mathcal{H}\left(\mathcal{A}_G\right)$. Thus, the adjoint operators fulfill the following equalities for all vertex functions $f \in \mathcal{H}\left(\mathcal{V}\right)$ and all arc functions $G \in \mathcal{H}\left(\mathcal{A}_G\right)$:\\
\begin{equation}
	{\langle G, \nabla_v f \rangle}_{\mathcal{H}\left(\mathcal{A}_G\right)} = {\langle f, \nabla^*_v G \rangle}_{\mathcal{H}\left(\mathcal{V}\right)}
\end{equation}
\begin{equation}
	{\langle f, \nabla_a G \rangle}_{\mathcal{H}\left(\mathcal{V}\right)} = {\langle G, \nabla^*_a f \rangle}_{\mathcal{H}\left(\mathcal{A}_G\right)}
\end{equation}

Based on the definitions of the adjoint operators $\nabla^*_v$ and $\nabla^*_a$, Laplacian and $p$-Laplacian operators for vertices and arcs can be defined later on.\\

\clearpage
\begin{definition}[\textbf{Vertex adjoint operator $\nabla^*_v$}]\label{Gnabla^*_v} \ \\
	For a weighted oriented normal graph $OG = \left(\mathcal{V}, \mathcal{A}_G, w, W\right)$ with symmetric arc weight functions $W_I$ and $W_G$, the vertex adjoint operator $\nabla^*_v$ is given by
	\begin{equation}
		\nabla^*_v: ~ \mathcal{H}\left(\mathcal{A}_G\right) \longrightarrow \mathcal{H}\left(\mathcal{V}\right) \qquad F \longmapsto \nabla^*_v F
	\end{equation}
	with the weighted adjoint of arc function $F \in \mathcal{H}\left(\mathcal{A}_G\right)$ at a vertex $v_i \in \mathcal{V}$ being defined as
	\begin{equation*}
		\nabla^*_v F: ~ \mathcal{V} \longrightarrow \mathbb{R} \qquad v_i \longmapsto \nabla^*_v F \left(v_i\right) =
	\end{equation*}
	\begin{equation}
		\sum_{v_j \in \mathcal{V}} \left(\delta\left(v_j, v_i\right) w_G \left(v_i\right)^\epsilon F\left(v_j, v_i\right) - \delta\left(v_i, v_j\right) w_G \left(v_i\right)^\eta F\left(v_i, v_j\right)\right) W_I \left(v_i, v_j\right)^\beta W_G \left(v_i, v_j\right)^\gamma.
	\end{equation}\vspace{0em}
\end{definition}

\begin{remark}[\textbf{Parameter choice for the vertex adjoint operator}]\label{Gnabla^*_vparam} \ \\
	Using the before defined parameter choices $\beta = 1$, $\gamma = \frac{1}{2}$, $\epsilon = 0$ and $\eta = 0$ together with $W_I \left(a_q\right) \equiv \frac{1}{2}$ from the inner product in $\mathcal{H}\left(\mathcal{A}_G\right)$ for all arcs $a_q \in \mathcal{A}_G$, results in the common definition of the vertex adjoint (Page 2419 in \cite{elmoataz2015p}: Adjoint operator $\mathcal{G}^*_w$) for all vertices $v_i \in \mathcal{V}$:
	\begin{equation}
		\nabla^*_v F \left(v_i\right) = \frac{1}{2} \sum_{v_j \in \mathcal{V}} \left(\delta\left(v_j, v_i\right) F\left(v_j, v_i\right) - \delta\left(v_i, v_j\right) F\left(v_i, v_j\right)\right) \sqrt{W\left(v_i, v_j\right)}.
	\end{equation}
\end{remark}

The following theorem proves that the vertex adjoint is well-defined with regard to the vertex gradient.\\

\begin{theorem}[\textbf{Connection vertex gradient $\nabla_v$ and vertex adjoint $\nabla^*_v$}]\label{Gvertexgradadj} \ \\
	For a weighted oriented normal graph $OG = \left(\mathcal{V}, \mathcal{A}_G, w, W\right)$ with symmetric arc weight functions $W_I$ and $W_G$, the vertex adjoint operator $\nabla^*_v$ fulfills the equality
	\begin{equation}
		{\langle G, \nabla_v f \rangle}_{\mathcal{H}\left(\mathcal{A}_G\right)} = {\langle f, \nabla^*_v G \rangle}_{\mathcal{H}\left(\mathcal{V}\right)}
	\end{equation}
	for all vertex functions $f \in \mathcal{H}\left(\mathcal{V}\right)$ and all arc functions $G \in \mathcal{H}\left(\mathcal{A}_G\right)$.\\
\end{theorem}

\begin{proof}\ \\
	Given a weighted oriented normal graph $OG = \left(\mathcal{V}, \mathcal{A}_G, w, W\right)$ with symmetric arc weight functions $W_I$ and $W_G$, a vertex function $f \in \mathcal{H}\left(\mathcal{V}\right)$ and an arc function $G \in \mathcal{H}\left(\mathcal{A}_G\right)$, using the definitions of the inner product in $\mathcal{H}\left(\mathcal{A}_G\right)$ and the vertex gradient operator $\nabla_v$ leads to the following:
	
	$\begin{aligned}[t]
		{\langle G, \nabla_v f \rangle}_{\mathcal{H}\left(\mathcal{A}_G\right)} = & \sum_{a_q \in \mathcal{A}_G} W_I \left(a_q\right)^\beta G\left(a_q\right) \nabla_v f \left(a_q\right) &\\
		= & \sum_{a_q = \left(v_i, v_j\right) \in \mathcal{A}_G} W_I \left(a_q\right)^\beta G\left(a_q\right) W_G \left(a_q\right)^\gamma &\\
		& \left(w_I \left(v_j\right)^\alpha w_G \left(v_j\right)^\epsilon f\left(v_j\right) - w_I \left(v_i\right)^\alpha w_G \left(v_i\right)^\eta f\left(v_i\right)\right) &\\
		= & \sum_{a_q = \left(v_i, v_j\right) \in \mathcal{A}_G} W_I \left(a_q\right)^\beta W_G \left(a_q\right)^\gamma G\left(a_q\right) &\\
		& \left(w_I \left(v_j\right)^\alpha w_G \left(v_j\right)^\epsilon f\left(v_j\right) - w_I \left(v_i\right)^\alpha w_G \left(v_i\right)^\eta f\left(v_i\right)\right) &\\
	\end{aligned}$\\
	
	Summing over all arcs $a_q = \left(v_i, v_j\right) \in \mathcal{A}_G$ is equivalent to firstly summing over all vertices $v_i \in \mathcal{V}$, then over all vertices $v_j \in \mathcal{V}$ and including $\delta\left(v_i, v_j\right)$:
	
	$\begin{aligned}[t]
		= & \sum_{v_i \in \mathcal{V}} \sum_{v_j \in \mathcal{V}} \delta\left(v_i, v_j\right) W_I \left(v_i, v_j\right)^\beta W_G \left(v_i, v_j\right)^\gamma G\left(v_i, v_j\right) &\\
		& \quad \left(w_I \left(v_j\right)^\alpha w_G \left(v_j\right)^\epsilon f\left(v_j\right) - w_I \left(v_i\right)^\alpha w_G \left(v_i\right)^\eta f\left(v_i\right)\right) &\\
		= & ~ \frac{1}{2} \sum_{v_i \in \mathcal{V}} \sum_{v_j \in \mathcal{V}} \delta\left(v_i, v_j\right) W_I \left(v_i, v_j\right)^\beta W_G \left(v_i, v_j\right)^\gamma G\left(v_i, v_j\right) &\\
		& \quad \left(w_I \left(v_j\right)^\alpha w_G \left(v_j\right)^\epsilon f\left(v_j\right) - w_I \left(v_i\right)^\alpha w_G \left(v_i\right)^\eta f\left(v_i\right)\right) + &\\
		& \quad \delta\left(v_j, v_i\right) W_I \left(v_j, v_i\right)^\beta W_G \left(v_j, v_i\right)^\gamma G\left(v_j, v_i\right) &\\
		& \quad \left(w_I \left(v_i\right)^\alpha w_G \left(v_i\right)^\epsilon f\left(v_i\right) - w_I \left(v_j\right)^\alpha w_G \left(v_j\right)^\eta f\left(v_j\right)\right) &\\
	\end{aligned}$\\

	Using the symmetry property of the arc weight functions $W_I$ and $W_G$ results in the following:
	
	$\begin{aligned}[t]
		= & ~ \frac{1}{2} \sum_{v_i \in \mathcal{V}} \sum_{v_j \in \mathcal{V}} \left(\delta\left(v_i, v_j\right) w_I \left(v_j\right)^\alpha w_G \left(v_j\right)^\epsilon G\left(v_i, v_j\right) f\left(v_j\right) - \right. &\\
		& \quad ~ \delta\left(v_i, v_j\right) w_I \left(v_i\right)^\alpha w_G \left(v_i\right)^\eta G\left(v_i, v_j\right) f\left(v_i\right) + &\\
		& \quad ~ \delta\left(v_j, v_i\right) w_I \left(v_i\right)^\alpha w_G \left(v_i\right)^\epsilon G\left(v_j, v_i\right) f\left(v_i\right) - &\\
		& \left. \quad \delta\left(v_j, v_i\right) w_I \left(v_j\right)^\alpha w_G \left(v_j\right)^\eta G\left(v_j, v_i\right) f\left(v_j\right)\right) W_I \left(v_i, v_j\right)^\beta W_G \left(v_i, v_j\right)^\gamma &\\
		= & ~ \frac{1}{2} \sum_{v_i \in \mathcal{V}} \sum_{v_j \in \mathcal{V}} \left(\delta\left(v_j, v_i\right) w_G \left(v_i\right)^\epsilon G\left(v_j, v_i\right) - \delta\left(v_i, v_j\right) w_G \left(v_i\right)^\eta G\left(v_i, v_j\right)\right) &\\
		& \quad W_I \left(v_i, v_j\right)^\beta W_G \left(v_i, v_j\right)^\gamma w_I \left(v_i\right)^\alpha f\left(v_i\right) + &\\
		& \quad \left(\delta\left(v_i, v_j\right) w_G \left(v_j\right)^\epsilon G\left(v_i, v_j\right) - \delta\left(v_j, v_i\right) w_G \left(v_j\right)^\eta G\left(v_j, v_i\right)\right) &\\
		& \quad W_I \left(v_i, v_j\right)^\beta W_G \left(v_i, v_j\right)^\gamma w_I \left(v_j\right)^\alpha f\left(v_j\right) &\\
	\end{aligned}$\\
	$\begin{aligned}[t]
		= & ~ \frac{1}{2} \sum_{v_i \in \mathcal{V}} \sum_{v_j \in \mathcal{V}} \left(\delta\left(v_j, v_i\right) w_G \left(v_i\right)^\epsilon G\left(v_j, v_i\right) - \delta\left(v_i, v_j\right) w_G \left(v_i\right)^\eta G\left(v_i, v_j\right)\right) &\\
		& \quad W_I \left(v_i, v_j\right)^\beta W_G \left(v_i, v_j\right)^\gamma w_I \left(v_i\right)^\alpha f\left(v_i\right) + &\\
		& \quad \left(\delta\left(v_i, v_j\right) w_G \left(v_j\right)^\epsilon G\left(v_i, v_j\right) - \delta\left(v_j, v_i\right) w_G \left(v_j\right)^\eta G\left(v_j, v_i\right)\right) &\\
		& \quad W_I \left(v_j, v_i\right)^\beta W_G \left(v_j, v_i\right)^\gamma w_I \left(v_j\right)^\alpha f\left(v_j\right) &\\
	\end{aligned}$\\
	
	Since every vertex combination of the double sum appears twice, once in the order $\left(v_k, v_l\right)$ and once in the order $\left(v_l, v_k\right)$ for $v_k, v_l \in \mathcal{V}$, the second summand together with the $\frac{1}{2}$ factor can be eliminated. Moreover, using the definitions for the inner product on the space of all vertex functions $\mathcal{H}\left(\mathcal{V}\right)$ and the vertex adjoint operator $\nabla^*_v$ yields the following:
	
	$\begin{aligned}[t]
		= & ~ \frac{1}{2} \sum_{v_i \in \mathcal{V}} \sum_{v_j \in \mathcal{V}} \left(\delta\left(v_j, v_i\right) w_G \left(v_i\right)^\epsilon G\left(v_j, v_i\right) - \delta\left(v_i, v_j\right) w_G \left(v_i\right)^\eta G\left(v_i, v_j\right)\right) &\\
		& ~ W_I \left(v_i, v_j\right)^\beta W_G \left(v_i, v_j\right)^\gamma w_I \left(v_i\right)^\alpha f\left(v_i\right) &\\
		= & ~ \frac{1}{2} \sum_{v_i \in \mathcal{V}} w_I \left(v_i\right)^\alpha f\left(v_i\right) \sum_{v_j \in \mathcal{V}} \left(\delta\left(v_j, v_i\right) w_G \left(v_i\right)^\epsilon G\left(v_j, v_i\right) - \delta\left(v_i, v_j\right) w_G \left(v_i\right)^\eta G\left(v_i, v_j\right)\right) &\\
		& ~ W_I \left(v_i, v_j\right)^\beta W_G \left(v_i, v_j\right)^\gamma  &\\
		= & \sum_{v_i \in \mathcal{V}} w_I \left(v_i\right)^\alpha f\left(v_i\right) \nabla^*_v G \left(v_i\right) = {\langle f, \nabla^*_v G \rangle}_{\mathcal{H}\left(\mathcal{V}\right)}
	\end{aligned}$\\
	
	Therefore, with the presented definition for the vertex adjoint operator $\nabla^*_v$, the equality ${\langle G, \nabla_v f \rangle}_{\mathcal{H}\left(\mathcal{A}_G\right)} = {\langle f, \nabla^*_v G \rangle}_{\mathcal{H}\left(\mathcal{V}\right)}$ holds true for all vertex functions $f \in \mathcal{H}\left(\mathcal{V}\right)$ and all arc functions $G \in \mathcal{H}\left(\mathcal{A}_G\right)$.\\
\end{proof}

Similarly to the vertex adjoint operator $\nabla^*_v$ for arc functions $F \in \mathcal{H}\left(\mathcal{A}_G\right)$, an arc adjoint operator $\nabla^*_a$ for vertex functions $f \in \mathcal{H}\left(\mathcal{V}\right)$ can be defined on the basis of the arc gradient operator $\nabla_a$.\\

\begin{definition}[\textbf{Arc adjoint operator $\nabla^*_a$}]\label{Gnabla^*_a} \ \\
	For a weighted oriented normal graph $OG = \left(\mathcal{V}, \mathcal{A}_G, w, W\right)$, the arc adjoint operator $\nabla^*_a$ is defined as
	\begin{equation}
		\nabla^*_a: ~ \mathcal{H}\left(\mathcal{V}\right) \longrightarrow \mathcal{H}\left(\mathcal{A}_G\right) \qquad f \longmapsto \nabla^*_a f
	\end{equation}
	with the weighted adjoint of vertex function $f \in \mathcal{H}\left(\mathcal{V}\right)$ at an arc $a_q = \left(v_i, v_j\right) \in \mathcal{A}_G$ being given by
	\begin{equation*}
		\nabla^*_a f: ~ \mathcal{A}_G \longrightarrow \mathbb{R} \qquad a_q = \left(v_i, v_j\right) \longmapsto \nabla^*_a f \left(a_q\right) = 
	\end{equation*}
	\begin{equation}
		W_G \left(v_i, v_j\right)^\theta \left(\frac{w_I \left(v_j\right)^\alpha w_G \left(v_j\right)^\zeta}{\deg_{in}\left(v_j\right)} f\left(v_j\right) - \frac{w_I \left(v_i\right)^\alpha w_G \left(v_i\right)^\zeta}{\deg_{out}\left(v_i\right)} f\left(v_i\right)\right).
	\end{equation}\vspace{0em}
\end{definition}

\begin{remark}[\textbf{Parameter choice for the arc adjoint operator}]\label{Gnabla^*_aparam} \ \\
	Using the previous parameter choices $\alpha = 0$, $\zeta = \frac{1}{2}$, $\theta = 0$ together with leaving out the degree functions $\deg_{in}\left(v_j\right)$ and $\deg_{out}\left(v_i\right)$ for all vertices $v_i, v_j \in \mathcal{V}$, results in a simplified definition of the arc adjoint, which matches the simplified version of the arc gradient $\nabla_a$, for all arcs $a_q = \left(v_i, v_j\right) \in \mathcal{A}_G$:
	\begin{equation}
		\nabla^*_a f \left(a_q\right) = \nabla^*_a f \left(v_i, v_j\right) = \sqrt{w_G \left(v_j\right)} f\left(v_j\right) - \sqrt{w_G \left(v_i\right)} f\left(v_i\right).
	\end{equation}\vspace{0em}
\end{remark}

As in the case of the vertex adjoint and vertex gradient, the following theorem proves that the arc adjoint is the well-defined dual counterpart to the arc gradient with respect to the inner products.\\

\begin{theorem}[\textbf{Connection arc gradient $\nabla_a$ and arc adjoint $\nabla^*_a$}]\label{Garcgradadj} \ \\
	For a weighted oriented normal graph $OG = \left(\mathcal{V}, \mathcal{A}_G, w, W\right)$, the arc adjoint operator $\nabla^*_a$ fulfills the equality
	\begin{equation}
		{\langle f, \nabla_a G \rangle}_{\mathcal{H}\left(\mathcal{V}\right)} = {\langle G, \nabla^*_a f \rangle}_{\mathcal{H}\left(\mathcal{A}_G\right)}
	\end{equation}
	for all vertex functions $f \in \mathcal{H}\left(\mathcal{V}\right)$ and all arc functions $G \in \mathcal{H}\left(\mathcal{A}_G\right)$.\\
\end{theorem}

\begin{proof}\ \\
	Given a weighted oriented normal graph $OG = \left(\mathcal{V}, \mathcal{A}_G, w, W\right)$, a vertex function $f \in \mathcal{H}\left(\mathcal{V}\right)$ and an arc function $G \in \mathcal{H}\left(\mathcal{A}_G\right)$, then using the definitions of the inner product in $\mathcal{H}\left(\mathcal{V}\right)$ and the arc gradient operator $\nabla_a$ results in the following:
	
	$\begin{aligned}[t]
		{\langle f, \nabla_a G \rangle}_{\mathcal{H}\left(\mathcal{V}\right)} = & \sum_{v_i \in \mathcal{V}} w_I \left(v_i\right)^\alpha f\left(v_i\right) \nabla_a G\left(v_i\right) &\\
		= & \sum_{v_i \in \mathcal{V}} w_I \left(v_i\right)^\alpha f\left(v_i\right) &\\
		& w_G \left(v_i\right)^\zeta \sum_{a_q \in \mathcal{A}_G} \left(\frac{\delta_{in}\left(v_i, a_q\right)}{\deg_{in}\left(v_i\right)} - \frac{\delta_{out}\left(v_i, a_q\right)}{\deg_{out}\left(v_i\right)}\right) W_I \left(a_q\right)^\beta W_G \left(a_q\right)^\theta G\left(a_q\right) &\\
		= & \sum_{v_i \in \mathcal{V}} \sum_{a_q \in \mathcal{A}_G} \left(\frac{\delta_{in}\left(v_i, a_q\right)}{\deg_{in}\left(v_i\right)} - \frac{\delta_{out}\left(v_i, a_q\right)}{\deg_{out}\left(v_i\right)}\right) &\\
		& W_I \left(a_q\right)^\beta W_G \left(a_q\right)^\theta w_I \left(v_i\right)^\alpha w_G \left(v_i\right)^\zeta f\left(v_i\right) G\left(a_q\right) &\\
	\end{aligned}$\\
	
	\clearpage
	Exchanging the two sums and writing every arc $a_q \in \mathcal{A}_G$ as $a_q = \left(v_j, v_k\right)$ yields:
	
	$\begin{aligned}[t]
		= & \sum_{a_q \in \mathcal{A}_G} \sum_{v_i \in \mathcal{V}} \left(\frac{\delta_{in}\left(v_i, a_q\right)}{\deg_{in}\left(v_i\right)} - \frac{\delta_{out}\left(v_i, a_q\right)}{\deg_{out}\left(v_i\right)}\right) W_I \left(a_q\right)^\beta W_G \left(a_q\right)^\theta w_I \left(v_i\right)^\alpha w_G \left(v_i\right)^\zeta f\left(v_i\right) G\left(a_q\right) &\\
		= & \sum_{a_q \in \mathcal{A}_G} W_I \left(a_q\right)^\beta G\left(a_q\right) W_G \left(a_q\right)^\theta \sum_{v_i \in \mathcal{V}} \left(\frac{\delta_{in}\left(v_i, a_q\right)}{\deg_{in}\left(v_i\right)} - \frac{\delta_{out}\left(v_i, a_q\right)}{\deg_{out}\left(v_i\right)}\right) w_I \left(v_i\right)^\alpha w_G \left(v_i\right)^\zeta f\left(v_i\right) &\\
		= & \sum_{a_q = \left(v_j, v_k\right) \in \mathcal{A}_G} W_I \left(v_j, v_k\right)^\beta G\left(v_j, v_k\right) W_G \left(v_j, v_k\right)^\theta &\\
		& \left(\frac{w_I \left(v_k\right)^\alpha w_G \left(v_k\right)^\zeta}{\deg_{in}\left(v_k\right)} f\left(v_k\right) - \frac{w_I \left(v_j\right)^\alpha w_G \left(v_j\right)^\zeta}{\deg_{out}\left(v_j\right)} f\left(v_j\right)\right) &\\
	\end{aligned}$\\
	
	Using the definitions for the inner product on the space of all arc functions $\mathcal{H}\left(\mathcal{A}_G\right)$ and the arc adjoint operator $\nabla^*_a$ leads to the following:
	
	$\begin{aligned}[t]	
		= & \sum_{a_q = \left(v_j, v_k\right) \in \mathcal{A}_G} W_I \left(v_j, v_k\right)^\beta G\left(v_j, v_k\right) \nabla^*_a f\left(v_j, v_k\right) &\\
		= & ~ {\langle G, \nabla^*_a f \rangle}_{\mathcal{H}\left(\mathcal{A}_G\right)} &\\
	\end{aligned}$\\

	Hence, with the introduced definition of the arc adjoint operator $\nabla^*_a$, the equality ${\langle f, \nabla_a G \rangle}_{\mathcal{H}\left(\mathcal{V}\right)} = {\langle G, \nabla^*_a f \rangle}_{\mathcal{H}\left(\mathcal{A}_G\right)}$ holds true for all vertex functions $f \in \mathcal{H}\left(\mathcal{V}\right)$ and all arc functions $G \in \mathcal{H}\left(\mathcal{A}_G\right)$.\\
\end{proof}
\clearpage
\section{Divergence and Laplacian operators on oriented normal graphs}\label{7}

The vertex adjoint and the arc adjoint can be used to calculate the vertex divergence and the arc divergence based on their relationship in the continuum case. The vertex and arc divergence together with the vertex gradient and the arc gradient are then necessary for the definition of the vertex Laplacian, the arc Laplacian, the vertex $p$-Laplacian and the arc $p$-Laplacian.\\

Thus, subsection (\ref{7.1}) gives a definition for the vertex and arc divergence. Moreover, subsection (\ref{7.2}) derives the vertex Laplacian and the arc Laplacian using the previous definitions and in subsection (\ref{7.3}) the vertex $p$-Laplacian and the arc $p$-Laplacian are introduced, also using the before defined vertex gradient, vertex divergence, arc gradient and arc divergence.\\

\subsection{Divergence operators on oriented normal graphs}\label{7.1}

Based on the underlying equivalence of the continuum setting
\begin{equation}
	\text{div}_v = - \nabla^*_v
\end{equation}
\begin{equation}
	\text{div}_a = - \nabla^*_a
\end{equation}
the vertex divergence operator $\text{div}_v$ for arc functions $F \in \mathcal{H}\left(\mathcal{A}_G\right)$ and the arc divergence operator $\text{div}_a$ for vertex functions $f \in \mathcal{H}\left(\mathcal{V}\right)$ can be defined.\\

\begin{definition}[\textbf{Vertex divergence operator $\text{div}_v$}]\label{Gdiv_v} \ \\
	For a weighted oriented normal graph $OG = \left(\mathcal{V}, \mathcal{A}_G, w, W\right)$ with symmetric arc weight functions $W_I$ and $W_G$, the vertex divergence operator $\text{div}_v$ is given by
	\begin{equation}
		\text{div}_v: ~ \mathcal{H}\left(\mathcal{A}_G\right) \longrightarrow \mathcal{H}\left(\mathcal{V}\right) \qquad F \longmapsto \text{div}_v F
	\end{equation}
	with the weighted divergence of arc function $F \in \mathcal{H}\left(\mathcal{A}_G\right)$ at a vertex $v_i \in \mathcal{V}$ being defined as
	\begin{equation*}
		\text{div}_v: ~ \mathcal{V} \longrightarrow \mathbb{R} \qquad v_i \longmapsto \text{div}_v \left(F\right) \left(v_i\right) = - \nabla^*_v F \left(v_i\right) =
	\end{equation*}
	\begin{equation}
		\sum_{v_j \in \mathcal{V}} \left(\delta\left(v_i, v_j\right) w_G \left(v_i\right)^\eta F\left(v_i, v_j\right) - \delta\left(v_j, v_i\right) w_G \left(v_i\right)^\epsilon F\left(v_j, v_i\right)\right) W_I \left(v_i, v_j\right)^\beta W_G \left(v_i, v_j\right)^\gamma.
	\end{equation}\vspace{0em}
\end{definition}

\clearpage
\begin{remark}[\textbf{Parameter choice for the vertex divergence operator}]\label{Gdiv_vparam} \ \\
	Using the parameters of remark (\ref{Gnabla^*_vparam}) $\beta = 1$, $\gamma = \frac{1}{2}$, $\epsilon = 0$ and $\eta = 0$ together with $W_I \left(a_q\right) \equiv \frac{1}{2}$ from the inner product in $\mathcal{H}\left(\mathcal{A}_G\right)$ for all arcs $a_q \in \mathcal{A}_G$, results in the common definition of the vertex divergence (Page 2419 in \cite{elmoataz2015p}: Divergence operator $D_w$) for all vertices $v_i \in \mathcal{V}$:
	\begin{equation}
		\text{div}_v \left(F\right) \left(v_i\right) = \frac{1}{2} \sum_{v_j \in \mathcal{V}} \left(\delta\left(v_i, v_j\right) F\left(v_i, v_j\right) - \delta\left(v_j, v_i\right) F\left(v_j, v_i\right)\right) \sqrt{W_G \left(v_i, v_j\right)}.
	\end{equation}\vspace{0em}
\end{remark}

The arc divergence is defined in the same manner as the vertex divergence, but based on the arc adjoint instead of the vertex adjoint.\\

\begin{definition}[\textbf{Arc divergence operator $\text{div}_a$}]\label{Gdiv_a} \ \\
	For a weighted oriented normal graph $OG = \left(\mathcal{V}, \mathcal{A}_G, w, W\right)$, the arc divergence operator $\text{div}_a$ is defined as
	\begin{equation}
		\text{div}_a: ~ \mathcal{H}\left(\mathcal{V}\right) \longrightarrow \mathcal{H}\left(\mathcal{A}_G\right) \qquad f \longmapsto \text{div}_a f
	\end{equation}
	with the weighted divergence of vertex function $f \in \mathcal{H}\left(\mathcal{V}\right)$ at an arc $a_q = \left(v_i, v_j\right) \in \mathcal{A}_G$ being given by
	\begin{equation*}
		\text{div}_a: ~ \mathcal{A}_G \longrightarrow \mathbb{R} \qquad a_q = \left(v_i, v_j\right) \longmapsto \text{div}_a \left(f\right) \left(a_q\right) = - \nabla^*_a f \left(a_q\right) =
	\end{equation*}
	\begin{equation}
		W_G \left(v_i, v_j\right)^\theta \left(\frac{w_I \left(v_i\right)^\alpha w_G \left(v_i\right)^\zeta}{\deg_{out}\left(v_i\right)} f\left(v_i\right) - \frac{w_I \left(v_j\right)^\alpha w_G \left(v_j\right)^\zeta}{\deg_{in}\left(v_j\right)} f\left(v_j\right)\right).
	\end{equation}\vspace{0em}
\end{definition}

\begin{remark}[\textbf{Parameter choice for the arc divergence operator}]\label{Gdiv_aparam} \ \\
	Using again the parameters of remark (\ref{Gnabla^*_aparam}) $\alpha = 0$, $\zeta = \frac{1}{2}$, $\theta = 0$ together with leaving out the degree functions $\deg_{in}\left(v_j\right)$ and $\deg_{out}\left(v_i\right)$ for all vertices $v_i, v_j \in \mathcal{V}$, results in the simplified definition of the arc divergence, which also matches the simplified version of the arc gradient $\nabla_a$, for all arcs $a_q = \left(v_i, v_j\right) \in \mathcal{A}_G$:
	\begin{equation}
		\text{div}_a \left(f\right) \left(a_q\right) = \text{div}_a \left(f\right) \left(v_i, v_j\right) = \sqrt{w_G \left(v_i\right)} f\left(v_i\right) - \sqrt{w_G \left(v_j\right)} f\left(v_j\right).
	\end{equation}\vspace{0em}
\end{remark}
\clearpage
\subsection{Laplacian operators on oriented normal graphs}\label{7.2}

As in the case of the divergence operators, the Laplacian operators can be both defined for vertices $v_i \in \mathcal{V}$ and for arcs $a_q \in \mathcal{A}_G$ on the basis of the relationship in the continuum setting\\
\begin{equation}
	\Delta_v f = \text{div}_v \left(\nabla_v f\right) 
\end{equation}
\begin{equation}
	\Delta_a F = \text{div}_a \left(\nabla_a F\right) 
\end{equation}

for any vertex function $f \in \mathcal{H}\left(\mathcal{V}\right)$ and any arc function $F \in \mathcal{H}\left(\mathcal{A}_G\right)$.\\

\begin{definition}[\textbf{Vertex Laplacian operator $\Delta_v$}]\label{GDelta_v} \ \\
	For a weighted oriented normal graph $OG = \left(\mathcal{V}, \mathcal{A}_G, w, W\right)$ with symmetric arc weight functions $W_I$ and $W_G$ the vertex Laplacian operator $\Delta_v$ is given by
	\begin{equation}
		\Delta_v: ~ \mathcal{H}\left(\mathcal{V}\right) \longrightarrow \mathcal{H}\left(\mathcal{V}\right) \qquad f \longmapsto \Delta_v f
	\end{equation}
	with the weighted Laplacian of vertex function $f \in \mathcal{H}\left(\mathcal{V}\right)$ at a vertex $v_i \in \mathcal{V}$ being defined as
	\begin{equation*}
		\Delta_v: ~ \mathcal{V} \longrightarrow \mathbb{R} \qquad v_i \longmapsto \Delta_v f \left(v_i\right) =
	\end{equation*}
	\begin{equation*}
		\sum_{v_j \in \mathcal{V}} \left(\left(\delta\left(v_i, v_j\right) w_G \left(v_i\right)^\eta w_I \left(v_j\right)^\alpha w_G \left(v_j\right)^\epsilon + \delta\left(v_j, v_i\right) w_I \left(v_j\right)^\alpha w_G \left(v_j\right)^\eta w_G \left(v_i\right)^\epsilon\right) f\left(v_j\right) \right. -
	\end{equation*}
	\begin{equation*}
		\left.\left(\delta\left(v_i, v_j\right) w_I \left(v_i\right)^\alpha w_G \left(v_i\right)^{2 \eta} + \delta\left(v_j, v_i\right) w_I \left(v_i\right)^\alpha w_G \left(v_i\right)^{2 \epsilon}\right) f\left(v_i\right)\right)
	\end{equation*}
	\begin{equation}
		W_I \left(v_i, v_j\right)^\beta W_G \left(v_i, v_j\right)^{2 \gamma}.
	\end{equation}\vspace{0em}
\end{definition}

\begin{remark}[\textbf{Parameter choice for the vertex Laplacian operator}]\label{GDelta_vparam} \ \\
	Using the previous parameters $\alpha = 0$, $\beta = 1$, $\gamma = \frac{1}{2}$, $\epsilon = 0$ and $\eta = 0$ together with $W_I \left(a_q\right)^\beta \equiv \frac{1}{2}$ from the inner product in $\mathcal{H}\left(\mathcal{A}_G\right)$ for all arcs $a_q \in \mathcal{A}_G$, leads to the common definition of the vertex Laplacian (Page 2421 in \cite{elmoataz2015p}: Classical unnormalized Laplacian) for all vertices $v_i \in \mathcal{V}$:
	\begin{equation}
		\Delta_v f \left(v_i\right) = \frac{1}{2}\sum_{v_j \in \mathcal{V}} \left(\delta\left(v_i, v_j\right) + \delta\left(v_j, v_i\right)\right) W_G \left(v_i, v_j\right) \left(f\left(v_j\right) - f\left(v_i\right)\right).
	\end{equation}\vspace{0em}
\end{remark}

The following theorem proves that the vertex Laplacian is well-defined based on the vertex gradient and vertex divergence definitions.\\

\begin{theorem}[\textbf{Connection vertex divergence $\text{div}_v$, vertex gradient $\nabla_v$, and vertex Laplacian $\Delta_v$}]\label{GvertexLap} \ \\
	On a weighted oriented normal graph $OG = \left(\mathcal{V}, \mathcal{A}_G, w, W\right)$ with symmetric arc weight functions $W_I$ and $W_G$, the vertex Laplacian operator $\Delta_v$ fulfills the equality
	\begin{equation}
		\Delta_v f = \text{div}_v \left(\nabla_v f\right)
	\end{equation}
	for all vertex functions $f \in \mathcal{H}\left(\mathcal{V}\right)$.\\
\end{theorem}

\begin{proof}\ \\
	Given a weighted oriented normal graph $OG = \left(\mathcal{V}, \mathcal{A}_G, w, W\right)$ with symmetric arc weight functions $W_I$ and $W_G$ and a vertex function $f \in \mathcal{H}\left(\mathcal{V}\right)$, then using the definitions of the vertex divergence operator $\text{div}_v$ and the vertex gradient operator $\nabla_v$ yields the following for all vertices $v_i \in \mathcal{V}$:
	
	$\begin{aligned}[t]	
		\text{div}_v \left(\nabla_v f\right)\left(v_i\right) = & \sum_{v_j \in \mathcal{V}} \left(\delta\left(v_i, v_j\right) w_G \left(v_i\right)^\eta \nabla_v f\left(v_i, v_j\right) - \delta\left(v_j, v_i\right) w_G \left(v_i\right)^\epsilon \nabla_v f\left(v_j, v_i\right)\right) &\\
		& \quad W_I \left(v_i, v_j\right)^\beta W_G \left(v_i, v_j\right)^\gamma &\\
		= & \sum_{v_j \in \mathcal{V}} \left(\delta\left(v_i, v_j\right) w_G \left(v_i\right)^\eta W_G \left(v_i, v_j\right)^\gamma \right. &\\
		& \quad \left(w_I \left(v_j\right)^\alpha w_G \left(v_j\right)^\epsilon f\left(v_j\right) - w_I \left(v_i\right)^\alpha w_G \left(v_i\right)^\eta f\left(v_i\right)\right) - &\\
		& \quad \delta\left(v_j, v_i\right) w_G \left(v_i\right)^\epsilon W_G \left(v_j, v_i\right)^\gamma &\\
		& \quad \left.\left(w_I \left(v_i\right)^\alpha w_G \left(v_i\right)^\epsilon f\left(v_i\right) - w_I \left(v_j\right)^\alpha w_G \left(v_j\right)^\eta f\left(v_j\right)\right)\right) &\\
		& \quad W_I \left(v_i, v_j\right)^\beta W_G \left(v_i, v_j\right)^\gamma &\\
	\end{aligned}$\\
	
	Using the symmetry property of the arc weight function $W$ and rearranging the summands leads to the vertex Laplacian definition:
	
	$\begin{aligned}[t]
		= & \sum_{v_j \in \mathcal{V}} \left(\delta\left(v_i, v_j\right) w_G \left(v_i\right)^\eta \left(w_I \left(v_j\right)^\alpha w_G \left(v_j\right)^\epsilon f\left(v_j\right) - w_I \left(v_i\right)^\alpha w_G \left(v_i\right)^\eta f\left(v_i\right)\right) - \right. &\\
		& \quad \left.\delta\left(v_j, v_i\right) w_G \left(v_i\right)^\epsilon \left(w_I \left(v_i\right)^\alpha w_G \left(v_i\right)^\epsilon f\left(v_i\right) - w_I \left(v_j\right)^\alpha w_G \left(v_j\right)^\eta f\left(v_j\right)\right)\right) &\\
		& \quad W_I \left(v_i, v_j\right)^\beta W_G \left(v_i, v_j\right)^{2 \gamma} &\\
		= & \sum_{v_j \in \mathcal{V}} \left(\delta\left(v_i, v_j\right) w_G \left(v_i\right)^\eta w_I \left(v_j\right)^\alpha w_G \left(v_j\right)^\epsilon f\left(v_j\right) - \right. &\\
		& \quad \delta\left(v_i, v_j\right) w_G \left(v_i\right)^\eta w_I \left(v_i\right)^\alpha w_G \left(v_i\right)^\eta f\left(v_i\right) - \delta\left(v_j, v_i\right) w_G \left(v_i\right)^\epsilon w_I \left(v_i\right)^\alpha w_G \left(v_i\right)^\epsilon f\left(v_i\right) + &\\
		& \quad \left.\delta\left(v_j, v_i\right) w_G \left(v_i\right)^\epsilon w_I \left(v_j\right)^\alpha w_G \left(v_j\right)^\eta f\left(v_j\right)\right) W_I \left(v_i, v_j\right)^\beta W_G \left(v_i, v_j\right)^{2 \gamma} &\\
	\end{aligned}$\\
	$\begin{aligned}[t]
		= & \sum_{v_j \in \mathcal{V}} \left(\left(\delta\left(v_i, v_j\right) w_G \left(v_i\right)^\eta w_I \left(v_j\right)^\alpha w_G \left(v_j\right)^\epsilon + \delta\left(v_j, v_i\right) w_G \left(v_i\right)^\epsilon w_I \left(v_j\right)^\alpha w_G \left(v_j\right)^\eta\right) f\left(v_j\right) - \right. &\\
		& \quad \left.\left(\delta\left(v_i, v_j\right) w_G \left(v_i\right)^\eta w_I \left(v_i\right)^\alpha w_G \left(v_i\right)^\eta + \delta\left(v_j, v_i\right) w_G \left(v_i\right)^\epsilon w_I \left(v_i\right)^\alpha w_G \left(v_i\right)^\epsilon\right) f\left(v_i\right)\right) &\\
		& \quad W_I \left(v_i, v_j\right)^\beta W_G \left(v_i, v_j\right)^{2 \gamma} &\\
		= & \sum_{v_j \in \mathcal{V}} \left(\left(\delta\left(v_i, v_j\right) w_G \left(v_i\right)^\eta w_I \left(v_j\right)^\alpha w_G \left(v_j\right)^\epsilon + \delta\left(v_j, v_i\right) w_G \left(v_i\right)^\epsilon w_I \left(v_j\right)^\alpha w_G \left(v_j\right)^\eta\right) f\left(v_j\right) - \right. &\\
		& \quad \left.\left(\delta\left(v_i, v_j\right) w_I \left(v_i\right)^\alpha w_G \left(v_i\right)^{2 \eta} + \delta\left(v_j, v_i\right) w_I \left(v_i\right)^\alpha w_G \left(v_i\right)^{2 \epsilon}\right) f\left(v_i\right)\right) &\\
		& \quad W_I \left(v_i, v_j\right)^\beta W_G \left(v_i, v_j\right)^{2 \gamma} &\\
		= & ~ \Delta_v f \left(v_i\right)
	\end{aligned}$\\
	
	Thus, with the presented definition for the vertex gradient $\nabla_v$, the vertex divergence $\text{div}_v$, and the vertex Laplacian $\Delta_v$, the equality $\Delta_v f\left(v_i\right) = \text{div}_v \left(\nabla_v f\right)\left(v_i\right)$ holds true for all vertices $v_i \in \mathcal{V}$ and for all vertex functions $f \in \mathcal{H}\left(\mathcal{V}\right)$.\\
\end{proof}

The Laplacian operator can be defined for arcs $a_q \in \mathcal{A}_G$ as well, using the arc gradient and the arc divergence.\\

\begin{definition}[\textbf{Arc Laplacian operator $\Delta_a$}]\label{GDelta_a} \ \\
	For a weighted oriented normal graph $OG = \left(\mathcal{V}, \mathcal{A}_G, w, W\right)$ the arc Laplacian operator $\Delta_a$ is defined as
	\begin{equation}
		\Delta_a: ~ \mathcal{H}\left(\mathcal{A}_G\right) \longrightarrow \mathcal{H}\left(\mathcal{A}_G\right) \qquad F \longmapsto \Delta_a F
	\end{equation}
	with the weighted Laplacian of arc function $F \in \mathcal{H}\left(\mathcal{A}_G\right)$ at an arc $a_q = \left(v_i, v_j\right) \in \mathcal{A}_G$ being given by
	\begin{equation*}
		\Delta_a: ~ \mathcal{A}_G \longrightarrow \mathbb{R} \qquad a_q = \left(v_i, v_j\right) \longmapsto \Delta_a F \left(a_q\right) =
	\end{equation*}
	\begin{equation*}
		W_G \left(v_i, v_j\right)^\theta \left(\frac{w_I \left(v_i\right)^\alpha w_G \left(v_i\right)^{2 \zeta}}{\deg_{out}\left(v_i\right)} \sum_{a_r \in \mathcal{A}_G} \left(\frac{\delta_{in}\left(v_i, a_r\right)}{\deg_{in}\left(v_i\right)} - \frac{\delta_{out}\left(v_i, a_r\right)}{\deg_{out}\left(v_i\right)}\right) W_I \left(a_r\right)^\beta W_G \left(a_r\right)^\theta F\left(a_r\right)\right.
	\end{equation*}
	\begin{equation}
		\left. - \frac{w_I \left(v_j\right)^\alpha w_G \left(v_j\right)^{2 \zeta}}{\deg_{in}\left(v_j\right)} \sum_{a_s \in \mathcal{A}_G} \left(\frac{\delta_{in}\left(v_j, a_s\right)}{\deg_{in}\left(v_j\right)} - \frac{\delta_{out}\left(v_j, a_s\right)}{\deg_{out}\left(v_j\right)}\right) W_I \left(a_s\right)^\beta W_G \left(a_s\right)^\theta F\left(a_s\right)\right).
	\end{equation}\vspace{0em}
\end{definition}

\clearpage
\begin{remark}[\textbf{Parameter choice for the arc Laplacian operator}]\label{GDelta_aparam} \ \\
	The previous parameter choices $\alpha = 0$, $\beta = 1$, $\zeta = \frac{1}{2}$, $\theta = 0$ together with leaving out $\deg_{in}\left(v_i\right)$ and $\deg_{out}\left(v_i\right)$ for all vertices $v_i \in \mathcal{V}$ and with $W_I \left(a_q\right) \equiv \frac{1}{2}$ from the inner product in $\mathcal{H}\left(\mathcal{A}_G\right)$ for all arcs $a_q \in \mathcal{A}_G$, results in the simplified definition of the arc Laplacian for all arcs $a_q = \left(v_i, v_j\right) \in \mathcal{A}_G$:
	\begin{equation*}
		\Delta_a F \left(a_q\right) = \Delta_a F \left(v_i, v_j\right) =
	\end{equation*}
	\begin{equation*}
		\frac{w_G \left(v_i\right)}{2} \sum_{a_r \in \mathcal{A}_G} \left(\delta_{in}\left(v_i, a_r\right) - \delta_{out}\left(v_i, a_r\right)\right) F\left(a_r\right) -
	\end{equation*}
	\begin{equation}
		\frac{w_G \left(v_j\right)}{2} \sum_{a_s \in \mathcal{A}_G} \left(\delta_{in}\left(v_j, a_s\right) - \delta_{out}\left(v_j, a_s\right)\right) F\left(a_s\right).
	\end{equation}\vspace{0em}
\end{remark}

The theorem below proves that the definition of the arc Laplacian is well-defined based on the previously introduced arc gradient and arc divergence.\\

\begin{theorem}[\textbf{Connection arc divergence $\text{div}_a$, arc gradient $\nabla_a$, and arc Laplacian $\Delta_a$}]\label{GarcLap} \ \\
	On a weighted oriented normal graph $OG = \left(\mathcal{V}, \mathcal{A}_G, w, W\right)$, the arc Laplacian operator $\Delta_a$ fulfills the equality
	\begin{equation}
		\Delta_a F = \text{div}_a \left(\nabla_a F\right)
	\end{equation}
	for all arc functions $F \in \mathcal{H}\left(\mathcal{A}_G\right)$.\\
\end{theorem}

\begin{proof}\ \\
	Given a weighted oriented normal graph $OG = \left(\mathcal{V}, \mathcal{A}_G, w, W\right)$ and an arc function $F \in \mathcal{H}\left(\mathcal{A}_G\right)$, then the definitions of the arc divergence operator $\text{div}_a$ and the arc gradient operator $\nabla_a$ results in the following for all arcs $a_q = \left(v_i, v_j\right) \in \mathcal{A}_G$:
	
	$\begin{aligned}[t]	
		\text{div}_a \left(\nabla_a F\right)\left(v_i, v_j\right) = & ~ W_G \left(v_i, v_j\right)^\theta &\\
		& \left(\frac{w_I \left(v_i\right)^\alpha w_G\left(v_i\right)^\zeta}{\deg_{out}\left(v_i\right)} \nabla_a F\left(v_i\right) - \frac{w_I \left(v_j\right)^\alpha w_G\left(v_j\right)^\zeta}{\deg_{in}\left(v_j\right)} \nabla_a F\left(v_j\right)\right) &\\
	\end{aligned}$\\
	$\begin{aligned}[t]	
		= & ~ W_G \left(v_i, v_j\right)^\theta \left(\frac{w_I \left(v_i\right)^\alpha w_G\left(v_i\right)^\zeta}{\deg_{out}\left(v_i\right)} w_G \left(v_i\right)^\zeta \right. &\\
		& \quad \sum_{a_r \in \mathcal{A}_G} \left(\frac{\delta_{in}\left(v_i, a_r\right)}{\deg_{in}\left(v_i\right)} - \frac{\delta_{out}\left(v_i, a_r\right)}{\deg_{out}\left(v_i\right)}\right) W_I \left(a_r\right)^\beta W_G \left(a_r\right)^\theta F\left(a_r\right) - &\\
		& \quad \frac{w_I \left(v_j\right)^\alpha w_G\left(v_j\right)^\zeta}{\deg_{in}\left(v_j\right)} w\left(v_j\right)^\zeta \left.\sum_{a_s \in \mathcal{A}_G} \left(\frac{\delta_{in}\left(v_j, a_s\right)}{\deg_{in}\left(v_j\right)} - \frac{\delta_{out}\left(v_j, a_s\right)}{\deg_{out}\left(v_j\right)}\right) W_I \left(a_s\right)^\beta W_G \left(a_s\right)^\theta F\left(a_s\right)\right) &\\
	\end{aligned}$\\
	$\begin{aligned}[t]	
		= & ~ W_G \left(v_i, v_j\right)^\theta \left(\frac{w_I \left(v_i\right)^\alpha w_G \left(v_i\right)^{2 \zeta}}{\deg_{out}\left(v_i\right)} \right. &\\
		& \quad \sum_{a_r \in \mathcal{A}_G} \left(\frac{\delta_{in}\left(v_i, a_r\right)}{\deg_{in}\left(v_i\right)} - \frac{\delta_{out}\left(v_i, a_r\right)}{\deg_{out}\left(v_i\right)}\right) W_I \left(a_r\right)^\beta W_G \left(a_r\right)^\theta F\left(a_r\right) - &\\
		& \quad \frac{w_I \left(v_j\right)^\alpha w_G \left(v_j\right)^{2 \zeta}}{\deg_{in}\left(v_j\right)} \left. \sum_{a_s \in \mathcal{A}_G} \left(\frac{\delta_{in}\left(v_j, a_s\right)}{\deg_{in}\left(v_j\right)} - \frac{\delta_{out}\left(v_j, a_s\right)}{\deg_{out}\left(v_j\right)}\right) W_I \left(a_s\right)^\beta W_G \left(a_s\right)^\theta F\left(a_s\right)\right) &\\
		= & ~ \Delta_a F \left(v_i, v_j\right)
	\end{aligned}$\\
	
	Hence, with the previously introduced definitions for the arc gradient $\nabla_a$, the arc divergence $\text{div}_a$, and the arc Laplacian $\Delta_a$, the equality $\Delta_a F\left(a_q\right) = \text{div}_a \left(\nabla_a F\right)\left(a_q\right)$ holds true for all arcs $a_q \in \mathcal{A}_G$ and for all arc functions $F \in \mathcal{H}\left(\mathcal{A}_G\right)$.\\
\end{proof}
\subsection{$p$-Laplacian operators on oriented normal graphs}\label{7.3}

The weighted $p$-Laplacian operators on oriented normal graphs are defined based on the continuum setting for all $p \in \left(1, \infty\right)$ with the aim of fulfilling the following two equalities\\
\begin{equation}
	\Delta_v^p f = \text{div}_v \left(\left\lvert \nabla_v f\right\rvert^{p - 2} \nabla_v f\right)
\end{equation}
\begin{equation}
	\Delta_a^p F = \text{div}_a \left(\left\lvert \nabla_a F\right\rvert^{p - 2} \nabla_a F\right)
\end{equation}

for any vertex function $f \in \mathcal{H}\left(\mathcal{V}\right)$ and any arc function $F \in \mathcal{H}\left(\mathcal{A}_G\right)$.\\

\begin{definition}[\textbf{Vertex $p$-Laplacian operator $\Delta_v^p$}]\label{GDelta_v^p} \ \\
	For a weighted oriented normal graph $OG = \left(\mathcal{V}, \mathcal{A}_G, w, W\right)$ with symmetric arc weight functions $W_I$ and $W_G$ the vertex $p$-Laplacian operator $\Delta_v^p$ is given by:
	\begin{equation}
		\Delta_v^p: ~ \mathcal{H}\left(\mathcal{V}\right) \longrightarrow \mathcal{H}\left(\mathcal{V}\right) \qquad f \longmapsto \Delta_v^p f
	\end{equation}
	with the weighted $p$-Laplacian of vertex function $f \in \mathcal{H}\left(\mathcal{V}\right)$ at a vertex $v_i \in \mathcal{V}$ being defined as:
	\begin{equation*}
		\Delta_v^p: ~ \mathcal{V} \longrightarrow \mathbb{R} \qquad v_i \longmapsto \Delta_v^p f \left(v_i\right) =
	\end{equation*}
	\begin{equation*}
		\sum_{v_j \in \mathcal{V}} \left(\delta\left(v_i, v_j\right) \left\lvert w_I \left(v_j\right)^\alpha w_G \left(v_j\right)^\epsilon f\left(v_j\right) - w_I \left(v_i\right)^\alpha w_G \left(v_i\right)^\eta f\left(v_i\right)\right\rvert^{p - 2} \right.
	\end{equation*}
	\begin{equation*}	
		\left(w_G \left(v_i\right)^\eta w_I \left(v_j\right)^\alpha w_G \left(v_j\right)^\epsilon f\left(v_j\right) - w_I \left(v_i\right)^\alpha w_G \left(v_i\right)^{2 \eta} f\left(v_i\right)\right)
	\end{equation*}
	\begin{equation*}
		- \delta\left(v_j, v_i\right) \left\lvert w_I \left(v_i\right)^\alpha w_G \left(v_i\right)^\epsilon f\left(v_i\right) - w_I \left(v_j\right)^\alpha w_G \left(v_j\right)^\eta f\left(v_j\right)\right\rvert^{p - 2}
	\end{equation*}
	\begin{equation}
		\left. \left(w_I \left(v_i\right)^\alpha w_G \left(v_i\right)^{2 \epsilon} f\left(v_i\right) - w_I \left(v_j\right)^\alpha w_G \left(v_j\right)^\eta w_G \left(v_i\right)^\epsilon  f\left(v_j\right)\right)\right) W_I \left(v_i, v_j\right)^\beta  W_G \left(v_i, v_j\right)^{p \gamma}.
	\end{equation}\vspace{0em}
\end{definition}

\begin{remark}[\textbf{Parameter choice for the vertex $p$-Laplacian operator}]\label{GDelta_v^pparam} \ \\
	The before defined parameter choices $\alpha = 0$, $\beta = 1$, $\gamma = \frac{1}{2}$, $\epsilon = 0$ and $\eta = 0$ together with $W_I \left(a_q\right) \equiv \frac{1}{2}$ from the inner product in $\mathcal{H}\left(\mathcal{A}_G\right)$ for all arcs $a_q \in \mathcal{A}_G$, lead to the common definition of the vertex $p$-Laplacian (Page 2421 in \cite{elmoataz2015p}: Anisotropic graph $p$-Laplacian) for all vertices $v_i \in \mathcal{V}$:
	\begin{equation*}
		\Delta_v^p f \left(v_i\right) = 
	\end{equation*}
	\begin{equation}
		\frac{1}{2} \sum_{v_j \in \mathcal{V}} \left(\delta\left(v_i, v_j\right) + \delta\left(v_j, v_i\right)\right) W_G \left(v_i, v_j\right)^{\frac{p}{2}} \left\lvert f\left(v_j\right) - f\left(v_i\right)\right\rvert^{p - 2}  \left(f\left(v_j\right) - f\left(v_i\right)\right).
	\end{equation}\vspace{0em}
\end{remark}

The following theorem proves that the $p$-Laplacian for vertices $v_i \in \mathcal{V}$ is well-defined with regard to the vertex gradient and the vertex divergence on normal graphs.\\

\begin{theorem}[\textbf{Connection vertex divergence $\text{div}_v$, vertex gradient $\nabla_v$, and vertex $p$-Laplacian $\Delta_v^p$}]\label{GvertexpLap} \ \\
	On a weighted oriented normal graph $OG = \left(\mathcal{V}, \mathcal{A}_G, w, W\right)$ with symmetric arc weight functions $W_I$ and $W_G$, the vertex $p$-Laplacian operator $\Delta_v^p$ fulfills the equality
	\begin{equation}
		\Delta_v^p f = \text{div}_v \left(\left\lvert \nabla_v f\right\rvert^{p - 2} \nabla_v f\right)
	\end{equation}
	for all vertex functions $f \in \mathcal{H}\left(\mathcal{V}\right)$.\\
\end{theorem}

\begin{proof}\ \\
	Given a weighted oriented normal graph $OG = \left(\mathcal{V}, \mathcal{A}_G, w, W\right)$ with symmetric arc weight functions $W_I$ and $W_G$ and a vertex function $f \in \mathcal{H}\left(\mathcal{V}\right)$, then using the definitions of the vertex divergence operator $\text{div}_v$ and the vertex gradient operator $\nabla_v$ yields the following for all vertices $v_i \in \mathcal{V}$:
	
	$\begin{aligned}[t]	
		\text{div}_v \left(\left\lvert \nabla_v f\right\rvert^{p - 2} \nabla_v f\right)\left(v_i\right) = & \sum_{v_j \in \mathcal{V}} \left(\delta\left(v_i, v_j\right) w_G \left(v_i\right)^\eta \left\lvert \nabla_v f\left(v_i, v_j\right)\right\rvert^{p - 2} \nabla_v f\left(v_i, v_j\right) - \right. &\\
		& \left. \delta\left(v_j, v_i\right) w_G \left(v_i\right)^\epsilon \left\lvert \nabla_v f\left(v_j, v_i\right)\right\rvert^{p - 2} \nabla_v f\left(v_j, v_i\right)\right) &\\
		& W_I \left(v_i, v_j\right)^\beta W_G \left(v_i, v_j\right)^\gamma &\\
	\end{aligned}$\\
	$\begin{aligned}[t]	
		= & \sum_{v_j \in \mathcal{V}} \left(\delta\left(v_i, v_j\right) w_G \left(v_i\right)^\eta \right. &\\
		& \quad ~ \left\lvert W_G \left(v_i, v_j\right)^\gamma \left(w_I \left(v_j\right)^\alpha w_G \left(v_j\right)^\epsilon f\left(v_j\right) - w_I \left(v_i\right)^\alpha w_G \left(v_i\right)^\eta f\left(v_i\right)\right)\right\rvert^{p - 2} &\\
		& \quad ~ W_G \left(v_i, v_j\right)^\gamma \left(w_I \left(v_j\right)^\alpha w_G \left(v_j\right)^\epsilon f\left(v_j\right) - w_I \left(v_i\right)^\alpha w_G \left(v_i\right)^\eta f\left(v_i\right)\right) - &\\
		& \quad ~ \delta\left(v_j, v_i\right) w_G\left(v_i\right)^\epsilon \left\lvert W_G \left(v_j, v_i\right)^\gamma \left(w_I \left(v_i\right)^\alpha w_G \left(v_i\right)^\epsilon f\left(v_i\right) - w_I \left(v_j\right)^\alpha w_G \left(v_j\right)^\eta f\left(v_j\right)\right)\right\rvert^{p - 2} &\\
		& \quad \left. W_G \left(v_j, v_i\right)^\gamma \left(w_I \left(v_i\right)^\alpha w_G \left(v_i\right)^\epsilon f\left(v_i\right) - w_I \left(v_j\right)^\alpha w_G \left(v_j\right)^\eta f\left(v_j\right)\right)\right) W_I \left(v_i, v_j\right)^\beta W_G \left(v_i, v_j\right)^\gamma &\\
	\end{aligned}$\\
	
	\clearpage
	Using the symmetry and the non-negativity properties of the arc weight functions $W_I$ and $W_G$ leads to the vertex $p$-Laplacian definition:
	
	$\begin{aligned}[t]	
		= & \sum_{v_j \in \mathcal{V}} \left(\delta\left(v_i, v_j\right) w_G \left(v_i\right)^\eta \left\lvert w_I \left(v_j\right)^\alpha w_G \left(v_j\right)^\epsilon f\left(v_j\right) - w_I \left(v_i\right)^\alpha w_G \left(v_i\right)^\eta f\left(v_i\right)\right\rvert^{p - 2} \right. &\\
		& \quad \left(w_I \left(v_j\right)^\alpha w_G \left(v_j\right)^\epsilon f\left(v_j\right) - w_I \left(v_i\right)^\alpha w_G \left(v_i\right)^\eta f\left(v_i\right)\right) - &\\
		& \quad ~ \delta\left(v_j, v_i\right) w_G \left(v_i\right)^\epsilon \left\lvert w_I \left(v_i\right)^\alpha w_G \left(v_i\right)^\epsilon f\left(v_i\right) - w_I \left(v_j\right)^\alpha w_G \left(v_j\right)^\eta f\left(v_j\right)\right\rvert^{p - 2} &\\
		& \quad \left. \left(w_I \left(v_i\right)^\alpha w_G \left(v_i\right)^\epsilon f\left(v_i\right) - w_I \left(v_j\right)^\alpha w_G \left(v_j\right)^\eta f\left(v_j\right)\right)\right) W_I \left(v_i, v_j\right)^\beta W_G \left(v_i, v_j\right)^{2 \gamma + \gamma \left(p - 2\right)} &\\
		= & \sum_{v_j \in \mathcal{V}} \left(\delta\left(v_i, v_j\right) \left\lvert w_I \left(v_j\right)^\alpha w_G \left(v_j\right)^\epsilon f\left(v_j\right) - w_I \left(v_i\right)^\alpha w_G \left(v_i\right)^\eta f\left(v_i\right)\right\rvert^{p - 2} \right. &\\
		& \quad \left(w_G \left(v_i\right)^\eta w_I \left(v_j\right)^\alpha w_G \left(v_j\right)^\epsilon f\left(v_j\right) - w_I \left(v_i\right)^\alpha w_G \left(v_i\right)^{2 \eta} f\left(v_i\right)\right) - &\\
		& \quad ~ \delta\left(v_j, v_i\right) \left\lvert w_I \left(v_i\right)^\alpha w_G \left(v_i\right)^\epsilon f\left(v_i\right) - w_I \left(v_j\right)^\alpha w_G \left(v_j\right)^\eta f\left(v_j\right) \right\rvert^{p - 2} &\\
		& \quad \left. \left(w_I \left(v_i\right)^\alpha w_G \left(v_i\right)^{2 \epsilon} f\left(v_i\right) - w_I \left(v_j\right)^\alpha w_G \left(v_j\right)^\eta w_G \left(v_i\right)^\epsilon f\left(v_j\right)\right)\right) W_I \left(v_i, v_j\right)^\beta W_G \left(v_i, v_j\right)^{p \gamma} &\\
		= & ~ \Delta_v^p f \left(v_i\right) &\\
	\end{aligned}$\\

	Therefore, with the definition for the vertex gradient $\nabla_v$, the vertex divergence $\text{div}_v$, and the vertex $p$-Laplacian $\Delta_v^p$, the equality $\Delta_v^p f\left(v_i\right) = \text{div}_v \left(\left\lvert \nabla_v f\right\rvert^{p - 2} \nabla_v f\right) \left(v_i\right)$ holds true for all vertices $v_i \in \mathcal{V}$ and for all vertex functions $f \in \mathcal{H}\left(\mathcal{V}\right)$.\\
\end{proof}

With the arc gradient and the arc divergence, the $p$-Laplacian operator can be defined for arcs $a_q \in \mathcal{A}_G$ as well.\\

\begin{definition}[\textbf{Arc $p$-Laplacian operator $\Delta_a^p$}]\label{GDelta_a^p} \ \\
	For a weighted oriented normal graph $OG = \left(\mathcal{V}, \mathcal{A}_G, w, W\right)$, the arc $p$-Laplacian operator $\Delta_a^p$ is defined as
	\begin{equation}
		\Delta_a^p: ~ \mathcal{H}\left(\mathcal{A}_G\right) \longrightarrow \mathcal{H}\left(\mathcal{A}_G\right) \qquad F \longmapsto \Delta_a^p F
	\end{equation}
	\clearpage
	with the weighted $p$-Laplacian of arc function $F \in \mathcal{H}\left(\mathcal{A}_G\right)$ at an arc $a_q = \left(v_i, v_j\right) \in \mathcal{A}_G$ being given by
	\begin{equation*}
		\Delta_a^p: ~ \mathcal{A}_G \longrightarrow \mathbb{R} \qquad a_q = \left(v_i, v_j\right) \longmapsto \Delta_a^p F \left(a_q\right) =
	\end{equation*}
	\begin{equation*}
		W_G \left(v_i, v_j\right)^\theta \left(\frac{w_I \left(v_i\right)^\alpha w_G \left(v_i\right)^{p \zeta}}{\deg_{out}\left(v_i\right)} \right.
	\end{equation*}
	\begin{equation*}
		\left\lvert \sum_{a_r \in \mathcal{A}_G} \left(\frac{\delta_{in}\left(v_i, a_r\right)}{\deg_{in}\left(v_i\right)} - \frac{\delta_{out}\left(v_i, a_r\right)}{\deg_{out}\left(v_i\right)}\right) W_I \left(a_r\right)^\beta W_G \left(a_r\right)^\theta F\left(a_r\right)\right\rvert^{p - 2}
	\end{equation*}
	\begin{equation*}
		\sum_{a_s \in \mathcal{A}_G} \left(\frac{\delta_{in}\left(v_i, a_s\right)}{\deg_{in}\left(v_i\right)} - \frac{\delta_{out}\left(v_i, a_s\right)}{\deg_{out}\left(v_i\right)}\right) W_I \left(a_s\right)^\beta W_G \left(a_s\right)^\theta F\left(a_s\right) -
	\end{equation*}
	\begin{equation*}
		\frac{w_I \left(v_j\right)^\alpha w_G \left(v_j\right)^{p \zeta}}{\deg_{in}\left(v_j\right)} \left\lvert \sum_{a_t \in \mathcal{A}_G} \left(\frac{\delta_{in}\left(v_j, a_t\right)}{\deg_{in}\left(v_j\right)} - \frac{\delta_{out}\left(v_j, a_t\right)}{\deg_{out}\left(v_j\right)}\right) W_I \left(a_t\right)^\beta W_G \left(a_t\right)^\theta F\left(a_t\right)\right\rvert^{p - 2}
	\end{equation*}
	\begin{equation}
		\left. \sum_{a_u \in \mathcal{A}_G} \left(\frac{\delta_{in}\left(v_j, a_u\right)}{\deg_{in}\left(v_j\right)} - \frac{\delta_{out}\left(v_j, a_u\right)}{\deg_{out}\left(v_j\right)}\right) W_I \left(a_u\right)^\beta W_G \left(a_u\right)^\theta F\left(a_u\right)\right).
	\end{equation}\vspace{0em}
\end{definition}

\begin{remark}[\textbf{Parameter choice for the arc $p$-Laplacian operator}]\label{GDelta_a^pparam} \ \\
	With the previous parameters $\alpha = 0$, $\beta = 1$, $\zeta = \frac{1}{2}$, $\theta = 0$ together with excluding the degree functions $\deg_{in}\left(v_i\right)$ and $\deg_{out}\left(v_i\right)$ for all vertices $v_i \in \mathcal{V}$ and with $W_I \left(a_q\right) \equiv \frac{1}{2}$ from the inner product in $\mathcal{H}\left(\mathcal{A}_G\right)$ for all arcs $a_q \in \mathcal{A}_G$, the arc $p$-Laplacian simplifies for all arcs $a_q = \left(v_i, v_j\right) \in \mathcal{V}$ to:
	\begin{equation*}
		\Delta_a^p F \left(a_q\right) = \Delta_a^p F \left(v_i, v_j\right) =
	\end{equation*}
	\begin{equation*}
		\frac{w_G\left(v_i\right)^{\frac{p}{2}}}{2^{p - 1}} \left\lvert \sum_{a_r \in \mathcal{A}_G} \left(\delta_{in}\left(v_i, a_r\right) - \delta_{out}\left(v_i, a_r\right)\right) F\left(a_r\right)\right\rvert^{p - 2}
	\end{equation*}
	\begin{equation*}
		\sum_{a_s \in \mathcal{A}_G} \left(\delta_{in}\left(v_i, a_s\right) - \delta_{out}\left(v_i, a_s\right)\right) F\left(a_s\right) -
	\end{equation*}
	\begin{equation*}
		\frac{w_G \left(v_j\right)^{\frac{p}{2}}}{2^{p - 1}} \left\lvert \sum_{a_t \in \mathcal{A}_G} \left(\delta_{in}\left(v_j, a_t\right) - \delta_{out}\left(v_j, a_t\right)\right) F\left(a_t\right)\right\rvert^{p - 2}
	\end{equation*}
	\begin{equation}
		\sum_{a_u \in \mathcal{A}_G} \left(\delta_{in}\left(v_j, a_u\right) - \delta_{out}\left(v_j, a_u\right)\right) F\left(a_u\right).
	\end{equation}\vspace{0em}
\end{remark}

The theorem below proves that the arc $p$-Laplacian is well defined by using the arc gradient and arc divergence for normal graphs.\\

\begin{theorem}[\textbf{Connection arc divergence $\text{div}_a$, arc gradient $\nabla_a$, and arc $p$-Laplacian $\Delta_a^p$}]\label{GarcpLap} \ \\
	On a weighted oriented normal graph $OG = \left(\mathcal{V}, \mathcal{A}_G, w, W\right)$, the arc $p$-Laplacian operator $\Delta_a^p$ fulfills the equality
	\begin{equation}
		\Delta_a^p F = \text{div}_a \left(\left\lvert \nabla_a F\right\rvert^{p - 2} \nabla_a F\right)
	\end{equation}
	for all arc functions $F \in \mathcal{H}\left(\mathcal{A}_G\right)$.\\
\end{theorem}

\begin{proof}\ \\
	Given a weighted oriented normal graph $OG = \left(\mathcal{V}, \mathcal{A}_G, w, W\right)$ and an arc function $F \in \mathcal{H}\left(\mathcal{A}_G\right)$, then the definitions of the arc divergence operator $\text{div}_a$ and the arc gradient operator $\nabla_a$ yields the following for all arcs $a_q = \left(v_i, v_j\right) \in \mathcal{A}_G$:
	
	$\begin{aligned}[t]	
		\text{div}_a \left(\left\lvert \nabla_a F\right\rvert^{p - 2} \nabla_a F\right)\left(v_i, v_j\right) = & ~ W_G \left(v_i, v_j\right)^\theta \left(\frac{w_I \left(v_i\right)^\alpha w_G \left(v_i\right)^\zeta}{\deg_{out}\left(v_i\right)} \left\lvert \nabla_a F \left(v_i\right)\right\rvert^{p - 2} \nabla_a F \left(v_i\right)\right. - &\\
		& \left. \frac{w_I \left(v_j\right)^\alpha w_G \left(v_j\right)^\zeta}{\deg_{in}\left(v_j\right)} \left\lvert \nabla_a F \left(v_j\right)\right\rvert^{p - 2} \nabla_a F \left(v_j\right)\right) &\\
	\end{aligned}$\\
	$\begin{aligned}[t]	
		= & ~ W_G \left(v_i, v_j\right)^\theta \left(\frac{w_I \left(v_i\right)^\alpha w_G \left(v_i\right)^\zeta}{\deg_{out}\left(v_i\right)} \right. &\\
		& \quad \left\lvert w_G \left(v_i\right)^\zeta \sum_{a_r \in \mathcal{A}_G} \left(\frac{\delta_{in}\left(v_i, a_r\right)}{\deg_{in}\left(v_i\right)} - \frac{\delta_{out}\left(v_i, a_r\right)}{\deg_{out}\left(v_i\right)}\right) W_I \left(a_r\right)^\beta W_G \left(a_r\right)^\theta F\left(a_r\right) \right\rvert^{p - 2} &\\
		& \quad w_G \left(v_i\right)^\zeta \sum_{a_s \in \mathcal{A}_G} \left(\frac{\delta_{in}\left(v_i, a_s\right)}{\deg_{in}\left(v_i\right)} - \frac{\delta_{out}\left(v_i, a_s\right)}{\deg_{out}\left(v_i\right)}\right) W_I \left(a_s\right)^\beta W_G \left(a_s\right)^\theta F\left(a_s\right) - &\\
		& \quad \frac{w_I \left(v_j\right)^\alpha w_G \left(v_j\right)^\zeta}{\deg_{in}\left(v_j\right)} &\\
		& \quad \left\lvert w_G \left(v_j\right)^\zeta \sum_{a_t \in \mathcal{A}_G} \left(\frac{\delta_{in}\left(v_j, a_t\right)}{\deg_{in}\left(v_j\right)} - \frac{\delta_{out}\left(v_j, a_t\right)}{\deg_{out}\left(v_j\right)}\right) W_I \left(a_t\right)^\beta W_G \left(a_t\right)^\theta F\left(a_t\right) \right\rvert^{p - 2} &\\
		& \quad \left. w_G \left(v_j\right)^\zeta \sum_{a_u \in \mathcal{A}_G} \left(\frac{\delta_{in}\left(v_j, a_u\right)}{\deg_{in}\left(v_j\right)} - \frac{\delta_{out}\left(v_j, a_u\right)}{\deg_{out}\left(v_j\right)}\right) W_I \left(a_u\right)^\beta W_G \left(a_u\right)^\theta F\left(a_u\right)\right) &\\
	\end{aligned}$\\
	$\begin{aligned}[t]
		= & ~ W_G \left(v_i, v_j\right)^\theta \left(\frac{w_I \left(v_i\right)^\alpha w_G \left(v_i\right)^{2 \zeta + \zeta \left(p - 2\right)}}{\deg_{out}\left(v_i\right)} \right. &\\
		& \quad \left\lvert \sum_{a_r \in \mathcal{A}_G} \left(\frac{\delta_{in}\left(v_i, a_r\right)}{\deg_{in}\left(v_i\right)} - \frac{\delta_{out}\left(v_i, a_r\right)}{\deg_{out}\left(v_i\right)}\right) W_I \left(a_r\right)^\beta W_G \left(a_r\right)^\theta F\left(a_r\right) \right\rvert^{p - 2} &\\
		& \quad \sum_{a_s \in \mathcal{A}_G} \left(\frac{\delta_{in}\left(v_i, a_s\right)}{\deg_{in}\left(v_i\right)} - \frac{\delta_{out}\left(v_i, a_s\right)}{\deg_{out}\left(v_i\right)}\right) W_I \left(a_s\right)^\beta W_G \left(a_s\right)^\theta F\left(a_s\right) - &\\
		& \quad \frac{w_I \left(v_j\right)^\alpha w_G \left(v_j\right)^{2 \zeta + \zeta \left(p - 2\right)}}{\deg_{in}\left(v_i\right)} &\\
		& \quad \left\lvert \sum_{a_t \in \mathcal{A}_G} \left(\frac{\delta_{in}\left(v_j, a_t\right)}{\deg_{in}\left(v_j\right)} - \frac{\delta_{out}\left(v_j, a_t\right)}{\deg_{out}\left(v_j\right)}\right) W_I \left(a_t\right)^\beta W_G \left(a_t\right)^\theta F\left(a_t\right) \right\rvert^{p - 2} &\\
		& \quad \left. \sum_{a_u \in \mathcal{A}_G} \left(\frac{\delta_{in}\left(v_j, a_u\right)}{\deg_{in}\left(v_j\right)} - \frac{\delta_{out}\left(v_j, a_u\right)}{\deg_{out}\left(v_j\right)}\right) W_I \left(a_u\right)^\beta W_G \left(a_u\right)^\theta F\left(a_u\right)\right) &\\
	\end{aligned}$\\
	$\begin{aligned}[t]
		= & ~ W_G \left(v_i, v_j\right)^\theta \left(\frac{w_I \left(v_i\right)^\alpha w_G \left(v_i\right)^{p \zeta}}{\deg_{out}\left(v_i\right)} \right. &\\
		& \quad \left\lvert \sum_{a_r \in \mathcal{A}_G} \left(\frac{\delta_{in}\left(v_i, a_r\right)}{\deg_{in}\left(v_i\right)} - \frac{\delta_{out}\left(v_i, a_r\right)}{\deg_{out}\left(v_i\right)}\right) W_I \left(a_r\right)^\beta W_G \left(a_r\right)^\theta F\left(a_r\right) \right\rvert^{p - 2} &\\
		& \quad \sum_{a_s \in \mathcal{A}_G} \left(\frac{\delta_{in}\left(v_i, a_s\right)}{\deg_{in}\left(v_i\right)} - \frac{\delta_{out}\left(v_i, a_s\right)}{\deg_{out}\left(v_i\right)}\right) W_I \left(a_s\right)^\beta W_G \left(a_s\right)^\theta F\left(a_s\right) - &\\
		& \quad \frac{w_I \left(v_j\right)^\alpha w_G \left(v_j\right)^{p \zeta}}{\deg_{in}\left(v_i\right)} &\\
		& \quad \left\lvert \sum_{a_t \in \mathcal{A}_G} \left(\frac{\delta_{in}\left(v_j, a_t\right)}{\deg_{in}\left(v_j\right)} - \frac{\delta_{out}\left(v_j, a_t\right)}{\deg_{out}\left(v_j\right)}\right) W_I \left(a_t\right)^\beta W_G \left(a_t\right)^\theta F\left(a_t\right) \right\rvert^{p - 2} &\\
		& \quad \left. \sum_{a_u \in \mathcal{A}_G} \left(\frac{\delta_{in}\left(v_j, a_u\right)}{\deg_{in}\left(v_j\right)} - \frac{\delta_{out}\left(v_j, a_u\right)}{\deg_{out}\left(v_j\right)}\right) W_I \left(a_u\right)^\beta W_G \left(a_u\right)^\theta F\left(a_u\right)\right) &\\
		= & ~ \Delta_a^p F \left(v_i, v_j\right)
	\end{aligned}$\\

	Thus, the equality $\Delta_a^p F\left(a_q\right) = \text{div}_a \left(\left\lvert \nabla_a F\right\rvert^{p - 2} \nabla_a F\right)\left(a_q\right)$ holds true for all arcs $a_q \in \mathcal{A}_G$ and for all arc functions $F \in \mathcal{H}\left(\mathcal{A}_G\right)$, when using the presented definitions of the arc gradient $\nabla_a$, the arc divergence $\text{div}_a$, and the arc $p$-Laplacian $\Delta_a^p$ .\\
\end{proof}
\clearpage
\section{Functions on hypergraphs}\label{8} 

Since hypergraphs are a valid generalization of normal graphs, the following sections aim at extending all previous definitions on normal grpahs to the hypergraph case. As on normal graphs, real functions can be defined both on the vertex set $\mathcal{V}$ and the hyperedge set $\mathcal{E}_H$ or the hyperarc set $\mathcal{A}_H$ of not oriented hypergraphs $NH$ or oriented hypergraphs $OH$.\\

Thus this section extends the definitions and concepts of section (\ref{5}) to hypergraphs, starting with vertex functions in subsection (\ref{8.1}). This subsection includes the introduction of the vertex weight function for hypergraphs and the definition of the space of all vertex functions as a Hilbert space with the corresponding inner product, absolute value and $\mathcal{L}^p$-norm of vertex functions.\\

Subsection (\ref{8.2}) then introduces hyperedge functions on not oriented hypergraphs and hyperarc functions on oriented hypergraphs in order to encode not only more information about pairwise connections between the data points $v_1, v_2, \dots v_N$, but about the relationship between an arbitrary subset of the data points. Moreover, hyperedge and hyperarc weight functions are defined and their symmetry property is analyzed. As in the previous subsection, the definition of the space of all hyperedge or hyperarc functions as a Hilbert space with the corresponding inner product, absolute value and $\mathcal{L}^p$-norm is given as well.\\

\subsection{Vertex functions}\label{8.1} 

Since vertex functions on normal graphs and hypergraphs are defined on the same set of vertices $\mathcal{V} = \left\{v_1, v_2, \dots v_N\right\}$, the corresponding real vertex function definitions and remarks for hypergraphs are only introduced briefly.\\

\begin{definition}[\textbf{Real vertex function $f$}]\label{Hf} \ \\	
	A real vertex function $f$ is defined on the domain of the vertices $\mathcal{V}$ of a not oriented hypergraph $NH = \left(\mathcal{V}, \mathcal{E}_H\right)$ or of an oriented hypergraph $OH = \left(\mathcal{V}, \mathcal{A}_H\right)$:
	\begin{equation}
		f: ~ \mathcal{V} \longrightarrow \mathbb{R} \qquad v_i \longmapsto f\left(v_i\right).
	\end{equation}\vspace{0em}
\end{definition}

An important example of real vertex functions are vertex weight functions, which assign a positive weight to each vertex $v_i \in \mathcal{V}$.\\

\clearpage
\begin{definition}[\textbf{Vertex weight function $w$}]\label{Hw} \ \\
	The vertex weight function $w$ can be defined for both a not oriented hypergraph $NH = \left(\mathcal{V}, \mathcal{E}_H\right)$ and an oriented hypergraph $OH = \left(\mathcal{V}, \mathcal{A}_H\right)$ in the same way:
	\begin{equation}
		w: ~ \mathcal{V} \longrightarrow \mathbb{R}_{> 0} \qquad v_i \longmapsto w\left(v_i\right)
	\end{equation}
	A not oriented hypergraph $NH = \left(\mathcal{V}, \mathcal{E}_H\right)$ or an oriented hypergraph $OH = \left(\mathcal{V}, \mathcal{A}_H\right)$ together with a vertex weight function $w$ can be written as $NH = \left(\mathcal{V}, \mathcal{E}_H, w\right)$ or $OH = \left(\mathcal{V}, \mathcal{A}_H, w\right)$.\\
\end{definition}

\begin{remark}[\textbf{Vertex weight functions $w_I$ and $w_G$}]\label{Hw_Iw_G} \ \\
	Different vertex weight functions $w$ of a weighted oriented hypergraph $OH = \left(\mathcal{V}, \mathcal{A}_H, w\right)$ will be used for the inner product of the space of real vertex functions and also for the vertex and hyperarc gradient definitions. Hence, different subscripts will indicate again if the vertex weight function stems from the inner product $w_I$ or the gradient definitions $w_G$.\\
\end{remark}

Other important examples for real vertex functions on not oriented hypergraphs $NH = \left(\mathcal{V}, \mathcal{E}_H\right)$ or on oriented hypergraphs $OH = \left(\mathcal{V}, \mathcal{A}_H\right)$ include the previously introduced vertex degree functions $\deg$, $\deg_{out}$, and $\deg_{in}$ of definition (\ref{degH}).\\

Since the number of vertices $|\mathcal{V}| = N \in \mathbb{N}$ is assumed to be finite, it is possible to associate the space of all real vertex functions $f$ with an $N-$dimensional Hilbert space.\\

\begin{definition}[\textbf{Space of real vertex functions $\mathcal{H}\left(\mathcal{V}\right)$}]\label{HH(V)} \ \\
	The space of real vertex functions $f$ is given by:
	\begin{equation}
		\mathcal{H}\left(\mathcal{V}\right) = \left\{f ~ \middle| ~ f: ~ \mathcal{V} \longrightarrow \mathbb{R}\right\}
	\end{equation}
	with the corresponding inner product:
	\begin{equation}
		{\langle f, g \rangle}_{\mathcal{H}\left(\mathcal{V}\right)} = \sum_{v_i \in \mathcal{V}} w_I \left(v_i\right)^\alpha f\left(v_i\right) g\left(v_i\right)
	\end{equation}
	for any two real vertex functions $f, g \in \mathcal{H}\left(\mathcal{V}\right)$ and parameter $\alpha \in \mathbb{R}$.\\
\end{definition}

\begin{remark}[\textbf{Parameter choice for the inner product of real vertex functions}]\label{HinproH(V)} \ \\
	Choosing $\alpha = 0$ results in a simplified definition of the inner product for real vertex functions, which matches the common inner product for real vertex functions on normal graphs (\ref{inproH(V)}) for any two real vertex functions $f, g \in \mathcal{H}\left(\mathcal{V}\right)$:
	\begin{equation}
		{\langle f, g \rangle}_{\mathcal{H}\left(\mathcal{V}\right)} = \sum_{v_i \in \mathcal{V}} f\left(v_i\right) g\left(v_i\right).
	\end{equation}\vspace{0em}
\end{remark}

The previous definitions of the absolute value and the $\mathcal{L}^p$-norm for real vertex functions on normal graphs in (\ref{Lpf}) can be generalized in order to also fit real vertex functions on hypergraphs.\\

\begin{definition}[\textbf{Absolute value and $\mathcal{L}^p$-norm on the space of real vertex functions}]\label{HLpf} \ \\
	For a real vertex function $f \in \mathcal{H}\left(\mathcal{V}\right)$ at a vertex $v_i \in \mathcal{V}$, the absolute value is given by:
	\begin{equation*}
		\left\lvert ~ \cdot ~ \right\rvert: ~ \mathbb{R} \longrightarrow \mathbb{R}_{\geq 0}
	\end{equation*}
	\begin{equation}
		f\left(v_i\right) \longmapsto \left\lvert f\left(v_i\right)\right\rvert = \left\{\begin{array}{ll}
			f\left(v_i\right) & \quad f\left(v_i\right) \geq 0\\
			-f\left(v_i\right) & \quad \text{otherwise}
		\end{array}\right..
	\end{equation}
	Moreover, the $\mathcal{L}^p$-norm of a real vertex function $f \in \mathcal{H}\left(\mathcal{V}\right)$ is defined as:
	\begin{equation*}
		\left\lvert \left\lvert ~ \cdot ~ \right\rvert\right\rvert_p: ~ \mathcal{H}\left(\mathcal{V}\right) \longrightarrow \mathbb{R}_{\geq 0}
	\end{equation*}
	\begin{equation}
		f \longmapsto \left\lvert \left\lvert f \right\rvert\right\rvert_p = \left\{\begin{array}{ll}
			\left(\sum_{v_i \in \mathcal{V}} \left\lvert f\left(v_i\right)\right\rvert^p\right)^{\frac{1}{p}} & \quad 1 \leq p < \infty\\
			\max_{v_i \in \mathcal{V}} \left(\left\lvert f\left(v_i\right) \right\rvert\right) & \quad p = \infty
		\end{array}\right..
	\end{equation}\vspace{0em}
\end{definition}
\subsection{Hyperedge and hyperarc functions}\label{8.2} 

In order to describe a relationship between two or more vertices, hyperedges in a not oriented hypergraph $NH = \left(\mathcal{V}, \mathcal{E}_H\right)$ or hyperarcs in an oriented hypergraph $OH = \left(\mathcal{V}, \mathcal{A}_H\right)$ can be used. Further information about such connections can be encoded with hyperedge functions defined on the hyperedge set $\mathcal{E}_H$ or with hyperarc functions defined on the hyperarc set $\mathcal{A}_H$.\\

\begin{definition}[\textbf{Real hyperedge or hyperarc function $F$}]\label{HF} \ \\
	A real hyperedge or hyperarc function $F$ is defined on the domain of the hyperedges $\mathcal{E}_H$ of a not oriented hypergraph $NH = \left(\mathcal{V}, \mathcal{E}_H\right)$ or on the domain of the hyperarcs $\mathcal{A}_H$ of an oriented hypergraph $OH = \left(\mathcal{V}, \mathcal{A}_H\right)$
	\begin{equation}
		F: ~ \mathcal{E}_H \longrightarrow \mathbb{R} \qquad e_q \longmapsto F\left(e_q\right),
	\end{equation}
	or alternatively in case of an oriented hypergraph
	\begin{equation}
		F: ~ \mathcal{A}_H \longrightarrow \mathbb{R} \qquad a_q \longmapsto F\left(a_q\right).
	\end{equation}\vspace{0em}
\end{definition}

Since the order of vertices in a hyperedge does not matter, any real hyperedge function is always symmetric.\\

\begin{remark}[\textbf{Symmetry of a real hyperedge function}]\label{HFsymm} \ \\
	For a not oriented hypergraph $NH = \left(\mathcal{V}, \mathcal{E}_H\right)$, a hyperedge function $F: ~ \mathcal{E}_H \longrightarrow \mathbb{R}$ is always symmetric since hyperedges have no orientation \big(Note: Any permutation of the vertices within the hyperedge describes the same hyperedge.\big).
	
	This property is generally not fulfilled by a real hyperarc function $F ~ \mathcal{A}_H \longrightarrow \mathbb{R}$, thus $F\left(a_q\right) = F\left(\widetilde{a_q}\right)$ can not be assumed to hold true for all hyperarcs $a_q \in \mathcal{A}_H$ and $\widetilde{a_q} \in \widetilde{\mathcal{A}_H}$, especially since $a_q = \left(a_q^{out}, a_q^{in}\right) \in \mathcal{A}_H$ does not imply $\widetilde{a_q} = \left(a_q^{in}, a_q^{out}\right) \in \mathcal{A}_H$.\\
\end{remark}

For the inner product on the space of real hyperedge or hyperarc functions and for the weighted gradients on oriented hypergraphs, the definition of a hyperedge and hyperarc weight function is necessary\\

\begin{definition}[\textbf{Hyperedge or hyperarc weight function $W$}]\label{HW} \ \\
	The hyperedge or hyperarc weight function $W$ can be defined for both a not oriented hypergraph $NH = \left(\mathcal{V}, \mathcal{E}_H\right)$ and an oriented hypergraph $OH = \left(\mathcal{V}, \mathcal{A}_H\right)$:
	\begin{itemize}
		\item[1)] In case of a not oriented hypergraph $NH$ the hyperedge weight function $W$ is given by:
		\begin{equation}
			W: ~ \mathcal{E}_H \longrightarrow \mathbb{R}_{> 0} \qquad e_q \longmapsto W\left(e_q\right).
		\end{equation}
		\item[2)] Similarly, in case of an oriented hypergraph $OH$, the weight function $W$ is defined on the domain of all hyperarcs $\mathcal{A}_H$ instead of all hyperedges $\mathcal{E}_H$:
		\begin{equation}
			W: ~ \mathcal{A}_H \longrightarrow \mathbb{R}_{> 0} \qquad a_q \longmapsto W\left(a_q\right).
		\end{equation}
	\end{itemize}
	A not oriented hypergraph $NH = \left(\mathcal{V}, \mathcal{E}_H\right)$ together with a hyperedge weight function $W$ or an oriented hypergraph $OH = \left(\mathcal{V}, \mathcal{A}_H\right)$ together with a hyperarc weight function $W$ are called weighted not oriented hypergraph $NH = \left(\mathcal{V}, \mathcal{E}_H, W\right)$ or weighted oriented hypergraph $OH = \left(\mathcal{V}, \mathcal{A}_H, W\right)$.\\
\end{definition}

\begin{remark}[\textbf{Hyperarc weight functions $W_I$ and $W_G$}]\label{HW_IW_G} \ \\
	Different hyperarc weight functions $W$ of a weighted oriented normal graph $OG = \left(\mathcal{V}, \mathcal{A}_G, W\right)$ will be used for the inner product of the space of real hyperarc functions and also for the vertex and hyperarc gradient definitions. Therefore, a subscript will indicate if the hyperarc weight function stems from the inner product $W_I$ or the gradient definitions $W_G$.\\
\end{remark}

As argued before, a real hyperarc function is not necessarily symmetric and therefore the symmetry of a hyperarc weight function $W$ on an oriented normal graph is a special property.\\

\begin{definition}[\textbf{Symmetric hyperarc weight function}]\label{HWsymm} \ \\
	A hyperarc weight function $W$ on an oriented hypergraph $OH = \left(\mathcal{V}, \mathcal{A}_H\right)$ is called symmetric if for all hyperarcs $a_q \in \mathcal{A}_H$ it holds true that:
	\begin{equation}
		W\left(a_q\right) = W\left(\widetilde{a_q}\right)
	\end{equation}
	with $\widetilde{a_q} = \left(a_q^{in}, a_q^{out}\right) \in \widetilde{\mathcal{A}_H}$ from the definition of the oriented hypergraph with switched orientation (\ref{switchedOH}).
	
	Note: This property either needs an extension of the domain of the hyperarc weight function $W$ from $\mathcal{A}_H$ to $\mathcal{A}_H \cup \widetilde{\mathcal{A}_H}$ or an oriented hypergraph with the property $\mathcal{A}_H = \widetilde{\mathcal{A}_H}$, which then implies for all hyperarcs $a_q \in \mathcal{A}_H$ that $\widetilde{a_q} \in \mathcal{A}_H$ also holds true.\\
\end{definition}

Every hyperedge or hyperarc weight function can be rescaled such that they only map to values less or equal to $1$.\\

\begin{theorem}[\textbf{Normalization of the hyperedge or hyperarc weight function}]\label{HWnorm} \ \\
	Every hyperedge or hyperarc weight function $W$ can be normalized in order to map to values in $\left(0, 1\right]$ instead of $\mathbb{R}_{> 0}$.\\
\end{theorem}

\begin{proof}\ \\
	All arguments of the following proof can be applied to both a not oriented hypergraph $NH$ and an oriented hypergraph $OH$, since the orientation of the hyperedges or hyperarcs and the symmetry of the weight function $W$ are not used in the argumentation. Thus, without loss of generality it is assumed that $W$ is defined on the domain of the hyperedge set $\mathcal{E}_H$ of a not oriented hypergraph $NH$.\\
	
	For every hyperedge $e_q \in \mathcal{E}_H$ it holds true that $W\left(e_q\right) \in \mathbb{R}_{> 0}$ and since the hyperedge set $\mathcal{E}_H$ is finite, the following maximum is well-defined:
	\begin{equation*}
		m := \max_{e_q \in \mathcal{E}_H} W\left(e_q\right) \in \mathbb{R}_{> 0}.
	\end{equation*}

	With the maximum weight $m$ a new hyperedge weight function $\overline{W}$ can be defined:
	\begin{equation}
		\overline{W}: ~ \mathcal{E}_H \longrightarrow \left(0, 1\right] \qquad e_q \longmapsto \frac{W\left(e_q\right)}{m}
	\end{equation}
	where $\frac{W\left(e_q\right)}{m} \in \left(0, 1\right]$ for all hyperedges $e_q \in \mathcal{E}_H$ holds true because:
	\begin{itemize}
		\item $W\left(e_q\right), m > 0 \quad \Longrightarrow \quad \overline{W}\left(e_q\right) = \frac{W\left(e_q\right)}{m} > 0$
		\item $W\left(e_q\right) \leq  m \quad \Longrightarrow \quad \overline{W}\left(e_q\right) = \frac{W\left(e_q\right)}{m} \leq 1$
	\end{itemize}
	Hence, any hyperedge or hyperarc weight function $W$ can be normalized in order to map to $\left(0, 1\right]$.\\
\end{proof}

The following theorem proves that the number of hyperedges in a not oriented hypergraph $NH = \left(\mathcal{V}, \mathcal{E}_H\right)$ and the the number of hyperarcs in an oriented hypergraph $OH = \left(\mathcal{V}, \mathcal{A}_H\right)$ are limited by $N^N$. Therefore, the space of real hyperedge or hyperarc functions can be described by an at most $N^N$-dimensional Hilbert space.\\

\begin{theorem}[\textbf{Number of hyperedges or hyperarcs in hypergraphs}]\label{numhypedgearc} \ \\
	With the assumption that no hyperedges or hyperarcs exist more than once, as argued in remark (\ref{Hfinite}), the number of hyperedges in a not oriented hypergraph $NH = \left(\mathcal{V}, \mathcal{E}_H\right)$ is limited by
	\begin{equation}
		\left\lvert \mathcal{E}_H \right\rvert \leq N^N
	\end{equation}
	and respectively, the number of hyperarcs in an oriented hypergraph $OH = \left(\mathcal{V}, \mathcal{A}_H\right)$ is also limited by
	\begin{equation}
		|\mathcal{A}_H| \leq N^N
	\end{equation}
	where $N \in \mathbb{N}_{\geq 1}$ is the cardinality of the vertex set $\mathcal{V}$.\\
\end{theorem}

\begin{proof}\ \\
	Given a not oriented hypergraph $NH = \left(\mathcal{V}, \mathcal{E}_H\right)$, there are at most:
	\begin{itemize}
		\item $\left(\begin{array}{c} N \\ 2 \end{array}\right) = \frac{N!}{2! \left(N - 2\right)!} \leq \frac{N!}{\left(N - 2\right)!} = N \cdot \left(N - 1\right)$ hyperedges connecting $2$ vertices
		\item $\left(\begin{array}{c} N \\ 3 \end{array}\right) = \frac{N!}{3! \left(N - 3\right)!} \leq \frac{N!}{\left(N - 3\right)!} = N \cdot \left(N - 1\right) \cdot \left(N - 2\right)$ hyperedges connecting $3$ vertices
		\item[] $\dots$
		\item $\left(\begin{array}{c} N \\ N \end{array}\right) = \frac{N!}{N! \left(N - N\right)!} \leq \frac{N!}{\left(N - N\right)!} = N \cdot \left(N - 1\right) \cdot ~ \dots ~ \cdot 1$ hyperedges connecting $N$ vertices
	\end{itemize}

	Therefore, the following estimation holds true:
	
	$\begin{aligned}[t]	
		\left\lvert \mathcal{E}_H \right\rvert \leq & \left(\begin{array}{c} N \\ 2 \end{array}\right) + \left(\begin{array}{c} N \\ 3 \end{array}\right) + \dots + \left(\begin{array}{c} N \\ N \end{array}\right) &\\
		\leq & ~ N \cdot \left(N - 1\right) + N \cdot \left(N - 1\right) \cdot \left(N - 2\right) + \dots + N \cdot \left(N - 1\right) \cdot ~ \dots ~ \cdot 1 &\\
		\leq & ~ N \cdot \left(N - 1\right) \cdot ~ \dots ~ \cdot 1 + N \cdot \left(N - 1\right) \cdot ~ \dots ~ \cdot 1 + \dots + N \cdot \left(N - 1\right) \cdot ~ \dots ~ \cdot 1 &\\
		\leq & \left(N - 1\right) N! \leq \left(N - 1\right) N^{\left(N - 1\right)} \leq N^N &\\
	\end{aligned}$\\
	
	\clearpage
	Similarly, given an oriented hypergraph $OH = \left(\mathcal{V}, \mathcal{A}_H\right)$, there are at most:
	\begin{itemize}
		\item $2^2 \left(\begin{array}{c} N \\ 2 \end{array}\right) = 2^2 \frac{N!}{2! \left(N - 2\right)!}$ hyperarcs connecting $2$ vertices
		\item $2^3 \left(\begin{array}{c} N \\ 3 \end{array}\right) = 2^3 \frac{N!}{3! \left(N - 3\right)!}$ hyperarcs connecting $3$ vertices
		\item[] $\dots$
		\item $2^N \left(\begin{array}{c} N \\ N \end{array}\right) = 2^N \frac{N!}{N! \left(N - N\right)!}$ hyperarcs connecting $N$ vertices
	\end{itemize}

	In comparison to the not oriented hypergraph, the differentiation between input and output vertices has to be considered as well. This leads to the multiplicative factor $2^k$ for $k  \in \left\{2, 3, \dots, N\right\}$, since for every vertex out of the $k$ vertices there are two possibilities: being an output vertex or being an input vertex.\\
	
	Mathematical induction over the number of vertices $N$ is used in order to prove:\\
	$\sum_{k = 2}^{N} 2^k \frac{N!}{k! \left(N - k\right)!} \leq N^N$ which then proves $|\mathcal{A}_H| \leq N^N$\\
	
	The theoretical base case $N = 1$ is proven by: 
	\begin{equation*}
		N = 1 \quad \Longrightarrow \quad \left\lvert A_H\right\rvert = 0 \quad \Longrightarrow \quad \left\lvert A_H\right\rvert = 0 \leq 1 = N^N
	\end{equation*}
	
	Furthermore, the base case of the mathematical induction $N = 2$ is proven by:
	\begin{equation*}
		\sum_{k = 2}^{N} 2^k \frac{N!}{k! \left(N - k\right)!} = 2^2 \frac{2!}{2! 0!} = 4 \leq 2^2 = N^N
	\end{equation*}
	
	Assume now that the inductive hypothesis is given by:
	\begin{equation*}
		\sum_{k = 2}^{N} 2^k \frac{N!}{k! \left(N - k\right)!} \leq N^N
	\end{equation*}
	and holds true for some $N \geq 2$. Then the following reformulations and arguments show that the inequality is also correct for $N + 1$ given the inductive hypothesis.
	
	$\begin{aligned}[t]	
		\sum_{k = 2}^{N + 1} 2^k \frac{\left(N + 1\right)!}{k! \left(N + 1 - k\right)!} & = 2^{N + 1} \frac{\left(N + 1\right)!}{\left(N + 1\right)! 0!} +
		\sum_{k = 2}^{N} 2^k \frac{\left(N + 1\right)!}{k! \left(N + 1 - k\right)!} &\\
		& = 2^{N + 1} + \left(N + 1\right) \sum_{k = 2}^{N} 2^k \frac{N!}{k! \left(N + 1 - k\right)!} &\\
		& \leq 2^{N + 1} + \left(N + 1\right) \sum_{k = 2}^{N} 2^k \frac{N!}{k! \left(N - k\right)!} &\\
		& \stackrel{\text{IH}}{\leq} 2^{N + 1} + \left(N + 1\right) N^N &\\
	\end{aligned}$\\

	Using $\left(N + 1\right)^{N + 1} > 0$ yields the following:
	
	$\sum_{k = 2}^{N + 1} 2^k \frac{\left(N + 1\right)!}{k! \left(N + 1 - k\right)!} \stackrel{\text{!}}{\leq} \left(N + 1\right)^{N + 1} \quad \Longleftarrow \quad 2^{N + 1} + \left(N + 1\right) N^N \leq \left(N + 1\right)^{N + 1} \quad \Longleftrightarrow$
	
	$\frac{2^{N + 1} + \left(N + 1\right) N^N}{\left(N + 1\right)^{N + 1}} \leq 1 \quad \Longleftrightarrow \quad \left(\frac{2}{N + 1}\right)^{N + 1} + \left(\frac{N}{N + 1}\right)^N \leq 1$\\
	
	In order to prove the last inequality helper functions
	\begin{equation*}
		h_1: ~ \mathbb{R}_{\geq 2} \longrightarrow \mathbb{R} \quad x \longmapsto h_1\left(x\right) = \left(\frac{2}{x + 1}\right)^{x + 1}
	\end{equation*}
	\begin{equation*}
		h_2: ~ \mathbb{R}_{\geq 2} \longrightarrow \mathbb{R} \quad x \longmapsto h_2\left(x\right) = \left(\frac{x}{x + 1}\right)^x
	\end{equation*}
	are defined in order to show:
	
	$h_1\left(2\right) + h_2\left(2\right) \leq 1$ and $h_1'\left(x\right), h_2'\left(x\right) \leq 0$ for all $x \geq 2$ which then implies:\\
	$\left(\frac{2}{N + 1}\right)^{N + 1} + \left(\frac{N}{N + 1}\right)^N \leq 1$ for all $N \in \mathbb{N}_{\geq 2}$
	\begin{itemize}
		\item
		$h_1\left(2\right) + h_2\left(2\right) = \left(\frac{2}{2 + 1}\right)^{2 + 1} + \left(\frac{2}{2 + 1}\right)^2 = \frac{8}{27} + \frac{4}{9} = \frac{20}{27} \leq 1$
		\item $h_1'\left(x\right) = \left(\frac{2}{x + 1}\right)^{x + 1} \left(\ln\left(\frac{2}{x + 1}\right) - 1\right) \leq \left(\frac{2}{x + 1}\right)^{x + 1} \left(\ln\left(1\right) - 1\right) = - \left(\frac{2}{x + 1}\right)^{x + 1} \leq 0$
		\item $h_2'\left(x\right) = \left(\frac{x}{x + 1}\right)^x \left(\ln\left(\frac{x}{x + 1}\right) + \frac{1}{x + 1}\right) \stackrel{\text{!}}{\leq} 0 \quad \Longleftrightarrow \quad \ln\left(\frac{x}{x + 1}\right) + \frac{1}{x + 1} \leq 0$
		
		Where for $\ln\left(\frac{x}{x + 1}\right) + \frac{1}{x + 1}$ it holds true that:
		\begin{itemize}
			\item[1)] $\ln\left(\frac{2}{2 + 1}\right) + \frac{1}{2 + 1} \leq -0.05 \leq 0$
			\item[2)] $\lim_{x \longrightarrow \infty} \left(\ln\left(\frac{x}{x + 1}\right) + \frac{1}{x + 1}\right) = \lim_{x \longrightarrow \infty} \left(\ln\left(\frac{1}{1 + \frac{1}{x}}\right) + \frac{1}{x + 1}\right) = \ln\left(1\right) + 0 = 0$
			\item[3)] $\frac{d}{dx}\left(\ln\left(\frac{x}{x + 1}\right) + \frac{1}{x + 1}\right) = \frac{1}{x} - \frac{1}{x + 1} - \frac{1}{\left(x + 1\right)^2} \stackrel{\text{!}}{\geq} 0 \quad \Longleftrightarrow \quad \frac{1}{x} \geq \frac{1}{x + 1} + \frac{1}{\left(x + 1\right)^2} \quad \Longleftrightarrow$\\
			$\left(x + 1\right)^2 \geq x \left(x + 1\right) + x \quad \Longleftrightarrow \quad x^2 + 2 x + 1 \geq x^2 + 2x \quad \Longleftrightarrow \quad 1 \geq 0$
		\end{itemize}
		Hence, the equality $\ln\left(\frac{x}{x + 1}\right) + \frac{1}{x + 1} \leq 0$ is correct for $x \geq 2$.		
	\end{itemize}

	This proves the inductive step and thus the upper bound for the number of hyperarcs $\left\lvert\mathcal{A}_H\right\rvert$:\\
	$\left(\frac{2}{N + 1}\right)^{N + 1} + \left(\frac{N}{N + 1}\right)^N \leq 1 \quad \Longrightarrow \quad \sum_{k = 2}^{N + 1} 2^k \frac{\left(N + 1\right)!}{k! \left(N + 1 - k\right)!} {\leq} \left(N + 1\right)^{N + 1} \quad \Longrightarrow \quad \left\lvert\mathcal{A}_H\right\rvert \leq N^N$\\
\end{proof}

With the extension of the domain of all real hyperedge functions from $\mathcal{E}_H$ to the set of all feasible hyperedges $2^\mathcal{V} \backslash \left\{\emptyset, \left\{v_1\right\}, \left\{v_2\right\}, \dots \left\{v_N\right\}\right\}$ and similarly the extension of the domain of all hyperarc functions from $\mathcal{A}_H$ to the set of all feasible hyperarcs $\left\{\left(a_q^{out}, a_q^{in}\right) ~ \middle| ~ a_q^{out}, a_q^{in} \in 2^\mathcal{V} \backslash \left\{\emptyset\right\} ~ \text{with} ~ a_q^{out} \cap a_q^{in} = \emptyset\right\}$ by setting $F\left(e_q\right) = F\left(a_q\right) = 0$ for any hyperedge $e_q \notin \mathcal{E}_H$ or any hyperarc $a_q \notin \mathcal{A}_H$, the space of all real hyperedge or hyperarc functions can be identified as an at most $N^N-$dimensional Hilbert space.\\

\begin{definition}[\textbf{Space of real hyperedge functions $\mathcal{H}\left(\mathcal{E}_H\right)$ or hyperarc functions $\mathcal{H}\left(\mathcal{A}_H\right)$}]\label{HH(E)(A)} \ \\
	The space of real hyperedge functions $F$ for a not oriented hypergraph $NH = \left(\mathcal{V}, \mathcal{E}_H, W\right)$ is given by
	\begin{equation}
		\mathcal{H}\left(\mathcal{E}_H\right) = \left\{F ~ \middle| ~ F: ~ \mathcal{E}_H \longrightarrow \mathbb{R}\right\}
	\end{equation}
	with the corresponding inner product
	\begin{equation}
		{\langle F, G \rangle}_{\mathcal{H}\left(\mathcal{E}_H\right)} = \sum_{e_q \in \mathcal{E}_H} W_I \left(e_q\right)^\beta F\left(e_q\right) G\left(e_q\right)
	\end{equation}
	for any two real hyperedge functions $F, G \in \mathcal{H}\left(\mathcal{E}_H\right)$ and parameter $\beta \in \mathbb{R}$. Similarly, in the case of an oriented hypergraph $OH = \left(\mathcal{V}, \mathcal{A}_H, W\right)$, the space of real hyperarc functions is defined as
	\begin{equation}
		\mathcal{H}\left(\mathcal{A}_H\right) = \left\{F ~ \middle| ~ F: ~ \mathcal{A}_H \longrightarrow \mathbb{R}\right\}
	\end{equation}
	with the corresponding inner product
	\begin{equation}
		{\langle F, G \rangle}_{\mathcal{H}\left(\mathcal{A}_H\right)} = \sum_{a_q \in \mathcal{A}_H} W_I \left(a_q\right)^\beta F\left(a_q\right) G\left(a_q\right)
	\end{equation}
	for any two real hyperarc functions $F, G \in \mathcal{H}\left(\mathcal{A}_H\right)$ and parameter $\beta \in \mathbb{R}$.\\
\end{definition}

\begin{remark}[\textbf{Parameter choice for the inner product of real hyperedge or hyperarc functions}]\label{HinproH(E)H(A)} \ \\
	By choosing $\beta = 1$, $W_I \left(e_q\right) \equiv \frac{1}{2}$ for all edges $e_q \in \mathcal{E}_H$ and $W_I \left(a_q\right) \equiv \frac{1}{2}$ for all arcs $a_q \in \mathcal{A}_H$, the inner product on the space of real hyperedge or hyperarc functions can be simplified in order to match the definition of the normal graph case (\ref{inproH(E)H(A)}) for any two real hyperedge functions $F, G \in \mathcal{H}\left(\mathcal{E}_H\right)$:
	\begin{equation}
		{\langle F, G \rangle}_{\mathcal{H}\left(\mathcal{E}_H\right)} = \frac{1}{2} \sum_{e_q \in \mathcal{E}_H} F\left(e_q\right) G\left(e_q\right)
	\end{equation}
	and similarly for any two real hyperarc functions $F, G \in \mathcal{H}\left(\mathcal{A}_H\right)$:
	\begin{equation}
		{\langle F, G \rangle}_{\mathcal{H}\left(\mathcal{A}_H\right)} = \frac{1}{2} \sum_{a_q \in \mathcal{A}_H} F\left(a_q\right) G\left(a_q\right).
	\end{equation}\vspace{0em}
\end{remark}

The previous definitions of the absolute value and the $\mathcal{L}^p$-norm for real edge and arc functions on normal graphs in (\ref{LpF}) can be generalized in order to also be applicable to hyperedge and hyperarc functions in the hypergraph case.\\

\clearpage
\begin{definition}[\textbf{Absolute value and $\mathcal{L}^p$-norm on the space of real hyperedge and hyperarc functions}]\label{HLpF} \ \\
	The absolute value and the $\mathcal{L}^p$-norm can be defined for real hyperedge functions $F \in \mathcal{H}\left(\mathcal{E}_H\right)$ and real hyperarc functions $F \in \mathcal{H}\left(\mathcal{A}_H\right)$:
	\begin{itemize}
		\item[1)] For a real hyperedge function $F \in \mathcal{H}\left(\mathcal{E}_H\right)$ at a hyperedge $e_q \in \mathcal{E}_H$, the absolute value is defined as:
		\begin{equation*}
			\left\lvert ~ \cdot ~ \right\rvert: ~ \mathbb{R} \longrightarrow \mathbb{R}_{\geq 0}
		\end{equation*}
		\begin{equation}
			F\left(e_q\right) \longmapsto \left\lvert F\left(e_q\right)\right\rvert = \left\{\begin{array}{ll}
				F\left(e_q\right) & \quad F\left(e_q\right) \geq 0\\
				-F\left(e_q\right) & \quad \text{otherwise}
			\end{array}\right..
		\end{equation}
		
		and the $\mathcal{L}^p$-norm of the hyperedge function $F$ is given by:
		\begin{equation*}
			\left\lvert \left\lvert ~ \cdot ~ \right\rvert\right\rvert_p: ~ \mathcal{H}\left(\mathcal{E}_H\right) \longrightarrow \mathbb{R}_{\geq 0}
		\end{equation*}
		\begin{equation}
			F \longmapsto \left\lvert \left\lvert F \right\rvert\right\rvert_p = \left\{\begin{array}{ll}
				\left(\frac{1}{2} \sum_{e_q \in \mathcal{E}_H} \left\lvert F\left(e_q\right)\right\rvert^p\right)^{\frac{1}{p}} & \quad 1 \leq p < \infty\\
				\max_{e_q \in \mathcal{E}_H} \left(\left\lvert F\left(e_q\right) \right\rvert\right) & \quad p = \infty
			\end{array}\right..
		\end{equation}\vspace{0em}
		
		\item[2)] Similarly, for a real hyperarc function $F \in \mathcal{H}\left(\mathcal{A}_H\right)$ at a hyperarc $a_q \in \mathcal{A}_H$, the absolute value is defined as:
		\begin{equation*}
			\left\lvert ~ \cdot ~ \right\rvert: ~ \mathbb{R} \longrightarrow \mathbb{R}_{\geq 0}
		\end{equation*}
		\begin{equation}
			F\left(a_q\right) \longmapsto \left\lvert F\left(a_q\right)\right\rvert = \left\{\begin{array}{ll}
				F\left(a_q\right) & \quad F\left(a_q\right) \geq 0\\
				-F\left(a_q\right) & \quad \text{otherwise}
			\end{array}\right..
		\end{equation}
		
		and the $\mathcal{L}^p$-norm of the hyperarc function $F$ is given by:
		\begin{equation*}
			\left\lvert \left\lvert ~ \cdot ~ \right\rvert\right\rvert_p: ~ \mathcal{H}\left(\mathcal{A}_H\right) \longrightarrow \mathbb{R}_{\geq 0}
		\end{equation*}
		\begin{equation}
			F \longmapsto \left\lvert \left\lvert F \right\rvert\right\rvert_p = \left\{\begin{array}{ll}
				\left(\frac{1}{2} \sum_{a_q \in \mathcal{A}_H} \left\lvert F\left(a_q\right)\right\rvert^p\right)^{\frac{1}{p}} & \quad 1 \leq p < \infty\\
				\max_{a_q \in \mathcal{A}_H} \left(\left\lvert F\left(a_q\right) \right\rvert\right) & \quad p = \infty
			\end{array}\right..
		\end{equation}\vspace{0em}
	\end{itemize}
\end{definition}
\clearpage
\section{Gradient and adjoint operators on oriented hypergraphs}\label{9} 

The definitions for the weighted gradient operators and the weighted adjoint operators for oriented normal graphs can be generalized to the setting of oriented hypergraphs, both for vertices and hyperarcs. Using gradients and adjoints on oriented hypergraphs then leads to feasible Laplacian and $p$-Laplacian definitions on oriented hypergraphs.\\

This section introduces gradient and adjoint operators on weighted oriented hypergraphs, which match the definitions for weighted oriented normal graphs from section (\ref{6}). Subsection (\ref{9.1}) defines the vertex gradient operator and the hyperarc gradient operator together with their specific properties. Furthermore, subsection (\ref{9.2}) derives the vertex adjoint operator and the hyperarc adjoint operator based on the previously defined vertex and hyperarc gradients. Both subsections also include proofs, which show that the hypergraph definitions for gradients and adjoints being applied to normal graphs results in the normal graph definitions from section (\ref{6}).\\

\subsection{Gradient operators on oriented hypergraphs}\label{9.1} 

As seen before, the gradient operator can be both defined for the set of all vertices $\mathcal{V}$ and on the set of all hyperarcs $\mathcal{A}_H$.\\

\begin{definition}[\textbf{Vertex gradient operator $\nabla_v$}]\label{Hnabla_v} \ \\
	For a weighted oriented hypergraph $OH = \left(\mathcal{V}, \mathcal{A}_H, w, W\right)$, the vertex gradient operator $\nabla_v$ is given by
	\begin{equation}
		\nabla_v: ~ \mathcal{H}\left(\mathcal{V}\right) \longrightarrow \mathcal{H}\left(\mathcal{A}_H\right) \qquad f \longmapsto \nabla_v f
	\end{equation}
	with the weighted difference of vertex function $f \in \mathcal{H}\left(\mathcal{V}\right)$ of a hyperarc $a_q \in \mathcal{A}_H$ being defined as
	\begin{equation*}
		\nabla_v f: ~ \mathcal{A}_H \longrightarrow \mathbb{R} \qquad a_q \longmapsto \nabla_v f \left(a_q\right) = 
	\end{equation*}
	\begin{equation}
		W_G \left(a_q\right)^\gamma \sum_{v_i \in \mathcal{V}} \left(\delta_{in}\left(v_i, a_q\right) \frac{w_I \left(v_i\right)^\alpha w_G \left(v_i\right)^\epsilon}{\left\lvert a_q^{in}\right\rvert} - \delta_{out}\left(v_i, a_q\right) \frac{w_I \left(v_i\right)^\alpha w_G \left(v_i\right)^\eta}{\left\lvert a_q^{out}\right\rvert}\right) f\left(v_i\right).
	\end{equation}\vspace{0em}
\end{definition}

The vertex gradient fulfills the two expected properties for a gradient, which are antisymmetry and the gradient of a constant function being equal to zero.\\

\clearpage
\begin{theorem}[\textbf{Vertex gradient operator properties}]\label{Hnabla_vprop} \ \\
	The weighted difference of a vertex function $f \in \mathcal{H}\left(\mathcal{V}\right)$ on a weighted oriented hypergraph $OH = \left(\mathcal{V}, \mathcal{A}_H, w, W\right)$ fulfills the following properties:
	\begin{itemize}
		\item[1)] Gradient of a constant vertex function is equal to zero:
		
		$f\left(v_i\right) \equiv \overline{f} \in \mathbb{R}$ for all vertices $v_i \in \mathcal{V}$ and $w_I \left(v_k\right)^\alpha w_G \left(v_k\right)^\epsilon = w_I \left(v_j\right)^\alpha w_G \left(v_j\right)^\eta$ for all vertex combinations $v_j, v_k \in \mathcal{V}$ with $v_j \in a_q^{out}$ and $v_k \in a_q^{in}$ for a hyperarc $a_q \in \mathcal{A}_H$
		
		$\Longrightarrow \quad \nabla_v f \left(a_q\right) = 0$ for all hyperarcs $a_q \in \mathcal{A}_H$
		
		\item[2)] Antisymmetry: 
		
		Symmetric hyperarc weight function $W_G$ and $\epsilon = \eta$
		
		$\Longrightarrow$ \quad $\nabla_v f \left(a_q\right) = - \nabla_v f \left(\widetilde{a_q}\right)$ for all hyperarcs $a_q \in \mathcal{A}_H$ \big(with the corresponding hyperarc $\widetilde{a_q} \in \widetilde{\mathcal{A}_H}$ of the hypergraph with switched orientation $\widetilde{OH}$ from definition (\ref{switchedOH})\big)\\
	\end{itemize}
\end{theorem}

\begin{proof}\ \\
	Given a real vertex function $f \in \mathcal{H}\left(\mathcal{V}\right)$ on a weighted oriented hypergraph $OH = \left(\mathcal{V}, \mathcal{A}_H, w, W\right)$, then it holds true that:
	\begin{itemize}
		\item[1)] The property $w_I \left(v_k\right)^\alpha w_G \left(v_k\right)^\epsilon = w_I \left(v_j\right)^\alpha w_G \left(v_j\right)^\eta$ for all vertex combinations $v_j, v_k \in \mathcal{V}$ with $v_j \in a_q^{out}$ and $v_k \in a_q^{in}$ for a hyperarc $a_q \in \mathcal{A}_H$ implies that there exists a constant $w_{a_q} \in \mathbb{R}_{> 0}$ such that:
		\begin{equation*}
			w_I \left(v_k\right)^\alpha w_G \left(v_k\right)^\epsilon = \left(v_j\right)^\alpha w_G \left(v_j\right)^\eta =: w_{a_q}
		\end{equation*}
		
		Thus, together with the additional property $f\left(v_i\right) \equiv \overline{f} \in \mathbb{R}$ for all vertices $v_i \in \mathcal{V}$ this results in:
		
		$\begin{aligned}[t]
			\nabla_v f \left(a_q\right) = & ~ W_G \left(a_q\right)^\gamma \sum_{v_i \in \mathcal{V}} \left(\delta_{in}\left(v_i, a_q\right) \frac{w_I \left(v_i\right)^\alpha w_G \left(v_i\right)^\epsilon}{\left\lvert a_q^{in}\right\rvert}\right. - &\\
			& \left.\delta_{out}\left(v_i, a_q\right) \frac{w_I \left(v_i\right)^\alpha w_G \left(v_i\right)^\eta}{\left\lvert a_q^{out}\right\rvert}\right) f\left(v_i\right) &\\
			= & ~ W_G \left(a_q\right)^\gamma \sum_{v_i \in \mathcal{V}} \left(\delta_{in}\left(v_i, a_q\right) \frac{w_{a_q}}{\left\lvert a_q^{in}\right\rvert} - \delta_{out}\left(v_i, a_q\right) \frac{w_{a_q}}{\left\lvert a_q^{out}\right\rvert}\right) \overline{f} &\\
		\end{aligned}$\\
		$\begin{aligned}[t]
			= & ~ W_G \left(a_q\right)^\gamma \left(\sum_{v_i \in \mathcal{V}} \delta_{in}\left(v_i, a_q\right) \frac{w_{a_q}}{\left\lvert a_q^{in}\right\rvert} - \sum_{v_j \in \mathcal{V}} \delta_{out}\left(v_j, a_q\right) \frac{w_{a_q}}{\left\lvert a_q^{out}\right\rvert}\right) \overline{f} &\\
			= & ~ W_G \left(a_q\right)^\gamma \left(\frac{w_{a_q}}{\left\lvert a_q^{in}\right\rvert} \sum_{v_i \in \mathcal{V}} \delta_{in}\left(v_i, a_q\right) - \frac{w_{a_q}}{\left\lvert a_q^{out}\right\rvert} \sum_{v_j \in \mathcal{V}} \delta_{out}\left(v_j, a_q\right)\right) \overline{f} &\\
			= & ~ W_G \left(a_q\right)^\gamma \left(\frac{w_{a_q}}{\left\lvert a_q^{in}\right\rvert} \left\lvert a_q^{in}\right\rvert - \frac{w_{a_q}}{\left\lvert a_q^{out}\right\rvert} \left\lvert a_q^{out}\right\rvert\right) \overline{f} &\\
			= & ~ W_G \left(a_q\right)^\gamma \cdot 0 \cdot \overline{f} &\\
			= & ~ 0 &\\
		\end{aligned}$\\
		
		Where the last equality is feasible due to the hyperarc weight function $W_G$ being a real function and therefore mapping to finite values.\\
		
		\item[2)] $W_G$ being a symmetric hyperarc weight function together with $\epsilon = \eta$ implies that for every hyperarc $a_q \in \mathcal{A}_H$ it holds true that:
		
		$\begin{aligned}[t]
			\nabla_v f \left(a_q\right) = & ~ W_G \left(a_q\right)^\gamma \sum_{v_i \in \mathcal{V}} \left(\delta_{in}\left(v_i, a_q\right) \frac{w_I \left(v_i\right)^\alpha w_G \left(v_i\right)^\epsilon}{\left\lvert a_q^{in}\right\rvert}\right. - &\\
			& \left.\delta_{out}\left(v_i, a_q\right) \frac{w_I \left(v_i\right)^\alpha w_G \left(v_i\right)^\eta}{\left\lvert a_q^{out}\right\rvert}\right) f\left(v_i\right) &\\
		\end{aligned}$\\
		$\begin{aligned}[t]
			& = - W_G \left(a_q\right)^\gamma \sum_{v_i \in \mathcal{V}} \left(\delta_{out}\left(v_i, a_q\right) \frac{w_I \left(v_i\right)^\alpha w_G \left(v_i\right)^\eta}{\left\lvert a_q^{out}\right\rvert} - \delta_{in}\left(v_i, a_q\right) \frac{w_I \left(v_i\right)^\alpha w_G \left(v_i\right)^\epsilon}{\left\lvert a_q^{in}\right\rvert}\right) f\left(v_i\right) &\\
			& = - W_G \left(a_q\right)^\gamma \sum_{v_i \in \mathcal{V}} \left(\delta_{out}\left(v_i, a_q\right) \frac{w_I \left(v_i\right)^\alpha w_G \left(v_i\right)^\epsilon}{\left\lvert a_q^{out}\right\rvert} - \delta_{in}\left(v_i, a_q\right) \frac{w_I \left(v_i\right)^\alpha w_G \left(v_i\right)^\eta}{\left\lvert a_q^{in}\right\rvert}\right) f\left(v_i\right) &\\
			& = - W_G \left(\widetilde{a_q}\right)^\gamma \sum_{v_i \in \mathcal{V}} \left(\delta_{out}\left(v_i, a_q\right) \frac{w_I \left(v_i\right)^\alpha w_G \left(v_i\right)^\epsilon}{\left\lvert a_q^{out}\right\rvert} - \delta_{in}\left(v_i, a_q\right) \frac{w_I \left(v_i\right)^\alpha w_G \left(v_i\right)^\eta}{\left\lvert a_q^{in}\right\rvert}\right) f\left(v_i\right) &\\
			& = - W_G \left(\widetilde{a_q}\right)^\gamma \sum_{v_i \in \mathcal{V}} \left(\delta_{in}\left(v_i, \widetilde{a_q}\right) \frac{w_I \left(v_i\right)^\alpha w_G \left(v_i\right)^\epsilon}{\left\lvert \widetilde{a_q^{in}}\right\rvert} - \delta_{out}\left(v_i, \widetilde{a_q}\right) \frac{w_I \left(v_i\right)^\alpha w_G \left(v_i\right)^\eta}{\left\lvert \widetilde{a_q^{out}}\right\rvert}\right) f\left(v_i\right) &\\
			& = - \nabla_v f \left(\widetilde{a_q}\right) &\\
		\end{aligned}$
		
		\clearpage
		Where the second to last equality results from $\delta_{out}\left(v_i, a_q\right) = \delta_{in}\left(v_i, \widetilde{a_q}\right)$,\\ $\delta_{in}\left(v_i, a_q\right) = \delta_{out}\left(v_i, \widetilde{a_q}\right)$, $\left\lvert a_q^{out}\right\rvert = \left\lvert \widetilde{a_q^{in}}\right\rvert$ and $\left\lvert a_q^{in}\right\rvert = \left\lvert \widetilde{a_q^{out}}\right\rvert$ since the output and input vertices of hyperarc $a_q \in \mathcal{A}_G$ are exchanged to get the hyperarc with switched orientation $\widetilde{a_q} \in \widetilde{\mathcal{A}_G}$.\\
	\end{itemize}
\end{proof}

Theorem (\ref{GaH}) shows that every oriented normal graph is a special case of an oriented hypergraph. Hence, the vertex gradient operator for hypergraphs can also be applied to normal graphs, resulting in a generalization of the vertex gradient operator for normal graphs.\\

\begin{theorem}[\textbf{Vertex gradient operator for hypergraphs a generalization of the normal graph case}]\label{nabla_v} \ \\
	Given a weighted oriented normal graph $OG = \left(\mathcal{V}, \mathcal{A}_G, w, W\right)$ and a vertex function $f \in \mathcal{H}\left(\mathcal{V}\right)$, then the vertex gradient operator $\nabla_v$ for weighted oriented hypergraphs matches the vertex gradient operator for normal graphs from definition (\ref{Gnabla_v}).\\
\end{theorem}

\begin{proof}\ \\
	Given a weighted oriented normal graph $OG = \left(\mathcal{V}, \mathcal{A}_G, w, W\right)$ together with a vertex function $f \in \mathcal{H}\left(\mathcal{V}\right)$, then applying the vertex gradient operator $\nabla_v$ for weighted oriented hypergraphs to any arc $a_q = \left(v_i, v_j\right) \in \mathcal{A}_G$ leads to:
	
	$\begin{aligned}[t]
		\nabla_v f \left(a_q\right) = & ~ \nabla_v f \left(v_i, v_j\right) &\\
		= & ~ W_G \left(v_i, v_j\right)^\gamma &\\
		& \sum_{v_k \in \mathcal{V}} \left(\delta_{in}\left(v_k, a_q\right) \frac{w_I \left(v_k\right)^\alpha w_G \left(v_k\right)^\epsilon}{\left\lvert a_q^{in}\right\rvert} - \delta_{out}\left(v_k, a_q\right) \frac{w_I \left(v_k\right)^\alpha w_G \left(v_k\right)^\eta}{\left\lvert a_q^{out}\right\rvert}\right) f\left(v_k\right) &\\
		= & ~ W_G \left(v_i, v_j\right)^\gamma &\\
		& \sum_{v_k \in \left\{v_i, v_j\right\}} \left(\delta_{in}\left(v_k, a_q\right) w_I \left(v_k\right)^\alpha w_G \left(v_k\right)^\epsilon - \delta_{out}\left(v_k, a_q\right) w_I \left(v_k\right)^\alpha w_G \left(v_k\right)^\eta\right) f\left(v_k\right) &\\
		= & ~ W_G \left(v_i, v_j\right)^\gamma \left(w_I \left(v_j\right)^\alpha w_G \left(v_j\right)^\epsilon f\left(v_j\right) - w_I \left(v_i\right)^\alpha w_G \left(v_i\right)^\eta f\left(v_i\right)\right) &\\
	\end{aligned}$\\

	These reformulations are based on the properties $a_q^{out} = \left\{v_i\right\}$ and $a_q^{in} = \left\{v_j\right\}$ for any arc $a_q = \left(v_i, v_j\right) \in \mathcal{A}_G$, which implies $\left\lvert a_q^{in} \right\rvert = \left\lvert a_q^{out} \right\rvert = 1$ and $\delta_{in}\left(v_j, a_q\right) = \delta_{out}\left(v_i, a_q\right) = 1$.\\
	
	Therefore, the vertex gradient operator $\nabla_v$ for hypergraphs is a valid generalization of the normal graph case.\\
\end{proof}

Similar to the vertex gradient operator $\nabla_v$, a linear gradient operator can be defined on the domain of hyperarc functions $\mathcal{H}\left(\mathcal{A}_H\right)$.\\

\begin{definition}[\textbf{Hyperarc gradient operator $\nabla_a$}]\label{Hnabla_a} \ \\
	For a weighted oriented hypergraph $OH = \left(\mathcal{V}, \mathcal{A}_H, w, W\right)$, the hyperarc gradient operator $\nabla_a$ is defined as
	\begin{equation}
		\nabla_a: ~ \mathcal{H}\left(\mathcal{A}_H\right) \longrightarrow \mathcal{H}\left(\mathcal{V}\right) \qquad F \longmapsto \nabla_a F
	\end{equation}
	with the weighted difference of hyperarc function $F \in \mathcal{H}\left(\mathcal{A}_H\right)$ at a vertex $v_i \in \mathcal{V}$ being being given by
	\begin{equation*}
		\nabla_a F: ~ \mathcal{V} \longrightarrow \mathbb{R} \qquad v_i \longmapsto \nabla_a F \left(v_i\right) =
	\end{equation*}
	\begin{equation}
		w_G \left(v_i\right)^\zeta \sum_{a_q \in \mathcal{A}_H} \left(\frac{\delta_{in}\left(v_i, a_q\right)}{\deg_{in}\left(v_i\right)} - \frac{\delta_{out}\left(v_i, a_q\right)}{\deg_{out}\left(v_i\right)}\right) W_I \left(a_q\right)^\beta W_G \left(a_q\right)^\theta F\left(a_q\right).
	\end{equation}\vspace{0em}
\end{definition}

The hyperarc gradient also fulfills the expected antisymmetry property and the gradient of a constant function being equal to zero.\\

\begin{theorem}[\textbf{Hyperarc gradient operator properties}]\label{Hnabla_aprop} \ \\
	The weighted difference of a hyperarc function $F \in \mathcal{H}\left(\mathcal{A}_H\right)$ on a weighted oriented hypergraph $OH = \left(\mathcal{V}, \mathcal{A}_H, w, W\right)$ fulfills the following properties:
	\begin{itemize}
		\item[1)] Gradient of a constant hyperarc function is equal to zero:
		
		$F\left(a_q\right) \equiv \overline{F} \in \mathbb{R}$ and $W_I \left(a_q\right)^\beta W_G \left(a_q\right)^\theta \equiv \overline{W} \in \mathbb{R}$ for all arcs $a_q \in \mathcal{A}_G$
		
		$\Longrightarrow$ \quad $\nabla_a F \left(v_i\right) = 0$ for all vertices $v_i \in \mathcal{V}$
		
		\item[2)] Antisymmetry:
		
		Symmetric hyperarc weight functions $W_I$ and $W_G$ and $F\left(a_q\right) = F\left(\widetilde{a_q}\right)$ for all hyperarcs $a_q \in \mathcal{A}_H$, $\widetilde{a_q} \in \widetilde{\mathcal{A}_H}$
		
		$\Longrightarrow$ \quad $\nabla_a F \left(v_i\right) = - \nabla_{\widetilde{a}} F \left(v_i\right)$ for all vertices $v_i \in \mathcal{V}$ \big(with $\nabla_{\widetilde{a}} F$ being the hyperarc gradient applied to the hypergraph with switched orientation $\widetilde{OH}$\big)\\
	\end{itemize}
\end{theorem}

\begin{proof}\ \\
	Given a real hyperarc function $F \in \mathcal{H}\left(\mathcal{A}_H\right)$ on a weighted oriented hypergraph $OH = \left(\mathcal{V}, \mathcal{A}_H, w, W\right)$, then it holds true that:
	\begin{itemize}
		\item[1)] $F\left(a_q\right) \equiv \overline{F} \in \mathbb{R}$ and $W_I \left(a_q\right)^\beta W_G \left(a_q\right)^\theta \equiv \overline{W} \in \mathbb{R}$ for all hyperarcs $a_q \in \mathcal{A}_H$ implies for all vertices $v_i \in \mathcal{V}$ that:
		
		$\begin{aligned}[t]
			\nabla_a F \left(v_i\right) = & ~ w_G \left(v_i\right)^\zeta \sum_{a_q \in \mathcal{A}_H} \left(\frac{\delta_{in}\left(v_i, a_q\right)}{\deg_{in}\left(v_i\right)} - \frac{\delta_{out}\left(v_i, a_q\right)}{\deg_{out}\left(v_i\right)}\right) W_I \left(a_q\right)^\beta W_G \left(a_q\right)^\theta F\left(a_q\right) &\\
			= & ~ w_G \left(v_i\right)^\zeta \sum_{a_q \in \mathcal{A}_H} \left(\frac{\delta_{in}\left(v_i, a_q\right)}{\deg_{in}\left(v_i\right)} - \frac{\delta_{out}\left(v_i, a_q\right)}{\deg_{out}\left(v_i\right)}\right) \overline{W} ~ \overline{F} &\\
			= & ~ w_G \left(v_i\right)^\zeta \left(\sum_{a_q \in \mathcal{A}_H} \frac{\delta_{in}\left(v_i, a_q\right)}{\deg_{in}\left(v_i\right)} - \sum_{a_r \in \mathcal{A}_H} \frac{\delta_{out}\left(v_i, a_r\right)}{\deg_{out}\left(v_i\right)}\right) \overline{W} ~ \overline{F} &\\
			= & ~ w\left(v_i\right)^\zeta &\\
			& \left( \frac{1}{\deg_{in}\left(v_i\right)} \sum_{a_q \in \mathcal{A}_H} \delta_{in}\left(v_i, a_q\right) - \frac{1}{\deg_{out}\left(v_i\right)} \sum_{a_r \in \mathcal{A}_H} \delta_{out}\left(v_i, a_r\right)\right) \overline{W} ~ \overline{F} &\\
			= & ~ w_G \left(v_i\right)^\zeta \left( \frac{1}{\deg_{in}\left(v_i\right)} \deg_{in}\left(v_i\right) - \frac{1}{\deg_{out}\left(v_i\right)} \deg_{out}\left(v_i\right)\right) \overline{W} ~ \overline{F} &\\
			= & ~ w_G \left(v_i\right)^\zeta \cdot 0 \cdot \overline{W} ~ \overline{F} &\\
			= & ~ 0 &\\
		\end{aligned}$
		
		Where the last equality is feasible due to the vertex weight function $w_G$ being a real function and thus mapping to finite values.\\
		
		\item[2)] $W_I$ and $W_G$ being symmetric hyperarc weight functions and $F\left(a_q\right) = F\left(\widetilde{a_q}\right)$ for all hyperarcs $a_q \in \mathcal{A}_H$, $\widetilde{a_q} \in \widetilde{\mathcal{A}_H}$ implies for all vertices $v_i \in \mathcal{V}$ that:
		
		$\begin{aligned}[t]
			\nabla_a F \left(v_i\right) & = w_G \left(v_i\right)^\zeta \sum_{a_q \in \mathcal{A}_H} \left(\frac{\delta_{in}\left(v_i, a_q\right)}{\deg_{in}\left(v_i\right)} - \frac{\delta_{out}\left(v_i, a_q\right)}{\deg_{out}\left(v_i\right)}\right)  W_I \left(a_q\right)^\beta W_G \left(a_q\right)^\theta F\left(a_q\right) &\\
			& = - w_G \left(v_i\right)^\zeta \sum_{a_q \in \mathcal{A}_H} \left(\frac{\delta_{out}\left(v_i, a_q\right)}{\deg_{out}\left(v_i\right)} - \frac{\delta_{in}\left(v_i, a_q\right)}{\deg_{in}\left(v_i\right)}\right) W_I \left(a_q\right)^\beta W_G \left(a_q\right)^\theta F\left(a_q\right) &\\
			& = - w_G \left(v_i\right)^\zeta \sum_{a_q \in \mathcal{A}_H} \left(\frac{\delta_{out}\left(v_i, a_q\right)}{\deg_{out}\left(v_i\right)} - \frac{\delta_{in}\left(v_i, a_q\right)}{\deg_{in}\left(v_i\right)}\right) W_I \left(\widetilde{a_q}\right)^\beta W_G \left(\widetilde{a_q}\right)^\theta F\left(\widetilde{a_q}\right) &\\
			& = - w_G \left(v_i\right)^\zeta \sum_{\widetilde{a_q} \in \widetilde{\mathcal{A}_H}} \left(\frac{\delta_{in}\left(v_i, \widetilde{a_q}\right)}{\widetilde{\deg_{in}}\left(v_i\right)} - \frac{\delta_{out}\left(v_i, \widetilde{a_q}\right)}{\widetilde{\deg_{out}}\left(v_i\right)}\right) W_I \left(\widetilde{a_q}\right)^\beta W_G \left(\widetilde{a_q}\right)^\theta F\left(\widetilde{a_q}\right) &\\
			& = - \nabla_{\widetilde{a}} F \left(v_i\right)\\
		\end{aligned}$\\
		
		Where the equality from the second to the third line is based on the symmetry of the hyperarc weight functions $W_I$ and $W_G$ and the reformulation from the third to the fourth line results from the fact that $\delta_{out}\left(v_i, a_q\right) = \delta_{in}\left(v_i, \widetilde{a_q}\right)$ and $\delta_{in}\left(v_i, a_q\right) = \delta_{out}\left(v_i, \widetilde{a_q}\right)$ for all vertices $v_i \in \mathcal{V}$. Moreover, it also holds true that $\deg_{out}\left(v_i\right)$ becomes $\widetilde{\deg_{in}}\left(v_i\right)$ and $\deg_{in}\left(v_i\right)$ becomes $\widetilde{\deg_{out}}\left(v_i\right)$ when considering the hypergraph with switched orientation.\\
	\end{itemize}
\end{proof}

The hyperarc gradient operator for oriented hypergraphs can also be applied to arcs of an oriented normal graph. The following theorem shows that the hyperarc gradient operator then results in a generalization of the arc gradient operator.\\

\begin{theorem}[\textbf{Hyperarc gradient operator a generalization of the arc gradient operator}]\label{nabla_a} \ \\
	Given a weighted oriented normal graph $OG = \left(\mathcal{V}, \mathcal{A}_G, w, W\right)$ and an arc function $F \in \mathcal{H}\left(\mathcal{A}_G\right)$, then the hyperarc gradient operator $\nabla_a$ for weighted oriented hypergraphs matches the arc gradient operator for normal graphs from definition (\ref{Gnabla_a}).\\
\end{theorem}

\begin{proof}\ \\
	Given a weighted oriented normal graph $OG = \left(\mathcal{V}, \mathcal{A}_G, w, W\right)$ together with an arc function $F \in \mathcal{H}\left(\mathcal{A}_G\right)$, then applying the hyperarc gradient operator $\nabla_a$ for weighted oriented hypergraphs to any vertex $v_i \in \mathcal{V}$ leads to the definition of the arc gradient operator:
	\begin{equation*}
		\nabla_a F \left(v_i\right) = w_G \left(v_i\right)^\zeta \sum_{a_q \in \mathcal{A}_G} \left(\frac{\delta_{in}\left(v_i, a_q\right)}{\deg_{in}\left(v_i\right)} - \frac{\delta_{out}\left(v_i, a_q\right)}{\deg_{out}\left(v_i\right)}\right) W_I \left(a_q\right)^\beta W_G \left(a_q\right)^\theta F\left(a_q\right).
	\end{equation*}
	Hence, the hyperarc gradient operator $\nabla_a$ for hypergraphs is a valid generalization of the arc gradient operator for normal graphs.\\
\end{proof}
\subsection{Adjoint operators on oriented hypergraphs}\label{9.2} 

On weighted oriented hypergraphs $OH = \left(\mathcal{V}, \mathcal{A}_H, w, W\right)$, the weighted adjoint operators $\nabla^*_v$ and $\nabla^*_a$ are defined as the dual counterpart to the weighted gradient operators $\nabla_v$ and $\nabla_a$ fulfilling again the following equalities for all vertex functions $f \in \mathcal{H}\left(\mathcal{V}\right)$ and all hyperarc functions $G \in \mathcal{H}\left(\mathcal{A}_H\right)$:\\
\begin{equation}
	{\langle G, \nabla_v f \rangle}_{\mathcal{H}\left(\mathcal{A}_H\right)} = {\langle f, \nabla^*_v G \rangle}_{\mathcal{H}\left(\mathcal{V}\right)}
\end{equation}
\begin{equation}
	{\langle f, \nabla_a G \rangle}_{\mathcal{H}\left(\mathcal{V}\right)} = {\langle G, \nabla^*_a f \rangle}_{\mathcal{H}\left(\mathcal{A}_H\right)}
\end{equation}

The vertex adjoint and the hyperarc adjoint are later on used to introduce Laplacian and $p$-Laplacian operators for vertices and hyperarcs.\\

\clearpage
\begin{definition}[\textbf{Vertex adjoint operator $\nabla^*_v$}]\label{Hnabla^*_v} \ \\
	For a weighted oriented hypergraph $OH = \left(\mathcal{V}, \mathcal{A}_H, w, W\right)$ with symmetric hyperarc weight functions $W_I$ and $W_G$, the vertex adjoint operator $\nabla^*_v$ is given by
	\begin{equation}
		\nabla^*_v: ~ \mathcal{H}\left(\mathcal{A}_H\right) \longrightarrow \mathcal{H}\left(\mathcal{V}\right) \qquad F \longmapsto \nabla^*_v F
	\end{equation}
	with the weighted adjoint of hyperarc function $F \in \mathcal{H}\left(\mathcal{A}_H\right)$ at a vertex $v_i \in \mathcal{V}$ being defined as
	\begin{equation*}
		\nabla^*_v F: ~ \mathcal{V} \longrightarrow \mathbb{R} \qquad v_i \longmapsto \nabla^*_v F \left(v_i\right) =
	\end{equation*}
	\begin{equation}
		\sum_{a_q \in \mathcal{A}_H} \left(\delta_{in}\left(v_i, a_q\right) \frac{w_G \left(v_i\right)^\epsilon}{\left\lvert a_q^{in}\right\rvert} - \delta_{out}\left(v_i, a_q\right) \frac{w_G \left(v_i\right)^\eta}{\left\lvert a_q^{out}\right\rvert}\right) W_I \left(a_q\right)^\beta W_G \left(a_q\right)^\gamma F\left(a_q\right).
	\end{equation}\vspace{0em}
\end{definition}

Applying the vertex adjoint operator for hypergraphs to arc functions on a normal graph results in the vertex adjoint operator for normal graphs.\\

\begin{theorem}[\textbf{Vertex adjoint operator for hypergraphs a generalization of the normal graph case}]\label{nabla^*_v} \ \\
	Given a weighted oriented normal graph $OG = \left(\mathcal{V}, \mathcal{A}_G, w, W\right)$ with symmetric arc weight functions $W_I$ and $W_G$ and an arc function $F \in \mathcal{H}\left(\mathcal{A}_G\right)$, then the vertex adjoint operator $\nabla_v^*$ for weighted oriented hypergraphs matches the vertex adjoint operator for normal graphs from definition (\ref{Gnabla^*_v}).\\
\end{theorem}

\begin{proof}\ \\
	Given a weighted oriented normal graph $OG = \left(\mathcal{V}, \mathcal{A}_G, w, W\right)$ with symmetric arc weight functions $W_I$ and $W_G$ together with an arc function $F \in \mathcal{H}\left(\mathcal{A}_G\right)$, then applying the vertex adjoint operator $\nabla_v^*$ for weighted oriented hypergraphs to any vertex $v_i \in \mathcal{V}$ leads to:
	
	$\begin{aligned}[t]
		\nabla^*_v F \left(v_i\right) = & \sum_{a_q \in \mathcal{A}_G} \left(\delta_{in}\left(v_i, a_q\right) \frac{w_G \left(v_i\right)^\epsilon}{\left\lvert a_q^{in}\right\rvert} - \delta_{out}\left(v_i, a_q\right) \frac{w_G \left(v_i\right)^\eta}{\left\lvert a_q^{out}\right\rvert}\right) W_I \left(a_q\right)^\beta W_G \left(a_q\right)^\gamma F\left(a_q\right) &\\
		= & \sum_{a_q \in \mathcal{A}_G} \left(\delta_{in}\left(v_i, a_q\right) w_G \left(v_i\right)^\epsilon - \delta_{out}\left(v_i, a_q\right) w_G \left(v_i\right)^\eta\right) W_I \left(a_q\right)^\beta W_G \left(a_q\right)^\gamma F\left(a_q\right) &\\
		= & \sum_{a_q \in \mathcal{A}_G} \delta_{in}\left(v_i, a_q\right) w_G \left(v_i\right)^\epsilon W_I \left(a_q\right)^\beta W_G \left(a_q\right)^\gamma F\left(a_q\right) - &\\
		& ~ \delta_{out}\left(v_i, a_q\right) w_G \left(v_i\right)^\eta W_I \left(a_q\right)^\beta W_G \left(a_q\right)^\gamma F\left(a_q\right) &\\
	\end{aligned}$\\
	$\begin{aligned}[t]
		= & \sum_{v_j \in \mathcal{V}} \delta\left(v_j, v_i\right) w_G \left(v_i\right)^\epsilon W_I \left(v_j, v_i\right)^\beta W_G \left(v_j, v_i\right)^\gamma F\left(v_j, v_i\right) - &\\
		& ~ \delta\left(v_i, v_j\right) w_G \left(v_i\right)^\eta W_I \left(v_i, v_j\right)^\beta W_G \left(v_i, v_j\right)^\gamma F\left(v_i, v_j\right) &\\
		= & \sum_{v_j \in \mathcal{V}} \left(\delta\left(v_j, v_i\right) w_G \left(v_i\right)^\epsilon F\left(v_j, v_i\right) - \delta\left(v_i, v_j\right) w_G \left(v_i\right)^\eta F\left(v_i, v_j\right)\right) &\\
		& W_I \left(v_i, v_j\right)^\beta W_G \left(v_i, v_j\right)^\gamma &\\
	\end{aligned}$\\
	
	These reformulations are based on the properties $a_q^{out} = \left\{v_j\right\}$ and $a_q^{in} = \left\{v_k\right\}$ for any arc $a_q = \left(v_j, v_k\right) \in \mathcal{A}_G$, which implies $\left\lvert a_q^{in} \right\rvert = \left\lvert a_q^{out} \right\rvert = 1$. Furthermore, the fact that summing over all arcs $a_q \in \mathcal{A}_G$ and only considering arcs $a_q$ with $\delta_{in}\left(v_i, a_q\right) = 1$ or $ \delta_{out}\left(v_i, a_q\right) = 1$ is the same as summing over all vertices $v_j \in \mathcal{V}$ and only considering vertices $v_j$ with $\delta\left(v_j, v_i\right) = 1$ or $\delta\left(v_i, v_j\right) = 1$, is utilized as well. Finally, the last equality holds true because of the the symmetry of the arc weight functions $W_I$ and $W_G$.\\
	
	Thus, the vertex adjoint operator $\nabla_v^*$ for hypergraphs is a valid generalization of the normal graph case.\\
\end{proof}

The theorem below proves that the vertex adjoint is the correct dual counterpart to the definition of the vertex gradient on hypergraphs.\\

\begin{theorem}[\textbf{Connection vertex gradient $\nabla_v$ and vertex adjoint $\nabla^*_v$}]\label{Hvertexgradadj} \ \\
	For a weighted oriented hypergraph $OH = \left(\mathcal{V}, \mathcal{A}_H, w, W\right)$ with symmetric hyperarc weight functions $W_I$ and $W_G$, the vertex adjoint operator $\nabla^*_v$ fulfills the equality
	\begin{equation}
		{\langle G, \nabla_v f \rangle}_{\mathcal{H}\left(\mathcal{A}_H\right)} = {\langle f, \nabla^*_v G \rangle}_{\mathcal{H}\left(\mathcal{V}\right)}
	\end{equation}
	for all vertex functions $f \in \mathcal{H}\left(\mathcal{V}\right)$ and all hyperarc functions $G \in \mathcal{H}\left(\mathcal{A}_H\right)$.\\
\end{theorem}

\begin{proof}\ \\
	Given a weighted oriented hypergraph $OH = \left(\mathcal{V}, \mathcal{A}_H, w, W\right)$ with symmetric hyperarc weight functions $W_I$ and $W_G$, a vertex function $f \in \mathcal{H}\left(\mathcal{V}\right)$ and a hyperarc function $G \in \mathcal{H}\left(\mathcal{A}_H\right)$, using the definitions of the inner product of $\mathcal{H}\left(\mathcal{A}_H\right)$ and the vertex gradient operator $\nabla_v$ yields the following:
	
	\clearpage
	$\begin{aligned}[t]
		{\langle G, \nabla_v f \rangle}_{\mathcal{H}\left(\mathcal{A}_H\right)} = & \sum_{a_q \in \mathcal{A}_H} W_I \left(a_q\right)^\beta G\left(a_q\right) \nabla_v f \left(a_q\right) &\\
	\end{aligned}$\\
	$\begin{aligned}[t]
		= & \sum_{a_q \in \mathcal{A}_H} W_I \left(a_q\right)^\beta G\left(a_q\right) W_G \left(a_q\right)^\gamma \\
		& \sum_{v_i \in \mathcal{V}} \left(\delta_{in}\left(v_i, a_q\right) \frac{w_I \left(v_i\right)^\alpha w_G \left(v_i\right)^\epsilon}{\left\lvert a_q^{in}\right\rvert} - \delta_{out}\left(v_i, a_q\right) \frac{w_I \left(v_i\right)^\alpha w_G \left(v_i\right)^\eta}{\left\lvert a_q^{out}\right\rvert}\right) f\left(v_i\right)&\\
		= & \sum_{a_q \in \mathcal{A}_H} \sum_{v_i \in \mathcal{V}} W_I \left(a_q\right)^\beta W_G \left(a_q\right)^\gamma G\left(a_q\right) &\\
		& \left(\delta_{in}\left(v_i, a_q\right) \frac{w_I \left(v_i\right)^\alpha w_G \left(v_i\right)^\epsilon}{\left\lvert a_q^{in}\right\rvert} - \delta_{out}\left(v_i, a_q\right) \frac{w_I \left(v_i\right)^\alpha w_G \left(v_i\right)^\eta}{\left\lvert a_q^{out}\right\rvert}\right) f\left(v_i\right)&\\
	\end{aligned}$\\

	Since both, the set of all hyperarcs $\mathcal{A}_H$ and the set of all vertices $\mathcal{V}$, are assumed to be finite, exchanging the sums is feasible:
	
	$\begin{aligned}[t]
		= & \sum_{v_i \in \mathcal{V}} \sum_{a_q \in \mathcal{A}_H} W_I \left(a_q\right)^\beta W_G \left(a_q\right)^\gamma G\left(a_q\right) &\\
		& \left(\delta_{in}\left(v_i, a_q\right) \frac{w_I \left(v_i\right)^\alpha w_G \left(v_i\right)^\epsilon}{\left\lvert a_q^{in}\right\rvert} - \delta_{out}\left(v_i, a_q\right) \frac{w_I \left(v_i\right)^\alpha w_G \left(v_i\right)^\eta}{\left\lvert a_q^{out}\right\rvert}\right) f\left(v_i\right)&\\
		= & \sum_{v_i \in \mathcal{V}} f\left(v_i\right) \sum_{a_q \in \mathcal{A}_H} \left(\delta_{in}\left(v_i, a_q\right) \frac{w\left(v_i\right)^{\alpha + \epsilon}}{\left\lvert a_q^{in}\right\rvert} - \delta_{out}\left(v_i, a_q\right) \frac{w\left(v_i\right)^{\alpha + \eta}}{\left\lvert a_q^{out}\right\rvert}\right) &\\
		& W_I \left(a_q\right)^\beta W_G \left(a_q\right)^\gamma G\left(a_q\right) &\\
	\end{aligned}$\\

	Using the definition of the vertex adjoint operator $\nabla^*_v$ and the definition of the inner product of $\mathcal{H}\left(\mathcal{V}\right)$ results in:
	
	$\begin{aligned}[t]
		= & \sum_{v_i \in \mathcal{V}} w_I \left(v_i\right)^\alpha f\left(v_i\right) &\\
		& \sum_{a_q \in \mathcal{A}_H} \left(\delta_{in}\left(v_i, a_q\right) \frac{w_G \left(v_i\right)^\epsilon}{\left\lvert a_q^{in}\right\rvert} - \delta_{out}\left(v_i, a_q\right) \frac{w_G \left(v_i\right)^\eta}{\left\lvert a_q^{out}\right\rvert}\right) W_I \left(a_q\right)^\beta W_G \left(a_q\right)^\gamma G\left(a_q\right) &\\
		= & \sum_{v_i \in \mathcal{V}} w_I \left(v_i\right)^\alpha f\left(v_i\right) \nabla_v^* G\left(a_q\right) = {\langle f, \nabla^*_v G \rangle}_{\mathcal{H}\left(\mathcal{V}\right)} &\\
	\end{aligned}$\\

	Therefore, with the previous definition for the vertex adjoint operator $\nabla^*_v$, the equality ${\langle G, \nabla_v f \rangle}_{\mathcal{H}\left(\mathcal{A}_H\right)} = {\langle f, \nabla^*_v G \rangle}_{\mathcal{H}\left(\mathcal{V}\right)}$ holds true for all vertex functions $f \in \mathcal{H}\left(\mathcal{V}\right)$ and all hyperarc functions $G \in \mathcal{H}\left(\mathcal{A}_H\right)$.\\
\end{proof}

Similarly to the vertex adjoint operator $\nabla^*_v$ for hyperarc functions $F \in \mathcal{H}\left(\mathcal{A}_H\right)$, the hyperarc ajoint operator $\nabla^*_a$ for vertex functions $f \in \mathcal{H}\left(\mathcal{V}\right)$ can be defined on the basis of the hyperarc gradient operator $\nabla_a$.\\

\begin{definition}[\textbf{Hyperarc adjoint operator $\nabla^*_a$}]\label{Hnabla^*_a} \ \\
	For a weighted oriented hypergraph $OH = \left(\mathcal{V}, \mathcal{A}_H, w, W\right)$, the hyperarc adjoint operator $\nabla^*_a$ is defined as
	\begin{equation}
		\nabla^*_a: ~ \mathcal{H}\left(\mathcal{V}\right) \longrightarrow \mathcal{H}\left(\mathcal{A}_H\right) \qquad f \longmapsto \nabla^*_a f
	\end{equation}
	with the weighted adjoint of vertex function $f \in \mathcal{H}\left(\mathcal{V}\right)$ at an hyperarc $a_q \in \mathcal{A}_H$ being given by
	\begin{equation*}
		\nabla^*_a f: ~ \mathcal{A}_H \longrightarrow \mathbb{R} \qquad a_q \longmapsto \nabla^*_a f \left(a_q\right) = 
	\end{equation*}
	\begin{equation}
		W_G \left(a_q\right)^\theta \sum_{v_i \in \mathcal{V}} \left(\frac{\delta_{in}\left(v_i, a_q\right)}{\deg_{in}\left(v_i\right)} - \frac{\delta_{out}\left(v_i, a_q\right)}{\deg_{out}\left(v_i\right)}\right) w_I \left(v_i\right)^\alpha w_G \left(v_i\right)^\zeta f\left(v_i\right).
	\end{equation}\vspace{0em}
\end{definition}

When the hyperarc adjoint operator for hypergraphs is used on normal graphs then this results in the arc adjoint operator definition for normal graphs.\\

\begin{theorem}[\textbf{Hyperarc adjoint operator a generalization of the arc adjoint operator}]\label{nabla^*_a} \ \\
	Given a weighted oriented normal graph $OG = \left(\mathcal{V}, \mathcal{A}_G, w, W\right)$ and a vertex function $f \in \mathcal{H}\left(\mathcal{V}\right)$, then the hyperarc adjoint operator $\nabla_a^*$ for weighted oriented hypergraphs matches the arc adjoint operator for normal graphs from definition (\ref{Gnabla^*_a}).\\
\end{theorem}

\begin{proof}\ \\
	Given a weighted oriented normal graph $OG = \left(\mathcal{V}, \mathcal{A}_G, w, W\right)$ together with a vertex function $f \in \mathcal{H}\left(\mathcal{V}\right)$, then applying the hyperarc adjoint operator $\nabla_a^*$ for weighted oriented hypergraphs to any arc $a_q = \left(v_i, v_j\right) \in \mathcal{A}_G$ leads to the definition of the arc adjoint operator:
	
	$\begin{aligned}[t]
		\nabla_a^* f \left(a_q\right) = & ~ \nabla_a^* f \left(v_i, v_j\right) &\\
		= & ~ W_G \left(v_i, v_j\right)^\theta \sum_{v_k \in \mathcal{V}} \left(\frac{\delta_{in}\left(v_k, a_q\right)}{\deg_{in}\left(v_k\right)} - \frac{\delta_{out}\left(v_k, a_q\right)}{\deg_{out}\left(v_k\right)}\right) w_I \left(v_k\right)^\alpha w_G \left(v_k\right)^\zeta f\left(v_k\right) &\\
		= & ~ W_G \left(v_i, v_j\right)^\theta \sum_{v_k \in \left\{v_i, v_j\right\}} \left(\frac{\delta_{in}\left(v_k, a_q\right)}{\deg_{in}\left(v_k\right)} w_I \left(v_k\right)^\alpha w_G \left(v_k\right)^\zeta f\left(v_k\right) \right. -&\\
		& \left. \frac{\delta_{out}\left(v_k, a_q\right)}{\deg_{out}\left(v_k\right)} w_I \left(v_k\right)^\alpha w_G \left(v_k\right)^\zeta f\left(v_k\right)\right) &\\
		= & ~ W_G \left(v_i, v_j\right)^\theta \left(\frac{w_I \left(v_j\right)^\alpha w_G \left(v_j\right)^\zeta}{\deg_{in}\left(v_j\right)} f\left(v_j\right) - \frac{w_I \left(v_i\right)^\alpha w_G \left(v_i\right)^\zeta}{\deg_{out}\left(v_i\right)} f\left(v_i\right)\right) &\\
	\end{aligned}$\\

	This reformulation is based on the properties $a_q^{out} = \left\{v_i\right\}$ and $a_q^{in} = \left\{v_j\right\}$ for any arc $a_q = \left(v_i, v_j\right) \in \mathcal{A}_G$, which implies $\delta_{in}\left(v_j, a_q\right) = \delta_{out}\left(v_i, a_q\right) = 1$.\\
	
	Hence, the hyperarc adjoint operator $\nabla_a^*$ for hypergraphs is a valid generalization of the arc adjoint operator for normal graphs.\\
\end{proof}

The following theorem proves that the hyperarc adjoint is well-defined with regard to the hyperarc gradient.\\

\begin{theorem}[\textbf{Connection hyperarc gradient $\nabla_a$ and hyperarc adjoint $\nabla^*_a$}]\label{Hhyperarcgradadj} \ \\
	For a weighted oriented hypergraph $OH = \left(\mathcal{V}, \mathcal{A}_H, w, W\right)$, the hyperarc adjoint operator $\nabla^*_a$ fulfills the equality
	\begin{equation}
		{\langle f, \nabla_a G \rangle}_{\mathcal{H}\left(\mathcal{V}\right)} = {\langle G, \nabla^*_a f \rangle}_{\mathcal{H}\left(\mathcal{A}_H\right)}
	\end{equation}
	for all vertex functions $f \in \mathcal{H}\left(\mathcal{V}\right)$ and all hyperarc functions $G \in \mathcal{H}\left(\mathcal{A}_H\right)$.\\
\end{theorem}

\begin{proof}\ \\
	Given a weighted oriented hypergraph $OH = \left(\mathcal{V}, \mathcal{A}_H, w, W\right)$, a vertex function $f \in \mathcal{H}\left(\mathcal{V}\right)$ and a hyperarc function $G \in \mathcal{H}\left(\mathcal{A}_H\right)$, then the definitions of the inner product in $\mathcal{H}\left(\mathcal{V}\right)$ and the hyperarc gradient operator $\nabla_a$ lead to the following:
	
	$\begin{aligned}[t]
		{\langle f, \nabla_a G \rangle}_{\mathcal{H}\left(\mathcal{V}\right)} = & \sum_{v_i \in \mathcal{V}} w_I \left(v_i\right)^\alpha f\left(v_i\right) \nabla_a G\left(v_i\right) &\\
		= & \sum_{v_i \in \mathcal{V}} w_I \left(v_i\right)^\alpha f\left(v_i\right) &\\
		& w_G \left(v_i\right)^\zeta \sum_{a_q \in \mathcal{A}_H} \left(\frac{\delta_{in}\left(v_i, a_q\right)}{\deg_{in}\left(v_i\right)} - \frac{\delta_{out}\left(v_i, a_q\right)}{\deg_{out}\left(v_i\right)}\right) W_I \left(a_q\right)^\beta W_G \left(a_q\right)^\theta G\left(a_q\right) &\\
		= & \sum_{v_i \in \mathcal{V}} \sum_{a_q \in \mathcal{A}_H} \left(\frac{\delta_{in}\left(v_i, a_q\right)}{\deg_{in}\left(v_i\right)} - \frac{\delta_{out}\left(v_i, a_q\right)}{\deg_{out}\left(v_i\right)}\right) w_I \left(v_i\right)^\alpha w_G \left(v_i\right)^\zeta &\\
		& W_I \left(a_q\right)^\beta W_G \left(a_q\right)^\theta f\left(v_i\right) G\left(a_q\right) &\\
	\end{aligned}$\\

	Exchanging the two sums is feasible since the set of all vertices $\mathcal{V}$ and the set of all hyperarcs $\mathcal{A}_H$ are assumed to be finite and this results in:
	
	$\begin{aligned}[t]
		= & \sum_{a_q \in \mathcal{A}_H} \sum_{v_i \in \mathcal{V}} \left(\frac{\delta_{in}\left(v_i, a_q\right)}{\deg_{in}\left(v_i\right)} - \frac{\delta_{out}\left(v_i, a_q\right)}{\deg_{out}\left(v_i\right)}\right) w_I \left(v_i\right)^\alpha w_G \left(v_i\right)^\zeta W_I \left(a_q\right)^\beta W_G \left(a_q\right)^\theta f\left(v_i\right) G\left(a_q\right) &\\
		= & \sum_{a_q \in \mathcal{A}_H} W_I \left(a_q\right)^\beta W_G \left(a_q\right)^\theta G\left(a_q\right) \sum_{v_i \in \mathcal{V}} \left(\frac{\delta_{in}\left(v_i, a_q\right)}{\deg_{in}\left(v_i\right)} - \frac{\delta_{out}\left(v_i, a_q\right)}{\deg_{out}\left(v_i\right)}\right) w_I \left(v_i\right)^\alpha w_G \left(v_i\right)^\zeta f\left(v_i\right) &\\
	\end{aligned}$\\
	
	Using the definitions for the inner product on the space of all real hyperarc funtions $\mathcal{H}\left(\mathcal{A}_H\right)$ and the hyperarc adjoint operator $\nabla^*_a$ yields the following:
	
	$\begin{aligned}[t]
		= & \sum_{a_q \in \mathcal{A}_H} W_I \left(a_q\right)^\beta G\left(a_q\right) W_G \left(a_q\right)^\theta \sum_{v_i \in \mathcal{V}} \left(\frac{\delta_{in}\left(v_i, a_q\right)}{\deg_{in}\left(v_i\right)} - \frac{\delta_{out}\left(v_i, a_q\right)}{\deg_{out}\left(v_i\right)}\right) w_I \left(v_i\right)^\alpha w_G \left(v_i\right)^\zeta f\left(v_i\right) &\\
		= & \sum_{a_q \in \mathcal{A}_H} W_I \left(a_q\right)^{\beta} G\left(a_q\right) \nabla_a^* f \left(a_q\right) = {\langle G, \nabla^*_a f \rangle}_{\mathcal{H}\left(\mathcal{A}_H\right)} &\\
	\end{aligned}$\\
	
	Therefore, with the presented definition for the hyperarc adjoint operator $\nabla^*_a$, the equality ${\langle f, \nabla_a G \rangle}_{\mathcal{H}\left(\mathcal{V}\right)} = {\langle G, \nabla^*_a f \rangle}_{\mathcal{H}\left(\mathcal{A}_H\right)}$ holds true for all vertex functions $f \in \mathcal{H}\left(\mathcal{V}\right)$ and all hyperarc functions $G \in \mathcal{H}\left(\mathcal{A}_H\right)$.\\
\end{proof}
\clearpage
\section{Divergence and Laplacian operators on oriented hypergraphs}\label{10} 

The following section uses the vertex gradient, vertex adjoint, hyperarc gradient and hyperarc adjoint to first of all calculate the vertex and hyperarc divergence for oriented hypergraphs, see subsection (\ref{10.1}), which are necessary for later definitions.\\

In subsection (\ref{10.2}), the vertex and hyperarc divergences are used to introduce a Laplacian operator for vertices and a Laplacian operator for hyperarcs. Furthermore, this subsection also proves that the resulting Laplacian operators on hypergraphs are valid generalizations of the normal graph case.\\

Subsection (\ref{10.3}) derives the vertex $p$-Laplacian and the hyperarc $p$-Laplacian with the before presented definitions. Moreover, this subsection includes a proof showing that the presented $p$-Laplacian operators match the definitions introduced in \cite{jost2021plaplaceoperators}, when choosing specific parameters.\\

\subsection{Divergence operators on oriented hypergraphs}\label{10.1} 

As in the case of oriented normal graphs, the vertex divergence operator $\text{div}_v$ for hyperarc functions $F \in \mathcal{H}\left(\mathcal{A}_H\right)$ and the hyperarc divergence operator $\text{div}_a$ for vertex functions $f \in \mathcal{H}\left(\mathcal{V}\right)$ can be defined based on the equivalence from the continuum setting:
\begin{equation}
	\text{div}_v = - \nabla^*_v
\end{equation}
\begin{equation}
	\text{div}_a = - \nabla^*_a
\end{equation}\vspace{0em}

\begin{definition}[\textbf{Vertex divergence operator $\text{div}_v$}]\label{Hdiv_v} \ \\
	For a weighted oriented hypergraph $OH = \left(\mathcal{V}, \mathcal{A}_H, w, W\right)$ with symmetric hyperarc weight functions $W_I$ and $W_G$, the vertex divergence operator $\text{div}_v$ is given by
	\begin{equation}
		\text{div}_v: ~ \mathcal{H}\left(\mathcal{A}_H\right) \longrightarrow \mathcal{H}\left(\mathcal{V}\right) \qquad F \longmapsto \text{div}_v F
	\end{equation}
	with the weighted divergence of hyperarc function $F \in \mathcal{H}\left(\mathcal{A}_H\right)$ at a vertex $v_i \in \mathcal{V}$ being defined as
	\begin{equation*}
		\text{div}_v: ~ \mathcal{V} \longrightarrow \mathbb{R} \qquad v_i \longmapsto \text{div}_v \left(F\right) \left(v_i\right) = - \nabla^*_v F \left(v_i\right) =
	\end{equation*}
	\begin{equation}
		\sum_{a_q \in \mathcal{A}_H} \left(\delta_{out}\left(v_i, a_q\right) \frac{w_G \left(v_i\right)^\eta}{\left\lvert a_q^{out}\right\rvert} - \delta_{in}\left(v_i, a_q\right) \frac{w_G \left(v_i\right)^\epsilon}{\left\lvert a_q^{in}\right\rvert}\right) W_I \left(a_q\right)^\beta W_G \left(a_q\right)^\gamma F\left(a_q\right).
	\end{equation}\vspace{0em}
\end{definition}

Since the vertex adjoint operator for hypergraphs is a valid generalization of the vertex adjoint operator for normal graphs, as proven in theorem (\ref{nabla^*_v}), the vertex divergence for hypergraphs is a valid generalization of the normal graph case as well.\\

\begin{remark}[\textbf{Vertex divergence operator for hypergraphs a generalization of the normal graph case}]\label{div_v} \ \\
	Given a weighted oriented normal graph $OG = \left(\mathcal{V}, \mathcal{A}_G, w, W\right)$ with symmetric arc weight functions $W_I$ and $W_G$ and an arc function $F \in \mathcal{H}\left(\mathcal{A}_G\right)$, then the vertex divergence operator $\text{div}_v$ for weighted oriented hypergraphs matches the vertex divergence operator for normal graphs from definition (\ref{Gdiv_v}).\\
\end{remark}

The hyperarc divergence is defined in the same way as the vertex divergence, but based on the hyperarc adjoint instead of the vertex adjoint.\\

\begin{definition}[\textbf{Hyperarc divergence operator $\text{div}_a$}]\label{Hdiv_a} \ \\
	For a weighted oriented hypergraph $OH = \left(\mathcal{V}, \mathcal{A}_H, w, W\right)$, the vertex divergence operator $\text{div}_v$ is defined as
	\begin{equation}
		\text{div}_a: ~ \mathcal{H}\left(\mathcal{V}\right) \longrightarrow \mathcal{H}\left(\mathcal{A}_H\right) \qquad f \longmapsto \text{div}_a f
	\end{equation}
	with the weighted divergence of vertex function $f \in \mathcal{H}\left(\mathcal{V}\right)$ at a hyperarc $a_q \in \mathcal{A}_H$ being given by
	\begin{equation*}
		\text{div}_a: ~ \mathcal{A}_H \longrightarrow \mathbb{R} \qquad a_q \longmapsto \text{div}_a \left(f\right) \left(a_q\right) = - \nabla^*_a f \left(a_q\right) =
	\end{equation*}
	\begin{equation}
		W_G \left(a_q\right)^\theta \sum_{v_i \in \mathcal{V}} \left(\frac{\delta_{out}\left(v_i, a_q\right)}{\deg_{out}\left(v_i\right)} - \frac{\delta_{in}\left(v_i, a_q\right)}{\deg_{in}\left(v_i\right)}\right) w_I \left(v_i\right)^\alpha w_G \left(v_i\right)^\zeta f\left(v_i\right).
	\end{equation}\vspace{0em}
\end{definition}

Theorem (\ref{nabla^*_a}) proves that the hyperarc adjoint operator for hypergraphs is a valid generalization of the arc adjoint operator for normal graphs and hence the hyperarc divergence for hypergraphs is a valid generalization of the normal graph case as well.\\

\begin{remark}[\textbf{Hyperarc divergence operator for hypergraphs a generalization of the normal graph case}]\label{div_a} \ \\
	Given a weighted oriented normal graph $OG = \left(\mathcal{V}, \mathcal{A}_G, w, W\right)$ and a vertex function $f \in \mathcal{H}\left(\mathcal{V}\right)$, then the hyperarc divergence operator $\text{div}_a$ for weighted oriented hypergraphs matches the arc divergence operator for normal graphs from definition (\ref{Gdiv_a}).\\
\end{remark}
\clearpage
\subsection{Laplacian operators on oriented hypergraphs}\label{10.2} 

Both, the vertex Laplacian operator and the hyperarc Laplacian operator for oriented hypergraphs, are based on the relationship in the continuum setting\\
\begin{equation}
	\Delta_v f = \text{div}_v \left(\nabla_v f\right) 
\end{equation}
\begin{equation}
	\Delta_a F = \text{div}_a \left(\nabla_a F\right) 
\end{equation}
for any vertex function $f \in \mathcal{H}\left(\mathcal{V}\right)$ and any hyperarc function $F \in \mathcal{H}\left(\mathcal{A}_H\right)$.\\

\begin{definition}[\textbf{Vertex Laplacian operator $\Delta_v$}]\label{HDelta_v} \ \\
	For a weighted oriented hypergraph $OH = \left(\mathcal{V}, \mathcal{A}_H, w, W\right)$ with symmetric hyperarc weight functions $W_I$ and $W_G$, the vertex Laplacian operator $\Delta_v$ is given by
	\begin{equation}
		\Delta_v: ~ \mathcal{H}\left(\mathcal{V}\right) \longrightarrow \mathcal{H}\left(\mathcal{V}\right) \qquad f \longmapsto \Delta_v f
	\end{equation}
	with the weighted Laplacian of vertex function $f \in \mathcal{H}\left(\mathcal{V}\right)$ at a vertex $v_i \in \mathcal{V}$ being defined as
	\begin{equation*}
		\Delta_v: ~ \mathcal{V} \longrightarrow \mathbb{R} \qquad v_i \longmapsto \Delta_v f \left(v_i\right) =
	\end{equation*}
	\begin{equation*}
		\sum_{a_q \in \mathcal{A}_H} \left(\delta_{out}\left(v_i, a_q\right) \frac{w_G \left(v_i\right)^\eta}{\left\lvert a_q^{out}\right\rvert} - \delta_{in}\left(v_i, a_q\right) \frac{w_G \left(v_i\right)^\epsilon}{\left\lvert a_q^{in}\right\rvert}\right) W_I \left(a_q\right)^\beta W_G \left(a_q\right)^{2 \gamma}
	\end{equation*}
	\begin{equation}
		\sum_{v_j \in \mathcal{V}} \left(\delta_{in}\left(v_j, a_q\right) \frac{w_I \left(v_j\right)^\alpha w_G \left(v_j\right)^\epsilon}{\left\lvert a_q^{in}\right\rvert} - \delta_{out}\left(v_j, a_q\right) \frac{w_I \left(v_j\right)^\alpha w_G \left(v_j\right)^\eta}{\left\lvert a_q^{out}\right\rvert}\right) f\left(v_j\right)
	\end{equation}\vspace{0em}
\end{definition}

The theorem below proves that applying the vertex Laplacian operator for hypergraphs to vertex functions on a normal graph results in the vertex Laplacian operator for normal graphs.\\

\begin{theorem}[\textbf{Vertex Laplacian operator for hypergraphs a generalization of the normal graph case}]\label{HDelta_vgen} \ \\
	Given a weighted oriented normal graph $OG = \left(\mathcal{V}, \mathcal{A}_G, w, W\right)$ with symmetric arc weight functions $W_I$ and $W_G$ and a vertex function $f \in \mathcal{H}\left(\mathcal{V}\right)$, then the vertex Laplacian operator $\Delta_v$ for weighted oriented hypergraphs matches the vertex Laplacian operator for normal graphs from definition (\ref{GDelta_v}).\\
\end{theorem}

\begin{proof}\ \\
	Given an oriented normal graph $OG = \left(\mathcal{V}, \mathcal{A}_G, w, W\right)$ with symmetric arc weight functions $W_I$ and $W_G$ together with a vertex function $f \in \mathcal{H}\left(\mathcal{V}\right)$, then applying the Laplacian gradient operator $\Delta_v$ for weighted oriented hypergraphs to any vertex $v_i \in \mathcal{V}$ results in:
	
	$\begin{aligned}[t]
		\Delta_v f \left(v_i\right) = & \sum_{a_q \in \mathcal{A}_H} \left(\delta_{out}\left(v_i, a_q\right) \frac{w_G \left(v_i\right)^\eta}{\left\lvert a_q^{out}\right\rvert} - \delta_{in}\left(v_i, a_q\right) \frac{w_G \left(v_i\right)^\epsilon}{\left\lvert a_q^{in}\right\rvert}\right) W_I \left(a_q\right)^\beta W_G \left(a_q\right)^{2 \gamma} &\\
		& \quad \sum_{v_j \in \mathcal{V}} \left(\delta_{in}\left(v_j, a_q\right) \frac{w_I \left(v_j\right)^\alpha w_G \left(v_j\right)^\epsilon}{\left\lvert a_q^{in}\right\rvert} - \delta_{out}\left(v_j, a_q\right) \frac{w_I \left(v_j\right)^\alpha w_G \left(v_j\right)^\eta}{\left\lvert a_q^{out}\right\rvert}\right) f\left(v_j\right) &\\
		= & \sum_{a_q \in \mathcal{A}_G} \left(\delta_{out}\left(v_i, a_q\right) w_G \left(v_i\right)^\eta - \delta_{in}\left(v_i, a_q\right) w_G \left(v_i\right)^\epsilon\right) W_I \left(a_q\right)^\beta W_G \left(a_q\right)^{2 \gamma} &\\
		& \quad \sum_{v_j \in \mathcal{V}} \left(\delta_{in}\left(v_j, a_q\right) w_I \left(v_j\right)^\alpha w_G \left(v_j\right)^\epsilon - \delta_{out}\left(v_j, a_q\right) w_I \left(v_j\right)^\alpha w_G \left(v_j\right)^\eta\right) f\left(v_j\right) &\\
		= & \sum_{a_q \in \mathcal{A}_G} \sum_{v_j \in \mathcal{V}} \left(\delta_{out}\left(v_i, a_q\right) w_G \left(v_i\right)^\eta - \delta_{in}\left(v_i, a_q\right) w_G \left(v_i\right)^\epsilon\right) W_I \left(a_q\right)^\beta W_G \left(a_q\right)^{2 \gamma} &\\
		& \quad \left(\delta_{in}\left(v_j, a_q\right) w_I \left(v_j\right)^\alpha w_G \left(v_j\right)^\epsilon - \delta_{out}\left(v_j, a_q\right) w_I \left(v_j\right)^\alpha w_G \left(v_j\right)^\eta\right) f\left(v_j\right) &\\
		= & \sum_{a_q \in \mathcal{A}_G} \sum_{v_j \in \mathcal{V}} \left(\delta_{out}\left(v_i, a_q\right) \delta_{in}\left(v_j, a_q\right) w_G \left(v_i\right)^\eta w_I \left(v_j\right)^\alpha w_G \left(v_j\right)^\epsilon \right. - &\\
		& \quad ~ \delta_{out}\left(v_i, a_q\right) \delta_{out}\left(v_j, a_q\right) w_G \left(v_i\right)^\eta w_I \left(v_j\right)^\alpha w_G \left(v_j\right)^\eta - &\\
		& \quad ~ \delta_{in}\left(v_i, a_q\right) \delta_{in}\left(v_j, a_q\right) w_G \left(v_i\right)^\epsilon w_I \left(v_j\right)^\alpha w_G \left(v_j\right)^\epsilon + &\\
		& \quad \left.\delta_{in}\left(v_i, a_q\right) \delta_{out}\left(v_j, a_q\right) w_G \left(v_i\right)^\epsilon w_I \left(v_j\right)^\alpha w_G \left(v_j\right)^\eta\right) W_I \left(a_q\right)^\beta W_G \left(a_q\right)^{2 \gamma} f\left(v_j\right) &\\
	\end{aligned}$\\
	
	This reformulation is based on the properties $a_q^{out} = \left\{v_k\right\}$ and $a_q^{in} = \left\{v_l\right\}$ for any arc $a_q = \left(v_k, v_l\right) \in \mathcal{A}_G$, which implies $\left\lvert a_q^{in} \right\rvert = \left\lvert a_q^{out} \right\rvert = 1$. Because of this, it is impossible, that $v_i \in \mathcal{V}$ and $v_j \in \mathcal{V} \backslash \left\{v_i\right\}$ are both output or input vertices of any arc $a_q \in \mathcal{A}_G$, which yields the following:
	
	$\begin{aligned}[t]
		= & \sum_{a_q \in \mathcal{A}_G} \left( - \delta_{out}\left(v_i, a_q\right) \delta_{out}\left(v_i, a_q\right) w_G \left(v_i\right)^\eta w_I \left(v_i\right)^\alpha w_G \left(v_i\right)^\eta \right. - &\\
		& \quad \left. \delta_{in}\left(v_i, a_q\right) \delta_{in}\left(v_i, a_q\right) w_G \left(v_i\right)^\epsilon w_I \left(v_i\right)^\alpha w_G \left(v_i\right)^\epsilon\right) W_I \left(a_q\right)^\beta W_G \left(a_q\right)^{2 \gamma} f\left(v_i\right) + &\\
		& \quad \sum_{v_j \in \mathcal{V} \backslash\{v_i\}} \left(\delta_{out}\left(v_i, a_q\right) \delta_{in}\left(v_j, a_q\right) w_G \left(v_i\right)^\eta w_I \left(v_j\right)^\alpha w_G \left(v_j\right)^\epsilon \right. + &\\
		& \quad \left. \delta_{out}\left(v_j, a_q\right) \delta_{in}\left(v_i, a_q\right) w_G \left(v_i\right)^\epsilon w_I \left(v_j\right)^\alpha w_G \left(v_j\right)^\eta\right) W_I \left(a_q\right)^\beta W_G \left(a_q\right)^{2 \gamma} f\left(v_j\right) &\\
		= & \sum_{a_q \in \mathcal{A}_G} \left( - \delta_{out}\left(v_i, a_q\right) w_i \left(v_i\right)^\alpha w_G \left(v_i\right)^{2 \eta} - \delta_{in}\left(v_i, a_q\right) w_I \left(v_i\right)^\alpha w_G \left(v_i\right)^{2 \epsilon}\right) &\\
		& \quad ~ W_I \left(a_q\right)^\beta W_G \left(a_q\right)^{2 \gamma} f\left(v_i\right) + &\\
		& \quad \sum_{v_j \in \mathcal{V} \backslash\{v_i\}} \left(\delta_{out}\left(v_i, a_q\right) \delta_{in}\left(v_j, a_q\right) w_G \left(v_i\right)^\eta w_I \left(v_j\right)^\alpha w_G \left(v_j\right)^\epsilon \right. + &\\
		& \quad \left. \delta_{out}\left(v_j, a_q\right) \delta_{in}\left(v_i, a_q\right) w_G \left(v_i\right)^\epsilon w_I \left(v_j\right)^\alpha w_G \left(v_j\right)^\eta\right) W_I \left(a_q\right)^\beta W_G \left(a_q\right)^{2 \gamma} f\left(v_j\right) &\\
	\end{aligned}$\\
	
	\clearpage
	Moreover, summing over all arcs $a_q \in \mathcal{A}_G$ which include vertex $v_i$, either as an output or as an input vertex \big(either $\delta_{out}\left(v_i, a_q\right) = 1$ or $\delta_{in}\left(v_i, a_q\right) = 1$\big), and summing over all vertices $v_j \in \mathcal{V}$ which are part of such arc $a_q$ \big(either $\delta_{in}\left(v_j, a_q\right) = 1$ or $\delta_{out}\left(v_j, a_q\right) = 1$\big) is the same as summing over all vertices $v_j \in \mathcal{V}$ and checking whether arc $\left(v_i, v_j\right) \in \mathcal{A}_G$ or arc $\left(v_j, v_i\right) \in \mathcal{A}_G$:
	
	$\begin{aligned}[t]
		= & \sum_{v_j \in \mathcal{V}}  \left( - \delta\left(v_i, v_j\right) w_I \left(v_i\right)^\alpha w_G \left(v_i\right)^{2 \eta} W_I \left(v_i, v_j\right)^\beta W_G \left(v_i, v_j\right)^{2 \gamma} \right. - &\\
		& \quad \left. \delta\left(v_j, v_i\right) w_I \left(v_i\right)^\alpha w_G \left(v_i\right)^{2 \epsilon} W_I \left(v_j, v_i\right)^\beta W_G \left(v_j, v_i\right)^{2 \gamma}\right) f\left(v_i\right) + &\\
		& \quad \sum_{v_j \in \mathcal{V} \backslash\{v_i\}} \left(\delta\left(v_i, v_j\right) w_G \left(v_i\right)^\eta w_I \left(v_j\right)^\alpha w_G \left(v_j\right)^\epsilon W_I \left(v_i, v_j\right)^\beta W_G \left(v_i, v_j\right)^{2 \gamma} \right. + &\\
		& \quad \left. \delta\left(v_j, v_i\right) w_I \left(v_j\right)^\alpha w_G \left(v_j\right)^\eta w_G \left(v_i\right)^\epsilon W_I \left(v_j, v_i\right)^\beta W_G \left(v_j, v_i\right)^{2 \gamma}\right) f\left(v_j\right) &\\
		= & \sum_{v_j \in \mathcal{V}}  \left( - \delta\left(v_i, v_j\right) w_I \left(v_i\right)^\alpha w_G \left(v_i\right)^{2 \eta} - \delta\left(v_j, v_i\right) w_I \left(v_i\right)^\alpha w_G \left(v_i\right)^{2 \epsilon} \right) &\\
		& \quad ~ W_I \left(v_i, v_j\right)^\beta W_G \left(v_i, v_j\right)^{2 \gamma} f\left(v_i\right) + &\\
		& \quad \left(\delta\left(v_i, v_j\right) w_G \left(v_i\right)^\eta w_I \left(v_j\right)^\alpha w_G \left(v_j\right)^\epsilon + \delta\left(v_j, v_i\right) w_I \left(v_j\right)^\alpha w_G \left(v_j\right)^\eta w_G \left(v_i\right)^\epsilon\right) &\\
		& \quad ~ W_I \left(v_i, v_j\right)^\beta W_G \left(v_i, v_j\right)^{2 \gamma} f\left(v_j\right) &\\
		= & \sum_{v_j \in \mathcal{V}}  \left(\left(\delta\left(v_i, v_j\right) w_G \left(v_i\right)^\eta w_I \left(v_j\right)^\alpha w_G \left(v_j\right)^\epsilon + \delta\left(v_j, v_i\right) w_I \left(v_j\right)^\alpha w_G \left(v_j\right)^\eta w_G \left(v_i\right)^\epsilon\right) f\left(v_j\right)\right. - &\\
		& \quad \left.\left(\delta\left(v_i, v_j\right) w_I \left(v_i\right)^\alpha w_G \left(v_i\right)^{2 \eta} + \delta\left(v_j, v_i\right) w_I \left(v_i\right)^\alpha w_G \left(v_i\right)^{2 \epsilon}\right) f\left(v_i\right)\right) &\\
		& \quad ~ W_I \left(v_i, v_j\right)^\beta W_G \left(v_i, v_j\right)^{2 \gamma}&\\
	\end{aligned}$\\
	
	Hence, the vertex Laplacian operator $\Delta_v$ for hypergraphs is a valid generalization of the normal graph case.\\
\end{proof}

The theorem below proves that the Laplacian for vertices $v_i \in \mathcal{V}$ is well-defined by using the previously introduced vertex gradient and vertex divergence on hypergraphs.\\

\begin{theorem}[\textbf{Connection vertex divergence $\text{div}_v$, vertex gradient $\nabla_v$, and vertex Laplacian $\Delta_v$}]\label{HvertexLap} \ \\
	On a weighted oriented hypergraph $OH = \left(\mathcal{V}, \mathcal{A}_H, w, W\right)$ with symmetric hyperarc weight functions $W_I$ and $W_G$, the vertex Laplacian operator $\Delta_v$ fulfills the equality
	\begin{equation}
		\Delta_v f = \text{div}_v \left(\nabla_v f\right)
	\end{equation}
	for all vertex functions $f \in \mathcal{H}\left(\mathcal{V}\right)$.\\
\end{theorem}

\clearpage
\begin{proof}\ \\
	Given a weighted oriented hypergraph $OH = \left(\mathcal{V}, \mathcal{A}_H, w, W\right)$ with symmetric hyperarc weight functions $W_I$ and $W_G$ and a vertex function $f \in \mathcal{H}\left(\mathcal{V}\right)$, then applying the definitions of the vertex divergence operator $\text{div}_v$ and the vertex gradient operator $\nabla_v$ yields the following for all vertices $v_i \in \mathcal{V}$:\\
	$\begin{aligned}[t]	
		\text{div}_v \left(\nabla_v f\right)\left(v_i\right) = & \sum_{a_q \in \mathcal{A}_H} \left(\delta_{out}\left(v_i, a_q\right) \frac{w_G \left(v_i\right)^\eta}{\left\lvert a_q^{out}\right\rvert} - \delta_{in}\left(v_i, a_q\right) \frac{w_G \left(v_i\right)^\epsilon}{\left\lvert a_q^{in}\right\rvert}\right) &\\
		& W_I \left(a_q\right)^\beta W_G \left(a_q\right)^\gamma \nabla_v f\left(a_q\right) &\\
	\end{aligned}$\\
	$\begin{aligned}[t]	
		= & \sum_{a_q \in \mathcal{A}_H} \left(\delta_{out}\left(v_i, a_q\right) \frac{w_G \left(v_i\right)^\eta}{\left\lvert a_q^{out}\right\rvert} - \delta_{in}\left(v_i, a_q\right) \frac{w_G \left(v_i\right)^\epsilon}{\left\lvert a_q^{in}\right\rvert}\right) W_I \left(a_q\right)^\beta W_G \left(a_q\right)^\gamma &\\
		& W_G \left(a_q\right)^\gamma \sum_{v_j \in \mathcal{V}} \left(\delta_{in}\left(v_j, a_q\right) \frac{w_I \left(v_j\right)^\alpha w_G \left(v_j\right)^\epsilon}{\left\lvert a_q^{in}\right\rvert} - \delta_{out}\left(v_j, a_q\right) \frac{w_I \left(v_j\right)^\alpha w_G \left(v_j\right)^\eta}{\left\lvert a_q^{out}\right\rvert}\right) f\left(v_j\right) &\\
		= & \sum_{a_q \in \mathcal{A}_H} \left(\delta_{out}\left(v_i, a_q\right) \frac{w_G \left(v_i\right)^\eta}{\left\lvert a_q^{out}\right\rvert} - \delta_{in}\left(v_i, a_q\right) \frac{w_G \left(v_i\right)^\epsilon}{\left\lvert a_q^{in}\right\rvert}\right) W_I \left(a_q\right)^\beta W_G \left(a_q\right)^{2 \gamma} &\\
		& \sum_{v_j \in \mathcal{V}} \left(\delta_{in}\left(v_j, a_q\right) \frac{w_I \left(v_j\right)^\alpha w_G \left(v_j\right)^\epsilon}{\left\lvert a_q^{in}\right\rvert} - \delta_{out}\left(v_j, a_q\right) \frac{w_I \left(v_j\right)^\alpha w_G \left(v_j\right)^\eta}{\left\lvert a_q^{out}\right\rvert}\right) f\left(v_j\right) &\\
		= & ~ \Delta_v f\left(v_i\right)
	\end{aligned}$\\

	Thus, with the previously introduced definitions of the vertex gradient $\nabla_v$, the vertex divergence $\text{div}_v$, and the vertex Laplacian $\Delta_v$, the equality $\Delta_v f\left(v_i\right) = \text{div}_v \left(\nabla_v f\right)\left(v_i\right)$ holds true for all vertices $v_i \in \mathcal{V}$ and for all vertex functions $f \in \mathcal{H}\left(\mathcal{V}\right)$.\\
\end{proof}

\begin{definition}[\textbf{Hyperarc Laplacian operator $\Delta_a$}]\label{HDelta_a} \ \\
	For a weighted oriented hypergraph $OH = \left(\mathcal{V}, \mathcal{A}_H, w, W\right)$, the hyperarc Laplacian operator $\Delta_a$ is defined as
	\begin{equation}
		\Delta_a: ~ \mathcal{H}\left(\mathcal{A}_H\right) \longrightarrow \mathcal{H}\left(\mathcal{A}_H\right) \qquad F \longmapsto \Delta_a F
	\end{equation}
	with the weighted Laplacian of hyperarc function $F \in \mathcal{H}\left(\mathcal{A}_H\right)$ at a hyperarc $a_q \in \mathcal{A}_H$ being given by
	\begin{equation*}
		\Delta_a: ~ \mathcal{A}_H \longrightarrow \mathbb{R} \qquad a_q \longmapsto \Delta_a F \left(a_q\right) =
	\end{equation*}
	\begin{equation*}
		W_G \left(a_q\right)^\theta \sum_{v_i \in \mathcal{V}} \left(\frac{\delta_{out}\left(v_i, a_q\right)}{\deg_{out}\left(v_i\right)} - \frac{\delta_{in}\left(v_i, a_q\right)}{\deg_{in}\left(v_i\right)}\right) w_I \left(v_i\right)^\alpha w_G \left(v_i\right)^{2 \zeta}
	\end{equation*}
	\begin{equation}
		\sum_{a_r \in \mathcal{A}_H} \left(\frac{\delta_{in}\left(v_i, a_r\right)}{\deg_{in}\left(v_i\right)} - \frac{\delta_{out}\left(v_i, a_r\right)}{\deg_{out}\left(v_i\right)}\right) W_I \left(a_r\right)^\beta W_G \left(a_r\right)^\theta F\left(a_r\right).
	\end{equation}\vspace{0em}
\end{definition}

Since the hyperarc Laplacian operator for oriented hypergraphs can be applied to arcs of an oriented normal graph, the next theorem shows that the hyperarc Laplacian operator matches the definition of the arc Laplacian operator for an oriented normal graph.\\

\begin{theorem}[\textbf{Hyperarc Laplacian operator a generalization of the arc Laplacian operator}]\label{HDelta_agen} \ \\
	Given a weighted oriented normal graph $OG = \left(\mathcal{V}, \mathcal{A}_G, w, W\right)$ and an arc function $F \in \mathcal{H}\left(\mathcal{A}_G\right)$, then the hyperarc Laplacian operator $\Delta_a$ for weighted oriented hypergraphs matches the arc Laplacian operator for normal graphs from definition (\ref{GDelta_a}).\\
\end{theorem}

\begin{proof}\ \\
	Given an oriented normal graph $OG = \left(\mathcal{V}, \mathcal{A}_G, w, W\right)$ together with an arc function $F \in \mathcal{H}\left(\mathcal{A}_G\right)$, then applying the hyperarc Laplacian operator $\nabla_a$ for weighted oriented hypergraphs to an arc $a_q = \left(v_i, v_j\right) \in \mathcal{A}_G$ leads to the definition of the arc Laplacian operator:
	
	$\begin{aligned}[t]
		\Delta_a F \left(a_q\right) = & ~ \Delta_a F \left(v_i, v_j\right) &\\
		= & ~ W_G \left(v_i, v_j\right)^\theta \sum_{v_k \in \mathcal{V}} \left(\frac{\delta_{out}\left(v_k, a_q\right)}{\deg_{out}\left(v_k\right)} - \frac{\delta_{in}\left(v_k, a_q\right)}{\deg_{in}\left(v_k\right)}\right) w_I \left(v_k\right)^\alpha w_G \left(v_k\right)^{2 \zeta} &\\
		& \quad \sum_{a_r \in \mathcal{A}_G} \left(\frac{\delta_{in}\left(v_k, a_r\right)}{\deg_{in}\left(v_k\right)} - \frac{\delta_{out}\left(v_k, a_r\right)}{\deg_{out}\left(v_k\right)}\right) W_I \left(a_r\right)^\beta W_G \left(a_r\right)^\theta F\left(a_r\right) &\\
		= & ~ W_G \left(v_i, v_j\right)^\theta \sum_{v_k \in \left\{v_i, v_j\right\}} \left(\frac{\delta_{out}\left(v_k, a_q\right)}{\deg_{out}\left(v_k\right)} - \frac{\delta_{in}\left(v_k, a_q\right)}{\deg_{in}\left(v_k\right)}\right) w_I \left(v_k\right)^\alpha w_G \left(v_k\right)^{2 \zeta} &\\
		& \quad \sum_{a_r \in \mathcal{A}_G} \left(\frac{\delta_{in}\left(v_k, a_r\right)}{\deg_{in}\left(v_k\right)} - \frac{\delta_{out}\left(v_k, a_r\right)}{\deg_{out}\left(v_k\right)}\right) W_I \left(a_r\right)^\beta W_G \left(a_r\right)^\theta F\left(a_r\right) &\\
		= & ~ W_G \left(v_i, v_j\right)^\theta \left(\frac{w_I \left(v_i\right)^\alpha w_G \left(v_i\right)^{2 \zeta}}{\deg_{out}\left(v_i\right)} \right. &\\
		& \quad \sum_{a_r \in \mathcal{A}_G} \left(\frac{\delta_{in}\left(v_i, a_r\right)}{\deg_{in}\left(v_i\right)} - \frac{\delta_{out}\left(v_i, a_r\right)}{\deg_{out}\left(v_i\right)}\right) W_I \left(a_s\right)^\beta W_G \left(a_s\right)^\theta F\left(a_r\right) - &\\
		& \quad \frac{w_I \left(v_j\right)^\alpha w_G \left(v_j\right)^{2 \zeta}}{\deg_{in}\left(v_j\right)} &\\
		& \quad \left. \sum_{a_s \in \mathcal{A}_G} \left(\frac{\delta_{in}\left(v_j, a_s\right)}{\deg_{in}\left(v_j\right)} - \frac{\delta_{out}\left(v_j, a_s\right)}{\deg_{out}\left(v_j\right)}\right) W_I \left(a_s\right)^\beta W_G \left(a_s\right)^\theta F\left(a_s\right)\right) &\\
	\end{aligned}$\\

	This reformulation is based on the property that for every arc $a_q = \left(v_i, v_j\right) \in \mathcal{A}_G$ it holds true that $a_q^{out} = \left\{v_i\right\}$ and $a_q^{in} = \left\{v_j\right\}$ which directly implies that $\delta_{out}\left(v_i, a_q\right) = 1$ and $\delta_{in}\left(v_j, a_q\right) = 1$.\\
	
	Therefore, the hyperarc Laplacian operator $\Delta_a$ for hypergraphs is a valid generalization of the arc Laplacian operator for normal graphs.\\
\end{proof}

The following theorem proves that the hyperarc Laplacian is well-defined based on the presented definitions of the hyperarc gradient and the hyperarc divergence.\\

\begin{theorem}[\textbf{Connection hyperarc divergence $\text{div}_a$, hyperarc gradient $\nabla_a$, and hyperarc Laplacian $\Delta_a$}]\label{HhyperarcLap} \ \\
	On a weighted oriented hypergraph $OH = \left(\mathcal{V}, \mathcal{A}_H, w, W\right)$, the hyperarc Laplacian operator $\Delta_a$ fulfills the equality
	\begin{equation}
		\Delta_a F = \text{div}_a \left(\nabla_a F\right)
	\end{equation}
	for all hyperarc functions $F \in \mathcal{H}\left(\mathcal{A}_H\right)$.\\
\end{theorem}

\begin{proof}\ \\
	Given a weighted oriented hypergraph $OH = \left(\mathcal{V}, \mathcal{A}_H, w, W\right)$ and a hyperarc function $F \in \mathcal{H}\left(\mathcal{A}_H\right)$, the definitions of the hyperarc divergence operator $\text{div}_a$ and the hyperarc gradient operator $\nabla_a$ result in the following for all hyperarcs $a_q \in \mathcal{A}_H$:
	
	$\begin{aligned}[t]	
		\text{div}_a \left(\nabla_a F\right)\left(a_q\right) = & ~ W_G \left(a_q\right)^\theta \sum_{v_i \in \mathcal{V}} \left(\frac{\delta_{out}\left(v_i, a_q\right)}{\deg_{out}\left(v_i\right)} - \frac{\delta_{in}\left(v_i, a_q\right)}{\deg_{in}\left(v_i\right)}\right) w_I \left(v_i\right)^\alpha w_G \left(v_i\right)^\zeta \nabla_a F\left(v_i\right) &\\
	\end{aligned}$\\
	$\begin{aligned}[t]	
		= & ~ W_G \left(a_q\right)^\theta \sum_{v_i \in \mathcal{V}} \left(\frac{\delta_{out}\left(v_i, a_q\right)}{\deg_{out}\left(v_i\right)} - \frac{\delta_{in}\left(v_i, a_q\right)}{\deg_{in}\left(v_i\right)}\right) w_I \left(v_i\right)^\alpha w_G \left(v_i\right)^\zeta w_G \left(v_i\right)^\zeta &\\
		& \sum_{a_r \in \mathcal{A}_H} \left(\frac{\delta_{in}\left(v_i, a_r\right)}{\deg_{in}\left(v_i\right)} - \frac{\delta_{out}\left(v_i, a_r\right)}{\deg_{out}\left(v_i\right)}\right) W_I \left(a_r\right)^\beta W_G \left(a_r\right)^\theta F\left(a_r\right) &\\
		= & ~ W_G \left(a_q\right)^\theta \sum_{v_i \in \mathcal{V}} \left(\frac{\delta_{out}\left(v_i, a_q\right)}{\deg_{out}\left(v_i\right)} - \frac{\delta_{in}\left(v_i, a_q\right)}{\deg_{in}\left(v_i\right)}\right) w_I \left(v_i\right)^\alpha w_G \left(v_i\right)^{2 \zeta} &\\
		& \sum_{a_r \in \mathcal{A}_H} \left(\frac{\delta_{in}\left(v_i, a_r\right)}{\deg_{in}\left(v_i\right)} - \frac{\delta_{out}\left(v_i, a_r\right)}{\deg_{out}\left(v_i\right)}\right) W_I \left(a_r\right)^\beta W_G \left(a_r\right)^\theta F\left(a_r\right) &\\
		= & ~ \Delta_a F \left(a_q\right) &\\
	\end{aligned}$\\
	
	Thus, with the presented definition for the hyperarc gradient $\nabla_a$, the hyperarc divergence $\text{div}_a$, and the hyperarc Laplacian $\Delta_a$, the equality $\Delta_a F\left(a_q\right) = \text{div}_a \left(\nabla_a F\right)\left(a_q\right)$ holds true for all hyperarcs $a_q \in \mathcal{A}_H$ and for all hyperarc functions $F \in \mathcal{H}\left(\mathcal{A}_H\right)$.\\
\end{proof}
\clearpage
\subsection{$p$-Laplacian operators on oriented hypergraphs}\label{10.3}

The definition of the weighted $p$-Laplacian operators extend the operators of the continuum setting which implies that for all $p \in \left(1, \infty\right)$ it must holds true that\\
\begin{equation}
	\Delta_v^p f = \text{div}_v \left(\left\lvert \nabla_v f\right\rvert^{p - 2} \nabla_v f\right)
\end{equation}
\begin{equation}
	\Delta_a^p F = \text{div}_a \left(\left\lvert \nabla_a F\right\rvert^{p - 2} \nabla_a F\right)
\end{equation}
for any vertex function $f \in \mathcal{H}\left(\mathcal{V}\right)$ and any hyperarc function $F \in \mathcal{H}\left(\mathcal{A}_H\right)$.\\

\begin{definition}[\textbf{Vertex $p$-Laplacian operator $\Delta_v^p$}]\ \\
	For a weighted oriented hypergraph $OH = \left(\mathcal{V}, \mathcal{A}_H, w, W\right)$ with symmetric hyperarc weight functions $W_I$ and $W_G$, the vertex $p$-Laplacian operator $\Delta_v^p$ is given by
	\begin{equation}
		\Delta_v^p: ~ \mathcal{H}\left(\mathcal{V}\right) \longrightarrow \mathcal{H}\left(\mathcal{V}\right) \qquad f \longmapsto \Delta_v^p f
	\end{equation}
	with the weighted $p$-Laplacian of vertex function $f \in \mathcal{H}\left(\mathcal{V}\right)$ at a vertex $v_i \in \mathcal{V}$ being defined as
	\begin{equation*}
		\Delta_v^p: ~ \mathcal{V} \longrightarrow \mathbb{R} \qquad v_i \longmapsto \Delta_v^p f \left(v_i\right) =
	\end{equation*}
	\begin{equation*}
		\sum_{a_q \in \mathcal{A}_H} \left(\delta_{out}\left(v_i, a_q\right) \frac{w_G \left(v_i\right)^\eta}{\left\lvert a_q^{out}\right\rvert} - \delta_{in}\left(v_i, a_q\right) \frac{w_G \left(v_i\right)^\epsilon}{\left\lvert a_q^{in}\right\rvert}\right) W_I \left(a_q\right)^\beta W_G \left(a_q\right)^{p \gamma} 
	\end{equation*}
	\begin{equation*}
		\left\lvert\sum_{v_j \in \mathcal{V}} \left(\delta_{in}\left(v_j, a_q\right) \frac{w_I \left(v_j\right)^\alpha w_G \left(v_j\right)^\epsilon}{\left\lvert a_q^{in}\right\rvert} - \delta_{out}\left(v_j, a_q\right) \frac{w_I \left(v_j\right)^\alpha w_G \left(v_j\right)^\eta}{\left\lvert a_q^{out}\right\rvert}\right) f\left(v_j\right)\right\rvert^{p - 2}
	\end{equation*}
	\begin{equation}
		\sum_{v_k \in \mathcal{V}} \left(\delta_{in}\left(v_k, a_q\right) \frac{w_I \left(v_k\right)^\alpha w_G \left(v_k\right)^\epsilon}{\left\lvert a_q^{in}\right\rvert} - \delta_{out}\left(v_k, a_q\right) \frac{w_I \left(v_k\right)^\alpha w_G \left(v_k\right)^\eta}{\left\lvert a_q^{out}\right\rvert}\right) f\left(v_k\right).
	\end{equation}\vspace{0em}
\end{definition}

\clearpage
\begin{theorem}[\textbf{Parameter choice for the vertex $p$-Laplacian operator}]\label{HDelta_v^pparam} \ \\
	A simplified definition of the vertex $p$-Laplacian is introduced in \cite{jost2021plaplaceoperators} and using the notation of this thesis, it can be written as the following:
	\begin{equation*}
		\Delta_p f \left(v_i\right) = \frac{1}{\deg\left(v_i\right)} \sum_{\substack{a_q \in \mathcal{A}_H: ~ \delta_{out}\left(v_i, a_q\right) = 1 \\ \text{or} ~ \delta_{in}\left(v_i, a_q\right) = 1}} \left\lvert\sum_{v_j \in a_q^{in}} f\left(v_j\right) - \sum_{v_j \in a_q^{out}} f\left(v_j\right)\right\lvert^{p - 2}
	\end{equation*}
	\begin{equation*}
		\left(\sum_{v_k \in \mathcal{V}} \left(\delta_{out}\left(v_i, a_q\right) \delta_{out}\left(v_k, a_q\right) + \delta_{in}\left(v_i, a_q\right) \delta_{in}\left(v_k, a_q\right)\right) f\left(v_k\right) - \right.
	\end{equation*}
	\begin{equation}
		\left.\sum_{v_k \in \mathcal{V}} \left(\delta_{out}\left(v_i, a_q\right) \delta_{in}\left(v_k, a_q\right) + \delta_{in}\left(v_i, a_q\right) \delta_{out}\left(v_k, a_q\right)\right) f\left(v_k\right)\right) 
	\end{equation}
	for any vertex function $f \in \mathcal{H}\left(\mathcal{V}\right)$ and for all vertices $v_i \in \mathcal{V}$.\\
	
	The factor $\left(\delta_{out}\left(v_i, a_q\right) \delta_{out}\left(v_k, a_q\right) + \delta_{in}\left(v_i, a_q\right) \delta_{in}\left(v_k, a_q\right)\right)$ is always equal to zero, unless $v_i, v_k \in a_q^{out}$ or $v_i, v_k \in a_q^{in}$, which means that the vertices $v_i$ and $v_k$ are co-oriented. Similarly, the factor $\left(\delta_{out}\left(v_i, a_q\right) \delta_{in}\left(v_k, a_q\right) + \delta_{in}\left(v_i, a_q\right) \delta_{out}\left(v_k, a_q\right)\right)$ ensures to only consider vertices $v_k \in \mathcal{V}$ which are anti-oriented compared to vertex $v_i$ and hence $v_i \in a_q^{out}, v_k \in a_q^{in}$ or $v_i \in a_q^{in}, v_k \in a_q^{out}$.\\
	
	Thus, choosing the parameters of the vertex $p$-Laplacian $\Delta_v^p$ as $\alpha = 0$, $\beta = 0$, $\gamma = 0$, $\epsilon = 0$ and $\eta = 0$ together with excluding the $\frac{1}{\left\lvert a_q^{out}\right\rvert}$ and $\frac{1}{\left\lvert a_q^{in}\right\rvert}$ multiplicative factors and including a new $- \frac{1}{\deg\left(v_i\right)}$ factor in the vertex adjoint and the vertex divergence, results in the simplified vertex $p$-Laplacian introduced in \cite{jost2021plaplaceoperators}.\\
	
	Furthermore, applying these parameter choices to the vertex gradient, the vertex adjoint and the vertex divergence leads to:
	\begin{itemize}
		\item $\nabla_v f \left(a_q\right) = \sum_{v_i \in \mathcal{V}} \left(\delta_{in}\left(v_i, a_q\right) - \delta_{out}\left(v_i, a_q\right)\right) f\left(v_i\right)$
		\item $\nabla^*_v F \left(v_i\right) = - \frac{1}{\deg\left(v_i\right)} \sum_{a_q \in \mathcal{A}_H} \left(\delta_{in}\left(v_i, a_q\right) - \delta_{out}\left(v_i, a_q\right)\right) F\left(a_q\right)$
		\item $\text{div}_v \left(F\right) \left(v_i\right) = - \frac{1}{\deg\left(v_i\right)} \sum_{a_q \in \mathcal{A}_H} \left(\delta_{out}\left(v_i, a_q\right) - \delta_{in}\left(v_i, a_q\right)\right) F\left(a_q\right)$
	\end{itemize}
	for any vertex function $f \in \mathcal{H}\left(\mathcal{V}\right)$, any hyperarc function $f \in \mathcal{H}\left(\mathcal{A}_H\right)$, and for all hyperarcs $a_q \in \mathcal{A}_H$ and all vertices $v_i \in \mathcal{V}$.\\
\end{theorem}

\clearpage
\begin{proof}\ \\
	Given an oriented hypergraph $OH = \left(\mathcal{V}, \mathcal{A}_H, w, W\right)$ with symmetric hyperarc weight functions $W_I$ and $W_G$ together with a vertex function $f \in \mathcal{H}\left(\mathcal{V}\right)$, then applying the described parameter choices, excluding the $\frac{1}{\left\lvert a_q^{out}\right\rvert}$ and $\frac{1}{\left\lvert a_q^{in}\right\rvert}$ multiplicative factors and including the new $- \frac{1}{\deg\left(v_i\right)}$ factor in the vertex divergence yields the following:
	
	$\begin{aligned}[t]
		\Delta_p f \left(v_i\right) = & ~ - \frac{1}{\deg\left(v_i\right)} \sum_{a_q \in \mathcal{A}_H} \left(\delta_{out}\left(v_i, a_q\right) - \delta_{in}\left(v_i, a_q\right)\right) &\\
		& \quad \left\lvert\sum_{v_j \in \mathcal{V}} \left(\delta_{in}\left(v_j, a_q\right) - \delta_{out}\left(v_j, a_q\right)\right) f\left(v_j\right)\right\rvert^{p - 2} &\\
		& \quad \sum_{v_k \in \mathcal{V}} \left(\delta_{in}\left(v_k, a_q\right) - \delta_{out}\left(v_k, a_q\right)\right) f\left(v_k\right) &\\
		= & ~ - \frac{1}{\deg\left(v_i\right)} \sum_{a_q \in \mathcal{A}_H} \left\lvert\sum_{v_j \in a_q^{in}} f\left(v_j\right) - \sum_{v_j \in a_q^{out}} f\left(v_j\right)\right\rvert^{p - 2} &\\
		& \quad \sum_{v_k \in \mathcal{V}} \left(\delta_{out}\left(v_i, a_q\right) \delta_{in}\left(v_k, a_q\right) - \delta_{out}\left(v_i, a_q\right) \delta_{out}\left(v_k, a_q\right) - \right. &\\
		& \quad \left. \delta_{in}\left(v_i, a_q\right) \delta_{in}\left(v_k, a_q\right) + \delta_{in}\left(v_i, a_q\right) \delta_{out}\left(v_k, a_q\right)\right) f\left(v_k\right) &\\
		= & ~ - \frac{1}{\deg\left(v_i\right)} \sum_{a_q \in \mathcal{A}_H} \left\lvert\sum_{v_j \in a_q^{in}} f\left(v_j\right) - \sum_{v_j \in a_q^{out}} f\left(v_j\right)\right\rvert^{p - 2} &\\
		& \quad \sum_{v_k \in \mathcal{V}} \left(- \delta_{out}\left(v_i, a_q\right) \delta_{out}\left(v_k, a_q\right) - \delta_{in}\left(v_i, a_q\right) \delta_{in}\left(v_k, a_q\right) + \right. &\\
		& \quad \left. \delta_{out}\left(v_i, a_q\right) \delta_{in}\left(v_k, a_q\right) + \delta_{in}\left(v_i, a_q\right) \delta_{out}\left(v_k, a_q\right)\right) f\left(v_k\right) &\\
	\end{aligned}$\\
	
	Moving the minus from the front to the last bracket and sorting the last bracket based on whether the vertices $v_i$ and $v_k$ are co-oriented or anti-oriented in hyperarc $a_q$ results in the simplified definition of the vertex $p$-Laplacian:
	
	$\begin{aligned}[t]
		= & ~ \frac{1}{\deg\left(v_i\right)} \sum_{a_q \in \mathcal{A}_H} \left\lvert\sum_{v_j \in a_q^{in}} f\left(v_j\right) - \sum_{v_j \in a_q^{out}} f\left(v_j\right)\right\rvert^{p - 2} &\\
		& \quad \sum_{v_k \in \mathcal{V}} \left(\delta_{out}\left(v_i, a_q\right) \delta_{out}\left(v_k, a_q\right) + \delta_{in}\left(v_i, a_q\right) \delta_{in}\left(v_k, a_q\right)\right. - &\\
		& \quad \left.\delta_{out}\left(v_i, a_q\right) \delta_{in}\left(v_k, a_q\right) - \delta_{in}\left(v_i, a_q\right) \delta_{out}\left(v_k, a_q\right)\right) f\left(v_k\right) &\\
	\end{aligned}$\\	
	$\begin{aligned}[t]
		= & ~ \frac{1}{\deg\left(v_i\right)} \sum_{a_q \in \mathcal{A}_H} \left\lvert\sum_{v_j \in a_q^{in}} f\left(v_j\right) - \sum_{v_j \in a_q^{out}} f\left(v_j\right)\right\rvert^{p - 2} &\\
		& \quad \left(\sum_{v_k \in \mathcal{V}} \left(\delta_{out}\left(v_i, a_q\right) \delta_{out}\left(v_k, a_q\right) + \delta_{in}\left(v_i, a_q\right) \delta_{in}\left(v_k, a_q\right) \right) f\left(v_k\right) \right. - &\\
		& \quad \left.\sum_{v_k \in \mathcal{V}} \left(\delta_{out}\left(v_i, a_q\right) \delta_{in}\left(v_k, a_q\right) + \delta_{in}\left(v_i, a_q\right) \delta_{out}\left(v_k, a_q\right)\right) f\left(v_k\right)\right) &\\
		= & ~ \frac{1}{\deg\left(v_i\right)} \sum_{\substack{a_q \in \mathcal{A}_H: ~ \delta_{out}\left(v_i, a_q\right) = 1 \\ \text{or} ~ \delta_{in}\left(v_i, a_q\right) = 1}} \left\lvert\sum_{v_j \in a_q^{in}} f\left(v_j\right) - \sum_{v_j \in a_q^{out}} f\left(v_j\right)\right\rvert^{p - 2} &\\
		& \quad \left(\sum_{v_k \in \mathcal{V}} \left(\delta_{out}\left(v_i, a_q\right) \delta_{out}\left(v_k, a_q\right) + \delta_{in}\left(v_i, a_q\right) \delta_{in}\left(v_k, a_q\right) \right) f\left(v_k\right) \right. - &\\
		& \quad \left.\sum_{v_k \in \mathcal{V}} \left(\delta_{out}\left(v_i, a_q\right) \delta_{in}\left(v_k, a_q\right) + \delta_{in}\left(v_i, a_q\right) \delta_{out}\left(v_k, a_q\right)\right) f\left(v_k\right)\right) &\\
	\end{aligned}$\\

	The change in the first sum \big(instead of summing over all hyperarcs $a_q \in \mathcal{A}_H$, now summing over all hyperarcs $a_q \in \mathcal{A}_H$ with $\delta_{out}\left(v_i, a_q\right) = 1$ or $\delta_{in}\left(v_i, a_q\right) = 1$\big) is no substantial change, since the last big bracket already ensures that only hyperarcs $a_q$ with $v_i \in a_q^{out}$ or $v_i \in a_q^{in}$ are considered, because otherwise the value of the big bracket would be equal to zero.\\
	
	Overall, this proof shows that the presented definition of the vertex $p$-Laplacian operator $\Delta_v^p$ is compatible with the definition of the vertex $p$-Laplacian in \cite{jost2021plaplaceoperators}.\\
\end{proof}

The following theorem proves that the vertex $p$-Laplacian on hypergraphs is well-defined based on the definitions for the vertex gradient and vertex divergence.\\

\begin{theorem}[\textbf{Connection vertex divergence $\text{div}_v$, vertex gradient $\nabla_v$, and vertex $p$-Laplacian $\Delta_v^p$}]\ \\
	On a weighted oriented hypergraph $OH = \left(\mathcal{V}, \mathcal{A}_H, w, W\right)$ with symmetric hyperarc weight functions $W_I$ and $W_G$, the vertex $p$-Laplacian operator $\Delta_v^p$ fulfills the equality
	\begin{equation}
		\Delta_v^p f = \text{div}_v \left(\left\lvert \nabla_v f\right\rvert^{p - 2} \nabla_v f\right)
	\end{equation}
	for all vertex functions $f \in \mathcal{H}\left(\mathcal{V}\right)$.\\
\end{theorem}

\clearpage
\begin{proof}\ \\
	Given a weighted oriented hypergraph $OH = \left(\mathcal{V}, \mathcal{A}_H, w, W\right)$ with symmetric hyperarc weight functions $W_I$ and $W_G$ and a vertex function $f \in \mathcal{H}\left(\mathcal{V}\right)$, then the definitions of the vertex divergence operator $\text{div}_v$ and the vertex gradient operator $\nabla_v$ lead to the following for all vertices $v_i \in \mathcal{V}$:
	
	$\begin{aligned}[t]	
		\text{div}_v \left(\left\lvert \nabla_v f\right\rvert^{p - 2} \nabla_v f\right)\left(v_i\right)	= & \sum_{a_q \in \mathcal{A}_H} \left(\delta_{out}\left(v_i, a_q\right) \frac{w_G \left(v_i\right)^\eta}{\left\lvert a_q^{out}\right\rvert} - \delta_{in}\left(v_i, a_q\right) \frac{w_G \left(v_i\right)^\epsilon}{\left\lvert a_q^{in}\right\rvert}\right) &\\
		& W_I \left(a_q\right)^\beta W_G \left(a_q\right)^\gamma \left\lvert \nabla_v f\left(a_q\right)\right\rvert^{p - 2} \nabla_v f\left(a_q\right) &\\
	\end{aligned}$\\
	$\begin{aligned}[t]	
		= & \sum_{a_q \in \mathcal{A}_H} \left(\delta_{out}\left(v_i, a_q\right) \frac{w_G \left(v_i\right)^\eta}{\left\lvert a_q^{out}\right\rvert} - \delta_{in}\left(v_i, a_q\right) \frac{w_G \left(v_i\right)^\epsilon}{\left\lvert a_q^{in}\right\rvert}\right) W_I \left(a_q\right)^\beta W_G \left(a_q\right)^\gamma &\\
		& \quad \left\lvert W_G \left(a_q\right)^\gamma \sum_{v_j \in \mathcal{V}} \left(\delta_{in}\left(v_j, a_q\right) \frac{w_I \left(v_j\right)^\alpha w_G \left(v_j\right)^\epsilon}{\left\lvert a_q^{in}\right\rvert} - \delta_{out}\left(v_j, a_q\right) \frac{w_I \left(v_j\right)^\alpha w_G \left(v_j\right)^\eta}{\left\lvert a_q^{out}\right\rvert}\right) f\left(v_j\right)\right\rvert^{p - 2} &\\
		& \quad W_G \left(a_q\right)^\gamma \sum_{v_k \in \mathcal{V}} \left(\delta_{in}\left(v_k, a_q\right) \frac{w_I \left(v_k\right)^\alpha w_G \left(v_k\right)^\epsilon}{\left\lvert a_q^{in}\right\rvert} - \delta_{out}\left(v_k, a_q\right) \frac{w_I \left(v_k\right)^\alpha w_G \left(v_k\right)^\eta}{\left\lvert a_q^{out}\right\rvert}\right) f\left(v_k\right) &\\
	\end{aligned}$\\
	
	Since the hyperarc weight function $W_G$ maps to positive values, the following reformulation holds true and leads to the vertex $p$-Laplacian definition:
	
	$\begin{aligned}[t]	
		= & \sum_{a_q \in \mathcal{A}_H} \left(\delta_{out}\left(v_i, a_q\right) \frac{w_G \left(v_i\right)^\eta}{\left\lvert a_q^{out}\right\rvert} - \delta_{in}\left(v_i, a_q\right) \frac{w_G \left(v_i\right)^\epsilon}{\left\lvert a_q^{in}\right\rvert}\right) W_I \left(a_q\right)^\beta W_G \left(a_q\right)^{2 \gamma + \gamma \left(p - 2\right)} &\\
		& \quad \left\lvert \sum_{v_j \in \mathcal{V}} \left(\delta_{in}\left(v_j, a_q\right) \frac{w_I \left(v_j\right)^\alpha w_G \left(v_j\right)^\epsilon}{\left\lvert a_q^{in}\right\rvert} - \delta_{out}\left(v_j, a_q\right) \frac{w_I \left(v_j\right)^\alpha w_G \left(v_j\right)^\eta}{\left\lvert a_q^{out}\right\rvert}\right) f\left(v_j\right)\right\rvert^{p - 2} &\\
		& \quad \sum_{v_k \in \mathcal{V}} \left(\delta_{in}\left(v_k, a_q\right) \frac{w_I \left(v_k\right)^\alpha w_G \left(v_k\right)^\epsilon}{\left\lvert a_q^{in}\right\rvert} - \delta_{out}\left(v_k, a_q\right) \frac{w_I \left(v_k\right)^\alpha w_G \left(v_k\right)^\eta}{\left\lvert a_q^{out}\right\rvert}\right) f\left(v_k\right) &\\
	\end{aligned}$\\
	$\begin{aligned}[t]	
		= & \sum_{a_q \in \mathcal{A}_H} \left(\delta_{out}\left(v_i, a_q\right) \frac{w_G \left(v_i\right)^\eta}{\left\lvert a_q^{out}\right\rvert} - \delta_{in}\left(v_i, a_q\right) \frac{w_G \left(v_i\right)^\epsilon}{\left\lvert a_q^{in}\right\rvert}\right) W_I \left(a_q\right)^\beta W_G \left(a_q\right)^{p \gamma} &\\
		& \quad \left\lvert \sum_{v_j \in \mathcal{V}} \left(\delta_{in}\left(v_j, a_q\right) \frac{w_I \left(v_j\right)^\alpha w_G \left(v_j\right)^\epsilon}{\left\lvert a_q^{in}\right\rvert} - \delta_{out}\left(v_j, a_q\right) \frac{w_I \left(v_j\right)^\alpha w_G \left(v_j\right)^\eta}{\left\lvert a_q^{out}\right\rvert}\right) f\left(v_j\right)\right\rvert^{p - 2} &\\
		& \quad \sum_{v_k \in \mathcal{V}} \left(\delta_{in}\left(v_k, a_q\right) \frac{w_I \left(v_k\right)^\alpha w_G \left(v_k\right)^\epsilon}{\left\lvert a_q^{in}\right\rvert} - \delta_{out}\left(v_k, a_q\right) \frac{w_I \left(v_k\right)^\alpha w_G \left(v_k\right)^\eta}{\left\lvert a_q^{out}\right\rvert}\right) f\left(v_k\right) &\\
		= & ~ \Delta_v^p f \left(v_i\right) &\\
	\end{aligned}$\\
	
	Therefore, using the previously introduced definitions for the vertex gradient $\nabla_v$, the vertex divergence $\text{div}_v$, and the vertex $p$-Laplacian $\Delta_v^p$ shows that the equality $\Delta_v^p f\left(v_i\right) = \text{div}_v \left(\left\lvert \nabla_v f\right\rvert^{p - 2} \nabla_v f\right) \left(v_i\right)$ holds true for all vertices $v_i \in \mathcal{V}$ and for all vertex functions $f \in \mathcal{H}\left(\mathcal{V}\right)$.\\
\end{proof}

\begin{definition}[\textbf{Hyperarc $p$-Laplacian operator $\Delta_a^p$}]\ \\
	For a weighted oriented hypergraph $OH = \left(\mathcal{V}, \mathcal{A}_H, w, W\right)$, the hyperarc $p$-Laplacian operator $\Delta_a^p$ is defined as
	\begin{equation}
		\Delta_a^p: ~ \mathcal{H}\left(\mathcal{A}_H\right) \longrightarrow \mathcal{H}\left(\mathcal{A}_H\right) \qquad F \longmapsto \Delta_a^p F
	\end{equation}
	with the weighted $p$-Laplacian of hyperarc function $F \in \mathcal{H}\left(\mathcal{A}_H\right)$ at a hyperarc $a_q \in \mathcal{A}_H$ being given by
	\begin{equation*}
		\Delta_a^p: ~ \mathcal{A}_H \longrightarrow \mathbb{R} \qquad a_q \longmapsto \Delta_a^p F \left(a_q\right) =
	\end{equation*}
	\begin{equation*}
		W_G \left(a_q\right)^{\theta} \sum_{v_i \in \mathcal{V}} \left(\frac{\delta_{out}\left(v_i, a_q\right)}{\deg_{out}\left(v_i\right)} - \frac{\delta_{in}\left(v_i, a_q\right)}{\deg_{in}\left(v_i\right)}\right) w_I \left(v_i\right)^\alpha w_G \left(v_i\right)^{2 \zeta}
	\end{equation*}
	\begin{equation*}
		\left\lvert\sum_{a_r \in \mathcal{A}_H} \left(\frac{\delta_{in}\left(v_i, a_r\right)}{\deg_{in}\left(v_i\right)} - \frac{\delta_{out}\left(v_i, a_r\right)}{\deg_{out}\left(v_i\right)}\right) W_I \left(a_r\right)^\beta W_G \left(a_r\right)^\theta F\left(a_r\right)\right\rvert^{p - 2}
	\end{equation*}
	\begin{equation}
		\sum_{a_s \in \mathcal{A}_H} \left(\frac{\delta_{in}\left(v_i, a_s\right)}{\deg_{in}\left(v_i\right)} - \frac{\delta_{out}\left(v_i, a_s\right)}{\deg_{out}\left(v_i\right)}\right) W_I \left(a_s\right)^\beta W_G \left(a_s\right)^\theta F\left(a_s\right).
	\end{equation}\vspace{0em}
\end{definition}

\clearpage
\begin{theorem}[\textbf{Parameter choice for the hyperarc $p$-Laplacian operator}]\label{HDelta_a^pparam} \ \\
	A simplified definition of the hyperarc $p$-Laplacian is introduced in \cite{jost2021plaplaceoperators} as hyperedge $p$-Laplacian and with the notation of this thesis, it can be written down as
	\begin{equation*}
		\Delta_p F \left(a_q\right) = \sum_{\substack{v_i \in \mathcal{V}: \\ v_i \in a_q^{out} \cup a_q^{in}}} \frac{1}{\deg\left(v_i\right)} \left\lvert\sum_{a_r \in \mathcal{A}_H} \delta_{in}\left(v_i, a_r\right) F\left(a_r\right) - \sum_{a_r \in \mathcal{A}_H} \delta_{out}\left(v_i, a_r\right) F\left(a_r\right)\right\lvert^{p - 2}
	\end{equation*}
	\begin{equation*}
		\left(\sum_{a_s \in \mathcal{A}_H} \left(\delta_{out}\left(v_i, a_q\right) \delta_{out}\left(v_i, a_s\right) + \delta_{in}\left(v_i, a_q\right) \delta_{in}\left(v_i, a_s\right)\right) F\left(a_s\right) - \right.
	\end{equation*}
	\begin{equation}
		\left.\sum_{a_s \in \mathcal{A}_H} \left(\delta_{out}\left(v_i, a_q\right) \delta_{in}\left(v_i, a_s\right) + \delta_{in}\left(v_i, a_q\right) \delta_{out}\left(v_i, a_s\right)\right) F\left(a_s\right)\right) 
	\end{equation}
	for any hyperarc function $F \in \mathcal{H}\left(\mathcal{A}_H\right)$ and for all hyperarcs $a_q \in \mathcal{A}_H$.\\
	
	The factor $\left(\delta_{out}\left(v_i, a_q\right) \delta_{out}\left(v_i, a_s\right) + \delta_{in}\left(v_i, a_q\right) \delta_{in}\left(v_i, a_s\right)\right)$ ensures that only hyperarcs $a_s \in \mathcal{A}_H$ are considered, where vertex $v_i$ is either in both hyperarcs $a_q$ and $a_s$ an output vertex or an input vertex and hence $a_q$ and $a_s$ have the same orientation concerning vertex $v_i$. Similarly, the factor $\left(\delta_{out}\left(v_i, a_q\right) \delta_{in}\left(v_k, a_q\right) + \delta_{in}\left(v_i, a_q\right) \delta_{out}\left(v_k, a_q\right)\right)$ is only nonzero if the hyperarcs $a_q$ and $a_s$ have the opposite orientation concerning vertex $v_i$, which means that either $v_i \in a_q^{out}, v_i \in a_s^{in}$ or $v_i \in a_q^{in}, v_i \in a_s^{out}$.\\
	
	Therefore, choosing the parameters of the hyperarc $p$-Laplacian $\Delta_a^p$ as $\alpha = 0$, $\beta = 0$, $\zeta = 0$ and $\theta = 0$ together with excluding the $\frac{1}{\deg_{out}\left(v_i\right)}$ and $\frac{1}{\deg_{in}\left(v_i\right)}$ multiplicative factors and including a new $- \frac{1}{\deg\left(v_i\right)}$ factor in the sums of the hyperarc adjoint and the hyperarc divergence, leads to the the simplified hyperarc $p$-Laplacian defined in \cite{jost2021plaplaceoperators}.\\
	
	Moreover, applying these parameter choices to the hyperarc gradient, the hyperarc adjoint and the hyperarc divergence leads to:
	\begin{itemize}
		\item $\nabla_a F \left(v_i\right) = \sum_{a_q \in \mathcal{A}_H} \left(\delta_{in}\left(v_i, a_q\right) - \delta_{out}\left(v_i, a_q\right)\right) F\left(a_q\right)$
		\item $\nabla^*_a f \left(a_q\right) = - \sum_{v_i \in \mathcal{V}} \frac{1}{\deg\left(v_i\right)}  \left(\delta_{in}\left(v_i, a_q\right) - \delta_{out}\left(v_i, a_q\right)\right) f\left(v_i\right)$
		\item $\text{div}_a \left(f\right) \left(a_q\right) = - \sum_{v_i \in \mathcal{V}} \frac{1}{\deg\left(v_i\right)} \left(\delta_{out}\left(v_i, a_q\right) - \delta_{in}\left(v_i, a_q\right)\right) f\left(v_i\right)$
	\end{itemize}
	for any vertex function $f \in \mathcal{H}\left(\mathcal{V}\right)$, any hyperarc function $f \in \mathcal{H}\left(\mathcal{A}_H\right)$, and for all hyperarcs $a_q \in \mathcal{A}_H$ and all vertices $v_i \in \mathcal{V}$.\\
\end{theorem}

\clearpage
\begin{proof}\ \\
	Given an oriented hypergraph $OH = \left(\mathcal{V}, \mathcal{A}_H, w, W\right)$ with symmetric hyperarc weight functions $W_I$ and $W_G$ together with a hyperarc function $F \in \mathcal{H}\left(\mathcal{A}_H\right)$, then implementing the described parameter choices, excluding the $\frac{1}{\deg_{out}\left(v_i\right)}$ and $\frac{1}{\deg_{in}\left(v_i\right)}$ multiplicative factors and including the new $- \frac{1}{\deg\left(v_i\right)}$ factor in the hyperarc divergence leads to the following:
	
	$\begin{aligned}[t]
		\Delta_p f \left(v_i\right) = & ~ - \sum_{v_i \in \mathcal{V}} \frac{1}{\deg\left(v_i\right)} \left(\delta_{out}\left(v_i, a_q\right) - \delta_{in}\left(v_i, a_q\right)\right) &\\
		& \left\lvert\sum_{a_r \in \mathcal{A}_H} \left(\delta_{in}\left(v_i, a_r\right) - \delta_{out}\left(v_i, a_r\right)\right) F\left(a_r\right)\right\rvert^{p - 2} &\\
		& \sum_{a_s \in \mathcal{A}_H} \left(\delta_{in}\left(v_i, a_s\right) - \delta_{out}\left(v_i, a_s\right)\right) F\left(a_s\right) &\\
	\end{aligned}$\\
	$\begin{aligned}[t]
		= & ~ - \sum_{v_i \in \mathcal{V}} \frac{1}{\deg\left(v_i\right)} \left\lvert\sum_{a_r \in \mathcal{A}_H} \delta_{in}\left(v_i, a_r\right) F\left(a_r\right) - \sum_{a_r \in \mathcal{A}_H} \delta_{out}\left(v_i, a_r\right) F\left(a_r\right)\right\lvert^{p - 2} &\\
		& \quad \sum_{a_s \in \mathcal{A}_H} \left(\delta_{out}\left(v_i, a_q\right) \delta_{in}\left(v_i, a_s\right) - \delta_{out}\left(v_i, a_q\right) \delta_{out}\left(v_i, a_s\right) \right. - &\\
		& \quad \left. \delta_{in}\left(v_i, a_q\right) \delta_{in}\left(v_i, a_s\right) + \delta_{in}\left(v_i, a_q\right) \delta_{out}\left(v_i, a_s\right)\right) F\left(a_s\right) &\\
		= & ~ - \sum_{v_i \in \mathcal{V}} \frac{1}{\deg\left(v_i\right)} \left\lvert\sum_{a_r \in \mathcal{A}_H} \delta_{in}\left(v_i, a_r\right) F\left(a_r\right) - \sum_{a_r \in \mathcal{A}_H} \delta_{out}\left(v_i, a_r\right) F\left(a_r\right)\right\lvert^{p - 2} &\\
		& \quad \sum_{a_s \in \mathcal{A}_H} \left(- \delta_{out}\left(v_i, a_q\right) \delta_{out}\left(v_i, a_s\right) - \delta_{in}\left(v_i, a_q\right) \delta_{in}\left(v_i, a_s\right) \right. + &\\
		& \quad \left. \delta_{out}\left(v_i, a_q\right) \delta_{in}\left(v_i, a_s\right) + \delta_{in}\left(v_i, a_q\right) \delta_{out}\left(v_i, a_s\right)\right) F\left(a_s\right) &\\
	\end{aligned}$\\
	
	Moving the minus from the front to the last bracket and sorting the last bracket based on whether the hyperarcs $a_q$ and $a_s$ have the same orientation or the opposite orientation concerning vertex $v_i$ results in the simplified definition of the hyperarc $p$-Laplacian:
	
	$\begin{aligned}[t]
		= & ~ \sum_{v_i \in \mathcal{V}} \frac{1}{\deg\left(v_i\right)} \left\lvert\sum_{a_r \in \mathcal{A}_H} \delta_{in}\left(v_i, a_r\right) F\left(a_r\right) - \sum_{a_r \in \mathcal{A}_H} \delta_{out}\left(v_i, a_r\right) F\left(a_r\right)\right\lvert^{p - 2} &\\
		& \quad \sum_{a_s \in \mathcal{A}_H} \left(\delta_{out}\left(v_i, a_q\right) \delta_{out}\left(v_i, a_s\right) + \delta_{in}\left(v_i, a_q\right) \delta_{in}\left(v_i, a_s\right) \right. - &\\
		& \quad \left. \delta_{out}\left(v_i, a_q\right) \delta_{in}\left(v_i, a_s\right) - \delta_{in}\left(v_i, a_q\right) \delta_{out}\left(v_i, a_s\right)\right) F\left(a_s\right) &\\
	\end{aligned}$\\
	$\begin{aligned}[t]
		= & ~ \sum_{v_i \in \mathcal{V}} \frac{1}{\deg\left(v_i\right)} \left\lvert\sum_{a_r \in \mathcal{A}_H} \delta_{in}\left(v_i, a_r\right) F\left(a_r\right) - \sum_{a_r \in \mathcal{A}_H} \delta_{out}\left(v_i, a_r\right) F\left(a_r\right)\right\lvert^{p - 2} &\\
		& \quad \sum_{a_s \in \mathcal{A}_H} \left(\delta_{out}\left(v_i, a_q\right) \delta_{out}\left(v_i, a_s\right) + \delta_{in}\left(v_i, a_q\right) \delta_{in}\left(v_i, a_s\right)\right) F\left(a_s\right) - &\\
		& \quad \sum_{a_s \in \mathcal{A}_H} \left(\delta_{out}\left(v_i, a_q\right) \delta_{in}\left(v_i, a_s\right) + \delta_{in}\left(v_i, a_q\right) \delta_{out}\left(v_i, a_s\right)\right) F\left(a_s\right) &\\
		= & ~ \sum_{\substack{v_i \in \mathcal{V}: \\ v_i \in a_q^{out} \cup a_q^{in}}} \frac{1}{\deg\left(v_i\right)} \left\lvert\sum_{a_r \in \mathcal{A}_H} \delta_{in}\left(v_i, a_r\right) F\left(a_r\right) - \sum_{a_r \in \mathcal{A}_H} \delta_{out}\left(v_i, a_r\right) F\left(a_r\right)\right\lvert^{p - 2} &\\
		& \quad \sum_{a_s \in \mathcal{A}_H} \left(\delta_{out}\left(v_i, a_q\right) \delta_{out}\left(v_i, a_s\right) + \delta_{in}\left(v_i, a_q\right) \delta_{in}\left(v_i, a_s\right)\right) F\left(a_s\right) - &\\
		& \quad \sum_{a_s \in \mathcal{A}_H} \left(\delta_{out}\left(v_i, a_q\right) \delta_{in}\left(v_i, a_s\right) + \delta_{in}\left(v_i, a_q\right) \delta_{out}\left(v_i, a_s\right)\right) F\left(a_s\right) &\\
	\end{aligned}$\\
	
	The change in the first sum \big(instead of summing over all vertices $v_i \in \mathcal{V}$, now summing over all vertices $v_i \in \mathcal{V}$ with $v_i \in a_q^{out} \cup a_q^{in}$\big) is no substantial change, since the last big bracket already ensures that only vertices $v_i$ with $v_i \in a_q^{out}$ or $v_i \in a_q^{in}$ are considered, because otherwise the value of the big bracket would be equal to zero.\\
	
	This proof shows that the previously introduced definition of the hyperarc $p$-Laplacian operator $\Delta_a^p$ is compatible with the definition of the hyperarc $p$-Laplacian in \cite{jost2021plaplaceoperators}.\\
\end{proof}

The theorem below proves that the $p$-Laplacian for hyperarcs $a_q \in \mathcal{H}\left(\mathcal{A}_H\right)$ is well-defined with regards to the hyperarc gradient and the hyperarc divergence.\\

\begin{theorem}[\textbf{Connection hyperarc divergence $\text{div}_a$, hyperarc gradient $\nabla_a$, and hyperarc $p$-Laplacian $\Delta_a^p$}]\ \\
	On a weighted oriented hypergraph $OH = \left(\mathcal{V}, \mathcal{A}_H, w, W\right)$, the hyperarc $p$-Laplacian operator $\Delta_a^p$ fulfills the equality
	\begin{equation}
		\Delta_a^p F = \text{div}_a \left(\left\lvert \nabla_a F\right\rvert^{p - 2} \nabla_a F\right)
	\end{equation}
	for all hyperarc functions $F \in \mathcal{H}\left(\mathcal{A}_H\right)$.\\
\end{theorem}

\clearpage
\begin{proof}\ \\
	Given a weighted oriented hypergraph $OH = \left(\mathcal{V}, \mathcal{A}_H, w, W\right)$ and a hyperarc function $F \in \mathcal{H}\left(\mathcal{A}_H\right)$, then using the definitions of the hyperarc divergence operator $\text{div}_a$ and the hyperarc gradient operator $\nabla_a$ yields the following for all hyperarcs $a_q \in \mathcal{A}_H$:
	
	$\begin{aligned}[t]	
		\text{div}_a \left(\left\lvert \nabla_a F\right\rvert^{p - 2} \nabla_a F\right)\left(v_i, v_j\right) = & ~ W_G \left(a_q\right)^\theta \sum_{v_i \in \mathcal{V}} \left(\frac{\delta_{out}\left(v_i, a_q\right)}{\deg_{out}\left(v_i\right)} - \frac{\delta_{in}\left(v_i, a_q\right)}{\deg_{in}\left(v_i\right)}\right) &\\
		& ~ w_I \left(v_i\right)^\alpha w_G \left(v_i\right)^\zeta \left\lvert \nabla_a F\left(v_i\right)\right\rvert^{p - 2} \nabla_a F\left(v_i\right) &\\
	\end{aligned}$\\
	$\begin{aligned}[t]	
		= & ~ W_G \left(a_q\right)^\theta \sum_{v_i \in \mathcal{V}} \left(\frac{\delta_{out}\left(v_i, a_q\right)}{\deg_{out}\left(v_i\right)} - \frac{\delta_{in}\left(v_i, a_q\right)}{\deg_{in}\left(v_i\right)}\right) w_I \left(v_i\right)^\alpha w_G \left(v_i\right)^\zeta &\\
		& \quad \left\lvert w_G \left(v_i\right)^\zeta \sum_{a_r \in \mathcal{A}_H} \left(\frac{\delta_{in}\left(v_i, a_r\right)}{\deg_{in}\left(v_i\right)} - \frac{\delta_{out}\left(v_i, a_r\right)}{\deg_{out}\left(v_i\right)}\right) W_I \left(a_r\right)^\beta W_G \left(a_r\right)^\theta F\left(a_r\right) \right\rvert^{p - 2} &\\
		& \quad w_G \left(v_i\right)^\zeta \sum_{a_s \in \mathcal{A}_H} \left(\frac{\delta_{in}\left(v_i, a_s\right)}{\deg_{in}\left(v_i\right)} - \frac{\delta_{out}\left(v_i, a_s\right)}{\deg_{out}\left(v_i\right)}\right) W_I \left(a_s\right)^\beta W_G \left(a_s\right)^\theta F\left(a_s\right) &\\
	\end{aligned}$\\

	Since the vertex weight function $w_G$ only maps to positive values, the following formulation is valid and ultimately result in the hyperarc $p$-Laplacian definition:
	
	$\begin{aligned}[t]
		= & ~ W_G \left(a_q\right)^\theta \sum_{v_i \in \mathcal{V}} \left(\frac{\delta_{out}\left(v_i, a_q\right)}{\deg_{out}\left(v_i\right)} - \frac{\delta_{in}\left(v_i, a_q\right)}{\deg_{in}\left(v_i\right)}\right) w_I \left(v_i\right)^\alpha w_G \left(v_i\right)^{2 \zeta + \zeta \left(p - 2\right)} &\\
		& \quad \left\lvert \sum_{a_r \in \mathcal{A}_H} \left(\frac{\delta_{in}\left(v_i, a_r\right)}{\deg_{in}\left(v_i\right)} - \frac{\delta_{out}\left(v_i, a_r\right)}{\deg_{out}\left(v_i\right)}\right) W_I \left(a_r\right)^\beta W_G \left(a_r\right)^\theta F\left(a_r\right) \right\rvert^{p - 2} &\\
		& \quad \sum_{a_s \in \mathcal{A}_H} \left(\frac{\delta_{in}\left(v_i, a_s\right)}{\deg_{in}\left(v_i\right)} - \frac{\delta_{out}\left(v_i, a_s\right)}{\deg_{out}\left(v_i\right)}\right) W_I \left(a_s\right)^\beta W_G \left(a_s\right)^\theta F\left(a_s\right) &\\
		= & ~ W_G \left(a_q\right)^\theta \sum_{v_i \in \mathcal{V}} \left(\frac{\delta_{out}\left(v_i, a_q\right)}{\deg_{out}\left(v_i\right)} - \frac{\delta_{in}\left(v_i, a_q\right)}{\deg_{in}\left(v_i\right)}\right) w_I \left(v_i\right)^\alpha w_G \left(v_i\right)^{p \zeta} &\\
		& \quad \left\lvert \sum_{a_r \in \mathcal{A}_H} \left(\frac{\delta_{in}\left(v_i, a_r\right)}{\deg_{in}\left(v_i\right)} - \frac{\delta_{out}\left(v_i, a_r\right)}{\deg_{out}\left(v_i\right)}\right) W_I \left(a_r\right)^\beta W_G \left(a_r\right)^\theta F\left(a_r\right) \right\rvert^{p - 2} &\\
		& \quad \sum_{a_s \in \mathcal{A}_H} \left(\frac{\delta_{in}\left(v_i, a_s\right)}{\deg_{in}\left(v_i\right)} - \frac{\delta_{out}\left(v_i, a_s\right)}{\deg_{out}\left(v_i\right)}\right) W_I \left(a_s\right)^\beta W_G \left(a_s\right)^\theta F\left(a_s\right) &\\
		= & ~ \Delta_a^p F \left(a_q\right)
	\end{aligned}$\\
	
	Thus, the equality $\Delta_a^p F\left(a_q\right) = \text{div}_a \left(\left\lvert \nabla_a F\right\rvert^{p - 2} \nabla_a F\right)\left(a_q\right)$ holds true for all hyperarcs $a_q \in \mathcal{A}_H$ and for all hyperarc functions $F \in \mathcal{H}\left(\mathcal{A}_H\right)$, when using the definitions of the hyperarc gradient $\nabla_a$, the hyperarc divergence $\text{div}_a$, and the hyperarc $p$-Laplacian $\Delta_a^p$.\\
\end{proof}
\clearpage
\section{Conclusion}

The most important contributions of this thesis can be summarized by the following five bullet points:
\begin{itemize}
	\item Section (\ref{2}):
	
	Introducing two adjacency tensors \big(a non-unique version with a dimension matching the not oriented hypergraph adjacency tensor and a unique version with bigger dimension\big) for oriented hypergraphs inspired by the adjacency tensor for not oriented normal graphs in \cite{zhang2019introducing}.
	
	\item Section (\ref{4}):
	
	Proving every normal graph is a hypergraph and analyzing the bipartite normal graph representation and the normal graph with complete subgraphs representation critically, which were briefly mentioned in \cite{mori2015peeling} and \cite{zhou2006learning}, but not explored in detail.
	
	\item Section (\ref{6}) and (\ref{7}):
	
	Generalizing the vertex gradient, adjoint and $p$-Laplacian operators for oriented normal graphs of \cite{elmoataz2015p} to include more parameters and weight functions and defining corresponding arc gradient, adjoint and $p$-Laplacian operators for normal graphs.
	
	\item Section (\ref{8}) and (\ref{9}):
	
	Extending the vertex gradient and adjoint operators and 
	arc gradient and adjoint operators for normal graphs to vertex gradient and adjoint operators and hyperarc gradient and adjoint operators for hypergraphs.
	
	\item Section (\ref{10}):
	
	Generalizing the vertex $p$-Laplacian and hyperarc $p$-Laplacian operators from \cite{jost2021plaplaceoperators} to include more parameters and weight functions and to match the definitions of the normal graph case.
\end{itemize}

\big(Note: If a definition, theorem, remark, lemma, corollary or example does not cite any source, then the author of this thesis is responsible for the content and could not find any prior publications about the content.\big)\\

Although the thesis is already very long, there are still many topics which have not been covered, such as deciding on an application and choosing suitable parameters and weight functions for the operators or doing a spectral analysis of the generalized vertex and (hyper)arc $1$-Laplacians. Unfortunately, it was not possible to include further content in this thesis, but at the same time this opens the possibility for future research and exploration.
\clearpage
\urlstyle{same}
\bibliographystyle{natdin}
\bibliography{bibliography}
\addcontentsline{toc}{section}{References}

\end{document}